\documentclass[11pt,twoside]{article}

\usepackage{fullpage}

\usepackage{epsf}
\usepackage{fancyhdr}
\usepackage{graphics}
\usepackage{graphicx}
\usepackage{psfrag}
\usepackage{algorithmicx}
\usepackage{multirow}
\usepackage{color}
\usepackage{subcaption}
\usepackage{graphicx}
\usepackage{amsthm}
\usepackage{amsfonts}
\usepackage{amsmath}
\usepackage{dsfont}
\usepackage{amssymb,bbm}
\usepackage{natbib}
\usepackage{hyperref}[] 
\hypersetup{
colorlinks = true,
linkcolor = blue, 
filecolor = magenta,
urlcolor = blue,
citecolor = blue 
}
\usepackage{nicefrac}
\usepackage{chngpage}
\usepackage{tabularx}
\usepackage{enumitem}
\usepackage{booktabs}
\usepackage{pbox}
\usepackage{caption,subcaption}
\usepackage{mathtools}
\usepackage{fullpage}

\newcommand{\dimension}{p}
\newcommand{\mP}{\mathbb{P}}
\newcommand{\mE}{\mathbb{E}}
\newcommand{\mV}{\text{Var}}

\newcommand{\ind}{\mathds{1}}

\DeclarePairedDelimiter\floor{\lfloor}{\rfloor}

\newtheorem{theorem}{Theorem}[section]
\newtheorem{definition}[theorem]{Definition}
\newtheorem{lemma}[theorem]{Lemma}

\newtheorem{proposition}[theorem]{Proposition}
\newtheorem{corollary}[theorem]{Corollary}
\newtheorem{remark}[theorem]{Remark}
\newtheorem{example}[theorem]{Example}

\begin{document}

\begin{center}
	
	{\bf{\LARGE{Minimax optimality of permutation tests}}}
	\vspace*{.2in}
	
	\begin{author}
		A
		Ilmun Kim$^{\dagger}$
		~ ~ Sivaraman Balakrishnan$^{\dagger}$ 
		~ ~  Larry Wasserman$^{\dagger}$ \\
		\texttt{\{ilmunk, siva, larry\}@stat.cmu.edu}
	\end{author}
	
	\vspace*{.2in}
	
	\begin{tabular}{c}
		Department of Statistics and Data Science$^\dagger$ \\
		Carnegie Mellon University\\
		Pittsburgh, PA 15213
	\end{tabular}

	\vspace*{.2in}

	\today

	\vspace*{.2in}
	
	\begin{abstract}
		Permutation tests are widely used in statistics, providing a finite-sample guarantee on the type I error rate whenever the distribution of the samples under the null hypothesis is invariant to some rearrangement.  Despite its increasing popularity and empirical success, theoretical properties of the permutation test, especially its power, have not been fully explored beyond simple cases. In this paper, we attempt to partly fill this gap by presenting a general non-asymptotic framework for analyzing the minimax power of the permutation test. The utility of our proposed framework is illustrated in the context of two-sample and independence testing under both discrete and continuous settings. In each setting, we introduce permutation tests based on $U$-statistics and study their minimax performance. We also develop exponential concentration bounds for permuted $U$-statistics based on a novel coupling idea, which may be of independent interest. 
		Building on these exponential bounds, we introduce permutation tests which are adaptive to unknown smoothness parameters without losing much power. The proposed framework is further illustrated using more sophisticated test statistics including weighted $U$-statistics for multinomial testing and Gaussian kernel-based statistics for density testing. Finally, we provide some simulation results that further justify the permutation approach.
	\end{abstract}
	
\end{center}

\vskip 1em

\section{Introduction}
A permutation test is a nonparametric approach to hypothesis testing routinely used in a variety of scientific applications such as astronomy~\citep{efron1992simple,freeman2017local}, biology~\citep{pesarin2001multivariate,blackford2009detecting}, neuroscience~\citep{hayasaka2004combining,maris2007nonparametric} and genomics~\citep{stranger2007population,maglietta2007selection}. The permutation test constructs the resampling distribution of a test statistic by permuting the labels of the observations. The resampling distribution, also called the permutation distribution, serves as a reference from which to assess the significance of the observed test statistic. A key property of the permutation test is that it provides exact control of the type I error rate for any test statistic whenever the labels are exchangeable under the null hypothesis~\citep[e.g.][]{hoeffding1952large}. Due to this attractive non-asymptotic property, the permutation test has received considerable attention and has been applied to a wide range of statistical tasks including testing independence, two-sample testing, change point detection, clustering, classification, principal component analysis \citep[see][]{anderson2001permutation, kirch2006permutation,park2009permutation,ojala2010permutation,zhou2018eigenvalue}. 

Once the type I error is controlled, the next concern is the type II error or equivalently the power of the resulting test. Despite its increasing popularity and empirical success, the power of the permutation test has yet to be fully understood. A major challenge in this regard is to control its random critical value that can have an intractable distribution. More specifically, in order to show that a test has high power, we need to ensure that the distribution of the test statistic is significantly distant from its critical value under the alternative. For the permutation test, this critical value is a random quantity, defined as a quantile of the permutation distribution. Studying this random value is in general difficult due to the combinatorial nature of the permutation distribution. While some progress has been made as we review in Section~\ref{Section: Challenges in power analysis}, our understanding of the permutation approach is still far from complete, especially in finite-sample scenarios. The purpose of this paper is to attempt to partly fill this gap by developing a general framework for analyzing the non-asymptotic minimax type II error of the permutation test. Using the tools developed in this paper, we demonstrate that permutation tests are minimax rate optimal in various scenarios.

\subsection{Alternative approaches and their limitations} \label{Section: Alternative approaches and their limitations}

\begin{figure}[t!]
	\begin{center}		
		\begin{minipage}[b]{0.325\textwidth}
			\includegraphics[width=\textwidth]{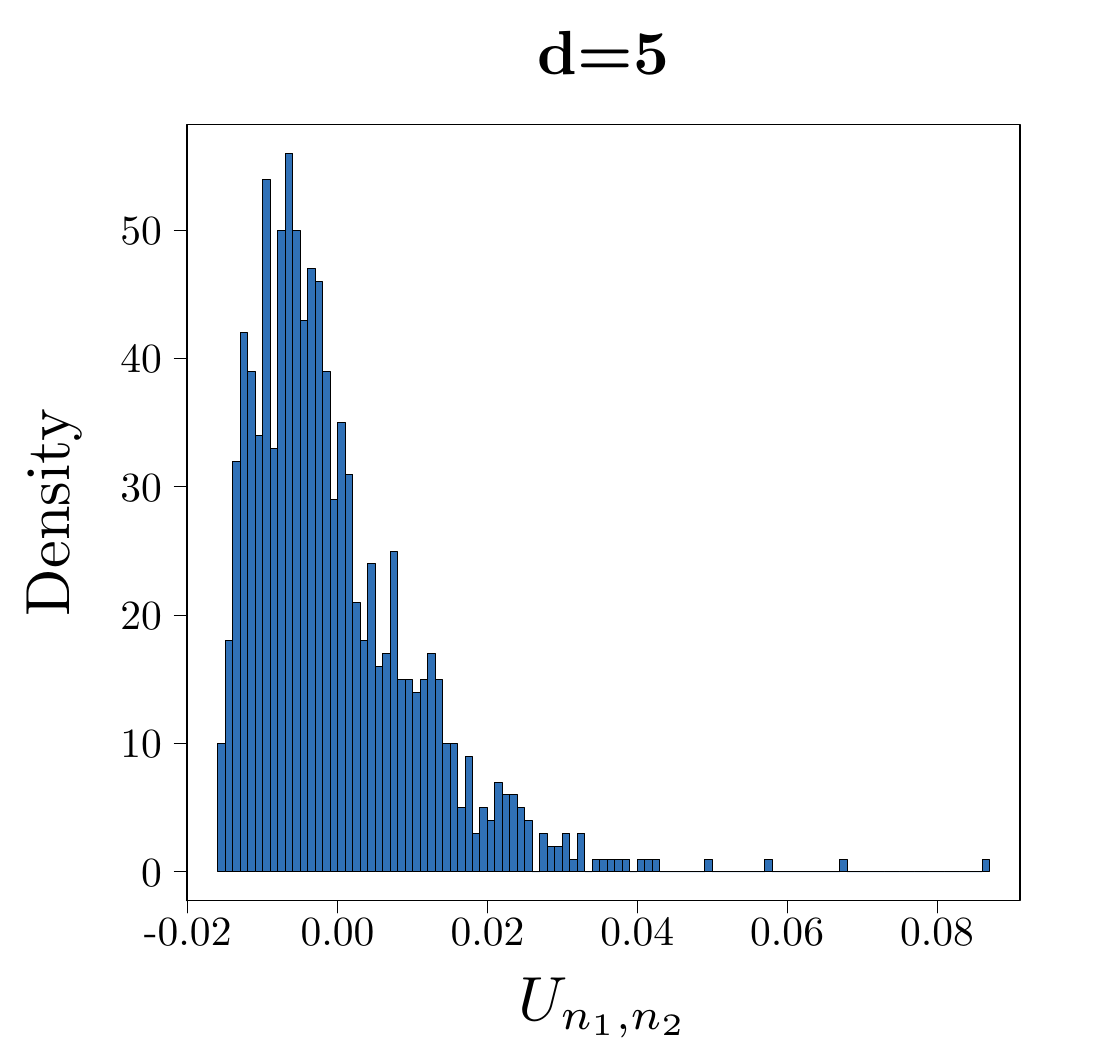}
		\end{minipage} 
		\begin{minipage}[b]{0.325\textwidth}
			\includegraphics[width=\textwidth]{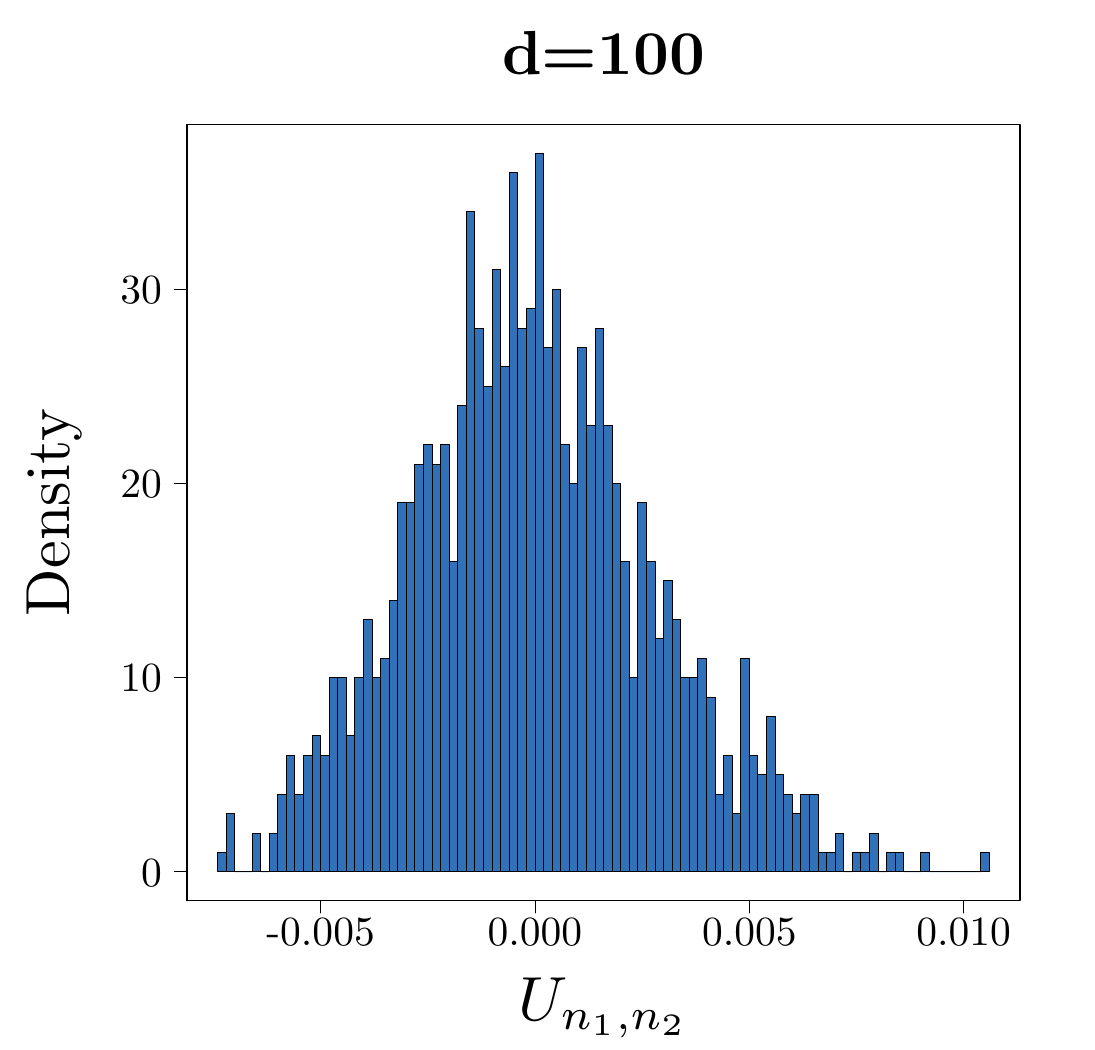}
		\end{minipage}
		\begin{minipage}[b]{0.325\textwidth}
			\includegraphics[width=\textwidth]{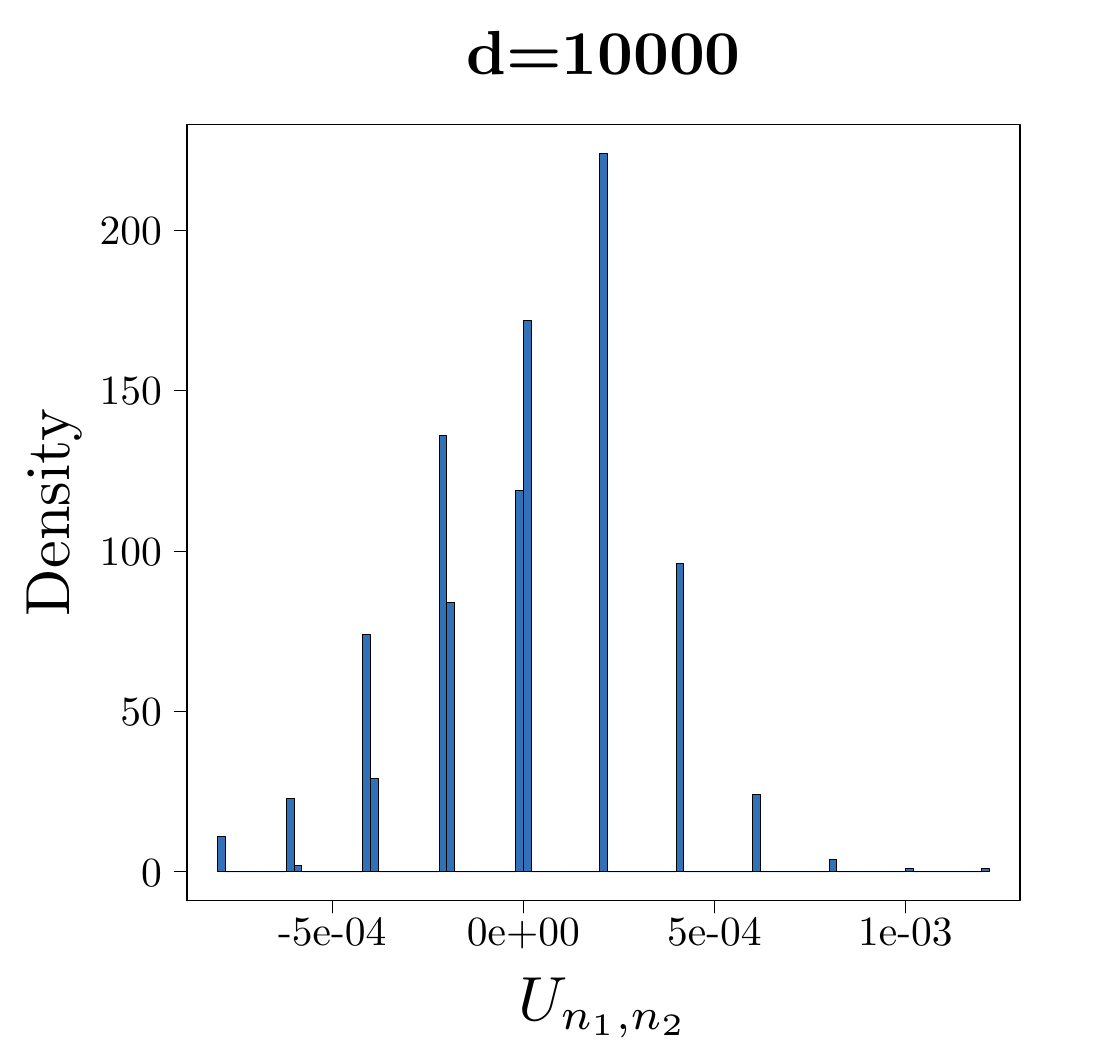}
		\end{minipage}
		\caption{\small Histograms of the $U$-statistic in Proposition~\ref{Proposition: Multinomial Two-Sample Testing} calculated under the uniform multinomial null by varying the number of bins $d$. The plots show that the shape of the null distribution is highly influenced by the number of bins and thus illustrate challenges of estimating the null distribution consistently over different scenarios. See Section~\ref{Section: Alternative approaches and their limitations} for details.} \label{Figure: histograms}
	\end{center}
\end{figure}

We first review a couple of other testing procedures and highlight the advantages of the permutation method. One common approach to determining the critical value of a test is based on the asymptotic null distribution of a test statistic. The validity of a test whose rejection region is calibrated using this asymptotic null distribution is well-studied in the classical regime where the number of parameters is held fixed and the sample size goes to infinity. However, it is no longer trivial to justify this asymptotic approach in a complex, high-dimensional setting where numerous parameters can interact in a non-trivial way and strongly influence the behavior of the test statistic. In such a case, the limiting null distribution is perhaps intractable without imposing stringent assumptions. To illustrate the challenge clearly, we consider the two-sample $U$-statistic $U_{n_1,n_2}$ defined later in Proposition~\ref{Proposition: Multinomial Two-Sample Testing} for multinomial testing. Here we compute $U_{n_1,n_2}$ based on samples from the multinomial distribution with uniform probabilities. To approximate the null distribution of $U_{n_1,n_2}$, we perform 1000 Monte-Carlo simulations for each number of bins $d \in \{5,100,10000\}$ while fixing the sample sizes as $n_1=n_2=100$. From the histograms in Figure~\ref{Figure: histograms}, we see that the shape of the null distribution heavily depends on the number of bins $d$ (and 
more generally on the probabilities of the null multinomial distribution). In particular, the null distribution tends to be more symmetric and sparser as $d$ increases. Since the underlying structure of the distribution is unknown beforehand, Figure~\ref{Figure: histograms} emphasizes difficulties of approximating the null distribution over different regimes. This in turn has led statisticians to 
impose 
stringent assumptions under which test statistics have simple, tractable limiting distributions. However, as noted in \cite{balakrishnan2018hypothesis}, this 
can exclude many high-dimensional cases when -- despite having non-normal null distributions -- carefully designed tests can have high (minimax) power.
We also note that the asymptotic approach does not have any finite sample guarantee, which is also true for other data-driven methods including bootstrapping~\citep{efron1994introduction} and subsampling~\citep{politis1999subsampling}. In sharp contrast, the permutation approach provides a valid test for any test statistic in any sample size under minimal assumptions. Furthermore, as we shall see, one can achieve minimax power through the permutation test even when a nice limiting null distribution is not available.

Another approach, that is commonly used in theoretical computer science, is based on concentration inequalities \citep[e.g.][]{chan2014optimal,acharya2014sublinear,bhattacharya2015testing,diakonikolas2016new,canonne2018testing}. In this approach the threshold of a test is determined using a tail bound of the test statistic under the null hypothesis. Then, owing to the non-asymptotic nature of the concentration bound, the resulting test can control the type I error rate in finite samples. This non-asymptotic approach is more robust to distributional assumptions than the previous asymptotic approach but comes with different challenges. For instance the resulting test tends to be too conservative as it depends on a loose tail bound.  
A more serious problem is that the threshold often relies on unspecified constants and even unknown parameters. By contrast, the permutation approach is entirely data-dependent and tightly controls the type I error rate.

\subsection{Challenges in power analysis and related work} \label{Section: Challenges in power analysis}
Having motivated the importance of the permutation approach, we now review the previous studies on permutation tests and also discuss challenges. The large sample power of the permutation test has been investigated by a number of authors including \cite{hoeffding1952large,robinson1973large,albers1976asymptotic,peter1978asymptotic,janssen1997studentized}. The main result in this line of research indicates that the permutation distribution of a certain test statistic (e.g.~Student's $t$-statistic and $F$-statistic) approximates its null distribution in large sample scenarios. Moreover this approximation is valid under both the null and local alternatives, which then guarantees that the permutation test is asymptotically as powerful as the test based on the asymptotic null distribution. In addition to these findings, power comparisons between permutation and bootstrap tests have been made by \cite{romano1989bootstrap,janssen2003bootstrap,janssen2005resampling} and among others. 
However these power analyses, which rely heavily on classical asymptotic theory, are not easily generalized to more complex settings.  In particular, they often require that alternate distributions satisfy certain regularity conditions under which the asymptotic power function is analytically tractable. Due to such restrictions, the focus has been on a limited class of test statistics applied to a relatively small set of distributions. Furthermore, most previous studies have studied the pointwise, instead of uniform, power that holds for any fixed sequence of alternatives but not uniformly over the class of alternatives.

Recently, there has been another line of research studying the power of the permutation test from a non-asymptotic point of view \citep[e.g.][]{albert2015tests,albert2019concentration,kim2018robust,kim2019global}. This framework, based on a concentration bound for a permuted test statistic, allows us to study the power in more general and complex settings than the asymptotic approach at the expense of being less precise (mainly in terms of constant factors). 
The main challenge in the non-asymptotic analysis, however, is to control the random critical value of the test. The distribution of this random critical value is in general difficult to study due to the non-i.i.d.~structure of the permuted test statistic. Several attempts have been made to overcome such difficulty focusing on linear-type statistics~\citep{albert2019concentration}, regressor-based statistics \citep{kim2019global}, the Cram\'{e}r--von Mises statistic~\citep{kim2018robust} and maximum-type kernel-based statistics~\citep{kim2019comparing}. Our work contributes to this line of research by developing some general tools for studying the finite-sample performance of permutation tests with a specific focus on degenerate $U$-statistics.

Concurrent with our work, and independently, \cite{berrett2020optimal} also develop results for the permutation test based on a degenerate $U$-statistic. While focusing on independence testing, \cite{berrett2020optimal} prove that one cannot hope to have a valid independence test that is uniformly powerful over alternatives in the $L_2$ distance. The authors then impose Sobolev-type smoothness conditions as well as boundedness conditions on density functions under which the proposed permutation test is minimax rate optimal in the $L_2$ distance\footnote{Throughout this paper, we distinguish the $L_p$ distance from the $\ell_p$ distance --- the former is defined with respect to Lebesgue measure and the latter is defined with respect to the counting measure.}. 

Finally, we also note 
that the robustness of permutation tests to violations of the exchangeability condition has been investigated by \cite{romano1990behavior,janssen1997studentized,chung2013exact,pauly2015asymptotic,chung2016multivariate,diciccio2017robust}.

\subsection{Overview of our results}
In this paper we take the non-asymptotic point of view as in \cite{albert2015tests} and establish general results to shed light on the power of permutation tests under a variety of scenarios. To concretely demonstrate our results, we focus on two canonical testing problems: 1) \emph{two-sample testing} and 2) \emph{independence testing}, for which the permutation approach rigorously controls the type I error rate (Section~\ref{Section: Background} for specific settings). These topics have been explored by a number of researchers across diverse fields including statistics and computer science and several optimal tests have been proposed in the minimax sense \citep[e.g.][]{chan2014optimal,bhattacharya2015testing,diakonikolas2016new,arias2018remember}. Nevertheless the existing optimal tests are mostly of theoretical interest, depending on loose or practically infeasible critical values. Motivated by this gap between theory and practice, the primary goal of this study is to introduce permutation tests that tightly control the type I error rate and have the same optimality guarantee as the existing optimal tests.

We summarize the major contributions of this paper and contrast them with the previous studies as follows:

\begin{itemize}
	\item \textbf{Two moments method (Lemma~\ref{Lemma: Two Moments Method}).} Leveraging the quantile approach introduced by \cite{fromont2013two} (see Section~\ref{Section: A general strategy with two moments} for details), we first present a general sufficient condition under which the permutation test has non-trivial power. This condition only involves the first two moments of a test statistic, hence called the \emph{two moments method}. To make this general condition more concrete, we consider degenerate $U$-statistics for two-sample testing and independence testing, respectively, and provide simple moment conditions that ensure that the resulting permutation test has non-trivial power for each testing problem. We then illustrate the efficacy of our results with concrete examples.
	\item \textbf{Multinomial testing (Proposition~\ref{Proposition: Multinomial Two-Sample Testing} and Proposition~\ref{Proposition: Multinomial Independence Testing}).} One example that we focus on is multinomial testing in the $\ell_2$ distance. \cite{chan2014optimal} study the multinomial two-sample problem in the $\ell_2$ distance but with some unnecessary conditions (e.g.~equal sample size, Poisson sampling, known squared norms etc). We remove these conditions and propose a permutation test that is minimax rate optimal for the two-sample problem. Similarly we introduce a minimax optimal test for independence testing in the $\ell_2$ distance based on the permutation procedure. 
	\item \textbf{Density testing (Proposition~\ref{Proposition: Two-Sample Testing for Holder Densities} and Proposition~\ref{Proposition: Minimum Separation for Independence Testing for Holder Densities}).} Another example that we focus on is density testing for H\"{o}lder classes. 
	Building on the work of \cite{ingster1987minimax}, the two-sample problem for H\"{o}lder densities has been recently studied by \cite{arias2018remember} and the authors propose an optimal test in the minimax sense. However their test depends on a loose critical value and also assumes equal sample sizes. We propose an alternative test based on the permutation procedure without such restrictions and show that it achieves the same minimax optimality. We also contribute to the literature by presenting an optimal permutation test for independence testing over H\"{o}lder classes. 
	\item \textbf{Combinatorial concentration inequalities (Theorem~\ref{Theorem: Two-Sample Concentration}, Theorem~\ref{Theorem: Concentration inequality for independence U-statistic} and Theorem~\ref{Theorem: Concentration inequality for independence U-statistic II}).} Although our two moments method is general, it might be sub-optimal in terms of the dependence on a nominal level $\alpha$. Focusing on degenerate $U$-statistics, we improve the dependence on $\alpha$ from polynomial to logarithmic with some extra assumptions. To do so, we develop combinatorial concentration inequalities inspired by the symmetrization trick \citep{duembgen1998symmetrization} and Hoeffding's average \citep{hoeffding1963probability}. We apply the developed inequalities to introduce adaptive tests to unknown smoothness parameters at the cost of $\log\log n$ factor. In contrast to the previous studies \citep[e.g.][]{chatterjee2007stein,bercu2015concentration,albert2019concentration} that are restricted to simple linear statistics, the proposed combinatorial inequalities are broadly applicable to the class of degenerate $U$-statistics. These results have potential applications beyond the problems considered in this paper (e.g.~providing concentration inequalities under sampling without replacement). 
\end{itemize}

In addition to the testing problems mentioned above, we also contribute to multinomial testing problems in the $\ell_1$ distance \citep[e.g.][]{chan2014optimal,bhattacharya2015testing,diakonikolas2016new}. First we revisit the chi-square test for multinomial two-sample testing considered in \cite{chan2014optimal} and show that the test based on the same test statistic but calibrated by the permutation procedure is also minimax rate optimal under Poisson sampling (Theorem~\ref{Theorem: Two-Sample testing under Poisson sampling}). Next, motivated by the flattening idea in \cite{diakonikolas2016new}, we introduce permutation tests based on weighted $U$-statistics and prove their minimax rate optimality for multinomial testing in the $\ell_1$ distance (Proposition~\ref{Proposition: Multinomial L1 Testing} and Proposition~\ref{Proposition: Multinomial independence testing in L1 distance}). Lastly, building on the recent work of \cite{meynaoui2019aggregated}, we analyze the permutation tests based on the maximum mean discrepancy \citep{gretton2012kernel} and the Hilbert--Schmidt independence criterion \citep{gretton2005measuring} for two-sample and independence testing, respectively, and illustrate their performance over certain Sobolev-type smooth function classes.

\subsection{Outline of the paper}
The remainder of the paper is organized as follows. Section~\ref{Section: Background} describes the problem setting and provides some background on the permutation procedure and minimax optimality. In Section~\ref{Section: A general strategy with two moments}, we give a general condition based on the first two moments of a test statistic under which the permutation test has non-trivial power. We concretely illustrate this condition using degenerate $U$-statistics for two-sample testing in Section~\ref{Section: The two moments method for two-sample testing} and for independence testing in Section~\ref{Section: The two moments method for independence testing}. Section~\ref{Section: Combinatorial concentration inequalities} is devoted to combinatorial concentration bounds for permuted $U$-statistics. Building on these results, we propose adaptive tests to unknown smoothness parameters in Section~\ref{Section: Adaptive tests to unknown parameters}. The proposed framework is further demonstrated using more sophisticated statistics in Section~\ref{Section: Further applications}. We present some simulation results that justify the permutation approach in Section~\ref{Section: Simulations} before concluding the paper in Section~\ref{Section: Discussion}. Additional results including concentration bounds for permuted linear statistics and the proofs omitted from the main text are provided in the appendices.

\paragraph{Notation.}
We use the notation $X \overset{d}{=} Y$ to denote that $X$ and $Y$ have the same distribution. The set of all possible permutations of $\{1,\ldots,n\}$ is denoted by $\mathbf{\Pi}_n$. For two deterministic sequences $a_n$ and $b_n$, we write $a_n \asymp b_n$ if $a_n/b_n$ is bounded away from zero and $\infty$ for large $n$. For integers $p,q$ such that $1 \leq q \leq p$, we let $(p)_q = p(p-1)\cdots{(p-q+1)}$. We use $\mathbf{i}_q^p$ to denote the set of all $q$-tuples drawn without replacement from the set $\{1,\ldots,p\}$. $C,C_1,C_2,\ldots,$ refer to positive absolute constants whose values may differ in different parts of the paper. We denote a constant that might depend on fixed parameters $\theta_1,\theta_2,\theta_3,\ldots$ by $C(\theta_1,\theta_2,\theta_3,\ldots)$. Given positive integers $p$ and $q$, we define $\mathbb{S}_p:=\{1,\ldots,p\}$ and similarly $\mathbb{S}_{p,q}:=\{1,\ldots,p\} \times \{1,\ldots,q\}$.

\section{Background} \label{Section: Background}
We start by formulating the problem of interest. Let $\mathcal{P}_0$ and $\mathcal{P}_1$ be two disjoint sets of distributions (or pairs of distributions) on a common measurable space. We are interested in testing whether the underlying data generating distributions belong to $\mathcal{P}_0$ or $\mathcal{P}_1$ based on mutually independent samples $\mathcal{X}_n:=\{X_1,\ldots,X_n \}$. Two specific examples of $\mathcal{P}_0$ and $\mathcal{P}_1$ are:
\begin{enumerate}
	\item \textbf{Two-sample testing.} Let $(P_Y,P_Z)$ be a pair of distributions that belongs to a certain family of pairs of distributions $\mathcal{P}$. Suppose we observe $\mathcal{Y}_{n_1} :=\{Y_1,\ldots,Y_{n_1}\} \overset{\text{i.i.d.}}{\sim} P_Y$ and, independently, $\mathcal{Z}_{n_2} := \{Z_1,\ldots,Z_{n_2}\} \overset{\text{i.i.d.}}{\sim} P_Z$ and denote the pooled samples by $\mathcal{X}_n := \mathcal{Y}_{n_1} \cup \mathcal{Z}_{n_2}$. Given the samples, two-sample testing is concerned with distinguishing the hypotheses:
	\begin{align*}
	H_0: P_Y = P_Z \quad \text{versus} \quad H_1: \delta(P_Y,P_Z) \geq \epsilon_{n_1,n_2},
	\end{align*}
	where $\delta(P_Y,P_Z)$ is a certain distance between $P_Y$ and $P_Z$ and $\epsilon_{n_1,n_2} > 0$. In this case, $\mathcal{P}_0$ is the set of $(P_Y,P_Z) \in \mathcal{P}$ such that $P_Y = P_Z$, whereas $\mathcal{P}_1:= \mathcal{P}_1(\epsilon_{n_1,n_2})$ is another set of $(P_Y,P_Z) \in \mathcal{P}$ such that $\delta(P_Y,P_Z) \geq \epsilon_{n_1,n_2}$.
	\item \textbf{Independence testing.} Let $P_{YZ}$ be a joint distribution of $Y$ and $Z$ that belongs to a certain family of distributions $\mathcal{P}$. Let $P_YP_Z$ denote the product of their marginal distributions. Suppose we observe $\mathcal{X}_n := ((Y_1,Z_1),\ldots,(Y_n,Z_n)) \overset{\text{i.i.d.}}{\sim} P_{YZ}$. Given the samples, the hypotheses for testing independence are
	\begin{align*}
	H_0: P_{YZ} = P_YP_Z \quad \text{versus} \quad H_1: \delta(P_{YZ}, P_YP_Z) \geq \epsilon_n,
	\end{align*}
	where $\delta(P_{YZ},P_YP_Z)$ is a certain distance between $P_{YZ}$ and $P_YP_Z$ and $\epsilon_n > 0$. In this case, $\mathcal{P}_0$ is the set of $P_{YZ} \in \mathcal{P}$ such that $P_{YZ} = P_YP_Z$, whereas $\mathcal{P}_1:= \mathcal{P}_1(\epsilon_n)$ is another set of $P_{YZ} \in \mathcal{P}$ such that $\delta(P_{YZ},P_YP_Z) \geq \epsilon_n$. 
\end{enumerate}
Let us consider a generic test statistic $T_n:=T_n(\mathcal{X}_n)$, which is designed to distinguish between the null and alternative hypotheses based on $\mathcal{X}_n$. Given a critical value $c_n$ and pre-specified constants $\alpha \in (0,1)$ and $\beta \in (0, 1-\alpha)$, the problem of interest is to find sufficient conditions on $\mathcal{P}_0$ and $\mathcal{P}_1$ under which the type I and II errors of the test $\ind(T_n > c_n)$ are uniformly bounded as
\begin{equation}
\begin{aligned} \label{Eq: uniform error control}
& \bullet ~ \text{Type I error:} ~ \sup_{P \in \mathcal{P}_0} \mP_P^{(n)} \left(T_n > c_{n}\right) \leq \alpha, \\
& \bullet ~ \text{Type II error:} ~ \sup_{P \in \mathcal{P}_1} \mP_P^{(n)} \left(T_n \leq c_{n} \right) \leq \beta.
\end{aligned}
\end{equation}
Our goal is to control these uniform (rather than pointwise) errors based on data-dependent critical values determined by the permutation procedure.

\subsection{Permutation procedure}
This section briefly overviews the permutation procedure and its well-known theoretical properties, referring readers to \cite{lehmann2006testing,pesarin2010permutation} for more details. Let us begin with some notation. Given a permutation $\pi := (\pi_1,\ldots,\pi_n) \in \mathbf{\Pi}_n$, we denote the permuted version of $\mathcal{X}_n$ by $\mathcal{X}_n^\pi$, that is, $\mathcal{X}_n^\pi:=\{X_{\pi_1},\ldots,X_{\pi_n}\}$. For the case of independence testing, $\mathcal{X}_n^{\pi}$ is defined by permuting the second variable $Z$, i.e.~$\mathcal{X}_n^\pi:= \{(Y_1,Z_{\pi_1}),\ldots,(Y_n,Z_{\pi_n})\}$. We write $T_n^{\pi}: = T_n(\mathcal{X}_n^{\pi})$ to denote the test statistic computed based on $\mathcal{X}_n^{\pi}$. Let $F_{T_n^{\pi}}(t)$ be the permutation distribution function of $T_n^{\pi}$ defined as
\begin{align*}
F_{T_n^{\pi}}(t) := M_n^{-1} \sum\nolimits_{\pi \in \mathbf{\Pi}_n} \ind\{T_n(\mathcal{X}_n^{\pi}) \leq t \}.
\end{align*}
Here $M_n$ denotes the cardinality of $\mathbf{\Pi}_n$. We write the $1-\alpha$ quantile of $F_{T_n^{\pi}}$ by $c_{1-\alpha,n}$ defined as
\begin{align} \label{Eq: permutation critical value}
c_{1-\alpha,n}  := \inf\{t: F_{T_n^{\pi}}(t) \geq 1 - \alpha \}.
\end{align}
Given the quantile $c_{1-\alpha,n}$, the permutation test rejects the null hypothesis when $T_n > c_{1-\alpha,n}$. This choice of the critical value provides finite-sample type I error control under the permutation-invariant assumption (or exchangeability). In more detail, the distribution of $\mathcal{X}_n$ is said to be permutation invariant if $\mathcal{X}_n$ and $\mathcal{X}_n^{\pi}$ have the same distribution whenever the null hypothesis is true. This permutation-invariance holds for two-sample and independence testing problems. When permutation-invariance holds, it is well-known that the permutation test $\ind(T_n > c_{1-\alpha,n} )$ has level $\alpha$, and by randomizing the test function we can also ensure it has size $\alpha$~\citep[see e.g.][]{hoeffding1952large,lehmann2006testing,hemerik2018exact}.

\vskip 1em

\begin{remark}[Computational aspects] \normalfont \label{Remark: Computational aspects}
	Exact calculation of the critical value~(\ref{Eq: permutation critical value}) is computationally prohibitive except for small sample sizes. Therefore it is common practice to use Monte Carlo simulations to approximate the critical value~\citep[e.g.][]{romano2005exact}. We note that this approximation error can be made arbitrary small by taking a sufficiently large number of Monte Carlo samples. We make this argument more rigorous in Appendix~\ref{Section: Monte Carlo-based permutation tests} and show that our master theorem for the permutation test (Lemma~\ref{Lemma: Two Moments Method}) also holds for its Monte Carlo counterpart as long as the number of Monte Carlo samples is greater than some constant that only depends on the pre-specified error rates $\alpha$ and $\beta$.
\end{remark}

\subsection{Minimax optimality} \label{Section: Minimax optimality}
A complementary aim of this paper is to show that the sufficient conditions for the error bounds in (\ref{Eq: uniform error control}) are indeed necessary in some applications. We approach this problem from the
minimax perspective pioneered by \cite{ingster1987minimax}, and further developed in subsequent works \citep{ingster2003nonparametric,ingster1993asymptotically,baraud2002non,lepski1999}. Let us define a test $\phi$, which is a Borel measurable map, $\phi: \mathcal{X}_n \mapsto \{0,1\}$. For a class of null distributions $\mathcal{P}_0$, we denote the set of all level $\alpha$ tests by
\begin{align*}
	\Phi_{n,\alpha} := \bigg\{ \phi : \sup_{P \in \mathcal{P}_0} \mP_P^{(n)} \left(\phi = 1 \right) \leq \alpha \bigg\}.
\end{align*}
Consider a class of alternative distributions $\mathcal{P}_1(\epsilon_n)$ associated with a positive sequence $\epsilon_n$. Two specific examples of this class of interest are $\displaystyle \mathcal{P}_1(\epsilon_{n_1,n_2}):=\{(P_Y,P_Z) \in \mathcal{P} : \delta(P_Y,P_Z) \geq \epsilon_{n_1,n_2} \}$ for two-sample testing and $\displaystyle \mathcal{P}_1(\epsilon_n) := \{P_{YZ} \in \mathcal{P}: \delta(P_{YZ},P_YP_Z) \geq \epsilon_n \}$ for independence testing. Given $\mathcal{P}_1(\epsilon_n)$, the maximum type II error of a test $\phi \in \Phi_{n,\alpha}$ is 
\begin{align*}
	R_{n,\epsilon_n}(\phi) := \sup_{P \in \mathcal{P}_1(\epsilon_n)} \mP_P^{(n)} (\phi = 0),
\end{align*}
and the minimax risk is defined as
\begin{align*}
	R_{n,\epsilon_n}^\dagger := \inf_{\phi \in \Phi_{n,\alpha}} R_{n,\epsilon_n}(\phi). 
\end{align*}
The minimax risk is frequently investigated via the minimum separation (or the critical radius), which is the smallest $\epsilon_n$ such that type II error becomes non-trivial. Formally, for some fixed $\beta \in (0,1-\alpha)$, the minimum separation is defined as
\begin{align*}
	\epsilon_n^\dagger := \inf \Big\{ \epsilon_n :  R_{n,\epsilon_n}^\dagger \leq \beta \Big\}.
\end{align*}
Of course, it is infeasible to obtain an optimal test that achieves the exact minimax risk in most realistic cases. We instead use the yardstick of minimax \emph{rate} optimality in order to assess the permutation procedure. Formally, consider a level $\alpha$ test $\phi \in \Phi_{n,\alpha}$. We call $\phi$ minimax rate optimal with separation rate $r_n \asymp \epsilon_n^\dagger$ if the following statement is true for every $n$: \emph{there exists an absolute constant $C>0$ such that, if $\epsilon_n \geq C r_n$, then $R_{n,\epsilon_n}(\phi) \leq \beta$.} This optimality guarantee is non-asymptotic and holds for every $n$.

Next we present a few remarks on the minimax framework and other testing procedures.

\begin{remark}  \leavevmode \normalfont
	\begin{itemize}
		\item In this paper, we focus on the minimax framework to evaluate the performance of tests. The minimax framework enables us to study the fundamental limits of hypothesis testing while ensuring 
		a (strong) uniform guarantee over a large class of (null and alternate) distributions. The notion of minimax performance is widely used to quantify the difficulty of a statistical problem.
		In the hypothesis testing literature, it is also common to study the performance of tests against fixed or directional alternatives. As discussed in past work \citep[see e.g.][for a discussion]{arias2018remember}, consistency against
		fixed or directional alternatives can be misleading and can conceal the curse-of-dimensionality that we would typically expect in non-parametric or high-dimensional settings. On the other hand, when we have additional prior information suggesting alternative hypotheses of importance, a na\"{i}ve use of the minimax framework to design or assess tests may not be adequate. It is certainly possible that a global minimax optimal test may not perform well against a particular fixed alternative~\citep[see e.g. our prior work][for empirical results]{balakrishnan2019hypothesis}. Although we restrict our attention to minimax optimality, it would also be interesting to consider other criteria (e.g.~robustness, computational costs, high power against pre-specified alternatives) to assess the performance of tests. We leave this direction to future work. 
		\item Besides the permutation procedure, bootstrap sampling and subsampling are two common practical ways of determining critical values. However, as mentioned earlier, they are only valid in a limited asymptotic setting where the test statistic has a nice limiting distribution. Therefore the tests based on these asymptotic methods do not belong to $\Phi_{n,\alpha}$ in general, and thus are often not minimax optimal in our framework.
	\end{itemize}
\end{remark}

With this background in place, we now focus on showing the minimax rate optimality of permutation tests in various settings.

\section{A general strategy with first two moments} \label{Section: A general strategy with two moments}
In this section, we discuss a general strategy for studying the testing errors of a permutation test based on the first two moments of a test statistic. As mentioned earlier, the permutation test is level $\alpha$ as long as permutation-invariance holds under the null hypothesis. Therefore we focus on the type II error rate and provide sufficient conditions under which the error bounds given in (\ref{Eq: uniform error control}) are fulfilled. 
Previous approaches to non-asymptotic minimax power analysis, reviewed in Section~\ref{Section: Alternative approaches and their limitations}, 
use a non-random critical value, often derived through upper bounds on the mean and variance of the test statistic under the null, 
and thus do not directly apply to the permutation test. To bridge the gap, we consider a deterministic quantile value that serves as a proxy for the permutation threshold $c_{1-\alpha,n}$. More precisely, let $q_{1-\gamma,n}$ be the $1-\gamma$ quantile of the distribution of the random critical value $c_{1-\alpha,n}$. Then by splitting the cases into $\{ c_{1-\alpha,n} \leq q_{1-\gamma,n}\}$ and $\{ c_{1-\alpha,n} > q_{1-\gamma,n}\}$ and using the definition of the quantile, it can be shown that the type II error of the permutation test is less than or equal to 
\begin{align*}
\sup_{P \in \mathcal{P}_1} \mP_P \left(T_n \leq c_{1-\alpha,n} \right) \leq \sup_{P \in \mathcal{P}_1} \mP_P \left(T_n \leq q_{1-\gamma,n} \right) + \gamma.
\end{align*}
Consequently, if one succeeds in showing that $\sup_{P \in \mathcal{P}_1} \mP_P \left(T_n \leq q_{1-\gamma,n} \right) \leq \gamma^\prime$ with $\gamma'$ such that $\gamma' + \gamma \leq \beta$, then the type II error of the permutation test is bounded by $\beta$ as desired. This quantile approach to dealing with a random threshold is not new and has been considered by \cite{fromont2013two} to study the power of a kernel-based test via a wild bootstrap method. In the next lemma, we build on this quantile approach and study the testing errors of the permutation test based on an arbitrary test statistic. Here and hereafter, we denote the expectation and variance of $T_n^\pi$ with respect to the permutation distribution by $\mE_{\pi} [T_n^{\pi} | \mathcal{X}_n]$ and $\text{Var}_{\pi} [T_n^{\pi} | \mathcal{X}_n]$, respectively.

\begin{lemma}[Two moments method] \label{Lemma: Two Moments Method}
	Suppose that for each permutation $\pi \in \mathbf{\Pi}_n$,  $T_n$ and $T_n^{\pi}$ have the same distribution under the null hypothesis. Given pre-specified error rates $\alpha \in (0,1)$ and $\beta \in (0,1 - \alpha)$, assume that for any $P \in \mathcal{P}_1$,
	\begin{equation}
	\begin{aligned} \label{Eq: Sufficient Condition}
	\mE_P[T_n]  ~\geq~ & \mE_P [\mE_{\pi} \{T_n^{\pi} | \mathcal{X}_n\}] + \sqrt{\frac{3\emph{\mV}_P[\mE_{\pi}\{T_n^{\pi}| \mathcal{X}_n \}]}{\beta}} \\
	&  + \sqrt{\frac{3\emph{\mV}_P[T_n]}{\beta}}  + \sqrt{\frac{3\mE_P[\emph{\mV}_{\pi}\{T_n^{\pi}| \mathcal{X}_n \}]}{\alpha\beta}}.
	\end{aligned}
	\end{equation}
	Then the permutation test $\ind(T_n > c_{1-\alpha,n})$ controls the type I and II error rates as in (\ref{Eq: uniform error control}).
\end{lemma}

\begin{figure}[t!]
	\begin{center}		
		\includegraphics[width=\textwidth]{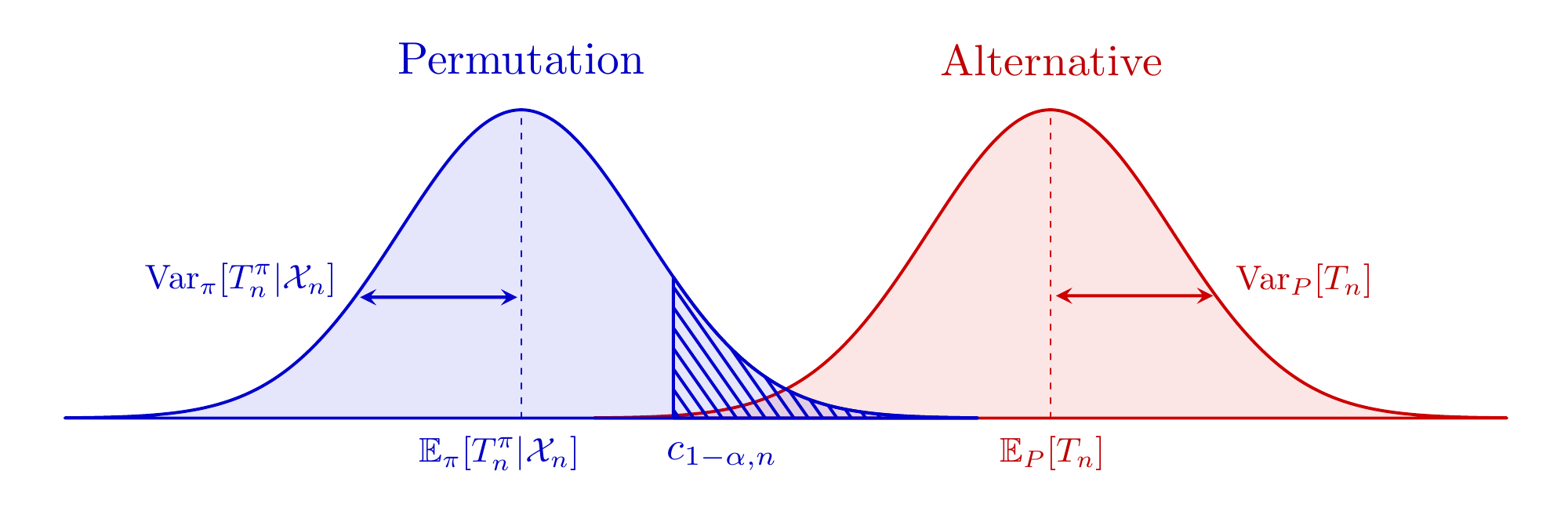}
		\caption{\small An illustration of Lemma~\ref{Lemma: Two Moments Method}. The lemma describes that the major components that determine the power of a permutation test are the mean and the variance of the alternative distribution as well as the permutation distribution. In particular, if the mean of the alternative distribution is sufficiently larger than the other components (on average since the permutation distribution is random), then the permutation test succeeds to reject the null with high probability.} \label{Figure: two moments mehtod}
	\end{center}
\end{figure}

The proof of this general statement follows by simple set algebra along with Markov and Chebyshev's inequalities. The details can be found in Appendix~\ref{Section: Proof of Lemma: Two Moments Method}. At a high-level, the sufficient condition~(\ref{Eq: Sufficient Condition}) roughly says that if the expected value of $T_n$ (say, signal) is much larger than the expected value of the permuted statistic $T_n^\pi$ (say, baseline) as well as the variances of $T_n$ and $T_n^{\pi}$ (say, noise), then the permutation test will have non-trivial power. We provide an illustration of Lemma~\ref{Lemma: Two Moments Method} in Figure~\ref{Figure: two moments mehtod}. Suppose further that $T_n^\pi$ is centered at zero under the permutation law, i.e.~$\mE_{\pi}[T_n^{\pi} | \mathcal{X}_n] = 0$. Then a modification of the proof of Lemma~\ref{Lemma: Two Moments Method} yields a simpler condition with improved constant factors. We show that if,
\begin{align}  \label{Eq: Sufficient Condition II}
\mE_P[T_n]  ~\geq~  \sqrt{\frac{2\mV_P[T_n]}{\beta}}  + \sqrt{\frac{2\mE_P[\mV_{\pi}\{T_n^{\pi}| \mathcal{X}_n \}]}{\alpha\beta}},
\end{align}
then the permutation test has type II error at most $\beta$.
In the following sections, we demonstrate the two moments method (Lemma~\ref{Lemma: Two Moments Method}) based on degenerate $U$-statistics for two-sample and independence testing. 

\section{The two moments method for two-sample testing} \label{Section: The two moments method for two-sample testing}
This section illustrates the two moments method given in Lemma~\ref{Lemma: Two Moments Method} for two-sample testing. By focusing on a $U$-statistic, we first present a general condition that ensures that the type I and II error rates of the permutation test are uniformly controlled (Theorem~\ref{Theorem: Two-Sample U-statistic}). We then turn to more specific cases of two-sample testing for multinomial distributions and H\"{o}lder densities. 

Let $g(x,y)$ be a bivariate function, which is symmetric in its arguments, i.e.~$g(x,y) = g(y,x)$. Based on this bivariate function, let us define a kernel for a two-sample $U$-statistic
\begin{equation}
\begin{aligned} \label{Eq: Two-Sample kernel}
h_{\text{ts}}(y_1,y_2;z_1,z_2)  :=  g(y_1,y_2) + g(z_1,z_2) - g(y_1,z_2) - g(y_2,z_1),
\end{aligned}
\end{equation}
and write the corresponding $U$-statistic as
\begin{align} \label{Eq: Two-Sample U-statistic}
U_{n_1,n_2} := \frac{1}{(n_1)_{(2)}(n_2)_{(2)}} \sum_{(i_1,i_2) \in \mathbf{i}_2^{n_1}} \sum_{(j_1,j_2) \in \mathbf{i}_2^{n_2}} h_{\text{ts}}(Y_{i_1},Y_{i_2}; Z_{j_1},Z_{j_2}).
\end{align}
Depending on the choice of kernel $h_{\text{ts}}$, the $U$-statistic includes frequently used two-sample test statistics in the literature such as the maximum mean discrepancy \citep{gretton2012kernel} and the energy statistic \citep{baringhaus2004new,szekely2004testing}. From the basic properties of $U$-statistics \citep[e.g.][]{lee1990u}, it is readily seen that $U_{n_1,n_2}$ is an unbiased estimator of $\mE_P[h_{\text{ts}}(Y_1,Y_2;Z_1,Z_2)]$. To describe the main result of this section, let us write the symmetrized kernel by
\begin{align} \label{Eq: symmetrized kernel}
\overline{h}_{\text{ts}} (y_1,y_2;z_1,z_2) := \frac{1}{2!2!} \sum_{(i_1,i_2) \in \mathbf{i}_2^{2}} \sum_{(j_1,j_2) \in \mathbf{i}_2^{2}} h_{\text{ts}}(y_{i_1},y_{i_2}; z_{j_1},z_{j_2}),
\end{align}
and define $\psi_{Y,1}(P),\psi_{Z,1}(P)$ and $\psi_{YZ,2}(P)$ by  
\begin{equation}
\begin{aligned} \label{Eq: definitions of psi functions}
\psi_{Y,1}(P) &:= \text{Var}_P[\mE_P\{\overline{h}_{\text{ts}}(Y_1,Y_2;Z_1,Z_2)|Y_1\}], \\[.5em]
\psi_{Z,1}(P) &:= \text{Var}_P[\mE_P\{\overline{h}_{\text{ts}}(Y_1,Y_2;Z_1,Z_2)|Z_1\}], \\[.5em]
\psi_{YZ,2}(P) &:= \max \{ \mE_P[g^2(Y_1,Y_2)], \  \mE_P[g^2(Y_1,Z_1)], \ \mE_P[g^2(Z_1,Z_2)] \}.
\end{aligned}
\end{equation}
As will be clearer in the sequel $\psi_{Y,1}(P),\psi_{Z,1}(P)$ and $\psi_{YZ,2}(P)$ are key quantities which we use to upper bound the variance of $U_{n_1,n_2}$. By leveraging Lemma~\ref{Lemma: Two Moments Method}, the next theorem presents a sufficient condition that guarantees that the type II error rate of the permutation test based on $U_{n_1,n_2}$ is uniformly bounded by $\beta$. 
\begin{theorem}[Two-sample $U$-statistic] \label{Theorem: Two-Sample U-statistic}
	Suppose that there is a sufficiently large constant $C>0$ such that
	\begin{align} \label{Eq: Two-Sample U-statistic Sufficient Condition}
	\mE_P[U_{n_1,n_2}] \geq C \sqrt{\max \Bigg\{ \frac{\psi_{Y,1}(P)}{\beta n_1}, \ \frac{\psi_{Z,1}(P)}{\beta n_2}, \  \frac{\psi_{YZ,2}(P)}{\alpha\beta}\left( \frac{1}{n_1} + \frac{1}{n_2}\right)^2 \Bigg\}},
	\end{align}
	 for all $P \in \mathcal{P}_1$. Then the type II error of the permutation test over $\mathcal{P}_1$ is uniformly bounded by $\beta$, that is
	\begin{align*}
	\sup_{P \in \mathcal{P}_1} \mP_P^{(n_1,n_2)} (U_{n_1,n_2} \leq c_{1-\alpha,n_1,n_2} )\leq \beta. 
	\end{align*}
\end{theorem}

\vskip 1em

\begin{remark} \leavevmode \normalfont
\begin{itemize}
\item This result, applicable broadly to degenerate $U$-statistics of the form in~\eqref{Eq: Two-Sample U-statistic}, simplifies the application of Lemma~\ref{Lemma: Two Moments Method}. The main difficulty in directly 
applying Lemma~\ref{Lemma: Two Moments Method} is that the sufficient condition depends on the conditional variance of the statistic under the permutation distribution which 
can be challenging to upper bound. On the other hand the sufficient condition of this theorem, only depends on the quantities in~\eqref{Eq: definitions of psi functions} which 
do not depend on the permutation distribution.
\item The main technical effort in establishing this result is in showing~\eqref{Eq: Variance upper bound under permutations}, which upper bounds the conditional variance of the test
statistic under the permutation distribution, as a function of the quantity $\psi_{YZ,2}(P)$ defined in~\eqref{Eq: definitions of psi functions}.
\item  In Theorem~\ref{Theorem: Two-Sample U-statistic} and other places, we say that $C$ is a universal constant when it does not depend on other quantities. When the constant is not 
universal but depends on problem-specific quantities (say $\theta_1,\theta_2,\ldots$), which are typically considered fixed (for instance, density smoothness parameters in non-parametric testing), 
we denote this via the explicit notation $C(\theta_1,\theta_2,\ldots)$. Importantly, these constants do not depend on the sample-size. 
\end{itemize}
\end{remark}

\vskip 1em

\begin{proof}[\textbf{\emph{Proof Sketch.}}] Let us give a high-level idea of the proof, while the details are deferred to Appendix~\ref{Section: Proof of Theorem: Two-Sample U-statistic}. First, by the linearity of expectation, it can be verified that the mean of the permuted $U$-statistic $U_{n_1,n_2}^\pi$ is zero. Therefore it suffices to check condition~(\ref{Eq: Sufficient Condition II}). By the well-known variance formula of a two-sample $U$-statistic \citep[e.g.~page 38 of][]{lee1990u}, we prove in Appendix~\ref{Section: Proof of Theorem: Two-Sample U-statistic} that
\begin{align} \label{Eq: Variance upper bound}
\mV_P[U_{n_1,n_2}] \leq C_1 \frac{\psi_{Y,1}(P)}{n_1}  +  C_2 \frac{\psi_{Z,1}(P)}{n_2} + C_3 \psi_{YZ,2}(P) \left( \frac{1}{n_1} + \frac{1}{n_2} \right)^2,
\end{align}
and this result can be used to bound the first term of condition~(\ref{Eq: Sufficient Condition II}). It is worth pointing out that the variance behaves differently under the null and alternative hypotheses. In particular, $\psi_{Y,1}$ and $\psi_{Z,1}$ are zero under the null hypothesis. Hence, in the null case, the third term dominates the variance of $U_{n_1,n_2}$ where we note that $\psi_{YZ,2}$ is a convenient upper bound for the variance of kernel $h_{\text{ts}}$.  Intuitively, the permuted $U$-statistic $U_{n_1,n_2}^\pi$ behaves similarly to $U_{n_1,n_2}$ computed based on samples from a certain null distribution (say a mixture of $P_Y$ and $P_Z$). This implies that the variance of $U_{n_1,n_2}^\pi$ is also dominated by the third term in the upper bound (\ref{Eq: Variance upper bound}). Having this intuition in mind, we use the symmetric structure of kernel $h_{\text{ts}}$ and prove that
\begin{align} \label{Eq: Variance upper bound under permutations}
\mE_P[\mV_{\pi}\{T_n^{\pi}| \mathcal{X}_n \}] \leq C_4 \psi_{YZ,2}(P) \left( \frac{1}{n_1} + \frac{1}{n_2} \right)^2,
\end{align}
which is one of our key technical contributions. Based on the previous two bounds in (\ref{Eq: Variance upper bound}) and (\ref{Eq: Variance upper bound under permutations}), we then complete the proof by verifying the sufficient condition~(\ref{Eq: Sufficient Condition II}). 
\end{proof}

The next two subsections focus on multinomial distributions and H\"{o}lder densities and give explicit expressions for condition (\ref{Eq: Two-Sample U-statistic Sufficient Condition}). We also demonstrate minimax optimality of permutation tests under the given scenarios.

\subsection{Two-sample testing for multinomials} \label{Section: Two-sample testing for multinomials}
Let $p_Y$ and $p_Z$ be multinomial distributions on a discrete domain $\mathbb{S}_d:=\{1,\ldots,d\}$. Throughout this subsection, we consider the kernel $h_{\text{ts}}(y_1,y_2;z_1,z_2)$ in (\ref{Eq: Two-Sample kernel}) defined with the following bivariate function:
\begin{align} \label{Eq: Two-Sample bivariate function}
g_{\text{Multi}}(x,y) := \sum_{k=1}^d \ind (x=k)\ind(y=k).
\end{align}
It is straightforward to see that the resulting $U$-statistic~(\ref{Eq: Two-Sample U-statistic}) is an unbiased estimator of $\|p_Y - p_Z\|_2^2$. Let us denote the maximum between the squared $\ell_2$ norms of $p_Y$ and $p_Z$ by
\begin{align}
b_{(1)} := \max \big\{\|p_Y\|_2^2,\|p_Z\|_2^2 \big\}.
\end{align}
Building on Theorem~\ref{Theorem: Two-Sample U-statistic}, the next result establishes a guarantee on the testing errors of the permutation test under the two-sample multinomial setting. 
\begin{proposition}[Multinomial two-sample testing in $\ell_2$ distance] \label{Proposition: Multinomial Two-Sample Testing}
	Let $\mathcal{P}_{\text{\emph{Multi}}}^{(d)}$ be the set of pairs of multinomial distributions defined on $\mathbb{S}_d$. Let $\mathcal{P}_0 = \{(p_Y,p_Z) \in \mathcal{P}_{\text{\emph{Multi}}}^{(d)} : p_Y = p_Z\}$ and $\mathcal{P}_1(\epsilon_{n_1,n_2}) = \{(p_Y,p_Z) \in \mathcal{P}_{\text{\emph{Multi}}}^{(d)} : \|p_Y - p_Z\|_2 \geq \epsilon_{n_1,n_2} \}$ where
	\begin{align*}
	\epsilon_{n_1,n_2} \geq C\frac{b_{(1)}^{1/4}}{\alpha^{1/4}\beta^{1/2}} \left( \frac{1}{n_1} + \frac{1}{n_2} \right)^{1/2},
	\end{align*}
	for a sufficiently large $C>0$. Consider the two-sample $U$-statistic $U_{n_1,n_2}$ defined with the bivariate function $g_{\text{\emph{Multi}}}$ given in (\ref{Eq: Two-Sample bivariate function}). Then the type I and II error rates of the resulting permutation test are uniformly bounded over the classes $\mathcal{P}_0$ and $\mathcal{P}_1$ as in (\ref{Eq: uniform error control}). 
\end{proposition}

\vskip 1em

\begin{proof}[\textbf{\emph{Proof Sketch.}}] We outline the proof of the result, while the details can be found in Appendix~\ref{Section: Proof of Proposition: Multinomial Two-Sample Testing}. Given the reduction
in Theorem~\ref{Theorem: Two-Sample U-statistic} the remaining technical effort is to show that there exist constants $C_1,C_2,C_3>0$ such that
\begin{equation}
\begin{aligned} \label{Eq: Three conditions for two-sample multinomials}
\psi_{Y,1}(P) & \leq C_1 \sqrt{b_{(1)}} \|p_Y - p_Z\|_2^2, \\[.5em]
\psi_{Z,1}(P) & \leq C_2 \sqrt{b_{(1)}} \|p_Y - p_Z\|_2^2, \\[.5em]
\psi_{YZ,2}(P) & \leq C_3 b_{(1)}.
\end{aligned}
\end{equation}
We note that the U-statistic with kernel in~\eqref{Eq: Two-Sample bivariate function} is similar (although not identical) to the statistic proposed by~\cite{chan2014optimal} for the case when the sample-sizes are equal, and they derive
similar bounds on the variance of their test statistic. These bounds together with Theorem~\ref{Theorem: Two-Sample U-statistic} imply that if there exists a sufficiently large $C_4>0$ such that 
\begin{align} \label{Eq: Two-Sample Bound}
\|p_Y - p_Z\|_2^2 \geq C_4 \frac{\sqrt{b_{(1)}}}{\alpha^{1/2} \beta} \left( \frac{1}{n_1} + \frac{1}{n_2} \right),
\end{align}
then the permutation test based on $U_{n_1,n_2}$ has non-trivial power as claimed. 
\end{proof}

For the balanced case where $n_1=n_2$, \cite{chan2014optimal} prove that no test can have uniform power if $\epsilon_{n_1,n_2}$ is of lower order than $b_{(1)}^{1/4} n_1^{-1/2}$. Hence the permutation test in Proposition~\ref{Proposition: Multinomial Two-Sample Testing} is minimax rate optimal in this balanced setting. The next proposition extends this result to the case of unequal sample sizes and shows that the permutation test is still optimal even for the unbalanced case. 

\begin{proposition}[Minimum separation for two-sample multinomial testing] \label{Proposition: Minimum Separation for Two-Sample Multinomial Testing}
	Consider the two-sample testing problem within the class of multinomial distributions $\mathcal{P}_{\text{\emph{Multi}}}^{(d)}$ where the null hypothesis and the alternative hypothesis are $H_0: p_Y = p_Z$ and $H_1: \|p_Y - p_Z\|_2 \geq \epsilon_{n_1,n_2}$. Under this setting and $n_1 \leq n_2$, the minimum separation satisfies $\epsilon_{n_1,n_2}^{\dagger} \asymp b_{(1)}^{1/4} n_1^{-1/2}$.
\end{proposition}

\vskip 0.5em

\begin{remark}[$\ell_1$- versus $\ell_2$-closeness testing] \normalfont
We note that the minimum separation strongly depends on the choice of metrics. As shown in \cite{bhattacharya2015testing} and \cite{diakonikolas2016new}, the minimum separation rate for two-sample testing in the $\ell_1$ distance is $\max \{ d^{1/2}n_2^{-1/4} n_1^{-1/2}, d^{1/4}n_1^{-1/2}\}$ for $n_1 \leq n_2$. This rate, in contrast to $b_{(1)}^{1/4} n_1^{-1/2}$, illustrates that the difficulty of $\ell_1$-closeness testing depends not only on the smaller sample size $n_1$ but also on the larger sample size $n_2$. In Section~\ref{Section: Two-sample testing via sample-splitting}, we provide a permutation test that is minimax rate optimal in the $\ell_1$ distance.
\end{remark}

\vskip 0.5em

\begin{proof}[\textbf{\emph{Proof Sketch.}}] We prove Proposition~\ref{Proposition: Minimum Separation for Two-Sample Multinomial Testing} indirectly by finding the minimum separation for one-sample multinomial testing. The goal of the one-sample problem is to test whether one set of samples is drawn from a known multinomial distribution. Intuitively the one-sample problem is no harder than the two-sample problem as the former can always be transformed into the latter by drawing another set of samples from the known distribution. This intuition was formalized by \cite{arias2018remember} in which they showed that the minimax risk of the one-sample problem is no larger than that of the two-sample problem (see their Lemma 1). We prove in Appendix~\ref{Section: Proposition: Minimum Separation for Two-Sample Multinomial Testing} that the minimum separation for the one-sample problem is of order $b_{(1)}^{1/4} n_1^{-1/2}$ and it thus follows that $b_{(1)}^{1/4} n_1^{-1/2} \lesssim \epsilon_{n_1,n_2}^{\dagger}$. The proof is completed by comparing this lower bound with the upper bound established in Proposition~\ref{Proposition: Multinomial Two-Sample Testing}. 
\end{proof}

\subsection{Two-sample testing for H\"{o}lder densities} \label{Section: Two-sample testing for Holder densities}
We next focus on testing for the equality between two density functions under H\"{o}lder's regularity condition. Adopting the notation used in \cite{arias2018remember}, let $\mathcal{H}^d_s(L)$ be the class of functions $f: [0,1]^d \mapsto \mathbb{R}$ such that
\begin{enumerate}
	\item $\big|f^{(\floor{s})}(x) - f^{(\floor{s})}(x') \big| \leq L \|x - x' \|^{s - \floor{s}}, \quad \forall x,x'\in [0,1]^d,$
	\item $|\!|\!|f^{(s')}|\!|\!|_{\infty} \leq L$ for each $s' \in \{1,\ldots,\floor{s}\}$,
\end{enumerate}
where $f^{(\floor{s})}$ denotes the $\floor{s}$-order derivative of $f$. Let us write the $L_2$ norm of $f \in \mathcal{H}^d_s(L)$ by $|\!|\!| f|\!|\!|_{L_2}^2:=\int f^2(x)dx$. By letting $f_Y$ and $f_Z$ be the density functions of $P_Y$ and $P_Z$ with respect to Lebesgue measure, we define the set of $(P_Y,P_Z)$, denoted by $\mathcal{P}_{\text{{H\"{o}lder}}}^{(d,s)}$, such that both $f_Y$ and $f_Z$ belong to $\mathcal{H}^d_s(L)$.  For this H\"{o}lder density class $\mathcal{P}_{\text{{H\"{o}lder}}}^{(d,s)}$ and $n_1 \leq n_2$, \cite{arias2018remember} establish that for testing $H_0: f_Y = f_Z$ against $H_1: |\!|\!|f_Y - f_Z|\!|\!|_{L_2} \geq \epsilon_{n_1,n_2}$, the minimum separation rate satisfies
\begin{align} \label{Eq: the minimum separation rate for two sample densities}
\epsilon_{n_1,n_2}^\dagger \asymp n_1^{-2s/(4s+d)}.
\end{align}
We note that this optimal testing rate is faster than the $n^{-s/(2s+d)}$ rate for estimating a H\"{o}lder density in the $L_2$ loss \citep[see for instance][]{tsybakov2009introduction}. It is further shown in \cite{arias2018remember} that the optimal rate (\ref{Eq: the minimum separation rate for two sample densities}) is achieved by the unnormalized chi-square test but with a somewhat loose threshold. Although they recommend a critical value calibrated by permutation in practice, it is unknown whether the resulting test has the same theoretical guarantees. We also note that their testing procedure discards $n_2-n_1$ observations to balance the sample sizes, which may lead to a less powerful test in practice. Motivated by these limitations, we propose an alternative test for H\"{o}lder densities, building on the multinomial permutation test in Proposition~\ref{Proposition: Multinomial Two-Sample Testing}. 

To implement the multinomial test for continuous data, we first need to discretize the support $[0,1]^d$. We follow the same strategy as in, for instance,~\cite{ingster1987minimax,arias2018remember,balakrishnan2019hypothesis} and consider bins of equal sizes that partition $[0,1]^d$. These bins are $d$-dimensional hypercubes whose length is set to $\kappa_{(1)}^{-1}$ where 
\begin{align*}
	\kappa_{(1)} := \floor{n_1^{2/(4s+d)}}.
\end{align*}	
Let $\{B_1,\ldots,B_{K}\}$ be an enumeration of such hypercubes. Next, consider a discretization function $Q: [0,1]^\dimension \mapsto \{1,\ldots,K\}$ such that $Q(x) = k$ if and only if $x \in B_k$. Then, based on the discretized data $\{Q(Y_1),\ldots,Q(Y_{n_1})\}$ and $\{Q(Z_1),\ldots,Q(Z_{n_2})\}$, we can now apply the multinomial test in Proposition~\ref{Proposition: Multinomial Two-Sample Testing} and the resulting test has the following theoretical guarantees for density testing.

\begin{proposition}[Two-sample testing for H\"{o}lder densities] \label{Proposition: Two-Sample Testing for Holder Densities}
	 Consider the multinomial test considered in Proposition~\ref{Proposition: Multinomial Two-Sample Testing} based on the equal-sized binned data described above. For a sufficiently large $C(s,d,L) >0$, consider $\epsilon_{n_1,n_2}$ such that
	\begin{align*}
	\epsilon_{n_1,n_2} \geq \frac{C(s,d,L)}{\alpha^{1/4}\beta^{1/2}} \left( \frac{1}{n_1} + \frac{1}{n_2}\right)^{\frac{2s}{4s+d}}. 
	\end{align*}
	Then for testing $\mathcal{P}_0 = \{(P_Y,P_Z) \in \mathcal{P}_{\text{\emph{H\"{o}lder}}}^{(d,s)} : f_Y=f_Z\}$ against $\mathcal{P}_1 =\{(P_Y,P_Z) \in \mathcal{P}_{\text{\emph{H\"{o}lder}}}^{(d,s)} : |\!|\!|f_Y - f_Z|\!|\!|_{L_2} \geq \epsilon_{n_1,n_2}\}$, the type I and II error rates of the resulting permutation test are uniformly controlled as in~(\ref{Eq: uniform error control}).
\end{proposition}

The proof of this result uses Proposition~\ref{Proposition: Multinomial Two-Sample Testing} along with careful analysis of the approximation errors from discretization, building on the analysis of \cite{ingster2000adaptive,arias2018remember}. The details can be found in Appendix~\ref{Section: Proposition: Two-Sample Testing for Holder Densities}. We remark that type I error control of the multinomial test follows clearly by the permutation principle, which is not affected by discretization. From the minimum separation rate given in (\ref{Eq: the minimum separation rate for two sample densities}), it is clear that the proposed test is minimax rate optimal for two-sample testing within H\"{o}lder class and it works for both equal and unequal sample sizes without discarding the data. However it is also important to note that the proposed test as well as the test introduced by \cite{arias2018remember} depend on knowledge of the smoothness parameter $s$, which is perhaps unrealistic in practice. To address this issue, \cite{arias2018remember} build upon the work of \cite{ingster2000adaptive} and propose a Bonferroni-type testing procedure that adapts to this unknown parameter at the cost of a $\log n$ factor. In Section~\ref{Section: Adaptive tests to unknown parameters}, we improve this logarithmic cost to an iterated logarithmic factor, leveraging combinatorial concentration inequalities developed in Section~\ref{Section: Combinatorial concentration inequalities}.

\section{The two moments method for independence testing} \label{Section: The two moments method for independence testing}
In this section we present analogous results to those in Section~\ref{Section: The two moments method for two-sample testing} for independence testing. We start by introducing a $U$-statistic for independence testing and establish a general condition under which the permutation test based on the $U$-statistic controls the type I and II error rates (Theorem~\ref{Theorem: U-statistic for independence testing}). We then move on to more specific cases of testing for multinomials and H\"{o}lder densities in Section~\ref{Section: Independence testing for multinomials} and Section~\ref{Section: Independence testing for Holder densities}, respectively.

Let us consider two bivariate functions $g_Y(y_1,y_2)$ and $g_Z(z_1,z_2)$, which are symmetric in their arguments. Define a product kernel associated with $g_Y(y_1,y_2)$ and $g_Z(z_1,z_2)$ by
\begin{equation}
\begin{aligned} \label{Eq: Independence kernel}
& h_{\text{in}}\{(y_1,z_1),(y_2,z_2),(y_3,z_3),(y_4,z_4)\} :=  \big\{g_Y(y_1,y_2) + g_Y(y_3,y_4) \\[.5em] 
- & g_Y(y_1,y_3) - g_Y(y_2,y_4) \big\} 
\cdot  \big\{g_Z(z_1,z_2) + g_Z(z_3,z_4) - g_Z(z_1,z_3) - g_Z(z_2,z_4) \big\}.
\end{aligned}
\end{equation}
For simplicity, we may also write $h_{\text{in}}\{(y_1,z_1),(y_2,z_2),(y_3,z_3),(y_4,z_4)\}$ as $h_{\text{in}}(x_1,x_2,x_3,x_4)$. Given this fourth order kernel, consider a $U$-statistic defined by
\begin{align} \label{Eq: U-statistic for independence testing}
U_n := \frac{1}{n_{(4)}} \sum_{(i_1,i_2,i_3,i_4) \in \mathbf{i}_4^n} h_{\text{in}}(X_{i_1},X_{i_2},X_{i_3},X_{i_4}). 
\end{align}
Again, by the unbiasedness property of $U$-statistics~\citep[e.g.][]{lee1990u}, it is clear that $U_n$ is an unbiased estimator of $\mE_P[h_{\text{in}}(X_{1},X_{2},X_{3},X_{4})]$. Depending on the choice of kernel $h_{\text{in}}$, the considered $U$-statistic covers numerous test statistics for independence testing including the Hilbert--Schmidt Independence Criterion (HSIC) \citep{gretton2005measuring} and distance covariance \citep{szekely2007measuring}. Let $\overline{h}_{\text{in}}(x_1,x_2,x_3,x_4)$ be the symmetrized version of $h_{\text{in}}(x_1,x_2,x_3,x_4)$ given by
\begin{align*} 
\overline{h}_{\text{in}}(x_1,x_2,x_3,x_4) := \frac{1}{4!} \sum_{(i_1,i_2,i_3,i_4) \in \mathbf{i}_4^{4}} h_{\text{in}}(x_{i_1},x_{i_2},x_{i_3},x_{i_4}).
\end{align*}
In a similar fashion to $\psi_{Y,1}(P),\psi_{Z,1}(P)$ and $\psi_{YZ,2}(P)$, we define $\psi_1'(P)$ and $\psi_2'(P)$ by  
\begin{equation}
\begin{aligned} \label{Eq: definition of psi prime functions}
\psi_1'(P) & := \text{Var}_P  [\mE_P\{ \overline{h}_{\text{in}}(X_1,X_2,X_3,X_4)|X_1\}], \\[.5em]
\psi_2'(P) & := \max \big\{ \mathbb{E}_P[g_Y^2(Y_1,Y_2)g_Z^2(Z_1,Z_2)], \ \mathbb{E}_P[g_Y^2(Y_1,Y_2)g_Z^2(Z_1,Z_3)],\\[.5em]
&~~~~~~~~~~~~~ \mathbb{E}_P[g_Y^2(Y_1,Y_2)g_Z^2(Z_3,Z_4)] \big\}.
\end{aligned}
\end{equation}
The following theorem studies the type II error of the permutation test based on $U_n$.
\begin{theorem}[$U$-statistic for independence testing] \label{Theorem: U-statistic for independence testing}
	Suppose that there is a sufficiently large constant $C>0$ such that  
	\begin{align*}
	\mE_P[U_{n}] \geq C \sqrt{\max \Bigg\{ \frac{\psi_{1}'(P)}{\beta n}, \  \frac{\psi_{2}'(P)}{\alpha\beta n^2} \Bigg\}},
	\end{align*}
	for all $P \in \mathcal{P}_1$. Then the type II error of the permutation test over $\mathcal{P}_1$ is uniformly bounded by $\beta$, that is
	\begin{align*}
	\sup_{P \in \mathcal{P}_1} \mP_P^{(n)} (U_{n} \leq c_{1-\alpha,n} )\leq \beta. 
	\end{align*}
\end{theorem}

\vskip 1em

\begin{remark} \normalfont
	Analogous to Theorem~\ref{Theorem: Two-Sample U-statistic}, for degenerate U-statistics of the form~\eqref{Eq: U-statistic for independence testing}, 
	this result simplifies the application of Lemma~\ref{Lemma: Two Moments Method} by reducing the sufficient condition to only depend on the quantities in~\eqref{Eq: definition of psi prime functions}. Importantly, these quantities do not depend on the permutation distribution of	the test statistic.
\end{remark}

\vskip 1em

\begin{proof}[\textbf{\emph{Proof Sketch.}}] The proof of Theorem~\ref{Theorem: U-statistic for independence testing} proceeds similarly as the proof of Theorem~\ref{Theorem: Two-Sample U-statistic}. Here we present a brief overview of the proof, while the details can be found in Appendix~\ref{Section: Theorem: U-statistic for independence testing}. First of all, the permuted $U$-statistic $U_{n}^\pi$ is centered and it suffices to verify the simplified condition~(\ref{Eq: Sufficient Condition II}). To this end, based on the explicit variance formula of a $U$-statistic \citep[e.g.~page 12 of][]{lee1990u}, we prove that 
\begin{align} \label{Eq: Variance upper bound II}
\mV_P[U_n] \leq C_1 \frac{\psi_1'(P)}{n} + C_2 \frac{\psi_2'(P)}{n^2}.
\end{align}
Analogous to the case of the two-sample $U$-statistic, the variance of $U_n$ behaves differently under the null and alternative hypotheses. In particular, under the null hypothesis, $\psi_1'(P)$ becomes zero and thus the second term dominates the upper bound (\ref{Eq: Variance upper bound II}). Since the permuted $U$-statistic $U_n^\pi$ mimics the behavior of $U_n$ under the null, the variance of $U_n^{\pi}$ is expected to be similarly bounded. We make this statement precise by proving the following result:
\begin{align} \label{Eq: Variance upper bound under permutations II}
\mE_P[\mV_{\pi}\{U_n^{\pi} | \mathcal{X}_n\}] \leq C_3 \frac{\psi_2'(P)}{n^2}.
\end{align}
Again, this part of the proof heavily relies on the symmetric structure of kernel $h_{\text{in}}$ and the details are deferred to Appendix~\ref{Section: Theorem: U-statistic for independence testing}. Now by combining the established bounds (\ref{Eq: Variance upper bound II}) and (\ref{Eq: Variance upper bound under permutations II}) together with the sufficient condition~(\ref{Eq: Sufficient Condition II}), we can conclude Theorem~\ref{Theorem: U-statistic for independence testing}. 
\end{proof}

In the following subsections, we illustrate Theorem~\ref{Theorem: U-statistic for independence testing} in the context of testing multinomial distributions and H\"{o}lder densities.

\subsection{Independence testing for multinomials} \label{Section: Independence testing for multinomials}
We begin with the case of multinomial distributions. Let $p_{YZ}$ denote a multinomial distribution on a product domain $\mathbb{S}_{d_1,d_2} := \{1,\ldots, d_1\} \times \{1,\ldots,d_2\}$ and $p_Y$ and $p_Z$ be its marginal distributions. Let us recall the kernel $h_{\text{in}}(x_1,x_2,x_3,x_4)$ in (\ref{Eq: Independence kernel}) and define it with the following bivariate functions:
\begin{equation}
\begin{aligned} \label{Eq: Independence bivariate functions}
& g_{\text{Multi},Y}(y_1,y_2) := \sum_{k=1}^{d_1} \ind(y_1=k)\ind(y_2=k) \quad \text{and} \\
& g_{\text{Multi},Z}(z_1,z_2) := \sum_{k=1}^{d_2} \ind(z_1=k)\ind(z_2=k).
\end{aligned}
\end{equation}
In this case, the expectation of the $U$-statistic is $4\|p_{YZ} - p_Yp_Z\|_2^2$. Analogous to the term $b_{(1)}$ in the two-sample case, let us define
\begin{align} \label{Eq: definition of b_2}
b_{(2)} := \max \big\{\|p_{YZ}\|_2^2, \|p_Yp_Z\|_2^2 \big\}.
\end{align}
Building on Theorem~\ref{Theorem: U-statistic for independence testing}, the next result establishes a guarantee on the testing errors of the permutation test for multinomial independence testing.

\begin{proposition}[Multinomial independence testing in $\ell_2$ distance] \label{Proposition: Multinomial Independence Testing}
	Let $\mathcal{P}_{\text{\emph{Multi}}}^{(d_1,d_2)}$ be the set of multinomial distributions defined on $\mathbb{S}_{d_1, d_2}$. Let $\mathcal{P}_0 = \{p_{YZ} \in \mathcal{P}_{\text{\emph{Multi}}}^{(d_1,d_2)}: p_{YZ} = p_Yp_Z \}$ and $\mathcal{P}_1(\epsilon_n) = \{p_{YZ} \in \mathcal{P}_{\text{\emph{Multi}}}^{(d_1,d_2)}: \|p_{YZ} - p_Yp_Z\|_2 \geq \epsilon_n \}$ where
	\begin{align*}
	\epsilon_{n} \geq \frac{C}{\alpha^{1/4}\beta^{1/2}} \frac{b_{(2)}^{1/4}}{n^{1/2}},
	\end{align*}
	for a sufficiently large $C>0$. Consider the $U$-statistic $U_{n}$ in (\ref{Eq: U-statistic for independence testing}) defined with the bivariate functions $g_{\text{\emph{Multi}},Y}$ and $g_{\text{\emph{Multi}},Z}$ given in (\ref{Eq: Independence bivariate functions}). Then, over the classes $\mathcal{P}_0$ and $\mathcal{P}_1$, the type I and II errors of the resulting permutation test are uniformly bounded as in (\ref{Eq: uniform error control}). 
\end{proposition}

\vskip 1em

\begin{proof}[\textbf{\emph{Proof Sketch.}}] We outline the proof of the result, while the details can be found in Appendix~\ref{Section: Proposition: Multinomial Independence Testing}. In the proof, we prove that there exist constants $C_1,C_2 > 0$ such that 
\begin{equation}
\begin{aligned} \label{Eq: Aim of the proof of multinomial independence testing}
& \psi_{1}'(P) \leq C_1 \sqrt{b_{(2)}}\|p_{YZ} - p_Yp_Z\|_2^2, \\[.5em]
& \psi_{2}'(P) \leq C_2 b_{(2)}.
\end{aligned}
\end{equation}
These bounds combined with Theorem~\ref{Theorem: U-statistic for independence testing} yields that if there exists a sufficiently large $C_3>0$ such that 
\begin{align} \label{Eq: Independence Bound}
\|p_{YZ} - p_Yp_Z\|_2^2 \geq \frac{C_3}{\alpha^{1/2}\beta} \frac{\sqrt{b_{(2)}}}{n},
\end{align}
then the type II error of the permutation test can be controlled by $\beta$ as desired. 
\end{proof}

The next proposition asserts that the minimum separation rate for independence testing in the $\ell_2$ distance is $\epsilon_{n}^{\dagger} \asymp b_{(2)}^{1/4} n^{-1/2}$. This implies that the permutation test based on $U_n$ in Proposition~\ref{Proposition: Multinomial Independence Testing} is minimax rate optimal in this scenario.

\begin{proposition}[Minimum separation for multinomial independence testing] \label{Proposition: Minimum Separation for Multinomial Independence Testing}
	Consider the independence testing problem within the class of multinomial distributions $\mathcal{P}_{\text{\emph{Multi}}}^{(d_1,d_2)}$ where the null hypothesis and the alternative hypothesis are $H_0: p_{YZ} = p_Yp_Z$ and $H_1: \|p_{YZ} - p_Yp_Z\|_2 \geq \epsilon_{n}$. Under this setting, the minimum separation satisfies $\epsilon_{n}^{\dagger} \asymp b_{(2)}^{1/4} n^{-1/2}$.
\end{proposition}

The proof of Proposition~\ref{Proposition: Minimum Separation for Multinomial Independence Testing} is based on the standard lower bound technique of \cite{ingster1987minimax} using a uniform mixture of alternative distributions. However, we remark that care is needed in order to ensure that alternative distributions are proper (normalized) multinomial distributions. To this end, we carefully perturb the uniform null distribution to generate a mixture of dependent alternative distributions, and use the property of negative association to deal with the dependency induced in ensuring the resulting distributions are normalized. The details can be found in Appendix~\ref{Section: Proposition: Minimum Separation for Multinomial Independence Testing}.
In the next subsection, we turn our attention to the class of H\"{o}lder densities and provide similar results of Section~\ref{Section: Two-sample testing for Holder densities} for independence testing.

\subsection{Independence testing for H\"{o}lder densities} \label{Section: Independence testing for Holder densities}
Turning to the case of H\"{o}lder densities, we leverage the previous multinomial result and establish the minimax rate for independence testing under the H\"{o}lder's regularity condition. As in Section~\ref{Section: Two-sample testing for Holder densities}, we restrict our attention to functions $f : [0,1]^{d_1+d_2} \mapsto \mathbb{R}$ that satisfy 
\begin{enumerate}
	\item $\big|f^{(\floor{s})}(x) - f^{(\floor{s})}(x') \big| \leq L \|x - x' \|^{s - \floor{s}}, \quad \forall x,x'\in [0,1]^{d_1 + d_2},$
	\item $|\!|\!|f^{(s')} |\!|\!|_{\infty} \leq L$ for each $s' \in \{1,\ldots,\floor{s}\}$.
\end{enumerate}
Let us write $\mathcal{H}_s^{d_1+d_2}(L)$ to denote the class of such functions. We further introduce the class of joint distributions, denoted by $\mathcal{P}_{\text{H\"{o}lder}}^{(d_1+d_2,s)}$, defined as follows. Let $f_{YZ}$ and $f_{Y}f_{Z}$ be the densities of $P_{YZ}$ and $P_YP_Z$ with respect to Lebesgue measure. Then $\mathcal{P}_{\text{{H\"{o}lder}}}^{(d_1+d_2,s)}$ is defined as the set of joint distributions $P_{YZ}$ such that both the joint density and the product density, $f_{YZ}$ and $f_Yf_Z$, belong to $\mathcal{H}^{d_1+d_2}_s(L)$. Consider partitions of $[0,1]^{d_1+d_2}$ into bins of equal size and set the bin size to be $\kappa_{(2)}^{-1}$ where $\kappa_{(2)} = \floor{n^{2/(4s+d_1+d_2)}}$. Based on these equal-sized partitions, one may apply the multinomial test for independence provided in Proposition~\ref{Proposition: Multinomial Independence Testing}. Despite discretization, the resulting test has valid level $\alpha$ due to the permutation principle and has the following theoretical guarantees for density testing over $\mathcal{P}_{\text{{H\"{o}lder}}}^{(d_1+d_2,s)}$. 
\begin{proposition}[Independence testing for H\"{o}lder densities] \label{Proposition: Independence Testing for Holder Densities}
	Consider the multinomial independence test considered in Proposition~\ref{Proposition: Multinomial Independence Testing} based on the binned data described above. For a sufficiently large $C(s,d_1,d_2,L) >0$, consider $\epsilon_{n}$ defined by 
	\begin{align*}
	\epsilon_{n} \geq \frac{C(s,d_1,d_2,L)}{\alpha^{1/4}\beta^{1/2}} \left(\frac{1}{n}\right)^{\frac{2s}{4s+d_1+d_2}}. 
	\end{align*}
	Then for testing $\mathcal{P}_0 = \{P_{YZ} \in \mathcal{P}_{\text{\emph{H\"{o}lder}}}^{(d_1+d_2,s)} : f_{YZ} = f_Yf_Z \}$ against $\mathcal{P}_1 =\{P_{YZ} \in \mathcal{P}_{\text{\emph{H\"{o}lder}}}^{(d_1+d_2,s)} : |\!|\!|f_{YZ} - f_Yf_Z|\!|\!|_{L_2} \geq \epsilon_{n}\}$, the type I and II errors of the resulting permutation test are uniformly controlled as in (\ref{Eq: uniform error control}).	
\end{proposition}

The proof of the above result follows similarly to the proof of Proposition~\ref{Proposition: Two-Sample Testing for Holder Densities} and can be found in Appendix~\ref{Section: Proposition: Independence Testing for Holder Densities}. Indeed, as shown in the next proposition, the proposed binning-based independence test is minimax rate optimal for the H\"{o}lder class density functions. That is, no test can have uniform power when the separation rate $\epsilon_n$ is of order smaller than $n^{-2s/(4s+d_1+d_2)}$. 
\begin{proposition}[Minimum separation for independence testing in H\"{o}lder class] \label{Proposition: Minimum Separation for Independence Testing for Holder Densities}
	Consider the independence testing problem within the class $\mathcal{P}_{\text{\emph{H\"{o}lder}}}^{(d_1+d_2,s)}$ in which the null hypothesis and the alternative hypothesis are $H_0: f_{YZ} = f_Yf_Z$ and $H_1: |\!|\!|f_{YZ} - f_Yf_Z|\!|\!|_{L_2} \geq \epsilon_{n}$. Under this setting, the minimum separation satisfies $\epsilon_{n}^{\dagger} \asymp n^{-2s/(4s+d_1+d_2)}$.
\end{proposition}

The proof of Proposition~\ref{Proposition: Minimum Separation for Independence Testing for Holder Densities} is again based on the standard lower bound technique by \cite{ingster1987minimax} and deferred to Appendix~\ref{Section: Proposition: Minimum Separation for Independence Testing for Holder Densities}. We note that the independence test in Proposition~\ref{Proposition: Independence Testing for Holder Densities} hinges on the assumption that the smoothness parameter $s$ is known. To avoid this assumption, we introduce an adaptive test to this smoothness parameter at the cost of $\log\log n$ factor in Section~\ref{Section: Adaptive tests to unknown parameters}. A building block for this adaptive result is combinatorial concentration inequalities developed in the next section.

\section{Combinatorial concentration inequalities} \label{Section: Combinatorial concentration inequalities}
Although the two moments method is broadly applicable, it may not yield sharp results when an extremely small significance level $\alpha$ is of interest (say, $\alpha$ shrinks to zero as $n$ increases). In particular, the sufficient condition (\ref{Eq: Sufficient Condition}) given by the two moments method has a polynomial dependency on $\alpha$. In this section, we develop exponential concentration inequalities for permuted $U$-statistics that allow us to improve this polynomial dependency. To this end, we introduce a novel strategy to couple a permuted $U$-statistic with i.i.d.~Bernoulli or Rademacher random variables, inspired by the symmetrization trick \citep{duembgen1998symmetrization} and Hoeffding's average \citep{hoeffding1963probability}. 

\vskip 1em

\noindent \textbf{Coupling with i.i.d.~random variables.} The core idea of our approach is fairly general and based on the following simple observation. Given a random permutation $\pi$ uniformly distributed over $\mathbf{\Pi}_n$, we randomly switch the order within $(\pi_{2i-1}, \pi_{2i})$ for $i=1,2,\ldots,\floor{n/2}$. We denote the resulting permutation by $\pi'$. It is clear that $\pi$ and $\pi'$ are dependent but identically distributed. The point of introducing this extra permutation $\pi'$ is that we are now able to associate $\pi'$ with i.i.d.~Bernoulli random variables without changing the distribution. To be more specific, let $\delta_1,\ldots,\delta_{\floor{n/2}}$ be i.i.d.~Bernoulli random variables with success probability $1/2$. Then $(\pi_{2i-1}',\pi_{2i}')$ can be written as
\begin{align*}
(\pi_{2i-1}',\pi_{2i}') = \big(\delta_i \pi_{2i-1} + (1-\delta_i) \pi_{2i}, (1-\delta_i) \pi_{2i-1} + \delta_i \pi_{2i}\big) \quad \text{for $i=1,2,\ldots,\floor{n/2}$}.
\end{align*}
Given that it is easier to work with i.i.d.~samples than permutations, the alternative representation of $\pi'$ gives a nice way to investigate a general permuted statistic. The next subsections provide concrete demonstrations of this coupling approach based on degenerate $U$-statistics.

\subsection{Degenerate two-sample $U$-statistics} \label{Section: Degenerate two-sample U-statistics}
We start with the two-sample $U$-statistic in (\ref{Eq: Two-Sample U-statistic}). Our strategy is outlined as follows. First, motivated by Hoeffding's average~\citep{hoeffding1963probability}, we express the permuted $U$-statistic as the average of more tractable statistics. We then link these tractable statistics to quadratic forms of i.i.d.~Rademacher random variables based on the coupling idea described before. Finally we apply existing concentration bounds for quadratic forms of i.i.d.~random variables to obtain the result in Theorem~\ref{Theorem: Two-Sample Concentration}.

Let us denote the permuted $U$-statistic associated with $\pi \in \mathbf{\Pi}_n$ by
\begin{align} \label{Eq: permuted two-sample U-statistic}
U_{n_1,n_2}^{\pi} := \frac{1}{(n_1)_{(2)}(n_2)_{(2)}} \sum_{(i_1,i_2) \in \mathbf{i}_2^{n_1}} \sum_{(j_1,j_2) \in \mathbf{i}_2^{n_2}} h_{\text{ts}}(X_{\pi_{i_1}},X_{\pi_{i_2}}; X_{\pi_{n_1+ j_1}},X_{\pi_{n_1+j_2}}).
\end{align}
By assuming $n_1 \leq n_2$, let $L:= \{\ell_1,\ldots, \ell_{n_1}\}$ be a $n_1$-tuple uniformly drawn without replacement from $\{1,\ldots,n_2\}$. Given $L$, we introduce another test statistic
\begin{align*}
\widetilde{U}_{n_1,n_2}^{\pi,L} := \frac{1}{(n_1)_{(2)}} \sum_{(k_1,k_2) \in \mathbf{i}_2^{n_1}} h_{\text{ts}}(X_{\pi_{k_1}},X_{\pi_{k_2}}; X_{\pi_{n_1+\ell_{k_1}}},X_{\pi_{n_1+\ell_{k_2}}}).
\end{align*}
By treating $L$ as a random quantity, $U_{n_1,n_2}^{\pi}$ can be viewed as the expected value of $\widetilde{U}_{n_1,n_2}^{\pi,L}$ with respect to $L$ (conditional on other random variables), that is,
\begin{align} \label{Eq: Equal in expectation}
U_{n_1,n_2}^{\pi} = \mE_L [\widetilde{U}_{n_1,n_2}^{\pi,L} | \mathcal{X}_n, \pi]. 
\end{align}
The idea of expressing a $U$-statistic as the average of more tractable statistics dates back to \cite{hoeffding1963probability}. The reason for introducing $\widetilde{U}_{n_1,n_2}^{\pi,L}$ is to connect $U_{n_1,n_2}^{\pi}$ with a Rademacher chaos. Recall that $\pi =(\pi_1,\ldots,\pi_n)$ is uniformly distributed over all possible permutations of $\{1,\ldots,n\}$. Therefore, as explained earlier, the distribution of $\widetilde{U}_{n_1,n_2}^{\pi,L}$ does not change even if we randomly switch the order between $X_{\pi_{k}}$ and $X_{\pi_{n_1+\ell_{k}}}$ for $k \in \{1,\ldots, n_1\}$. More formally, recall that $\delta_1,\ldots,\delta_{n_1}$ are i.i.d.~Bernoulli random variables with success probability $1/2$. For $k=1,\ldots,n_1$, define
\begin{align}
\widetilde{X}_{\pi_{k}} := \delta_k X_{\pi_{k}} + (1-\delta_k) X_{\pi_{n_1 + \ell_k}} \quad \text{and} \quad \widetilde{X}_{\pi_{n_1 + \ell_k}}  := (1-\delta_k) X_{\pi_{k}}  + \delta_k X_{\pi_{n_1 + \ell_k}}.
\end{align}
Then it can be seen that $\widetilde{U}_{n_1,n_2}^{\pi,L}$ is equal in distribution to
\begin{align*}
\widetilde{U}_{n_1,n_2}^{\pi,L,\delta} := \frac{1}{(n_1)_{(2)}} \sum_{(k_1,k_2) \in \mathbf{i}_2^{n_1}} h_{\text{ts}}(\widetilde{X}_{\pi_{k_1}},\widetilde{X}_{\pi_{k_2}}; \widetilde{X}_{\pi_{n_1+\ell_{k_1}}},\widetilde{X}_{\pi_{n_1+\ell_{k_2}}}).
\end{align*}
In other words, we link $\widetilde{U}_{n_1,n_2}^{\pi,L}$ to i.i.d.~Bernoulli random variables, which are easier to work with. Furthermore, by the symmetry of $g(x,y)$ in its arguments and letting $\zeta_1,\ldots,\zeta_{n}$ be i.i.d.~Rademacher random variables, one can observe that $\widetilde{U}_{n_1,n_2}^{\pi,L,\delta}$ is equal in distribution to the following Rademacher chaos:
\begin{align*}
\widetilde{U}_{n_1,n_2}^{\pi,L,\zeta} := \frac{1}{(n_1)_{(2)}} \sum_{(k_1,k_2) \in \mathbf{i}_2^{n_1}} \zeta_{k_1} \zeta_{k_2} h_{\text{ts}}(X_{\pi_{k_1}},X_{\pi_{k_2}}; X_{\pi_{n_1+\ell_{k_1}}},X_{\pi_{n_1+\ell_{k_2}}}).
\end{align*}
Consequently, we observe that $\widetilde{U}_{n_1,n_2}^{\pi,L}$ and $\widetilde{U}_{n_1,n_2}^{\pi,L,\zeta}$ are equal in distribution, i.e.
\begin{align} \label{Eq: Equal in distribution}
\widetilde{U}_{n_1,n_2}^{\pi,L}  \overset{d}{=} \widetilde{U}_{n_1,n_2}^{\pi,L,\zeta}.
\end{align}
We now have all the ingredients ready for obtaining an exponential bound for $U_{n_1,n_2}^{\pi}$. By the Chernoff bound~\cite[e.g.][]{boucheron2013concentration}, for any $\lambda >0$, 
\begin{equation}
\begin{aligned} \label{Eq: Chernoff Bound}
\mP_{\pi} \left( U_{n_1,n_2}^{\pi} > t | \mathcal{X}_n \right) ~ \leq ~ &  e^{-\lambda t} \mathbb{E}_{\pi} \big[ \exp \big( \lambda U_{n_1,n_2}^{\pi} \big) |  \mathcal{X}_n \big] \\[.5em]
\overset{(i)}{\leq} ~ &  e^{-\lambda t} \mathbb{E}_{\pi,L} \big[ \exp \big( \lambda \widetilde{U}_{n_1,n_2}^{\pi,L} \big) |  \mathcal{X}_n \big] \\[.5em]
\overset{(ii)}{=} ~ &  e^{-\lambda t} \mathbb{E}_{\pi,L,\zeta} \big[ \exp \big( \lambda \widetilde{U}_{n_1,n_2}^{\pi,L,\zeta} \big) |  \mathcal{X}_n \big]
\end{aligned}
\end{equation}
where step~$(i)$ uses Jensen's inequality together with (\ref{Eq: Equal in expectation}) and step~$(ii)$ holds from (\ref{Eq: Equal in distribution}). Finally, conditional on $\pi$ and $L$, we can associate the last equation with the moment generating function of a quadratic form of i.i.d.~Rademacher random variables. This quadratic form has been well-studied in the literature through a decoupling argument \citep[e.g. Chapter 6 of][]{vershynin2018high}, which leads to the following theorem. The remainder of the proof of Theorem~\ref{Theorem: Two-Sample Concentration} can be found in Section~\ref{Section: Theorem: Two-Sample Concentration}. 

\begin{theorem}[Concentration of $U_{n_1,n_2}^{\pi}$] \label{Theorem: Two-Sample Concentration}
	Consider the permuted two-sample $U$-statistic $U_{n_1,n_2}^{\pi}$ (\ref{Eq: permuted two-sample U-statistic}) and define 
	\begin{align*}
	\Sigma_{n_1,n_2}^2:= \frac{1}{n_1^2(n_1-1)^2} \sup_{\pi \in \mathbf{\Pi}_n} \Bigg\{ \sum_{(i_1,i_2) \in \mathbf{i}_2^{n_1}} g^2(X_{\pi_{i_1}},X_{\pi_{i_2}}) \Bigg\}.
	\end{align*}
	Then, for every $t>0$ and some constant $C>0$, we have 
	\begin{align*}
	\mP_{\pi} \left( U_{n_1,n_2}^{\pi} \geq t ~| \mathcal{X}_n \right) \leq \exp \bigg\{ - C \min \left( \frac{t^2}{\Sigma_{n_1,n_2}^2}, \ \frac{t}{\Sigma_{n_1,n_2}} \right)  \bigg\}.
	\end{align*}
\end{theorem}
In our application, it is convenient to have an upper bound for $\Sigma_{n_1,n_2}$ without involving the supremum operator. One trivial bound, suitable for our purpose, is given by
\begin{align} \label{Eq: trivial bound}
\Sigma_{n_1,n_2}^2 \leq  \frac{1}{n_1^2(n_1-1)^2} \sum_{(i_1,i_2) \in \mathbf{i}_2^{n}} g^2(X_{i_1},X_{i_2}).
\end{align}
The next subsection presents an analogous result for degenerate $U$-statistics in the context of independence testing.

\subsection{Degenerate $U$-statistics for independence testing} \label{Section: Degenerate U-statistics for independence testing}
Let us recall the $U$-statistic for independence testing in (\ref{Eq: U-statistic for independence testing}) and denote the permuted version by
\begin{align} \label{Eq: permuted U-statistic}
U_n^{\pi} := \frac{1}{n_{(4)}} \sum_{(i_1,i_2,i_3,i_4) \in \mathbf{i}_4^n} h_{\text{in}} \big\{(Y_{i_1},Z_{\pi_{i_1}}),(Y_{i_2},Z_{\pi_{i_2}}),(Y_{i_3},Z_{\pi_{i_3}}),(Y_{i_4},Z_{\pi_{i_4}}) \big\}. 
\end{align}
We follow a similar strategy taken in the previous subsection to obtain an exponential bound for $U_n^{\pi}$. To this end, we first introduce some notation. Let $L:=\{\ell_1,\ldots, \ell_{\floor{n/2}} \}$ be a $\floor{n/2}$-tuple uniformly sampled without replacement from $\{1,\ldots,n\}$ and similarly $L':= \{\ell_1',\ldots,\ell_{\floor{n/2}}' \}$ be another $\floor{n/2}$-tuple uniformly sampled without replacement from $\{1,\ldots,n\} \setminus L$. By construction, $L$ and $L'$ are disjoint. Given $L$ and $L'$, we define another test statistic $\widetilde{U}_{n}^{\pi,L,L'}$ as
\begin{align*}
\widetilde{U}_{n}^{\pi,L,L'} := \frac{1}{\floor{n/2}_{(2)}} \sum_{(i_1,i_2) \in \mathbf{i}_2^{\floor{n/2}}} h_{\text{in}}  \big\{\big(Y_{\ell_{i_1}},Z_{\pi_{\ell_{i_1}}}\big),\big(Y_{\ell_{i_2}},Z_{\pi_{\ell_{i_2}}}\big),\big(Y_{\ell_{i_2}'},Z_{\pi_{\ell_{i_2}'}}\big),\big(Y_{\ell_{i_1}'},Z_{\pi_{\ell_{i_1}'}}\big)\big\}.
\end{align*}
By treating $L$ and $L'$ as random quantities, $U_n^{\pi}$ can be viewed as the expected value of $\widetilde{U}_{n}^{\pi,L,L'}$ with respect to $L$ and $L'$, i.e.
\begin{align} \label{Eq: Equal in Expectation 2}
U_n^{\pi} = \mE_{L,L'}[\widetilde{U}_{n}^{\pi,L,L'} | \mathcal{X}_n, \pi].
\end{align} 
From the same reasoning as before, the distribution of $\widetilde{U}_{n}^{\pi,L,L'}$ does not change even if we randomly switch the order between $Z_{\pi_{\ell_k}}$ and $Z_{\pi_{\ell_k'}}$ for $k=1,\ldots,\floor{n/2}$, which allows us to introduce i.i.d.~Bernoulli random variables with success probability $1/2$. By the symmetry of $g_Y(y_1,y_2)$ and $g_Z(z_1,z_2)$, we may further observe that $\widetilde{U}_{n}^{\pi,L,L'}$ is equal in distribution to
\begin{equation}
\begin{aligned} \label{Eq: U with L zeta}
\widetilde{U}_{n}^{\pi,L,L',\zeta} :=~ & \frac{1}{\floor{n/2}_{(2)}} \sum_{(i_1,i_2) \in \mathbf{i}_2^{\floor{n/2}}} \zeta_{i_1} \zeta_{i_2} \times \\
& h_{\text{in}}  \big\{\big(Y_{\ell_{i_1}},Z_{\pi_{\ell_{i_1}}}\big),\big(Y_{\ell_{i_2}},Z_{\pi_{\ell_{i_2}}}\big),\big(Y_{\ell_{i_2}'},Z_{\pi_{\ell_{i_2}'}}\big),\big(Y_{\ell_{i_1}'},Z_{\pi_{\ell_{i_1}'}}\big)\big\}.
\end{aligned} 
\end{equation}
Thus, based on the alternative expression of $U_n^{\pi}$ in (\ref{Eq: Equal in Expectation 2}) along with the relationship $\widetilde{U}_{n}^{\pi,L,L'} \overset{d}{=} \widetilde{U}_{n}^{\pi,L,L',\zeta}$, we can establish a similar exponential tail bound as in Theorem~\ref{Theorem: Two-Sample Concentration} for $U_n^{\pi}$ as follows.

\begin{theorem}[Concentration I of $U_{n}^{\pi}$] \label{Theorem: Concentration inequality for independence U-statistic}
	Consider the permuted $U$-statistic $U_{n}^{\pi}$ (\ref{Eq: permuted U-statistic}) and define 
	\begin{align} \label{Eq: definition of Sigma}
	\Sigma_n^2:= \frac{1}{n^2(n-1)^2} \sup_{\pi \in \mathbf{\Pi}_n} \Bigg\{ \sum_{(i_1,i_2) \in \mathbf{i}_2^{n}} g_Y^2(Y_{i_1},Y_{i_2})g_Z^2(Z_{\pi_{i_1}},Z_{\pi_{i_2}}) \Bigg\}.
	\end{align}
	Then, for every $t>0$ and some constant $C>0$, we have 
	\begin{align*}
	\mP_{\pi} \left( U_{n}^{\pi} \geq t ~| \mathcal{X}_n \right) \leq \exp \bigg\{ - C \min \left( \frac{t^2}{\Sigma_n^2}, \ \frac{t}{\Sigma_n} \right)  \bigg\}.
	\end{align*}
\end{theorem}
We omit the proof of the result as it follows exactly the same line of the proof of Theorem~\ref{Theorem: Two-Sample Concentration}. Similar to the upper bound~(\ref{Eq: trivial bound}), H\"{o}lder's inequality yields two convenient bounds for $\Sigma_n^2$ as
\begin{align} \label{Eq: A bound of Sigma}
& \Sigma_n^2 \leq \frac{1}{n^2(n-1)^2} |\!|\!|g_Z^2 |\!|\!|_{\infty} \sum_{(i_1,i_2) \in \mathbf{i}_2^{n}} g_Y^2(Y_{i_1},Y_{i_2}) \quad \text{and} \\[.5em] \nonumber
& \Sigma_n^2 \leq \frac{1}{n^2(n-1)^2}  \sqrt{\sum_{(i_1,i_2) \in \mathbf{i}_2^{n}} g_Y^4(Y_{i_1},Y_{i_2})} \sqrt{\sum_{(i_1,i_2) \in \mathbf{i}_2^{n}} g_Z^4(Z_{i_1},Z_{i_2})}.
\end{align}
At the end of this subsection, we provide an application of Theorem~\ref{Theorem: Concentration inequality for independence U-statistic} to a dependent Rademacher chaos. 

\vskip 1em

\noindent \textbf{A refined version}. Although Theorem~\ref{Theorem: Concentration inequality for independence U-statistic} presents a fairly strong exponential concentration of $U_n^\pi$, it may lead to a sub-optimal result for independence testing. Indeed, for the minimax result, we want to obtain a similar bound but by replacing the supremum with the average over $\pi \in \mathbf{\Pi}_n$ in (\ref{Eq: definition of Sigma}). To this end, we borrow decoupling ideas from \cite{duembgen1998symmetrization} and \cite{de1999decoupling} and present a refined concentration inequality in Theorem~\ref{Theorem: Concentration inequality for independence U-statistic II}. The proposed bound~(\ref{Eq: Concentration inequality for independence U-statistic}) can be viewed as Bernstein-type inequality in a sense that it contains the variance term $\Lambda_n$ (not depending on the supremum) and maximum term $M_n$ defined as
\begin{equation}
\begin{aligned} \label{Eq: definition of lambda and M} 
& \Lambda_n^2:= \frac{1}{n^4} \sum_{1 \leq i_1,i_2 \leq n} \sum_{1 \leq j_1,j_2 \leq n} g_Y^2(Y_{i_1},Y_{i_2})g_Z^2(Z_{j_1},Z_{j_2}) \quad \text{and} \\[.5em]
& M_n: = \max_{1 \leq i_1,i_2,j_1,j_2 \leq n} |g_Y(Y_{i_1},Y_{i_2}) g_Z(Z_{j_1},Z_{j_2})|.
\end{aligned}
\end{equation}
In particular, the revised inequality would be sharper than the one in Theorem~\ref{Theorem: Concentration inequality for independence U-statistic} especially when $\Lambda_n$ is much smaller than $n \Sigma_n$.
\begin{theorem}[Concentration II of $U_{n}^{\pi}$] \label{Theorem: Concentration inequality for independence U-statistic II}
	Consider the permuted $U$-statistic $U_{n}^{\pi}$ (\ref{Eq: permuted U-statistic}) and recall $\Lambda_n$ and $M_n$ from (\ref{Eq: definition of lambda and M}). Then, for every $t>0$ and some constant $C_1,C_2>0$, we have 
	\begin{align}  \label{Eq: Concentration inequality for independence U-statistic}
	\mP_{\pi} \left( U_{n}^{\pi} \geq t ~| \mathcal{X}_n \right) \leq C_1 \exp \Bigg\{ - C_2 \min \left( \frac{nt}{\Lambda_n}, \ \frac{nt^{2/3}}{M_n^{2/3}} \right)  \Bigg\}.
	\end{align}
\end{theorem}

\vskip 1em

\begin{proof}[\textbf{\emph{Proof Sketch.}}] Here we sketch the proof of the result while the details are deferred to Appendix~\ref{Section: Theorem: Concentration inequality for independence U-statistic II}. Let $\psi(\cdot)$ be a nondecreasing convex function on $[0,\infty)$ and $\Psi(x) = \psi(|x|)$. Based on the equality in (\ref{Eq: Equal in Expectation 2}), Jensen's inequality yields
\begin{align*}
\mE_{\pi} [ \Psi(\lambda U_n^{\pi}) | \mathcal{X}_n] \leq \mE_{\pi,L,L',\zeta} \big[\Psi \big(\lambda \widetilde{U}_n^{\pi,L,L',\zeta} \big) | \mathcal{X}_n\big], 
\end{align*}
where $\widetilde{U}_n^{\pi,L,L',\zeta}$ can be recalled from (\ref{Eq: U with L zeta}). Let $\pi'$ be i.i.d.~copy of permutation $\pi$. Then, by letting $m = \floor{n/2}$ and observing that $(i)$ $\{ \zeta_i \}_{i=1}^m \overset{d}{=} \{ \zeta_i \zeta_{i+m}\}_{i=1}^m$ and $(ii)$ $\{L,L'\} \overset{d}{=} \{\pi_1',\ldots,\pi_{2m}'\}$, we have
\begin{align*}
\widetilde{U}_{n}^{\pi,L,L',\zeta} \overset{d}{=} ~ & \widetilde{U}_{n}^{\pi,\pi',\zeta} := ~ \frac{1}{m_{(2)}} \sum_{(i_1,i_2) \in \mathbf{i}_2^{m}} \zeta_{i_1} \zeta_{i_2} \zeta_{i_1+m} \zeta_{i_2+m} \times \\
&~~~~~~~~~~~~~~~ h_{\text{in}} \{(Y_{\pi_{i_1}'},Z_{\pi_{i_1}}),(Y_{\pi_{i_2}'},Z_{\pi_{i_2}}),(Y_{\pi_{i_2+m}'},Z_{\pi_{\pi_{i_2+m}}}),(Y_{\pi_{i_1+m}'},Z_{\pi_{i_1+m}})\}.
\end{align*}
Next denote the decoupled version of $\pi$ by $\widetilde{\pi} := (\widetilde{\pi}_1,\ldots,\widetilde{\pi}_n)$ whose components are independent and identically distributed as $\pi_1$. Let $\widetilde{\pi}'$ be i.i.d.~copy of $\widetilde{\pi}$. Building on the decoupling idea of \cite{duembgen1998symmetrization}, our proof proceeds by replacing $\pi,\pi'$ in $\widetilde{U}_{n}^{\pi,\pi',\zeta}$ with $\widetilde{\pi},\widetilde{\pi}'$. If this decoupling step succeeds, then we can view the corresponding $U$-statistic as a second order degenerate $U$-statistic of i.i.d.~random variables (conditional on $\mathcal{X}_n$). We are then able to apply concentration inequalities for degenerate $U$-statistics in \cite{de1999decoupling} to finish the proof. 
\end{proof}

\vskip .5em

\noindent \textbf{Dependent Rademacher chaos.} To illustrate the efficacy of Theorem~\ref{Theorem: Concentration inequality for independence U-statistic}, let us consider a Rademacher chaos under sampling without replacement, which has been recently studied by \cite{hodara2019exponential}. To describe the problem, let $\widetilde{\zeta}_1,\ldots,\widetilde{\zeta}_n$ be dependent Rademacher random variables such that $\sum_{i=1}^n \widetilde{\zeta}_i = 0$ where $n$ is assumed to be even. For real numbers $\{a_{i,j}\}_{i,j=1}^n$, the Rademacher chaos under sampling without replacement is given by
\begin{align*}
T_{\text{Rad}} := \sum_{(i_1,i_2) \in \mathbf{i}_2^{n}} \widetilde{\zeta}_{i_1} \widetilde{\zeta}_{i_2} a_{i_1,i_2}.
\end{align*}
\cite{hodara2019exponential} present two exponential concentration inequalities for $T_{\text{Rad}}$ based on the coupling argument introduced by \cite{chung2013exact}. Intuitively, $T_{\text{Rad}}$ should behave like i.i.d.~Rademacher chaos, replacing $\{\widetilde{\zeta}_i\}_{i=1}^n$ with $\{\zeta_i\}_{i=1}^n$, at least in the large sample size. Both of their results, however, do not fully recover a well-known concentration bound for i.i.d.~Rademacher chaos \citep[e.g.~Corollary 3.2.6 of][]{de1999decoupling}; namely,
\begin{align} \label{Eq: i.i.d. Rademacher chaos}
\mP \Big\{ \Big| \sum\nolimits_{(i_1,i_2) \in \mathbf{i}_2^{n}} \zeta_{i_1} \zeta_{i_2} a_{i_1,i_2}  \Big| \geq t \Big\} \leq 2 \exp\left( -C t / A_n \right),
\end{align}
where $ A_n^2 := \sum\nolimits_{(i_1,i_2) \in \mathbf{i}_2^{n}} a_{i_1,i_2}^2$. In the next corollary, we leverage Theorem~\ref{Theorem: Concentration inequality for independence U-statistic} and present an alternative tail bound for $T_{\text{Rad}}$ that precisely captures the tail bound~(\ref{Eq: i.i.d. Rademacher chaos}) for large $t$. Note that, unlike i.i.d.~Rademacher chaos, $T_{\text{Rad}}$ has a non-zero expectation. Hence we construct a tail bound for the chaos statistic centered by $\overline{a} := n_{(2)}^{-1} \sum_{(i_1,i_2) \in \mathbf{i}_2^{n}} a_{i_1,i_2}$. The proof of the result can be found in Appendix~\ref{Section: Corollary: Dependent Rademacher Chaos}.

\begin{corollary}[Dependent Rademacher chaos] \label{Corollary: Dependent Rademacher Chaos}
	For every $t>0$ and some constant $C>0$, the dependent Rademacher chaos is concentrated as
	\begin{align*}
	\mP \Big\{ \Big| \sum\nolimits_{(i_1,i_2) \in \mathbf{i}_2^{n}} \widetilde{\zeta}_{i_1} \widetilde{\zeta}_{i_2}\left(a_{i_1,i_2}-\overline{a}\right) \Big| \geq t \Big\} \leq 2  \exp \bigg\{\!\! - C \min \left( \frac{t^2}{A_n^2}, \ \frac{t}{A_n} \right)  \bigg\}.
	\end{align*}
\end{corollary}

The next section studies adaptive tests based on the combinatorial concentration bounds provided in this section.

\section{Adaptive tests} \label{Section: Adaptive tests to unknown parameters}
In this section, we revisit two-sample testing and independence testing for H\"{o}lder densities considered in Section~\ref{Section: Two-sample testing for Holder densities} and Section~\ref{Proposition: Independence Testing for Holder Densities}, respectively. As mentioned earlier, minimax optimality of the multinomial tests for H\"{o}lder densities depends on an unknown smoothness parameter (see Proposition~\ref{Proposition: Two-Sample Testing for Holder Densities} and Proposition~\ref{Proposition: Minimum Separation for Independence Testing for Holder Densities}). The aim of this section is to introduce adaptive permutation tests to this unknown parameter at the expense of an iterated logarithm factor. To this end, we generally follow the Bonferroni-type approach in \cite{ingster2000adaptive} combined with the exponential concentration bounds in Section~\ref{Section: Degenerate U-statistics for independence testing}. Here and hereafter, we restrict our attention to the nominal level $\alpha$ less than $e^{-1} \approx 0.368$, for which $\log(1/\alpha)$ is larger than $\sqrt{\log(1/\alpha)}$, to simplify our results.

\vskip 1em

\noindent \textbf{Two-sample testing.} Let us start with the two-sample problem. Without loss of generality, assume that $n_1 \leq n_2$ and consider a set of integers such that $\mathsf{K} := \{ 2^j : j = 1,\ldots, \gamma_{\text{max}}\}$ where
\begin{align*}
\gamma_{\text{max}} := \bigg\lceil\frac{2}{d} \log_2 \left( \frac{n_1}{\log \log n_1}\right)\bigg\rceil. 
\end{align*}
For each $\kappa \in \mathsf{K}$, we denote by $\phi_{\kappa, \alpha/\gamma_{\text{max}}}:= \ind(U_{n_1,n_2} > c_{1-\alpha/\gamma_{\text{max}},n})$, the multinomial two-sample test in Proposition~\ref{Proposition: Two-Sample Testing for Holder Densities} with the bin size $\kappa^{-1}$. We note that the type I error of an individual test is controlled at $\alpha /\gamma_{\text{max}}$ instead of $\alpha$. By taking the maximum of the resulting tests, we introduce an adaptive test for two-sample testing as follows:
\begin{align*}
\phi_{\text{adapt}} := \max_{\kappa \in \mathsf{K}}  \phi_{\kappa, \alpha/\gamma_{\text{max}}} .
\end{align*}
This adaptive test does not require knowledge on the smoothness parameter. 
We describe this result in the following proposition.

\begin{proposition}[Adaptive two-sample test] \label{Proposition: Adaptive two-sample test}
	Consider the same problem setting in Proposition~\ref{Proposition: Two-Sample Testing for Holder Densities} with an additional assumption that $n_1 \asymp n_2$. For a sufficiently large $C(s,d,L,\alpha,\beta) >0$, consider $\epsilon_{n_1,n_2}$ such that
	\begin{align*}
	\epsilon_{n_1,n_2} \geq C(s,d,L,\alpha,\beta)  \left( \frac{\log \log n_1}{n_1} \right)^{\frac{2s}{4s+d}}. 
	\end{align*}
	Then for testing $\mathcal{P}_0 = \{(P_Y,P_Z) \in \mathcal{P}_{\text{\emph{H\"{o}lder}}}^{(d,s)} : f_Y=f_Z\}$ against $\mathcal{P}_1 =\{(P_Y,P_Z) \in \mathcal{P}_{\text{\emph{H\"{o}lder}}}^{(d,s)} : |\!|\!|f_Y - f_Z|\!|\!|_{L_2} \geq \epsilon_{n_1,n_2}\}$, the type I and II errors of the adaptive test $\phi_{\text{\emph{adapt}}}$ are uniformly controlled as in (\ref{Eq: uniform error control}).
\end{proposition}
The test we propose is adaptive to an unknown smoothness parameter, but we pay a price of a factor of $(\log \log n_1)^{2s/(4s + d)}$ in the scaling of the critical radius. 
By the results of \cite{ingster2000adaptive}, we would expect the price for adaptation to scale as $\sqrt{\log \log n_1}^{2s/(4s + d)}$, and we hope to develop a more precise understanding of this gap 
in future work.

Type I error control of the adaptive test is trivial via the union bound. The proof of the type II error control is an application of Theorem~\ref{Theorem: Two-Sample Concentration} and can be found in Appendix~\ref{Section: Proposition: Adaptive two-sample test}. We note that the assumption $n_1 \asymp n_2$ is necessary to apply the concentration result in Theorem~\ref{Theorem: Two-Sample Concentration}, and it remains an open question whether the same result can be established without $n_1 \asymp n_2$.

\vskip 1em

\noindent \textbf{Independence testing.}  Let us now turn to the independence testing problem. Similarly as before, we define a set of integers by $\mathsf{K}^\dagger := \{ 2^j : j =1,\ldots, \gamma_{\text{max}}^\ast \}$ where 
\begin{align*}
\gamma_{\text{max}}^\ast := \bigg\lceil\frac{2}{d_1+ d_2} \log_2 \left( \frac{n}{\log \log n}\right)\bigg\rceil. 
\end{align*}
For each $\kappa \in \mathsf{K}^\dagger$, we use the notation $\phi_{\kappa, \alpha/\gamma_{\text{max}}^\ast}^\dagger:= \ind(U_n > c_{1-\alpha/\gamma_{\text{max}}^\ast,n})$ to denote the multinomial independence test in Proposition~\ref{Proposition: Minimum Separation for Independence Testing for Holder Densities} with the bin size $\kappa^{-1}$. Again we note that the type I error of an individual test is controlled at $\alpha / \gamma_{\text{max}}^\ast$ instead of $\alpha$. We then introduce an adaptive test for independence testing by taking the maximum of individual tests as
\begin{align*}
	\phi_{\text{adapt}}^\dagger := \max_{\kappa \in \mathsf{K}^\dagger} \phi_{\kappa, \alpha/\gamma_{\text{max}}^\ast}^\dagger.
\end{align*}
As in the two-sample case, the adaptive test does not depend on the smoothness parameter. In addition, when densities are smooth enough such that $4s > d_1 + d_2$, the adaptive test is minimax rate optimal up to an iterated logarithm factor as shown in the next proposition.

\begin{proposition}[Adaptive independence test] \label{Proposition: Adaptive independence test}
	Consider the same problem setting in Proposition~\ref{Proposition: Independence Testing for Holder Densities} and suppose that $4s > d_1 + d_2$. For a sufficiently large $C(s,d_1,d_2,L,\alpha,\beta) >0$, consider $\epsilon_{n}$ such that
	\begin{align*}
	\epsilon_{n} \geq C(s,d_1,d_2,L,\alpha,\beta) \left(\frac{\log \log n}{n}\right)^{\frac{2s}{4s+d_1+d_2}}. 
	\end{align*}
	Then for testing $\mathcal{P}_0 = \{P_{YZ} \in \mathcal{P}_{\text{\emph{H\"{o}lder}}}^{(d_1+d_2,s)} : f_{YZ} = f_Yf_Z \}$ against $\mathcal{P}_1 =\{P_{YZ} \in \mathcal{P}_{\text{\emph{H\"{o}lder}}}^{(d_1+d_2,s)} : |\!|\!|f_{YZ} - f_Yf_Z|\!|\!|_{L_2} \geq \epsilon_{n}\}$, the type I and II errors of the resulting permutation test are uniformly controlled as in (\ref{Eq: uniform error control}).
\end{proposition}

The proof of this result relies on Theorem~\ref{Theorem: Concentration inequality for independence U-statistic II} and is similar to that of Proposition~\ref{Proposition: Adaptive two-sample test}. The details can be found in Appendix~\ref{Section: Proposition: Adaptive independence test}. The restriction $4s > d_1 + d_2$ is imposed to guarantee that the first term $nt \Lambda_n^{-1}$ is smaller than the second term $n t^{2/3} M_n^{-3/2}$ in the tail bound~(\ref{Eq: Concentration inequality for independence U-statistic}) with high probability. Although it seems difficult, we believe that this restriction can be dropped with a more careful analysis. Alternatively one can convert independence testing to two-sample testing via sample-splitting (see Section~\ref{Section: Independence testing via sample-splitting} for details) and then apply the adaptive two-sample test in Proposition~\ref{Proposition: Adaptive two-sample test}. The resulting test has the same theoretical guarantee as in Proposition~\ref{Proposition: Adaptive independence test} without this restriction. However the sample-splitting approach should be considered with caution as it only uses a fraction of the data, which may result in a loss of power in practice.

\begin{remark}[Comparison to the two moments method] \normalfont
	While the exponential inequalities in Section~\ref{Section: Combinatorial concentration inequalities} lead to the adaptivity at the cost of $\log \log n$ factor, they are limited to degenerate $U$-statistics and require additional assumptions such as $n_1 \asymp n_2$ and $4s > d_1 + d_2$ to yield minimax rates. On the other hand, the two moments method is applicable beyond $U$-statistics and yields minimax rates without these extra assumptions. However we highlight that this generality comes at the cost of $\log n$ factor rather than $\log \log n$ to obtain the same adaptivity results. 
\end{remark}

\section{Further applications} \label{Section: Further applications}
In this section, we further investigate the power performance of permutation tests under different problem settings. One specific problem that we focus on is testing for multinomial distributions in the $\ell_1$ distance. The $\ell_1$ distance has an intuitive interpretation in terms of the total variation distance and has been considered as a metric for multinomial distribution testing \citep[see e.g.][and also references therein]{paninski2008coincidence,chan2014optimal,diakonikolas2016new,balakrishnan2019hypothesis}. Unlike the previous work, we approach this problem using the permutation procedure and study its minimax rate optimality in the $\ell_1$ distance. We also consider the problem of testing for continuous distributions and demonstrate the performance of the permutation tests based on reproducing kernel-based test statistics in Section~\ref{Section: Gaussian MMD} and Section~\ref{Section: Gaussian HSIC}.

\subsection{Two-sample testing under Poisson sampling with equal sample sizes} \label{Section: Two-sample testing under Poisson sampling}
Let $p_Y$ and $p_Z$ be multinomial distributions defined on $\mathbb{S}_d$. Suppose that we observe samples from Poisson distributions as $\{Y_{1,k},\ldots,Y_{n,k} \} \overset{\text{i.i.d.}}{\sim} \text{Poisson}\{p_Y(k)\}$ and $\{Z_{1,k},\ldots,Z_{n,k}\} \overset{\text{i.i.d.}}{\sim} \text{Poisson}\{p_Z(k)\}$ for each $k \in \{1,\ldots,d\}$. Assume that all these samples are mutually independent. Let us write $V_k := \sum_{i=1}^n Y_{i,k}$ and $W_k := \sum_{i=1}^n Z_{i,k}$ where $V_k$ and $W_k$ have Poisson distributions with parameters $np_Y(k)$ and $np_Z(k)$, respectively. Under this distributional assumption, \cite{chan2014optimal} consider a centered chi-square test statistic given by
\begin{align} \label{Eq: Chi-square statistic}
T_{\chi^2}:= \sum_{k=1}^d \frac{(V_k - W_k)^2 - V_k - W_k}{V_k + W_k} \ind(V_k+W_k > 0).
\end{align}
Based on this statistic, they show that if one rejects the null $H_0: p_Y = p_Z$ when $T_{\chi^2}$ is greater than $C \sqrt{\min\{n,d\}}$ for some constant $C$, then the resulting test is minimax rate optimal for the class of alternatives determined by the $\ell_1$ distance. In particular, the minimax rate is shown to be
\begin{align} \label{Eq: Minimax Separation for equal sized two-sample testing}
\epsilon_n^\dagger \asymp \max \Bigg\{ \frac{d^{1/2}}{n^{3/4}}, \ \frac{d^{1/4}}{n^{1/2}} \Bigg\}.
\end{align}
However, in their test, the choice of $C$ is implicit and based on a loose concentration inequality. Here, by letting $\{X_{i,k}\}_{i=1}^{2n}$ be the pooled samples of $\{Y_{i,k}\}_{i=1}^n$ and $\{Z_{i,k}\}_{i=1}^n$, we instead determine the critical value via the permutation procedure. In this setting the permuted test statistic is 
\begin{align*} 
T_{\chi^2}^{\pi} := \sum_{k=1}^d \frac{(\sum_{i=1}^n X_{\pi_i,k} - \sum_{i=1}^n X_{\pi_{i+n},k} )^2 - V_k - W_k}{V_k + W_k} \ind(V_k+W_k > 0).
\end{align*}
The next theorem shows that the resulting permutation test is also minimax rate optimal. 
\begin{theorem}[Two-sample testing under Poisson sampling] \label{Theorem: Two-Sample testing under Poisson sampling}
	Consider the distributional setting described above. For a sufficiently large $C >0$, let us consider a positive sequence $\epsilon_n$ such that 
	\begin{align*}
	\epsilon_n \geq \frac{C}{\beta}\sqrt{\log \left( \frac{1}{\alpha} \right)} \cdot \max \Bigg\{ \frac{d^{1/2}}{n^{3/4}}, \ \frac{d^{1/4}}{n^{1/2}} \Bigg\}.
	\end{align*}
	Then for testing $\mathcal{P}_0 = \{(p_Y,p_Z) : p_Y = p_Z\}$ against $\mathcal{P}_1 = \{(p_Y,p_Z) : \|p_Y - p_Z\|_1 \geq \epsilon_n \}$, the type I and II errors of the permutation test based on $T_{\chi^2}$ are uniformly controlled as in (\ref{Eq: uniform error control}).	
\end{theorem}
It is worth noting that $\sqrt{\log(1/\alpha)}$ factor in Theorem~\ref{Theorem: Two-Sample testing under Poisson sampling} is a consequence of applying the exponential concentration inequality in Section~\ref{Section: Combinatorial concentration inequalities}. We also note that this logarithmic factor cannot be obtained by the technique used in \cite{chan2014optimal} which only bounds the mean and variance of the test statistic. On the other hand, the dependency on $\beta$ may be sub-optimal and may be improved via a more sophisticated analysis. We leave this direction to future work.

\subsection{Two-sample testing via sample-splitting} \label{Section: Two-sample testing via sample-splitting}
Although the chi-square two-sample test in Theorem~\ref{Theorem: Two-Sample testing under Poisson sampling} is simple and comes with a theoretical guarantee of minimax optimality, it is only valid in the setting of equal sample sizes. The goal of this subsection is to provide an alternative permutation test via sample-splitting which is minimax rate optimal regardless of the sample size ratio. When the two sample sizes are different, \cite{bhattacharya2015testing} modify the chi-square statistic (\ref{Eq: Chi-square statistic}) and propose an optimal test but with the additional assumption that $\epsilon_{n_1,n_2} \geq d^{-1/12}$. \cite{diakonikolas2016new} remove this extra assumption and introduce another test with the same statistical guarantee. Their test is based on the flattening idea that artificially transforms the probability distributions to be roughly uniform. The same idea is considered in \cite{canonne2018testing} for conditional independence testing. 

Despite their optimality, neither \cite{bhattacharya2015testing} nor \cite{diakonikolas2016new} presents a concrete way of choosing the critical value that leads to a level $\alpha$ test. Here we address this issue based on the permutation procedure.

Suppose that we observe $\mathcal{Y}_{2n_1}$ and $\mathcal{Z}_{2n_2}$ samples from two multinomial distributions $p_Y$ and $p_Z$ defined on $\mathbb{S}_d$, respectively. Without loss of generality, we assume that $n_1 \leq n_2$. Let us define $m := \min\{n_2,d\}$ and denote data-dependent weights, computed based on $\{Z_{n_2+1},\ldots, Z_{n_2+m} \}$, by
\begin{align*}
w_k := \frac{1}{2d} + \frac{1}{2m} \sum_{i=1}^m \ind(Z_{i+n_2} = k) \quad \text{for $k=1,\ldots,d$.} 
\end{align*}
Under the given scenario, we consider the two-sample $U$-statistic~(\ref{Eq: Two-Sample U-statistic}) defined with the following bivariate function:
\begin{align} \label{Eq: weighted kernel}
g_{\text{Multi},w}(x,y) := \sum_{k=1}^d w_k^{-1} \ind(x= k) \ind(y=k).
\end{align}
We emphasize that the considered $U$-statistic is evaluated based on the first $n_1$ observations from each group, i.e.~$\mathcal{X}_{2n_1}^{\text{split}}:= \{Y_1,\ldots,Y_{n_1},Z_1,\ldots,Z_{n_1}\}$, which are clearly independent of weights $\{w_1,\ldots,w_d\}$. Let us denote the $U$-statistic computed in this way by $U_{n_1,n_2}^{\text{split}}$. Let us consider the critical value of a permutation test obtained by permuting the labels within $\mathcal{X}_{2n_1}^\text{split}$. Then the resulting permutation test via sample-splitting has the following theoretical guarantee. 
\begin{proposition}[Multinomial two-sample testing in the $\ell_1$ distance]\label{Proposition: Multinomial L1 Testing}
	Let $\mathcal{P}_{\text{\emph{Multi}}}^{(d)}$ be the set of pairs of multinomial distributions defined on $\mathbb{S}_d$. Let $\mathcal{P}_0 = \{(p_Y,p_Z) \in \mathcal{P}_{\text{\emph{Multi}}}^{(d)} : p_Y = p_Z\}$ and $\mathcal{P}_1(\epsilon_{n_1,n_2}) = \{(p_Y,p_Z) \in \mathcal{P}_{\text{\emph{Multi}}}^{(d)} : \|p_Y - p_Z\|_1 \geq \epsilon_{n_1,n_2} \}$ where
	\begin{align} \label{Eq: sufficient condition for L1 two-sample testing}
	\epsilon_{n_1,n_2} \geq \frac{C}{\beta^{3/4}} \sqrt{\log\left( \frac{1}{\alpha} \right)} \cdot \max \Bigg\{ \frac{d^{1/2}}{n_1^{1/2}n_2^{1/4}}, \ \frac{d^{1/4}}{n_1^{1/2}} \Bigg\}, 
	\end{align}
	for a sufficiently large $C>0$. Consider the two-sample $U$-statistic $U_{n_1,n_2}^{\text{\emph{split}}}$ described above. Then, over the classes $\mathcal{P}_0$ and $\mathcal{P}_1$, the type I and II errors of the resulting permutation test via sample-splitting are uniformly bounded as in (\ref{Eq: uniform error control}). 
\end{proposition}

\vskip 1em

\begin{proof}[\textbf{\emph{Proof Sketch.}}] The proof of this result can be found in Appendix~\ref{Section: Proof of Proposition: Multinomial L1 Testing}. To sketch the proof, 
conditional on weights $w_1,\ldots,w_d$, the problem of interest is essentially the same as that of Proposition~\ref{Proposition: Multinomial Two-Sample Testing}. One difference is that $U_{n_1,n_2}^{\text{split}}$ is not an unbiased estimator of $\|p_Y - p_Z\|_1$. However, by noting that $\sum_{k=1}^d w_i =1$, one can lower bound the expected value in terms of the $\ell_1$ distance by Cauchy-Schwarz inequality as
\begin{align*}
\mE_P \big[U_{n_1,n_2}^{\text{split}}|w_1,\ldots,w_n \big] = \sum_{k=1}^d \frac{\{ p_Y(k) - p_Z(k) \}^2}{w_k} \geq \| p_Y - p_Z\|_1^2.
\end{align*}
The conditional variance can be similarly bounded as in Proposition~\ref{Proposition: Multinomial Two-Sample Testing} and we use Theorem~\ref{Theorem: Two-Sample Concentration} to study the critical value of the permutation test. Finally, we remove the randomness from the weights $w_1,\ldots,w_d$ via Markov's inequality to complete the proof.
\end{proof}

The results of \cite{bhattacharya2015testing} and \cite{diakonikolas2016new} show that the minimum separation for $\ell_1$-closeness testing satisfies
\begin{align*} 
\epsilon_{n_1,n_2}^\dagger \asymp \max \Bigg\{ \frac{d^{1/2}}{n_2^{1/4} n_1^{1/2}}, \ \frac{d^{1/4}}{n_1^{1/2}} \Bigg\}.
\end{align*}
This means that the proposed permutation test is minimax rate optimal for multinomial testing in the $\ell_1$ distance. On the other hand the procedure depends on sample-splitting which may result in a loss of practical power.
Indeed all of the previous approaches~\citep{acharya2014sublinear,bhattacharya2015testing,diakonikolas2016new} also depend on sample-splitting, which leaves the important question as to whether it is possible to obtain the same minimax guarantee without sample-splitting.

\subsection{Independence testing via sample-splitting} \label{Section: Independence testing via sample-splitting}
We now turn to independence testing for multinomial distributions in the $\ell_1$ distance. To take full advantage of the two-sample test developed in the previous subsection, we follow the idea of \cite{diakonikolas2016new} in which the independence testing problem is converted into the two-sample problem as follows. Suppose that we sample $\mathcal{X}_{3n}$ observations from a joint multinomial distribution $p_{YZ}$ on $\mathbb{S}_{d_1, d_2}$. We then take the first one-third of the data and denote it by $\widetilde{Y}_{n} := \{(Y_1,Z_1),\ldots,(Y_n,Z_n) \}$. Using the remaining data, we define another set of samples $\widetilde{Z}_{n} := \{(Y_{n+1},Z_{2n+1}),\ldots, (Y_{2n},Z_{3n})\}$. By construction, it is clear that $\widetilde{Y}_{n}$ consists of samples from the joint distribution $p_{YZ}$ whereas $\widetilde{Z}_{n}$ consists of samples from the product distribution $p_Yp_Z$. In other words, we have a fresh dataset $\widetilde{\mathcal{X}}_n:= \widetilde{Y}_n \cup \widetilde{Z}_n$ for two-sample testing. It is interesting to mention, however, that the direct application of the two-sample test in Proposition~\ref{Proposition: Multinomial L1 Testing} to $\widetilde{\mathcal{X}}_n$ does not guarantee optimality. In particular, by replacing $d$ with $d_1d_2$ and letting $n_1=n_2=n$ in condition~(\ref{Eq: sufficient condition for L1 two-sample testing}), we see that the permutation test has power when $\epsilon_{n_1,n_2}$ is sufficiently larger than $\max \big\{d_1^{1/2}d_2^{1/2}n^{-3/4}, d_1^{1/4}d_2^{1/4}n^{-1/2} \big\}$, whereas by assuming $d_1 \leq d_2$, the minimum separation for independence testing in the $\ell_1$ distance \citep{diakonikolas2016new} is given by
\begin{align} \label{Eq: minimum separation for multinomial independence testing}
\epsilon_n^{\dagger} \asymp \max \Bigg\{ \frac{d_1^{1/4}d_2^{1/2}}{n^{3/4}}, \ \frac{d_1^{1/4}d_2^{1/4} }{n^{1/2}} \Bigg\}.
\end{align} 
The main reason is that, unlike the original two-sample problem where two distributions can be arbitrary different, we have the further restriction that the marginal distributions of $p_{YZ}$ are the same as those of $p_Yp_Z$. Therefore we need to consider a more refined weight function for independence testing to derive an optimal test. To this end, for each $(k_1,k_2) \in \mathbb{S}_{d_1,d_2}$, we define a product weight by
\begin{align*}
w_{k_1,k_2} := \Bigg[ \frac{1}{2d_1} + \frac{1}{2m_1} \sum_{i=1}^{m_1} \ind(Y_{3n/2+i} = k_1)  \Bigg] \times \Bigg[ \frac{1}{2d_2} + \frac{1}{2m_2} \sum_{j=1}^{m_2} \ind(Z_{5n/2+i} = k_2) \Bigg],  
\end{align*}
where $m_1 := \min\{n/2,d_1\}$ and $m_2 := \min\{n/2,d_2\}$ and we assume $n$ is even. Notice that by construction, the given product weights are independent of the first half of $\widetilde{\mathcal{X}}_{n}$, denoted by $\widetilde{\mathcal{X}}_{n/2}^{\text{split}}$. Similarly as before, we use $\widetilde{\mathcal{X}}_{n/2}^{\text{split}}$ to compute the two-sample $U$-statistic~(\ref{Eq: Two-Sample U-statistic}) defined with the following bivariate function:
\begin{align*} 
g_{\text{Multi},w}^{\ast}\{(x_1,y_1),(x_2,y_2)\} := \sum_{k_1=1}^{d_1} \sum_{k_2=1}^{d_2}  w_{k_1,k_2}^{-1} \ind(x_1=k_1,y_1=k_2)  \ind(x_2=k_1,y_2=k_2),
\end{align*}
and denote the resulting test statistic by $U_{n_1,n_2}^{\text{split} \ast}$. The critical value is determined by permuting the labels within $\widetilde{\mathcal{X}}_{n/2}^{\text{split}}$ and the resulting test has the following theoretical guarantee. 
\begin{proposition}[Multinomial independence testing in $\ell_1$ distance]  \label{Proposition: Multinomial independence testing in L1 distance}
	Let $\mathcal{P}_{\text{\emph{Multi}}}^{(d_1,d_2)}$ be the set of multinomial distributions defined on $\mathbb{S}_{d_1, d_2}$. Let $\mathcal{P}_0 = \{p_{YZ} \in \mathcal{P}_{\text{\emph{Multi}}}^{(d_1,d_2)}: p_{YZ} = p_Yp_Z \}$ and $\mathcal{P}_1(\epsilon_n) = \{p_{YZ} \in \mathcal{P}_{\text{\emph{Multi}}}^{(d_1,d_2)}: \|p_{YZ} - p_Yp_Z\|_2 \geq \epsilon_n \}$ where
	\begin{align*}
	\epsilon_{n} \geq \frac{C}{\beta^{3/4}}  \sqrt{\log\left( \frac{1}{\alpha} \right)} \cdot \max \Bigg\{ \frac{d_1^{1/4}d_2^{1/2}}{n^{3/4}}, \ \frac{d_1^{1/4}d_2^{1/4} }{n^{1/2}} \Bigg\},
	\end{align*}
	for a sufficiently large $C>0$ and $d_1 \leq d_2$. Consider the two-sample $U$-statistic $U_{n_1,n_2}^{\text{\emph{split}}\ast}$ described above. Then, over the classes $\mathcal{P}_0$ and $\mathcal{P}_1$, the type I and II errors of the resulting permutation test via sample-splitting are uniformly bounded as in (\ref{Eq: uniform error control}). 
\end{proposition}
In view of the minimum separation~(\ref{Eq: minimum separation for multinomial independence testing}), the proposed test above is minimax rate optimal for multinomial independence testing in the $\ell_1$ distance.  We note again that sample-splitting is mainly for technical convenience and it might result in a loss of efficiency in practice. In order to use the data more efficiently and remove the randomness from a single split, we recommend that one considers many splits and uses the average of the resulting statistics as their final test statistic. However this strategy requires a more delicate power analysis, which is beyond the scope of the present paper. An interesting direction of future work is therefore to see whether one can obtain the same minimax guarantee by using multiple splits or developing other statistics that avoid sample-splitting.

\subsection{Gaussian MMD} \label{Section: Gaussian MMD}
In this subsection we switch gears to continuous distributions and focus on the two-sample $U$-statistic with a Gaussian kernel. For $x,y \in \mathbb{R}^d$ and $\lambda_1,\ldots,\lambda_d >0$, the Gaussian kernel is defined by
\begin{align} \label{Eq: Gaussian kernel}
g_{\text{Gau}}(x,y):= K_{\lambda_1,\ldots,\lambda_d,d}(x-y) = \frac{1}{(2\pi)^{d/2} \lambda_1 \cdots \lambda_d} \exp \Bigg\{ - \frac{1}{2} \sum_{i=1}^d \frac{(x_i - y_i)^2}{\lambda_i^2} \Bigg\}.
\end{align}
The two-sample $U$-statistic defined with this Gaussian kernel is known as the Gaussian maximum mean discrepancy (MMD) statistic due to \cite{gretton2012kernel} and is also related to the test statistic considered in \cite{anderson1994two}. The Gaussian MMD statistic has a nice property that its expectation becomes zero if and only if $P_Y = P_Z$. Given the $U$-statistic with the Gaussian kernel, we want to find a sufficient condition under which the resulting permutation test has non-trivial power against alternatives determined with respect to the $L_2$ distance. In detail, by letting $f_Y$ and $f_Z$ be the density functions of $P_Y$ and $P_Z$ with respect to Lebesgue measure, consider the set of paired distributions $(P_Y,P_Z)$ such that the infinity norms of their densities are uniformly bounded, i.e.~$\max\{|\!|\!|f_Y|\!|\!|_\infty, |\!|\!|f_Z|\!|\!|_\infty \} \leq M_{f,d} < \infty$. We denote such a set by $\mathcal{P}_{\infty}^d$. Then for the class of alternatives $\mathcal{P}_1(\epsilon_{n_1,n_2}) = \{ (P_Y,P_Z) \in \mathcal{P}_{\infty}^d: |\!|\!|f_Y - f_Z|\!|\!|_{L_2} \geq \epsilon_{n_1,n_2} \}$, the following proposition gives a sufficient condition on $\epsilon_{n_1,n_2}$ under which the permutation-based MMD test has non-trivial power. It is worth noting that a similar result exists in \cite{fromont2013two} where they study the two-sample problem for Poisson processes using a wild bootstrap method. The next proposition differs from their result in three different ways: (1) we consider the usual i.i.d.~sampling scheme, (2) we do not assume that $n_1$ and $n_2$ are the same and (3) we use the permutation procedure, which is more generally applicable than the wild bootstrap procedure. 

\begin{proposition}[Gaussian MMD] \label{Proposition: Gaussian MMD}
	Consider the permutation test based on the two-sample $U$-statistic $U_{n_1,n_2}$ with the Gaussian kernel where we assume $\prod_{i=1}^d \lambda_i \leq 1$ and $n_1 \asymp n_2$. For a sufficiently large $C(M_{f,d},d) >0$, consider $\epsilon_{n_1,n_2}$ such that 
	\begin{equation}
	\begin{aligned} \label{Eq: Gaussian MMD epsilon}
	\epsilon_{n_1,n_2}^2 ~\geq~ & |\!|\!| (f_Y - f_Z)  - (f_Y - f_Z) \ast K_{\lambda,d} |\!|\!|_{L_2}^2 \\[.5em]
	 + ~& \frac{C(M_{f,d},d)}{\beta \sqrt{\lambda_1 \cdots \lambda_d}}  \log \left( \frac{1}{\alpha} \right) \cdot \left( \frac{1}{n_1} + \frac{1}{n_2} \right), 
	\end{aligned}
	\end{equation}
	where $\ast$ is the convolution operator with respect to Lebesgue measure. Then for testing $\mathcal{P}_0 = \{(P_Y,P_Z) \in \mathcal{P}_\infty^d : f_Y=f_Z\}$ against $\mathcal{P}_1 =\{(P_Y,P_Z) \in \mathcal{P}_\infty^d : |\!|\!|f_Y - f_Z|\!|\!|_{L_2} \geq \epsilon_{n_1,n_2}\}$, the type I and II errors of the resulting permutation test are uniformly controlled as in (\ref{Eq: uniform error control}).
\end{proposition}
The proof of this result is based on the exponential concentration inequality in Theorem~\ref{Theorem: Two-Sample Concentration} and the details are deferred to Appendix~\ref{Section: Proof of Proposition: Gaussian MMD}. One can remove the assumption that $n_1 \asymp n_2$ using the two moment method in Theorem~\ref{Theorem: Two-Sample U-statistic} but in this case, the result relies on a polynomial dependence on $\alpha$. The first term on the right-hand side of condition~(\ref{Eq: Gaussian MMD epsilon}) can be interpreted as a bias term, which measures a difference between the $L_2$ distance and the Gaussian MMD. The second term is related to the variance of the test statistic. We note that there is a certain trade-off between the bias and the variance, depending on the choice of tuning parameters $\{\lambda_i\}_{i=1}^d$. To make the bias term more explicit, we make some regularity assumptions on densities, following \cite{fromont2013two} and \cite{meynaoui2019aggregated}, and discuss the optimal choice of $\{\lambda_i\}_{i=1}^d$ under each condition.

\begin{example}[Sobolev ball] \normalfont  \label{Example: Sobolev}
	For $s, R >0$, the Sobolev ball $\mathcal{S}_d^s(R)$ is defined as
	\begin{align*}
	\mathcal{S}_d^s(R) := \bigg\{ q:  \mathbb{R}^d \mapsto \mathbb{R} \bigg/ q \in L^1\big(\mathbb{R}^d\big) \cap L^2\big(\mathbb{R}^d\big), \ \int_{\mathbb{R}^d} \|u\|^{2s} | \widehat{q}(u)|^2 du \leq (2\pi)^dR^2  \bigg\},
	\end{align*} 
	where $\widehat{q}$ is the Fourier transform of $q$, i.e.~$\widehat{q}(u) := \int_{\mathbb{R}^d} q(x) e^{i \langle x, u \rangle } dx$ and $\langle x, u \rangle$ is the scalar product in $\mathbb{R}^d$. Suppose that $f_Y - f_Z \in \mathcal{S}_d^s(R)$ where $s \in (0, 2]$. Then following Lemma~3 of \cite{meynaoui2019aggregated}, it can be seen that the bias term is bounded by
	\begin{align*}
	|\!|\!| (f_Y - f_Z)  - (f_Y - f_Z) \ast K_{\lambda,d} |\!|\!|_{L_2}^2 \leq C(R,s,d) \sum_{k=1}^d \lambda_k^{2s}.
	\end{align*} 
	Now we further upper bound the right-hand side of condition~(\ref{Eq: Gaussian MMD epsilon}) using the above result and then optimize it over $\lambda_1,\ldots,\lambda_d$. This can be done by putting $\lambda_1 = \cdots = \lambda_d = (n_1^{-1}  + n_2^{-1})^{2/(4s + d)}$, which in turn yields
	\begin{align} \label{Eq: Sobolev ball condition}
	\epsilon_{n_1,n_2} \geq C(M_{f,d},R,s,d,\alpha,\beta) \left( \frac{1}{n_1} + \frac{1}{n_2} \right)^{\frac{2s}{4s + d}}.
	\end{align} 
	In other words, Proposition~\ref{Proposition: Gaussian MMD} holds over the Sobolev ball as long as condition~(\ref{Eq: Sobolev ball condition}) is satisfied. 
\end{example}

By leveraging the minimax lower bound result in \cite{meynaoui2019aggregated} and the proof of Proposition~\ref{Proposition: Minimum Separation for Two-Sample Multinomial Testing}, it is straightforward to prove that the minimum separation rate for two-sample testing over the Sobolev ball is $n_1^{-2s/(4s+d)}$ for $n_1 \leq n_2$. This means that the permutation-based MMD test is minimax rate optimal over the Sobolev ball. In the next example, we consider an anisotropic Nikol'skii-Besov ball that can have different regularity conditions over $\mathbb{R}^d$.

\begin{example}[Nikol'skii-Besov ball] \normalfont \label{Example: Nikolskii}
	For $\boldsymbol{s}:= (s_1,\ldots,s_d) \in (0,\infty)^d$ and $R>0$, the anisotropic Nikol'skii-Besov ball $\mathcal{N}_{2,d}^{\boldsymbol{s}}(R)$ defined by 
	\begin{align*}
	\mathcal{N}_{2,d}^{\boldsymbol{s}}(R) := \bigg\{ & q: \mathbb{R}^d \mapsto \mathbb{R} \bigg/  \text{$q$ has continuous partial derivatives $D_i^{\floor{s_i}}$ of order $\floor{s_i}$} \\
	&\text{with respect to $u_i$ and for all $i=1,\ldots,d$, $u_1,\ldots,u_d,v \in \mathbb{R}$,} \\
	& \big|\!\big|\!\big|D_i^{\floor{s_i}} q(u_1,\ldots,u_i+v,\ldots,u_d) - D_i^{\floor{s_i}} q(u_1,\ldots,u_d) \big|\!\big|\!\big|_{L_2} \leq R|v|^{s_i - \floor{s_i}} \bigg\}.
	\end{align*}
	Suppose that $f_Y - f_Z \in \mathcal{N}_{2,d}^{\boldsymbol{s}}(R)$ where $\boldsymbol{s} \in (0,2]^d$. Then similarly to Lemma~4 of \cite{meynaoui2019aggregated}, it can be shown that the bias term is bounded by
	\begin{align*}
	|\!|\!| (f_Y - f_Z)  - (f_Y - f_Z) \ast K_{\lambda,d} |\!|\!|_{L_2}^2 \leq C(R,\boldsymbol{s},d) \sum_{k=1}^d \lambda_k^{2s_k}.
	\end{align*}
	Again we further upper bound the right-hand side of condition~(\ref{Eq: Gaussian MMD epsilon}) using the above result and then minimize it over $\lambda_1,\ldots,\lambda_d$. Letting $\eta^{-1} = \sum_{k=1}^d s_k^{-1}$, the minimum (up to a constant factor) can be achieved when $\lambda_k = (n_1^{-1} + n_2^{-1})^{2\eta/\{s_k(1 + 4\eta) \}}$ for $k=1,\ldots,d$, which yields 
	\begin{align} \label{Eq: Nikolskii-Besov ball condition}
	\epsilon_{n_1,n_2} \geq C(M_{f,d},R,\boldsymbol{s},d,\alpha,\beta) \left( \frac{1}{n_1} + \frac{1}{n_2} \right)^{\frac{2\eta}{1 + 4\eta}}.
	\end{align} 
	Therefore we are guaranteed that Proposition~\ref{Proposition: Gaussian MMD} holds over the Nikol'skii-Besov ball as long as condition~(\ref{Eq: Nikolskii-Besov ball condition}) is satisfied. 
\end{example}

\subsection{Gaussian HSIC} \label{Section: Gaussian HSIC}
We now focus on independence testing for continuous distributions. In particular we study the performance of the permutation test using the $U$-statistic~(\ref{Eq: U-statistic for independence testing}) defined with Gaussian kernels. For $y_1,y_2 \in \mathbb{R}^{d_1}$, $z_1,z_2 \in \mathbb{R}^{d_2}$ and $\lambda_1,\ldots,\lambda_{d_1}, \gamma_1,\ldots, \gamma_{d_2} > 0$, let us recall the definition of a Gaussian kernel~(\ref{Eq: Gaussian kernel}) and similarly write
\begin{align} \label{Eq: Gaussian kernels for independence testing}
g_{\text{Gau},Y}(y_1,y_2):= K_{\lambda_1,\ldots,\lambda_{d_1},d_1}(y_1-y_2)  \quad \text{and} \quad g_{\text{Gau},Z}(z_1,z_2):= K_{\gamma_1,\ldots,\gamma_{d_2},d_2}(z_1-z_2).
\end{align}
The $U$-statistic~(\ref{Eq: U-statistic for independence testing}) defined with these Gaussian kernels is known as the Hilbert--Schmidt independence criterion (HSIC) statistic \citep{gretton2005measuring}. As in the case of the Gaussian MMD, it is well-known that the expected value of the Gaussian HSIC statistic becomes zero if and only if $P_{YZ} = P_YP_Z$. Using this property, the resulting test can be consistent against any fixed alternative. \cite{meynaoui2019aggregated} consider the same statistic and study the power of a HSIC-based test over Sobolev and Nikol'skii-Besov balls. It is important to note, however, that the critical value of their test is calculated based on the (theoretical) null distribution of the test statistic, which is unknown in general. The aim of this subsection is to extend their results to the permutation test that does not require knowledge of the null distribution. To describe the main result, let us write the density functions of $P_{YZ}$ and $P_YP_Z$ with respect to Lebesgue measure by $f_{YZ}$ and $f_Yf_Z$. As in Section~\ref{Section: Gaussian MMD}, we use $\mathcal{P}_{\infty}^{d_1,d_2}$ to denote the set of distributions $P_{YZ}$ whose joint and product densities are uniformly bounded, i.e.~$\max\{|\!|\!|f_{YZ}|\!|\!|_{\infty}, |\!|\!|f_Yf_Z|\!|\!|_{\infty} \} \leq M_{f,d_1,d_2} < \infty$. Then the following proposition presents a theoretical guarantee for the permutation-based HSIC test.

\begin{proposition}[Gaussian HSIC] \label{Proposition: Gaussian HSIC}
	Consider the permutation test based on the $U$-statistic $U_{n}$ with the Gaussian kernels~(\ref{Eq: Gaussian kernels for independence testing}) where we assume $\prod_{i=1}^{d_1} \lambda_i \leq 1$ and $\prod_{i=1}^{d_2} \gamma_i \leq 1$. For a sufficiently large $C(M_{f,d_1,d_2},d_1,d_2) >0$, consider $\epsilon_{n}$ such that
	\begin{equation}
	\begin{aligned} \label{Eq: Gaussian HSIC epsilon}
	\epsilon_{n}^2 ~\geq ~& |\!|\!| (f_{YZ} - f_Yf_Z)  - (f_{YZ} - f_Yf_Z) \ast (K_{\lambda,d_1}  K_{\gamma,d_2}) |\!|\!|_{L_2}^2 \\[.5em]
	 +~ &\frac{C(M_{f,d_1,d_2},d_1,d_2)}{\alpha^{1/2} \beta n\sqrt{\lambda_1 \cdots \lambda_{d_1} \gamma_1 \cdots \gamma_{d_2}}}, 
	\end{aligned}
	\end{equation}
	where $\ast$ is the convolution operator with respect to Lebesgue measure. Then for testing $\mathcal{P}_0 = \{ P_{YZ} \in \mathcal{P}_{\infty}^{d_1,d_2} : f_{YZ}=f_Yf_Z\}$ against $\mathcal{P}_1 =\{ P_{YZ} \in \mathcal{P}_{\infty}^{d_1,d_2} : |\!|\!|f_{YZ} - f_Yf_Z|\!|\!|_{L_2} \geq \epsilon_{n}\}$, the type I and II errors of the resulting permutation test are uniformly controlled as in (\ref{Eq: uniform error control}).
\end{proposition}
The proof of this result is based on the two moments method in Proposition~\ref{Theorem: U-statistic for independence testing}. We omit the proof of this result since it is very similar to that of Proposition~\ref{Proposition: Gaussian MMD} and Theorem 1 of \cite{meynaoui2019aggregated}. As before, the first term on the right-hand side of condition~(\ref{Eq: Gaussian HSIC epsilon}) can be viewed as a bias, which measures a difference between the $L_2$ distance and the Gaussian HSIC. To make this bias term more tractable, we now consider Sobolev and Nikol'skii-Besov balls and further illustrate Proposition~\ref{Proposition: Gaussian HSIC}. The following two examples correspond Corollary 2 and Corollary 3 of \cite{meynaoui2019aggregated} but based on the permutation test.

\begin{example}[Sobolev ball] \normalfont
	Recall the definition of the Sobolev ball from Example~\ref{Example: Sobolev} and assume that $f_{YZ} - f_Yf_Z \in \mathcal{S}_{d_1+d_2}^s(R)$ where $s \in (0,2]$. Then from Lemma~3 of \cite{meynaoui2019aggregated}, the bias term in condition~(\ref{Eq: Gaussian HSIC epsilon}) can be bounded by 
	\begin{align*}
	|\!|\!| (f_{YZ} - f_Yf_Z)  - (f_{YZ} - f_Yf_Z) \ast (K_{\lambda,d_1}  K_{\gamma,d_2}) |\!|\!|_{L_2}^2 \leq C(R,s,d_1,d_2) \Bigg\{ \sum_{i=1}^{d_1} \lambda_i^{2s} +  \sum_{j=1}^{d_2} \gamma_j^{2s} \Bigg\}.
	\end{align*}
	For each $i \in \{1,\ldots,d_1\}$ and $j \in \{1,\ldots,d_2\}$, we choose $\lambda_i = \gamma_j = n^{-2/(4s+d_1+d_2)}$ such that the lower bound of $\epsilon_n$ in condition~(\ref{Eq: Gaussian HSIC epsilon}) is minimized. Then by plugging these parameters, it can be seen that Proposition~\ref{Proposition: Gaussian HSIC} holds as long as $\epsilon_n \geq C(M_{f,d_1,d_2},s,R,d_1,d_2,\alpha,\beta) n^{-\frac{2s}{4s + d_1 + d_2}}$. Furthermore, this rate matches with the lower bound given in \cite{meynaoui2019aggregated}.
\end{example}

\begin{example}[Nikol'skii-Besov ball] \normalfont
	Recall the definition of the Nikol'skii-Besov ball from Example~\ref{Example: Nikolskii} and assume that $f_{YZ} - f_Yf_Z \in \mathcal{N}_{2,d_1+d_2}^{\boldsymbol{s}}(R)$ where $\boldsymbol{s} \in (0,2]^{d_1+d_2}$. Then followed by Lemma 4 of \cite{meynaoui2019aggregated}, the bias term in condition~(\ref{Eq: Gaussian HSIC epsilon}) can be bounded by 
	\begin{align*}
	|\!|\!| (f_{YZ} - f_Yf_Z)  - (f_{YZ} - f_Yf_Z) \ast (K_{\lambda,d_1}  K_{\gamma,d_2}) |\!|\!|_{L_2}^2 \leq C(R,\boldsymbol{s},d_1,d_2) \Bigg\{ \sum_{i=1}^{d_1} \lambda_i^{2s_i} +  \sum_{j=1}^{d_2} \gamma_j^{2s_{j+d_1}} \Bigg\}.
	\end{align*}
	Let us write $\eta^{-1} := \sum_{i=1}^{d_1+d_2} s_i^{-1}$. Then by minimizing the lower bound of $\epsilon_n$ in condition~(\ref{Eq: Gaussian HSIC epsilon}) using the above result with $\lambda_i = n^{-\frac{2\eta}{s_i(1+4\eta)}}$ for $i=1,\ldots,d_1$ and $\gamma_i = n^{-\frac{2\eta}{s_{i+d_1}(1+4\eta)}}$ for $i=1,\ldots,d_2$, it can be seen that the same conclusion of Proposition~\ref{Proposition: Gaussian HSIC} holds as long as $\epsilon_n \geq C(M_{f,d_1,d_2},\boldsymbol{s},R,d_1,d_2,\alpha,\beta) n^{-\frac{2\eta}{1 + 4 \eta}}$.
\end{example}

From the above two examples, we see that the permutation-based HSIC test has the same power guarantee as the theoretical test considered in \cite{meynaoui2019aggregated}. However, our results do not fully recover those in \cite{meynaoui2019aggregated} in terms of $\alpha$. It remains an open question as to whether Proposition~\ref{Proposition: Gaussian HSIC} continues to hold when $\alpha^{-1/2}$ is replaced by $\log(1/\alpha)$. Alternatively one can employ the sample-splitting idea in Section~\ref{Section: Independence testing via sample-splitting} and apply the permutation-based MMD test in Proposition~\ref{Proposition: Gaussian MMD} for independence testing. The result of Proposition~\ref{Proposition: Gaussian MMD} then guarantees that the MMD test achieves the same rate of the power as the permutation-based HSIC test but it improves the dependency on $\alpha$ in condition~(\ref{Eq: Gaussian HSIC epsilon}) to a logarithmic factor.

\section{Simulations} \label{Section: Simulations}
\begin{figure}[t!]
	\begin{center}		
		\begin{minipage}[b]{0.49\textwidth}
			\includegraphics[width=\textwidth]{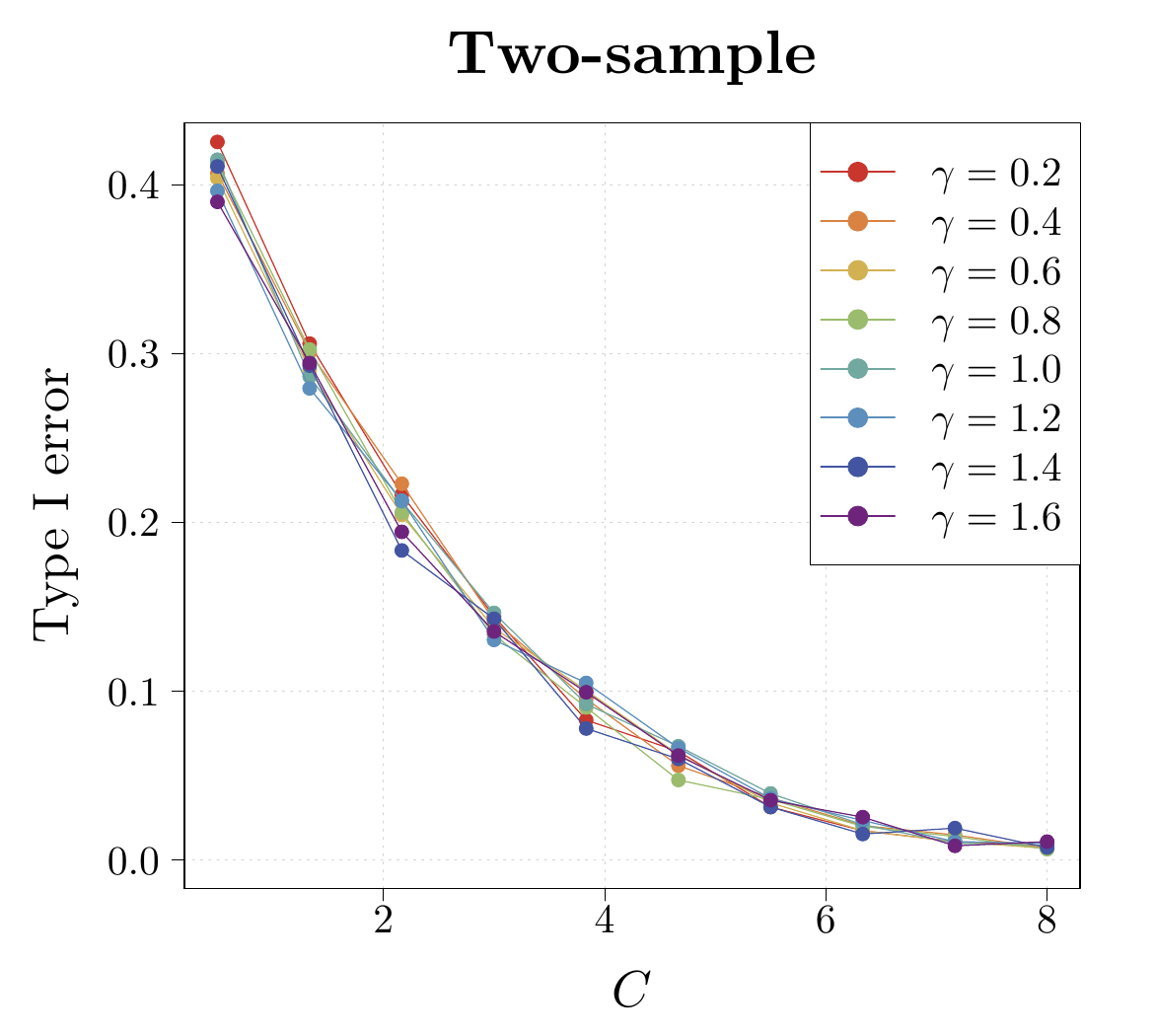}
		\end{minipage} 
		\hskip 0.0em
		\begin{minipage}[b]{0.49\textwidth}
			\includegraphics[width=\textwidth]{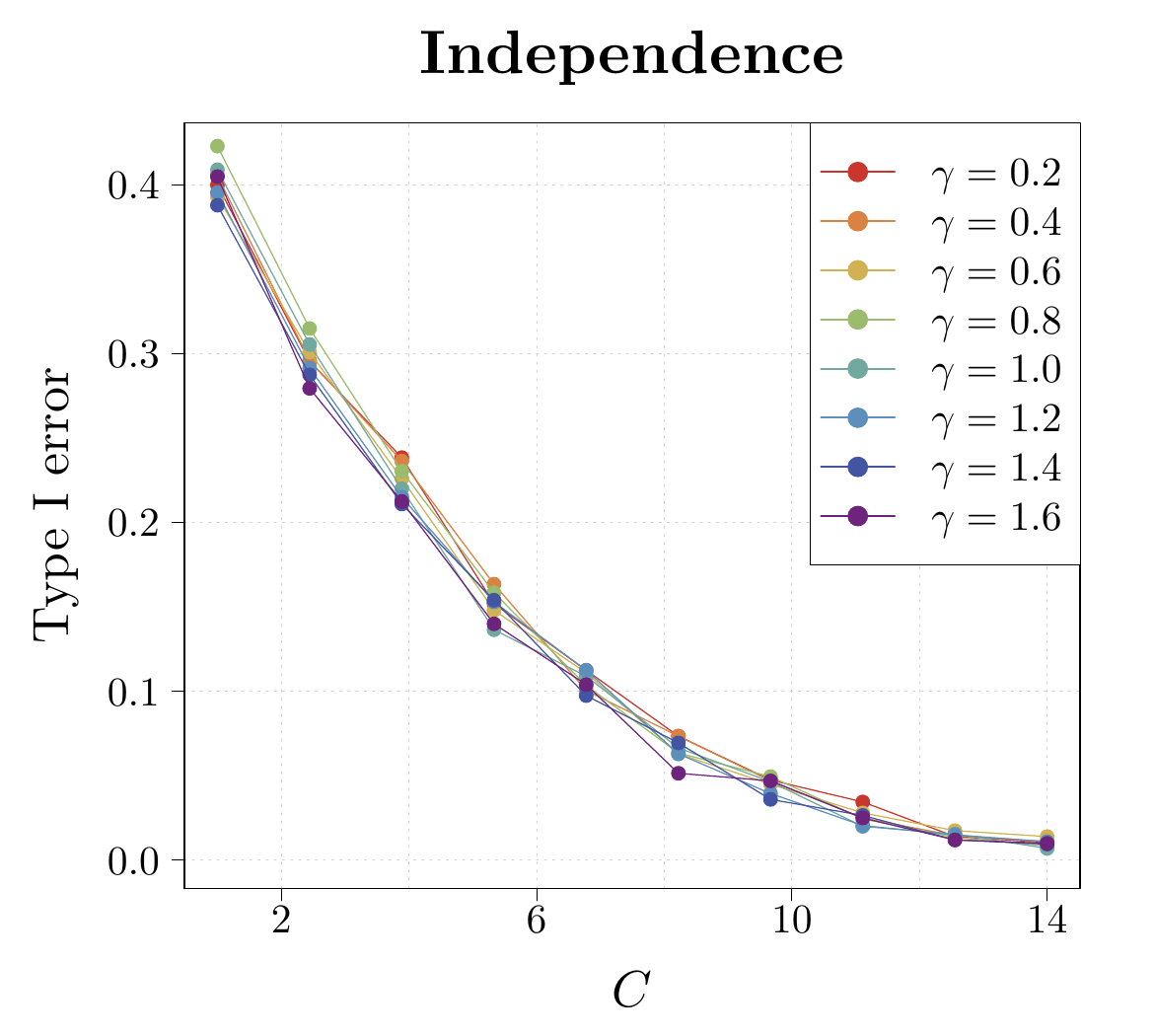}
		\end{minipage}
		\caption{\small Type I error rates of the tests based on concentration bounds by varying constant $C$ in their thresholds. Here we approximated the type I error rates via Monte-Carlo simulations under different power law distributions with parameter $\gamma$. The results show that the error rates vary considerably depending on the choice of $C$.} \label{Figure: loose thresholds}
	\end{center}
\end{figure}

This section provides empirical results to further justify the permutation approach. As emphasized before, the most significant feature of the permutation test is that it tightly controls the type I error rate for any sample size. This is in sharp contrast to non-asymptotic tests based on concentration bounds. The latter tests are typically conservative as they depend on a loose threshold. More seriously it is often the case that this threshold depends on a number of unspecified constants or even unknown parameters which raises the issue of practicality. In the first part of the simulation study, we demonstrate the sensitivity of the latter approach to the choice of constants in terms of type I error control. For this purpose, we focus on the problems of multinomial two-sample and independence testing and the simulation settings are described below. 
\begin{enumerate}
	\item \textbf{Two-sample testing.} We consider various power law multinomial distributions under the two-sample null hypothesis. Specifically the probability of each bin is defined to be $p_Y(k) = p_Z(k) \propto k^\gamma$ for $k \in \{1,\ldots,d\}$ and $\gamma \in \{0.2,\ldots,1.6\}$. We let the sample sizes be $n_1=n_2 = 50$ and the bin size be $d=50$. Following \cite{chan2014optimal} and \cite{diakonikolas2016new}, we use the threshold $C \|p_Y\|_2 n_1^{-1}$ for some constant $C$ and reject the null when $U_{n_1,n_2} > C \|p_Y\|_2 n_1^{-1}$ where $U_{n_1,n_2}$ is the $U$-statistic considered in Proposition~\ref{Proposition: Multinomial Two-Sample Testing}.
	\item \textbf{Independence testing.} We again consider power law multinomial distributions under the independence null hypothesis. In particular the probability of each bin is defined to be $p_{YZ}(k_1,k_2) = p_Y(k_1)p_Z(k_2) \propto k_1^\gamma k_2^\gamma$ for $k_1,k_2 \in \{1,\ldots,d\}$ and $\gamma \in \{0.2,\ldots,1.6\}$. We let the sample size be $n=100$ and the bin sizes be $d_1=d_2=20$. Similarly as before, we use the threshold $C \|p_Yp_Z\|_2 n^{-1}$ for some constant $C$ and reject the null when $U_n > C \|p_Yp_Z\|_2 n^{-1}$ where $U_n$ is the $U$-statistic considered in Proposition~\ref{Proposition: Multinomial Independence Testing}.
\end{enumerate}
The simulations were repeated 2000 times to approximate the type I error rate of the tests as a function of $C$. The results are presented in Figure~\ref{Figure: loose thresholds}. One notable aspect of the results is that, in both two-sample and independence cases, the error rates are fairly stable over different null scenarios for each fixed $C$. However these error rates vary a lot over different $C$, which clearly shows the sensitivity of the non-asymptotic approach to the choice of $C$. Furthermore it should be emphasized that both tests are not practical as they depend on unknown parameters $\|p_Y\|_2$ and $\|p_Yp_Z\|_2$, respectively.

\begin{figure}[t!]
	\begin{center}		
		\begin{minipage}[b]{0.32\textwidth}
			\includegraphics[width=\textwidth]{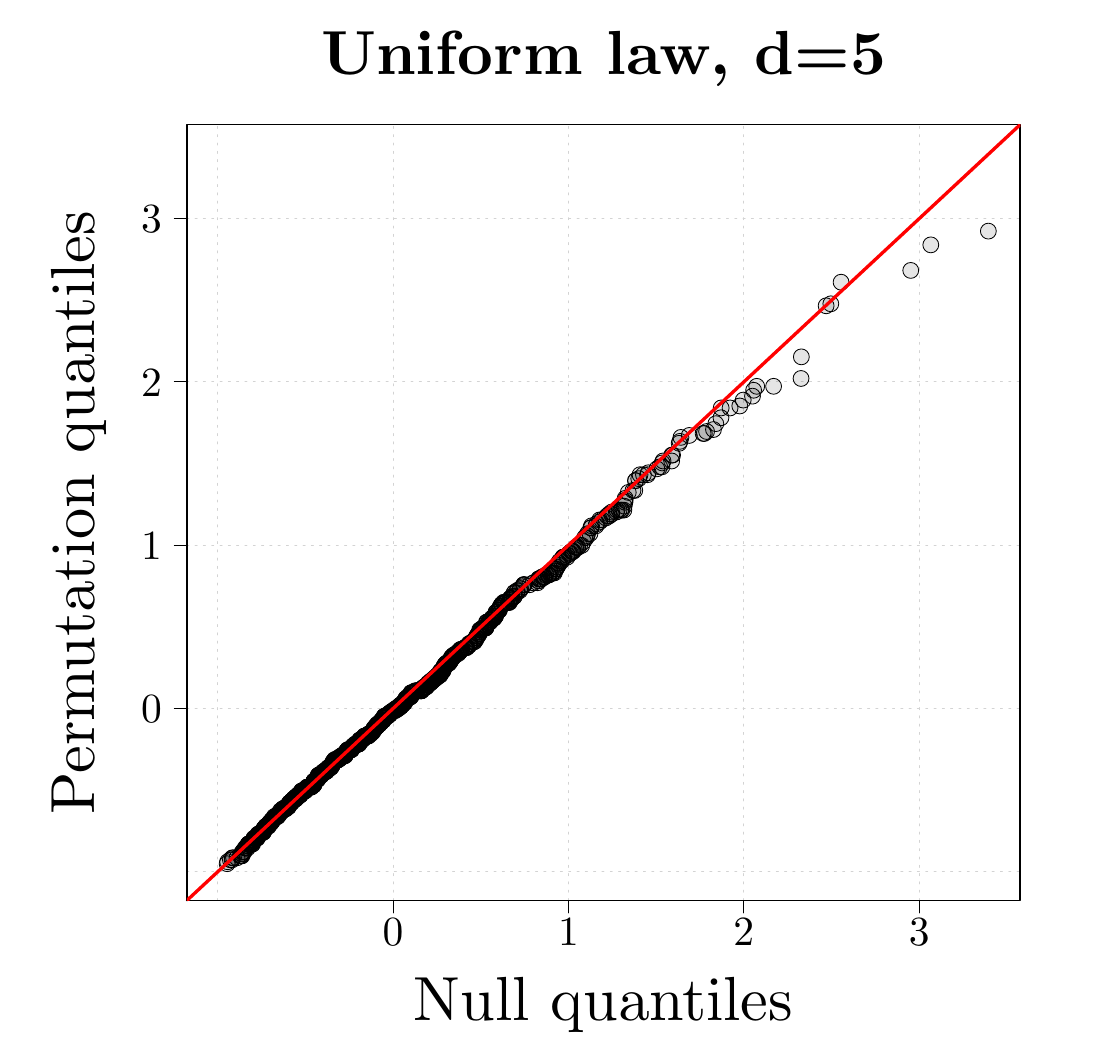}
		\end{minipage} 
		\begin{minipage}[b]{0.32\textwidth}
			\includegraphics[width=\textwidth]{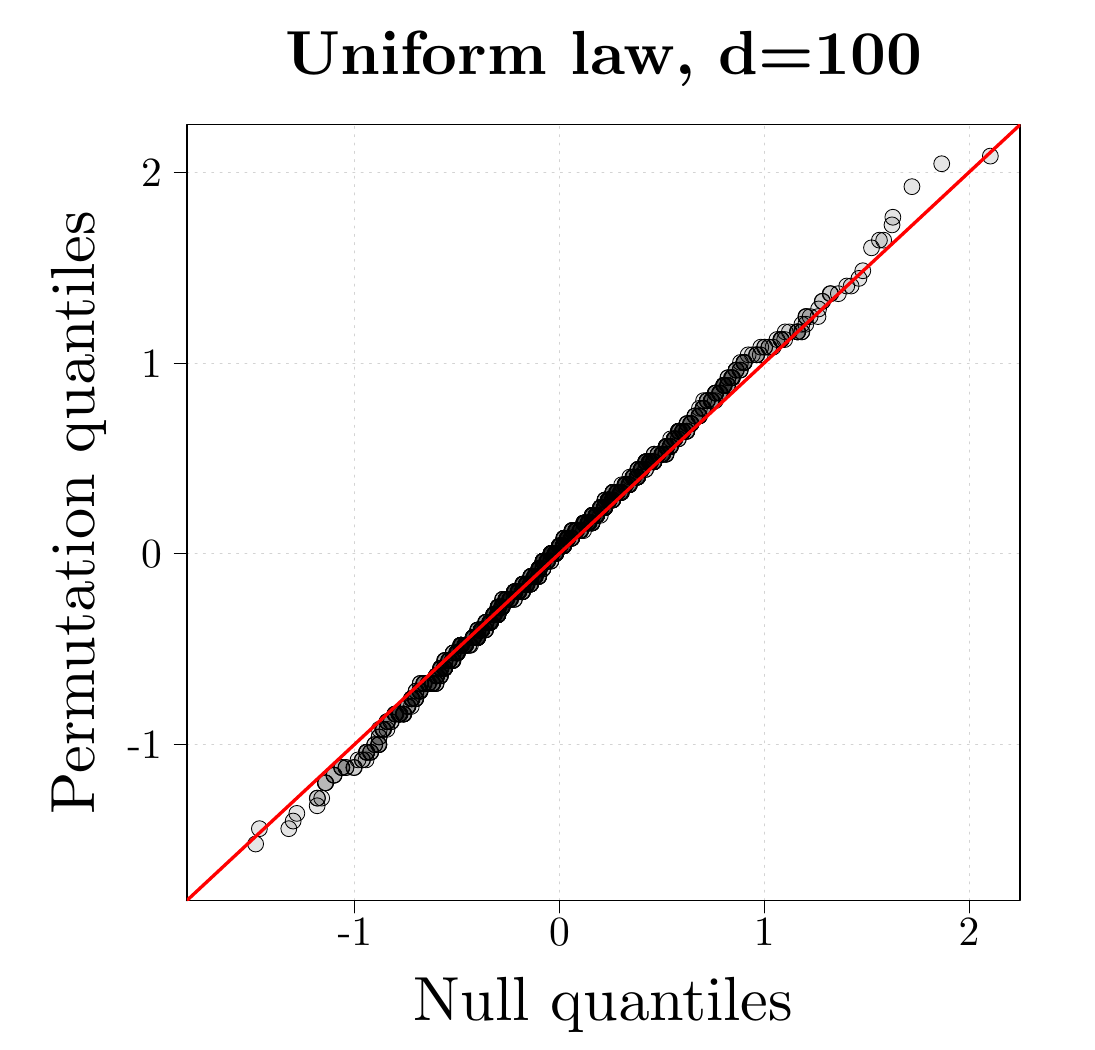}
		\end{minipage}
		\begin{minipage}[b]{0.32\textwidth}
			\includegraphics[width=\textwidth]{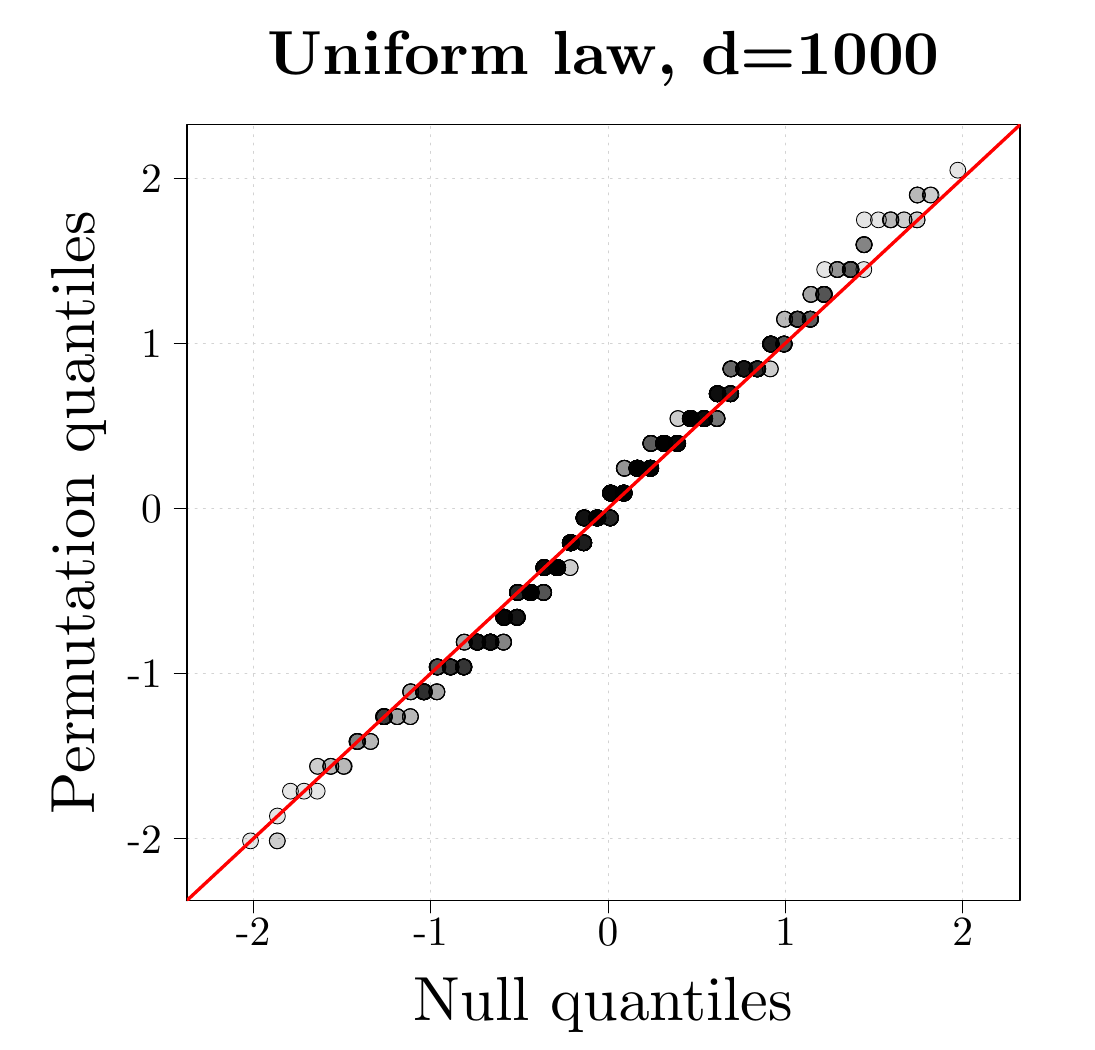}
		\end{minipage}
		\vskip 0.8em
		\begin{minipage}[b]{0.32\textwidth}
			\includegraphics[width=\textwidth]{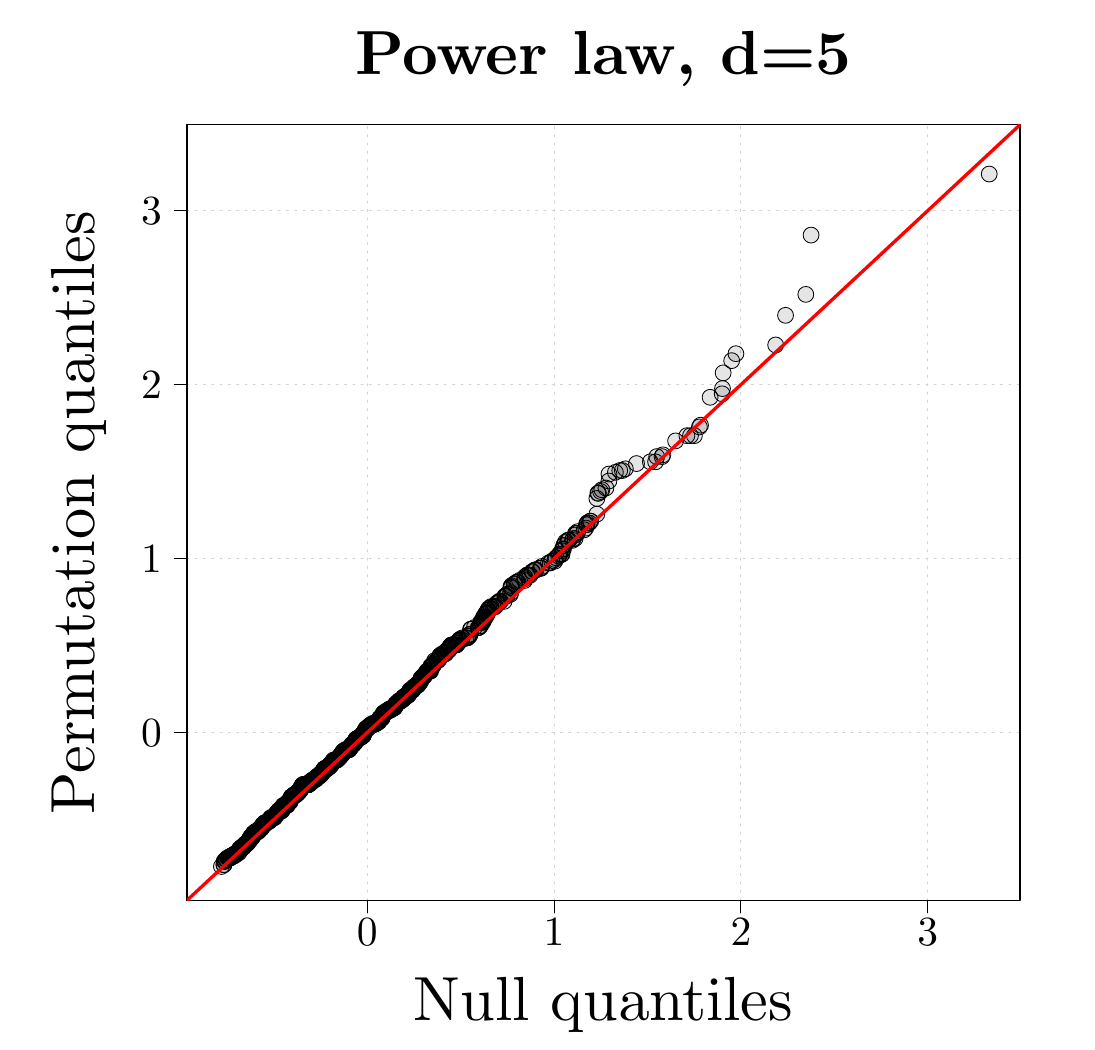}
		\end{minipage} 
		\begin{minipage}[b]{0.32\textwidth}
			\includegraphics[width=\textwidth]{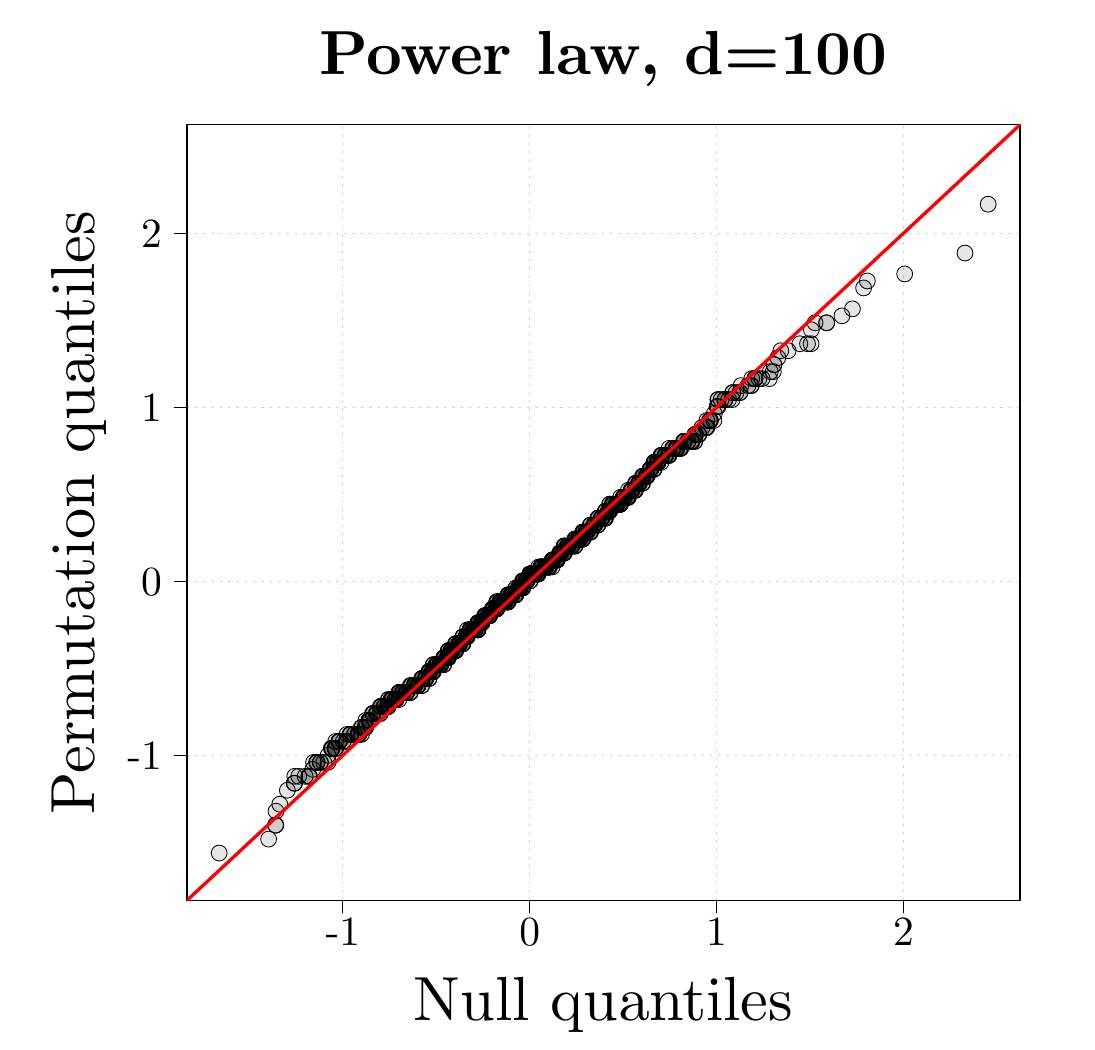}
		\end{minipage}
		\begin{minipage}[b]{0.32\textwidth}
			\includegraphics[width=\textwidth]{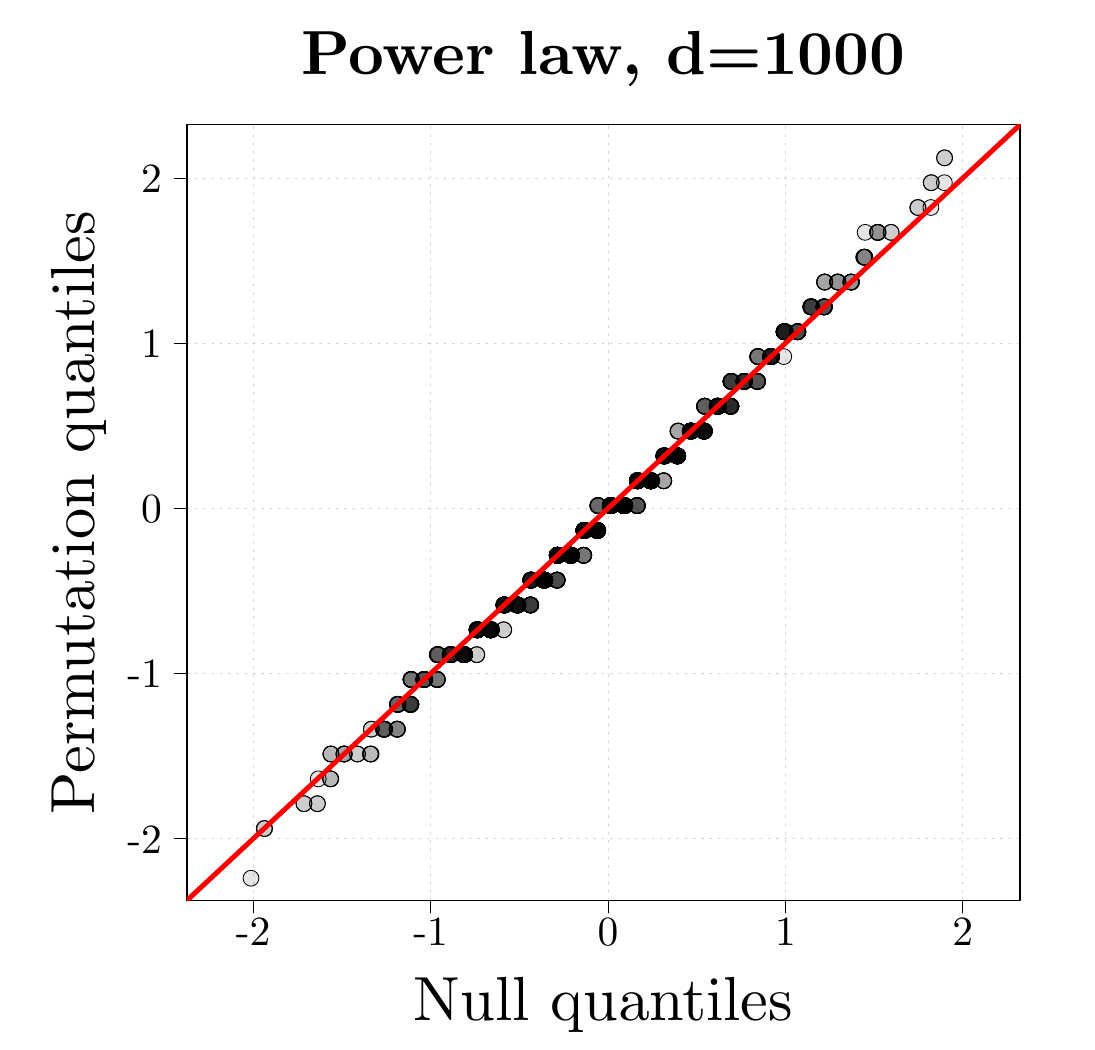}
		\end{minipage}
		\caption{\small Q-Q plots between the null distribution and the permutation distribution of the two-sample $U$-statistic. The quantiles of the two distributions approximately lie on the straight line $y=x$ in all cases, which demonstrates the similarity of the two distributions. Here we rescaled the test statistic by an appropriate constant for display purpose only.} \label{Figure: QQ plot}
	\end{center}
\end{figure}

It has been demonstrated by several authors \citep[e.g.][]{hoeffding1952large} that the permutation distribution of a test statistic mimics the underlying null distribution of the same test statistic in low-dimensional settings. In the next simulation, we provide empirical evidence that the same conclusion still holds in high-dimensional settings. This may further imply that the power of the permutation test approximates that of the theoretical test based on the null distribution of the test statistic. To illustrate, we focus on the two-sample $U$-statistic for multinomial testing in Proposition~\ref{Proposition: Multinomial Two-Sample Testing} and consider two different scenarios as follows.
\begin{enumerate}
	\item \textbf{Uniform law under the null.} We simulate $n_1=n_2=200$ samples from the uniform multinomial distributions under the null such that $p_Y(k) = p_Z(k) = 1/d$ for $k=1,\ldots,d$ where $d \in \{5, 100, 1000\}$. Conditional on these samples, we compute the permutation distribution of the test statistic. On the other hand, the null distribution of the test statistic is estimated based on $n_1=n_2=200$ samples from the uniform distribution by running a Monte-Carlo simulation with 2000 repetitions.
	\item \textbf{Power law under the alternative.} In order to argue that the power of the permutation test is similar to that of the theoretical test, we need to study the behavior of the permutation distribution under the alternative. For this reason, we simulate $n_1=200$ samples from the uniform distribution $p_Y(k)=1/d$ and $n_2=200$ samples from the power law distribution $p_Z(k) \propto k$ for $k = 1,\ldots, d$ where $d \in \{5,100,1000\}$. Conditional on these samples, we compute the permutation distribution of the test statistic. On the other hand, the null distribution of the test statistic is estimated based on $n_1=n_2=200$ samples from the mixture distribution $1/2 \times p_Y + 1/2 \times p_Z$ with $2000$ repetitions.   
\end{enumerate}
In the simulation study, due to the computational difficulty of considering all possible permutations, we approximated the original permutation distribution using the Monte-Carlo method. Nevertheless the difference between the original permutation distribution and its Monte-Carlo counterpart can be made arbitrary small uniformly over the entire real line, which can be shown by using Dvoretzky-Kiefer-Wolfowitz inequality \citep{dvoretzky1956asymptotic}. In our simulations, we randomly sampled $2000$ permutations from the entire permutations, based on which we computed the empirical distribution of the permuted test statistic. 

We recall from Figure~\ref{Figure: histograms} that the null distribution of the test statistic changes a lot depending on the size of $d$. In particular, it tends to be skewed to the right (similar to a $\chi^2$ distribution) when $d$ is small and tends to be symmetric (similar to a normal distribution) when $d$ is large. Also note that the null distribution tends to be more discrete when $d$ is large relative to the sample size. In Figure~\ref{Figure: QQ plot}, we present the Q-Q plots of the null and (approximate) permutation distributions. It is apparent from the figure that the quantiles of these two distributions approximately lie along the straight line $y=x$ in all the scenarios. In other words, the permutation distribution closely follows the null distribution, regardless of the size of $d$, from which we conjecture that the null distribution and the permutation distribution might have the same even in high-dimensional settings.

\section{Discussion} \label{Section: Discussion}
In this work we presented a general framework for analyzing the type II error rate of the permutation test based on the first two moments of the test statistic. We illustrated the utility of the proposed framework in the context of two-sample testing and independence testing in both discrete and continuous cases. In particular, we introduced the permutation tests based on degenerate $U$-statistics and explored their minimax optimality for multinomial testing as well as density testing. To improve a polynomial dependency on the nominal level $\alpha$, we developed exponential concentration inequalities for permuted $U$-statistics based on an idea that links permutations to i.i.d.~Bernoulli random variables. The utility of the exponential bounds was highlighted by introducing adaptive tests to unknown parameters and also providing a concentration bound for Rademacher chaos under sampling without replacement.

Our work motivates several lines of future directions. First, while this paper considered the problem of 
unconditional independence testing, it would be interesting to extend our results to the problem of conditional independence testing. When the conditional variable is discrete, one can apply unconditional independence tests within categories and combine them, in a suitable way, to test for conditional independence. When the conditional variable is continuous, however, this strategy does not work and this has led several authors to use ``local permutation'' heuristics~\citep[see for instance][]{zhang2014permutation,fukumizu2008kernel,neykov2019minimax}. In contrast to the two-sample and independence testing problems 
we have considered in our paper, even justifying the type I error control of these methods is not straightforward. Second, based on the coupling idea in Section~\ref{Section: Combinatorial concentration inequalities}, further work can be done to develop combinatorial concentration inequalities for other statistics. It would also be interesting to see whether one can obtain tighter concentration bounds, especially for $U_{n}^\pi$ in (\ref{Eq: permuted U-statistic}). 
Finally, identifying settings in which we can improve the dependence on the type II error rate $\beta$ for the two-moment method is another interesting direction for future research. 

\section{Acknowledgments}
The authors are grateful to Arthur Gretton for pointing us to the work of \cite{meynaoui2019aggregated} and for helpful discussions. We would like thank Antonin Schrab for pointing out several typos in the original version of this manuscript. We thank Richard J. Samworth, Thomas Berrett, Tracy Ke, Bodhisattva Sen for their helpful comments. Finally, we thank three anonymous reviewers and an associate editor for their constructive comments. This work was partially supported by the NSF grant DMS-17130003. 
 
\bibliographystyle{apalike}
\bibliography{reference}

\appendix

\section{Overview of Appendix}
In this supplementary material, we provide some additional results and the technical proofs omitted in the main text. The remainder of this material is organized as follows.  
\begin{itemize}
	\item In Appendix~\ref{Section: Exponential inequalities for permuted linear statistics}, we develop exponential inequalities for permuted linear statistics, building on the concept of negative association. 
	\item In Appendix~\ref{Section: Improved version of Theorem}, we provide the result that improves Theorem~\ref{Theorem: Two-Sample U-statistic} based on the exponential bound in Theorem~\ref{Theorem: Two-Sample Concentration} with an extra assumption that $n_1 \asymp n_2$. 
	\item In Appendix~\ref{Section: Monte Carlo-based permutation tests}, we show that Monte Carlo-based permutation tests can achieve the same error guarantee as the corresponding full permutation tests.
	\item The proof of Lemma~\ref{Lemma: Two Moments Method} on the two moments method is provided in Appendix~\ref{Section: Proof of Lemma: Two Moments Method}.
	\item The proofs of the results on two-sample testing in Section~\ref{Section: The two moments method for two-sample testing} are presented in Appendix~\ref{Section: Proof of Theorem: Two-Sample U-statistic}, \ref{Section: Proof of Proposition: Multinomial Two-Sample Testing}, \ref{Section: Proposition: Minimum Separation for Two-Sample Multinomial Testing} and \ref{Section: Proposition: Two-Sample Testing for Holder Densities}.
	\item The proofs of the results on independence testing in Section~\ref{Section: The two moments method for independence testing} are presented in Appendix~\ref{Section: Theorem: U-statistic for independence testing}, \ref{Section: Proposition: Multinomial Independence Testing}, \ref{Section: Proposition: Minimum Separation for Multinomial Independence Testing}, \ref{Section: Proposition: Independence Testing for Holder Densities} and \ref{Section: Proposition: Minimum Separation for Independence Testing for Holder Densities}.
	\item The proofs of the results on combinatorial concentration inequalities in Section~\ref{Section: Combinatorial concentration inequalities} are presented in Appendix~\ref{Section: Corollary: Dependent Rademacher Chaos} and \ref{Section: Theorem: Concentration inequality for independence U-statistic II}.
	\item The proofs of the results on adaptive tests in Section~\ref{Section: Adaptive tests to unknown parameters} are presented in Appendix~\ref{Section: Proposition: Adaptive two-sample test} and \ref{Section: Proposition: Adaptive independence test}.
	\item The proofs of the results on multinomial tests and Gaussian kernel tests in Section~\ref{Section: Further applications} are presented in Appendix~\ref{Section: Theorem: Two-Sample testing under Poisson sampling}, \ref{Section: Proof of Proposition: Multinomial L1 Testing}, \ref{Section: Proof of Proposition: Multinomial independence testing in L1 distance} and \ref{Section: Proof of Proposition: Gaussian MMD}.
	\item Lastly, the proof of the results in Appendix~\ref{Section: Monte Carlo-based permutation tests} can be found in Appendix~\ref{Section: Proof of Proposition: Two-moment method for MC-based tests}. 
\end{itemize}

\section{Exponential inequalities for permuted linear statistics} \label{Section: Exponential inequalities for permuted linear statistics}

Suppose that $\mathcal{X}_n=\{(Y_1,Z_1),\ldots,(Y_n,Z_n)\}$ is a set of bivariate random variables where $Y_i \in \mathbb{R}$ and $Z_i \in \mathbb{R}$. Following the convention, let us write the sample means of $Y$ and $Z$ by $\overline{Y} := n^{-1} \sum_{i=1}^n Y_i$ and $\overline{Z} := n^{-1} \sum_{i=1}^n Z_i$, respectively. The sample covariance, which measures a linear relationship between $Y$ and $Z$, is given by
\begin{align*}
L_n := \frac{1}{n} \sum_{i=1}^n (Y_i - \overline{Y})(Z_{i} - \overline{Z}).
\end{align*}
We also call $L_n$ as a linear statistic as opposed to quadratic statistics or degenerate $U$-statistics considered in the main text. Let us denote the permuted linear statistic, associated with a permutation $\pi$ of $\{1,\ldots,n\}$, by
\begin{align*}
L_n^\pi = \frac{1}{n} \sum_{i=1}^n (Y_i - \overline{Y})(Z_{\pi_i} - \overline{Z}).
\end{align*}
In this section, we provide two exponential concentration bounds for $L_n^{\pi}$ conditional on $\mathcal{X}_n$; namely Hoeffding-type inequality (Proposition~\ref{Proposition: Hoeffding bound}) and Bernstein-type inequality (Proposition~\ref{Proposition: Bernstein bound}). These results have potential applications in studying the power of the permutation test based on $L_n$ and also concentration inequalities for sampling without replacement. We describe the second application in more detail in Appendix~\ref{Section: Concentration inequalities for sampling without replacement} after we develop the results.

\vskip .8em

\noindent \textbf{Related work and negative association.} We should note that the same problem has been considered by several authors using Stein's method \citep{chatterjee2007stein}, a martingale method \citep[Chapter 4.2 of][]{bercu2015concentration} and Talagrand's inequality \citep{albert2019concentration}. In fact they consider a more general linear statistic which has the form of $\sum_{i=1}^n d_{i,\pi_i}$ where $\{d_{i,j}\}_{i,j=1}^n$ is an arbitrary bivariate sequence. Thus their statistic includes $L_n$ as a special case by letting $d_{i,\pi_i} =(Y_i - \overline{Y})(Z_{\pi_i} - \overline{Z})$. However their proofs are quite involved at the expense of being more general. Here we provide a much simpler proof with sharper constant factors by taking advantage of the decomposability of $d_{i,j}$. To this end, we utilize the concept of negative association \citep[e.g.][]{joag1983negative,dubhashi1998balls}, defined as follows. 

\begin{definition}[Negative association]
	Random variables $X_1,\ldots,X_n$ are negatively associated (NA) if for every two disjoint index sets $\mathcal{I},\mathcal{J} \subseteq \{1,\ldots,n\}$,
	\begin{align*}
	\mE[f(X_i, i \in \mathcal{I}) g(X_j, j \in \mathcal{J})] \leq \mE [f(X_i, i \in \mathcal{I})] \mE[g(X_j, j \in \mathcal{J})]
	\end{align*}
	for all functions $f: \mathbb{R}^{|\mathcal{I}|} \mapsto \mathbb{R}$ and $g : \mathbb{R}^{|\mathcal{J}|} \mapsto \mathbb{R}$ that are both non-decreasing or both non-increasing.  
\end{definition}
Let us state several useful facts about negatively associated random variables that we shall leverage to prove the main results of this section. The proofs of the given facts can be found in \cite{joag1983negative} and \cite{dubhashi1998balls}. 

\vskip 0.8em 

\begin{itemize}
	\item \textbf{Fact 1.} Let $\{ x_1,\ldots,x_n \}$ be a set of $n$ real values. Suppose that $\{X_1,\ldots,X_n\}$ are random variables with the probability such that 
	\begin{align*}
	\mP(X_1= x_{\pi_1}, \ldots, X_n = x_{\pi_n}) = \frac{1}{n!} \quad \text{for any permutation $\pi$ of $\{1,\ldots,n\}$.}
	\end{align*}
	Then $\{X_1,\ldots,X_n\}$ are negatively associated. 
	\item  \textbf{Fact 2.} Let $\{X_1,\ldots,X_n\}$ be negatively associated. Let $\mathcal{I}_1,\ldots, \mathcal{I}_k \subseteq \{1,\ldots,n \}$ be disjoint index sets, for some positive integer $k$. For $j \in \{1,\ldots,n\}$, let $h_j : \mathbb{R}^{|\mathcal{I}_k|} \mapsto \mathbb{R}$ be functions that are all non-decreasing or all non-increasing and define $Y_j = h_j(X_i, i \in \mathcal{I}_j)$. Then $\{Y_1,\ldots, Y_k\}$ are also negatively associated.   
	\item \textbf{Fact 3.} Let $\{X_1,\ldots,X_n\}$ be negatively associated. Then for any non-decreasing functions $f_i$, $i \in \{1,\ldots,n\}$, we have that
	\begin{align} \label{Eq: NA property}
	\mE \bigg[ \prod_{i=1}^n f_i(X_i) \bigg] \leq \prod_{i=1}^n \mE [f_i(X_i)].
	\end{align}
\end{itemize}

\noindent \textbf{Description of the main idea.} Notice that $L_n^\pi$ is a function of non-i.i.d.~random variables for which standard techniques relying on i.i.d.~assumptions do not work directly. We avoid this difficulty by connecting $L_n^\pi$ with negatively associated random variables and then applying Chernoff bound combined with the inequality~(\ref{Eq: NA property}). The details are as follows. 
For notational simplicity, let us denote 
\begin{align*}
& \{a_{1}, \ldots, a_{n}\} = \{Y_{1} - \overline{Y}, \ldots, Y_{n} - \overline{Y} \} \quad \text{and} \\[.5em]
& \{b_{\pi_1}, \ldots, b_{\pi_n}\} = \{Z_{\pi_1} - \overline{Z}, \ldots, Z_{\pi_n} - \overline{Z} \}.
\end{align*}
To proceed, we make several important observations. 
\begin{itemize}
	\item \textbf{Observation 1.} First, since $\{ b_{\pi_1}, \ldots, b_{\pi_n} \}$ has a permutation distribution, we can use Fact 1 and conclude that $\{ b_{\pi_1}, \ldots, b_{\pi_n} \}$ are negatively associated. 
	\item \textbf{Observation 2.} Second, let $\mathcal{I}_+$ be the set of indices such that $a_i > 0$ and similarly $\mathcal{I}_-$ be the set of indices such that $a_i < 0$. Since $h_i(X_i, i \in \mathcal{I}_+) = a_i X_i$ is non-decreasing function and $h_i(X_i, i \in \mathcal{I}_-) = a_i X_i$ is non-increasing functions, it can be seen that 
	$\{a_ib_{\pi_i} \}_{i \in \mathcal{I}_+}$ and $\{a_ib_{\pi_i} \}_{i \in \mathcal{I}_-}$ are negatively associated by Fact 2. Using this notation, the linear statistic can be written as 
	\begin{align*}
	L_n^\pi = \frac{1}{n} \sum_{i \in \mathcal{I}_+} a_i b_{\pi_i} + \frac{1}{n} \sum_{i \in \mathcal{I}_-} a_i b_{\pi_i}. 
	\end{align*}
	It can be easily seen that $\mE_{\pi}[b_{\pi_i}|\mathcal{X}_n] =0$ for each $i$ and thus $\mE_{\pi}[L_n^\pi | \mathcal{X}_n] = 0$ by linearity of expectation. Hence, for $\lambda > 0$, applying the Chernoff bound yields
	\begin{align*}
	& \mP_{\pi}(L_n^\pi \geq t | \mathcal{X}_n) \\[.5em]
	\leq ~ & e^{-\lambda t}  \mE_{\pi} \Bigg[ \exp \Bigg( \lambda n^{-1} \sum_{i \in \mathcal{I}_+} a_i b_{\pi_i} + \lambda n^{-1} \sum_{i \in \mathcal{I}_-} a_i b_{\pi_i} \Bigg) \Bigg| \mathcal{X}_n  \Bigg] \\[.5em]
	\leq ~ & \frac{e^{-\lambda t}}{2}  \mE_{\pi} \Bigg[ \exp \Bigg( 2\lambda n^{-1} \sum_{i \in \mathcal{I}_+} a_i b_{\pi_i} \Bigg)  \Bigg| \mathcal{X}_n \Bigg] + \frac{e^{-\lambda t}}{2}  \mE_{\pi} \Bigg[ \exp \Bigg( 2\lambda n^{-1} \sum_{i \in \mathcal{I}_-} a_i b_{\pi_i} \Bigg)  \Bigg| \mathcal{X}_n \Bigg] \\[.5em]
	:= ~ & (I) + (II),
	\end{align*}
	where the last inequality uses the elementary inequality $xy \leq x^2/2 + y^2/2$. 
	\item \textbf{Observation 3.} Third, based on fact that $\{a_ib_{\pi_i} \}_{i \in \mathcal{I}_+}$ and $\{a_ib_{\pi_i} \}_{i \in \mathcal{I}_-}$ are negatively associated, we may apply Fact 3 to have that 
	\begin{equation}
	\begin{aligned} \label{Eq: upper bounds for I and II}
	& (I) ~\leq~  \frac{e^{-\lambda t}}{2} \prod_{i \in \mathcal{I}_+} \mE_{\widetilde{b}} \big[ \exp \big( 2\lambda n^{-1} a_i \widetilde{b}_{i} \big)  \big| \mathcal{X}_n \big]  =  \frac{e^{-\lambda t}}{2} \prod_{i = 1}^n \mE_{\widetilde{b}} \big[ \exp \big( 2\lambda n^{-1} a_i^+ \widetilde{b}_{i} \big)  \big| \mathcal{X}_n \big]   \quad \text{and} \\[.5em]
	& (II) ~\leq~  \frac{e^{-\lambda t}}{2} \prod_{i \in \mathcal{I}_-} \mE_{\widetilde{b}} \big[ \exp \big( 2\lambda n^{-1} a_i \widetilde{b}_{i} \big)  \big| \mathcal{X}_n \big] = \frac{e^{-\lambda t}}{2} \prod_{i = 1}^n \mE_{\widetilde{b}} \big[ \exp \big( -2\lambda n^{-1} a_i^- \widetilde{b}_{i} \big)  \big| \mathcal{X}_n \big],
	\end{aligned}
	\end{equation}
	where $\widetilde{b}_1,\ldots,\widetilde{b}_n$ are i.i.d.~random variables uniformly distributed over $\{b_1,\ldots,b_n\}$. Here $a_i^+$ and $a_i^-$ represent $a_i^+ = a_i \ind(a_i \geq 0)$ and $a_i^- = - a_i \ind(a_i \leq 0)$ respectively. 
\end{itemize}
With these upper bounds for $(I)$ and $(II)$ in place, we are now ready to present the main results of this section. The first one is a Hoeffding-type bound which provides a sharper constant factor than  \cite{duembgen1998symmetrization}.

\begin{proposition}[Hoeffding-type bound] \label{Proposition: Hoeffding bound}
	Let us define $a_{\text{\emph{range}}} := Y_n - Y_1$ and $b_{\text{\emph{range}}} := Z_n - Z_1$. Then
	\begin{align*}
	\mP_{\pi} (L_n^\pi \geq t | \mathcal{X}_n)  \leq \exp \Bigg[ - \max \Bigg\{ \frac{n^2t^2}{a_{\text{\emph{range}}}^2\sum_{i=1}^n b_i^2}, ~  \frac{n^2t^2}{b_{\text{\emph{range}}}^2\sum_{i=1}^n a_i^2} \Bigg\} \Bigg].
	\end{align*}
	\vskip .5em
	\begin{proof}
		The proof directly follows by applying Hoeffding's lemma \citep{hoeffding1963probability}, which states that when $Z$ has zero mean and $a \leq Z \leq b$, 
		\begin{align*}
		\mE[e^{\lambda Z}] \leq e^{\lambda^2 (b-a)^2/8}.
		\end{align*}
		Notice that Hoeffding's lemma yields
		\begin{align*}
		 & \prod_{i = 1}^n \mE_{\widetilde{b}} \big[ \exp \big( 2\lambda n^{-1} a_i^+ \widetilde{b}_{i} \big)  \big| \mathcal{X}_n \big] \leq \exp \Bigg\{ \frac{\lambda^2b_{\text{range}}^2}{2n^2} \sum_{i=1}^n (a_i^+)^2 \Bigg\} \leq \exp \Bigg\{ \frac{\lambda^2b_{\text{range}}^2}{2n^2} \sum_{i=1}^n a_i^2 \Bigg\} \quad \text{and} \\[.5em]
		 & \prod_{i = 1}^n \mE_{\widetilde{b}} \big[ \exp \big( 2\lambda n^{-1} a_i^- \widetilde{b}_{i} \big)  \big| \mathcal{X}_n \big] \leq \exp \Bigg\{ \frac{\lambda^2b_{\text{range}}^2}{2n^2} \sum_{i=1}^n (a_i^-)^2 \Bigg\} \leq \exp \Bigg\{ \frac{\lambda^2b_{\text{range}}^2}{2n^2} \sum_{i=1}^n a_i^2 \Bigg\}. 
		\end{align*}
		Thus combining the above with the upper bounds for $(I)$ and $(II)$ in (\ref{Eq: upper bounds for I and II}) yields
		\begin{align*}
		\mP_{\pi} (L_n^\pi \geq t | \mathcal{X}_n) \leq \exp\Bigg\{ -\lambda t + \frac{\lambda^2b_{\text{range}}^2}{2n^2} \sum_{i=1}^n a_i^2  \Bigg\}.
		\end{align*}
		By optimizing over $\lambda$ on the right-hand side, we obtain that
		\begin{align} \label{Eq: Bound 1}
		\mP_{\pi} (L_n^\pi \geq t | \mathcal{X}_n) \leq \exp\Bigg\{ - \frac{n^2t^2}{b_{\text{range}}^2\sum_{i=1}^n a_i^2} \Bigg\}.
		\end{align}
		Since $\sum_{i=1}^n a_{\pi_i} b_i$ and $\sum_{i=1}^n a_i b_{\pi_i}$ have the same permutation distribution, it also holds that 		\begin{align} \label{Eq: Bound 2}
		\mP_{\pi} (L_n^\pi \geq t | \mathcal{X}_n) \leq \exp\Bigg\{ - \frac{n^2t^2}{a_{\text{range}}^2\sum_{i=1}^n b_i^2} \Bigg\}.
		\end{align}
		Then putting together these two bounds (\ref{Eq: Bound 1}) and (\ref{Eq: Bound 2}) gives the desired result. 
	\end{proof}
\end{proposition}

Note that Proposition~\ref{Proposition: Hoeffding bound} depends on the variance of either $\{a_i\}_{i=1}^n$ or $\{b_i\}_{i=1}^n$. In the next proposition, we provide a Bernstein-type bound which depends on the variance of the bivariate sequence $\{a_ib_j\}_{i,j=1}^n$. Similar results can be found in  \cite{bercu2015concentration} and \cite{albert2019concentration} but in terms of constants, the bound below is much shaper than the previous ones.

\begin{proposition}[Bernstein-type bound] \label{Proposition: Bernstein bound}
	Based on the same notation in Proposition~\ref{Proposition: Hoeffding bound}, a Bernstein-type bound is provided by
	\begin{align*}
	\mP_{\pi}(L_n^\pi \geq t | \mathcal{X}_n) \leq \exp \Bigg\{ - \frac{nt^2}{2n^{-2} \sum_{i,j=1}^n a_i^2 b_j^2 + \frac{2}{3} t \max_{1 \leq i,j \leq n} |a_ib_j|}  \Bigg\}.
	\end{align*}
	\begin{proof}
		Once we have the upper bounds for $(I)$ and $(II)$ in (\ref{Eq: upper bounds for I and II}), the remainder of the proof is routine. First it is straightforward to verify that 
		for $|Z| \leq c$, $\mathbb{E}[Z]=0$ and $\mathbb{E}[Z^2] = \sigma^2$, we have that 
		\begin{align*}
		\mathbb{E}[e^{\lambda Z}] = 1 + \sum_{k=2}^\infty \frac{\mathbb{E}[(\lambda Z)^k]}{k!} \leq 1 + \frac{\sigma^2}{c^2} \sum_{k=2}^\infty \frac{\lambda^k c^{k}}{k!} \leq \exp \bigg\{ \frac{\sigma^2}{c^2} \left( e^{\lambda c} - 1 - \lambda c \right) \bigg\}.
		\end{align*}
		Let us write $\widehat{\sigma}_i^2 = n^{-3} a_i^2\sum_{j=1}^n b_i^2$ and $M = n^{-1} \max_{1\leq i,j \leq n} |a_ib_j|$. Then based on the above inequality, we can obtain that
		\begin{align*}
		& \prod_{i = 1}^n \mE_{\widetilde{b}} \big[ \exp \big( 2\lambda n^{-1} a_i^+ \widetilde{b}_{i} \big)  \big| \mathcal{X}_n \big]  \leq \exp \Bigg\{ \frac{\sum_{i=1}^n \widehat{\sigma}_i^2}{M^2} \left( e^{\lambda M} - 1 - \lambda M \right)  \Bigg\} \quad \text{and} \\[.5em]
		& \prod_{i = 1}^n \mE_{\widetilde{b}} \big[ \exp \big( 2\lambda n^{-1} a_i^- \widetilde{b}_{i} \big)  \big| \mathcal{X}_n \big]  \leq \exp \Bigg\{ \frac{\sum_{i=1}^n \widehat{\sigma}_i^2}{M^2} \left( e^{\lambda M} - 1 - \lambda M \right)  \Bigg\}. 
		\end{align*}
		Combining these two upper bounds with the result in (\ref{Eq: upper bounds for I and II}) yields
		\begin{align*}
		\mathbb{P}_{\pi} \left( L_{n}^{\pi} \geq t | \mathcal{X}_n \right) \leq e^{-\lambda t} \exp \Bigg\{ \frac{\sum_{i=1}^n \widehat{\sigma}_i^2}{M^2} \left( e^{\lambda M} - 1 - \lambda M \right)  \Bigg\}.
		\end{align*}
		By optimizing the right-hand side in terms of $\lambda$, we obtain a Bennett-type inequality
		\begin{align*}
		\mathbb{P}_{\pi} \left( L_{n}^{\pi} \geq t | \mathcal{X}_n \right)  \leq \exp \Bigg\{ - \frac{\sum_{i=1}^n \widehat{\sigma}_i^2}{M^2} h\left( \frac{tM}{\sum_{i=1}^n \widehat{\sigma}_i^2} \right)  \Bigg\},
		\end{align*}
		where $h(x) = (1+x) \log(1+x) -x$. Then the result follows by noting that $h(x) \geq x^2 /(2 + 2x/3)$. 
	\end{proof}
\end{proposition}

In the next subsection, we apply our results to derive concentration inequalities for sampling without replacement.

\subsection{Concentration inequalities for sampling without replacement} \label{Section: Concentration inequalities for sampling without replacement}
To establish the explicit connection to sampling without replacement, we focus on the case where $Z_i$ is binary, say $Z_i \in \{-a, a\}$. Then the linear statistic $L_n$ is related to the unscaled two-sample $t$-statistic. More specifically, let us write $n_1 = \sum_{i=1}^n \ind(Z_i = a)$ and $n_2= \sum_{i=1}^n \ind(Z_i = -a)$. Additionally we use the notation $\overline{Y}_1 = n_1^{-1} \sum_{i=1}^{n} Y_i \ind(Z_i = a)$ and $\overline{Y}_2 = n_2^{-1}\sum_{i=1}^n Y_i \ind(Z_i = -a)$. Then some algebra shows that the sample covariance $L_n$ is exactly the form of 
\begin{align*}
L_n = 2a \frac{n_1n_2}{n^2} (\overline{Y}_1 - \overline{Y}_2).
\end{align*}
Without loss of generality, we assume $a=1$, i.e. $Z_i \in \{-1,1\}$. Then Proposition~\ref{Proposition: Hoeffding bound} gives a concentration inequality for the unscaled $t$-statistic as
\begin{align*}
\mathbb{P}_\pi \Bigg\{ \frac{2n_1n_2}{n^2}  \left( \overline{Y}_{1,\pi} - \overline{Y}_{2,\pi} \right)  \geq t \Bigg| \mathcal{X}_n \Bigg\}  ~= ~ & \mathbb{P}_\pi \Bigg\{  \left( \overline{Y}_{1,\pi} - \overline{Y}_{2,\pi} \right)  \geq \frac{tn^2}{2n_1n_2}  \Bigg| \mathcal{X}_n  \Bigg\} \\[.5em]
\leq ~ & \exp \Bigg\{ - \frac{n^2 t^2}{4 \sum_{i=1}^n (Y_i - \overline{Y})^2} \Bigg\},
\end{align*}
where $\overline{Y}_{1,\pi} = n_1^{-1} \sum_{i=1}^{n} Y_{\pi_i} \ind(Z_i = 1)$ and $\overline{Y}_{2,\pi} = n_2^{-1}\sum_{i=1}^n Y_{\pi_i} \ind(Z_i = -1)$. This implies that
\begin{align*}
\mathbb{P}_\pi \left( \overline{Y}_{1,\pi} - \overline{Y}_{2,\pi} \geq t  \big| \mathcal{X}_n \right) \leq \exp \left( - \frac{n_1^2 n_2^2 t^2}{n^3 \widehat{\sigma}_{\text{lin}}^2} \right),
\end{align*}
where $\widehat{\sigma}_{\text{lin}}^2 = n^{-1} \sum_{i=1}^n (Z_i - \overline{Z})^2$. By symmetry, it also holds that 
\begin{align*}
\mathbb{P}_\pi \left( |\overline{Y}_{1,\pi} - \overline{Y}_{2,\pi}| \geq t  \big| \mathcal{X}_n \right) \leq 2\exp \left( - \frac{n_1^2 n_2^2 t^2}{n^3 \widehat{\sigma}_{\text{lin}}^2} \right).
\end{align*}
Let us denote the sample mean of the entire samples by $\overline{Y} = n^{-1} \sum_{i=1}^n Y_i$. Then using the exact relationship
\begin{align} \label{Eq: exact relationship}
|\overline{Y}_{1,\pi} - \overline{Y}_{2,\pi}| = \frac{n}{n_2} |\overline{Y}_{1,\pi} - \overline{Y}|,
\end{align}
the above inequality is equivalent to 
\begin{align} \label{Eq: Hoeffding type tail bound}
\mathbb{P}_\pi \left( |\overline{Y}_{1,\pi} - \overline{Y}| \geq t  \big| \mathcal{X}_n \right) \leq 2\exp \left( - \frac{n_1^2 t^2}{n \widehat{\sigma}_{\text{lin}}^2} \right).
\end{align}
Notice that $\overline{Y}_{1,\pi}$ is the sample mean of $n_1$ observations sampled without replacement from $\{Y_1,\ldots, Y_n\}$. This implies that the permutation law of the sample mean is equivalent to the probability law under sampling without replacement. The same result (including the constant factor) exists in \cite{massart1986rates} (see Lemma 3.1 therein). However, the result given there only holds when $n=n_1 \times m$ where $m$ is a positive integer whereas our result does not require such restriction. 

\vskip 1em

\noindent \textbf{An improvement via Berstein-type bound.} Although the tail bound (\ref{Eq: Hoeffding type tail bound}) is simple depending only on the variance term $\widehat{\sigma}_{\text{lin}}^2$, it may not be effective when $n_1$ is much smaller than $n$ (e.g.~$n_1^2/n \rightarrow 0$ as $n_1 \rightarrow \infty$). In such case, Proposition~\ref{Proposition: Bernstein bound} gives a tighter bound. More specifically, following the same steps as before, Proposition~\ref{Proposition: Bernstein bound} presents a concentration inequality for the two-sample (unscaled) $t$-statistic as
\begin{align*}
\mP_{\pi}(L_n^\pi \geq t | \mathcal{X}_n) ~ = ~& \mathbb{P}_\pi \Bigg\{  \left( \overline{Y}_{1,\pi} - \overline{Y}_{2,\pi} \right)  \geq \frac{tn^2}{2n_1n_2}  \Bigg| \mathcal{X}_n  \Bigg\} \\[.5em] 
~ \leq ~ & \exp \Bigg\{ - \frac{nt^2}{ 8\frac{n_1n_2}{n^2} \widehat{\sigma}_{\text{lin}}^2 + \frac{4}{3} t \cdot \max \left( \frac{n_1}{n}, \frac{n_2}{n} \right) \cdot M_Z }  \Bigg\},
\end{align*}
where $M_Z :=\max_{1 \leq i \leq n} |Z_i - \overline{Z}|$. Furthermore, using the relationship (\ref{Eq: exact relationship}) and by symmetry, 
\begin{align} \label{Eq: Bernstein type tail bound}
\mathbb{P}_\pi \left( |\overline{Y}_{1,\pi} - \overline{Y}| \geq t  \big| \mathcal{X}_n \right) \leq 2\exp \Bigg\{ - \frac{12n_1 t^2}{24 \frac{n_2}{n} \widehat{\sigma}_{\text{lin}}^2 + 8\frac{n_2}{n} M_Z t } \Bigg\},
\end{align}
where we assumed $n_1 \leq n_2$. 
\begin{remark} \normalfont
	We remark that the bounds in (\ref{Eq: Hoeffding type tail bound}) and (\ref{Eq: Bernstein type tail bound}) are byproducts of more general bounds and are not necessary the sharpest ones in the context of sampling without replacement. We refer to \cite{bardenet2015concentration} and among others for some recent developments of concentration bounds for sampling without replacement. 
\end{remark}

\section{Improved version of Theorem~\ref{Theorem: Two-Sample U-statistic}} \label{Section: Improved version of Theorem}
In this section, we improve the result of Theorem~\ref{Theorem: Two-Sample U-statistic} based on the exponential bound in Theorem~\ref{Theorem: Two-Sample Concentration}. In particular we replace the dependency on $\alpha^{-1}$ there with $\log(1/\alpha)$ by adding an extra assumption that $n_1 \asymp n _2$ as follows.
\begin{lemma}[Two-sample $U$-statistic] \label{Lemma: Two-Sample U-statistic Improved Version}
	For $0 < \alpha < e^{-1}$, suppose that there is a sufficiently large constant $C>0$ such that
	\begin{align} \label{Eq: Two-Sample U-statistic Sufficient Condition 2}
	\mE_P[U_{n_1,n_2}] \geq C \max \Bigg\{ \sqrt{\frac{\psi_{Y,1}(P)}{\beta n_1}}, \ \sqrt{\frac{\psi_{Z,1}(P)}{\beta n_2}}, \  \sqrt{\frac{\psi_{YZ,2}(P)}{\beta}} \log \left(\frac{1}{\alpha}\right) \cdot \left( \frac{1}{n_1} + \frac{1}{n_2}\right) \Bigg\},
	\end{align}
	for all $P \in \mathcal{P}_1 \subset \mathcal{P}_{h_\text{\emph{ts}}}$. Then under the assumptions that $n_1 \asymp n_2$, the type II error of the permutation test over $\mathcal{P}_1$ is uniformly bounded by $\beta$, that is
	\begin{align*}
	\sup_{P \in \mathcal{P}_1} \mP_P^{(n_1,n_2)} (U_{n_1,n_2} \leq c_{1-\alpha,n_1,n_2} )\leq \beta. 
	\end{align*}
\begin{proof}
To prove the above lemma, we employ the quantile approach described in Section~\ref{Section: A general strategy with two moments} \citep[see also][]{fromont2013two}. More specifically we let $q_{1-\beta/2,n}$ denote the quantile of the permutation critical value $c_{1-\alpha,n}$ of $U_{n_1,n_2}$. Then as shown in the proof of Lemma~\ref{Lemma: Two Moments Method}, if 
\begin{align*} 
\mE_P[U_{n_1,n_2}] ~ \geq ~ q_{1-\beta/2,n} + \sqrt{\frac{2\mV_P[U_{n_1,n_2}]}{\beta}},
\end{align*}
then the type II error of the permutation test is controlled as 
\begin{align*}
\sup_{P \in \mathcal{P}_1} \mP_P( U_{n_1,n_2} \leq c_{1-\alpha,n}) ~\leq~ & \sup_{P \in \mathcal{P}_1} \mP_P( U_{n_1,n_2} \leq q_{1-\beta/2,n})  + \sup_{P \in \mathcal{P}_1} \mP_P( q_{1-\beta/2,n} < c_{1-\alpha,n})  \\[.5em]
\leq ~ &\beta. 
\end{align*}
Therefore it is enough to verify that the right-hand side of (\ref{Eq: Two-Sample U-statistic Sufficient Condition 2}) is lower bounded by $ q_{1-\beta/2,n} + \sqrt{2\mV_P[U_{n_1,n_2}]/\beta}$. As shown in the proof of Theorem~\ref{Theorem: Two-Sample U-statistic}, the variance is bounded by
\begin{align} \label{Eq: Variance Bound}
\mV_P[U_{n_1,n_2}] \leq C_1 \frac{\psi_{Y,1}(P)}{n_1}  +  C_2 \frac{\psi_{Z,1}(P)}{n_2} + C_3 \psi_{YZ,2}(P) \left( \frac{1}{n_1} + \frac{1}{n_2} \right)^2.
\end{align}
Moving onto an upper bound for $q_{1-\beta/2,n}$, let us denote
\begin{align*}
\Sigma_{n_1,n_2}^\dagger := \frac{1}{n_1^2(n_1-1)^2} \sum_{(i_1,i_2) \in \mathbf{i}_2^{n}} g^2(X_{i_1},X_{i_2}).
\end{align*}
From Theorem~\ref{Theorem: Two-Sample Concentration} together with the trivial bound (\ref{Eq: trivial bound}), we know that $c_{1-\alpha,n}$ is bounded by 
\begin{equation}
\begin{aligned} \label{Eq: quantile bound}
c_{1-\alpha,n} & ~\leq~ \max \Bigg\{ \sqrt{\frac{\Sigma_{n_1,n_2}^{\dagger 2} }{C_4} \log \left( \frac{1}{\alpha} \right)}, \  \frac{\Sigma_{n_1,n_2}^\dagger}{C_4} \log \left( \frac{1}{\alpha} \right) \Bigg\}  \\[.5em]
& ~\leq ~ C_5 \Sigma_{n_1,n_2}^\dagger \log \left( \frac{1}{\alpha} \right),
\end{aligned}
\end{equation}
where the last inequality uses the assumption that $\alpha < e^{-1}$. Now applying Markov's inequality yields
\begin{align*}
\mP_P \left( \Sigma_{n_1,n_2}^\dagger  \geq t \right) \leq \frac{\mE_P[\Sigma_{n_1,n_2}^{\dagger 2}]}{t^2} \leq C_6 \frac{\psi_{YZ,2}(P)}{t^2 n_1^2}.
\end{align*}
By setting the right-hand side to be $\beta/2$, we can find an upper bound for the $1-\beta/2$ quantile of $\Sigma_{n_1,n_2}^\dagger$. Combining this observation with inequality (\ref{Eq: quantile bound}) yields
\begin{align*}
q_{1-\beta/2,n} \leq \frac{C_7}{\beta^{1/2}} \log \left( \frac{1}{\alpha} \right) \frac{\sqrt{\psi_{YZ,2}(P)}}{n_1}.
\end{align*}
Therefore, from the above bound and (\ref{Eq: Variance Bound}), 
\begin{align*}
&q_{1-\beta/2,n} + \sqrt{\frac{2\mV_P[U_{n_1,n_2}]}{\beta}} \\[.5em]
\leq ~ & C \sqrt{\max \Bigg\{ \frac{\psi_{Y,1}(P)}{\beta n_1}, \ \frac{\psi_{Z,1}(P)}{\beta n_2}, \  \frac{\psi_{YZ,2}(P)}{\beta}\log^2 \left(\frac{1}{\alpha}\right) \cdot \left( \frac{1}{n_1} + \frac{1}{n_2}\right)^2 \Bigg\}}.
\end{align*}
This completes the proof of Lemma~\ref{Lemma: Two-Sample U-statistic Improved Version}.
\end{proof}
\end{lemma}

\section{Monte Carlo-based permutation tests} \label{Section: Monte Carlo-based permutation tests}
As mentioned in Remark~\ref{Remark: Computational aspects}, exact calculation of a permutation test is computationally infeasible for large sample sizes. To mitigate this computational issue, it is now standard to consider a Monte Carlo-based (MC-based) permutation test in practical applications. The goal of this section is to demonstrate that this MC-based permutation test can achieve the same minimax rate optimality as the original permutation test, while controlling the type I error rate. Of course, this guarantee requires a sufficiently large number of Monte Carlo samples and we show that this number depends on the pre-specified type I and II error rates $\alpha$ and $\beta$, but not on the sample size.

To fix ideas, let $\pi_1,\ldots,\pi_B$ be random permutations uniformly sampled with replacement from $\mathbf{\Pi}_n$, where we recall that $\mathbf{\Pi}_n$ is the set of all possible permutations of $\{1,\ldots,n\}$. Then the MC-based permutation $p$-value is given by
\begin{align} \label{Eq: MC-based p-value}
	\widehat{p} := \frac{1}{B+1} \bigg\{ 1 + \sum_{i=1}^B \mathds{1} \big(T_n^{\pi_i} \geq T_n \big)  \bigg\} = \frac{1}{B+1} \sum_{i=1}^{B+1}  \mathds{1} \big(T_n^{\pi_i} \geq T_n \big),
\end{align}
where we define $T_n^{\pi_{B+1}} := T_n$. It is worth pointing out that adding the identity permutation is required to control the type I error rate in finite sample sizes~\citep[e.g.][]{phipson2010permutation}. In particular, we have $\mP(\widehat{p} \leq \alpha) \leq \alpha$ whenever $\mathcal{X}_n$ and $\mathcal{X}_n^{\pi}$ have the same distribution for every $\pi \in \mathbf{\Pi}_n$. A formal proof of this statement can be found, for example, in \cite{hemerik2018exact}. Let $\widehat{c}_{1-\alpha,n}$ be the $1-\alpha$ quantile of the distribution with random permutations defined as
\begin{equation}
	\begin{aligned} \label{Eq: MC critical value}
		\widehat{c}_{1-\alpha,n} ~:= ~& \inf\bigg\{t: \frac{1}{B+1} \sum_{i=1}^{B+1}  \mathds{1} \big(T_n^{\pi_i} \leq t \big) \geq 1-\alpha \bigg\}. 
	\end{aligned}
\end{equation}
Then it can be verified that the MC-based permutation test can be equivalently written as $\mathds{1}(\widehat{p} \leq \alpha)  = \mathds{1}(T_n > \widehat{c}_{1-\alpha,n})$. Intuitively, when $B$ is sufficiently large, $ \widehat{c}_{1-\alpha,n}$ becomes close to $c_{1-\alpha,n}$, the $1-\alpha$ quantile of the original permutation distribution. Consequently, the type II error of the MC-based permutation test should approximate that of the original permutation test $\mathds{1}(T_n > c_{1-\alpha,n})$. The next proposition formalizes this intuition by reproducing Lemma~\ref{Lemma: Two Moments Method} for the MC-based permutation test.

\begin{proposition}[Two-moment method for MC-based tests] \label{Proposition: Two-moment method for MC-based tests}
	Suppose that for each permutation $\pi \in \mathbf{\Pi}_n$,  $T_n$ and $T_n^{\pi}$ have the same distribution under the null hypothesis. Given pre-specified error rates $\alpha \in (0,1)$ and $\beta \in (0,1 - \alpha)$, assume that for any $P \in \mathcal{P}_1$,
	\begin{equation}
		\begin{aligned} \label{Eq: Sufficient Condition for MC}
			\mE_P[T_n]  ~\geq~ & \mE_P [\mE_{\pi} \{T_n^{\pi} | \mathcal{X}_n\}] + \sqrt{\frac{6\emph{\mV}_P[\mE_{\pi}\{T_n^{\pi}| \mathcal{X}_n \}]}{\beta}} \\
			&  + \sqrt{\frac{6\emph{\mV}_P[T_n]}{\beta}}  + \sqrt{\frac{12\mE_P[\emph{\mV}_{\pi}\{T_n^{\pi}| \mathcal{X}_n \}]}{\alpha\beta}}.
		\end{aligned}
	\end{equation}
	Suppose further that the number of Monte Carlo simulations $B$ is greater than $B \geq 8 \alpha^{-2} \log(4/\beta)$. Then the MC-based permutation test $\ind(\widehat{p} \leq \alpha)$ controls the type I and II error rates as in (\ref{Eq: uniform error control}) where $\widehat{p}$ is defined in (\ref{Eq: MC-based p-value}).
\end{proposition}

We present several remarks on Proposition~\ref{Proposition: Two-moment method for MC-based tests}. 
\begin{remark} \leavevmode \normalfont
	\begin{itemize}
		\item The new condition~(\ref{Eq: Sufficient Condition for MC}) is the same as the condition for the full permutation test in Lemma~\ref{Lemma: Two Moments Method} up to constant factors. In other words, the MC-based permutation test can achieve the same error guarantee as the original permutation test by adjusting constant factors. \\[-0.8em]
		\item Using Proposition~\ref{Proposition: Two-moment method for MC-based tests} and its modification, we can reproduce the rate optimality results of the permutation test based on the Monte Carlo counterpart. \\[-0.8em]
		\item The proof of Proposition~\ref{Proposition: Two-moment method for MC-based tests} is based on Dvoretzky--Kiefer--Wolfowitz inequality. Using this, we verify the type II error can be bounded by
		\begin{align*}
			\mP_P^{(n)}(\widehat{p} > \alpha) ~\leq~ \mP_P^{(n)}(T_n \leq c_{1-\alpha^\ast}) + \beta/2,
		\end{align*}
		where $1-\alpha^\ast = \frac{B+1}{B}(1-\alpha) + \sqrt{\frac{1}{2B} \log \left( 4/\beta \right)}$. We then build on Lemma~\ref{Lemma: Two Moments Method} and show that $\mP_P(T_n \leq c_{1-\alpha^\ast}) \leq \beta/2$ under the given condition to complete the proof. The details can be found in Appendix~\ref{Section: Proof of Proposition: Two-moment method for MC-based tests}.   \\[-0.8em]
		\item It is interesting to point out that the condition for $B$ is independent of the sample size $n$. This may be explained by the fact that the distance between $\widehat{c}_{1-\alpha,n}$ and $c_{1-\alpha,n}$ can be made small independent of $n$. 
	\end{itemize}
\end{remark}

\section{Proof of Lemma~\ref{Lemma: Two Moments Method}} \label{Section: Proof of Lemma: Two Moments Method}
As discussed in the main text, the key difficulty of studying the type II error of the permutation test lies in the fact that its critical value is data-dependent and thereby random. Our strategy to overcome this problem is to bound the random critical value by a quantile value with high probability \citep[see also][]{fromont2013two}. We split the proof of Lemma~\ref{Lemma: Two Moments Method} into three steps. In the first step, we present a sufficient condition under which the type II error of the test with a non-random cutoff value is small. In the second step, we provide a non-random upper bound for the permutation critical value, which holds with high probability. In the last step, we combine the results and complete the proof.

\paragraph{Step 1.} For a given $P \in \mathcal{P}_1$, let $\omega(P)$ be any constant depending on $P$ such that  
\begin{align} \label{Eq: Condition of a critical value}
\mE_P[T_n] ~ \geq ~ \omega(P) + \sqrt{\frac{3\mV_P[T_n]}{\beta}}.
\end{align}
Based on such $\omega(P)$, we define a test $\ind\{T_n > \omega(P)\}$, which controls the type II error by $\beta/3$. To see this, let us apply Chebyshev's inequality 
\begin{align*}
\beta/3 ~ \geq ~ & \mP_P \big(  \big|T_n - \mE_P[T_n] \big| \geq \sqrt{3\beta^{-1} \mV_P[T_n]} \big) \\[.5em]
\geq ~ & \mP_P \big(  -T_n + \mE_P[T_n]  \geq \sqrt{3\beta^{-1} \mV_P[T_n]} \big) \\[.5em]
\geq ~ & \mP_P \big( \omega(P) \geq T_n \big),
\end{align*}
where the last inequality uses the condition of $\omega(P)$ in (\ref{Eq: Condition of a critical value}). In other words, the type II error of the test $\ind\{T_n > \omega(P)\}$ is less than or equal to $\beta/3$ as desired. 

\paragraph{Step 2.} In this step, we provide an upper bound for $c_{1-\alpha,n}$, which may hold with high probability. First, applying Chebyshev's inequality yields 
\begin{align*}
\mP_{\pi} \left( \big|T_n^{\pi} - \mE_{\pi}[T_n^{\pi} | \mathcal{X}_n] \big| \geq \sqrt{\alpha^{-1}\mV_{\pi}[T_n^{\pi} | \mathcal{X}_n]} ~\big| \mathcal{X}_n \right) \leq \alpha.
\end{align*}
Therefore, by the definition of the quantile, we see that $c_{1-\alpha,n}$ satisfies
\begin{align} \label{Eq: Quantile Bound}
c_{1-\alpha,n}  ~\leq~  \mE_{\pi}[T_n^{\pi} | \mathcal{X}_n] + \sqrt{\alpha^{-1}\mV_{\pi}[T_n^{\pi} | \mathcal{X}_n]}.
\end{align}
Note that the two terms on the right-hand side are random variables depending on $\mathcal{X}_n$. In order to use the result from the first step, we want to further upper bound these two terms by some constants. To this end, let us define two good events:
\begin{align*}
& \mathcal{A}_1 := \Big\{ \mE_{\pi}[T_n^{\pi} | \mathcal{X}_n] <  \mE_P[\mE_{\pi}\{T_n^{\pi} | \mathcal{X}_n\}] + \sqrt{3\beta^{-1}\mV_P [ \mE_{\pi}\{T_n^{\pi} | \mathcal{X}_n \}] }  \Big\}, \\[.5em]
& \mathcal{A}_2 := \Big\{ \sqrt{\alpha^{-1}\mV_{\pi}[T_n^{\pi} | \mathcal{X}_n]} <  \sqrt{3\alpha^{-1}\beta^{-1} \mE_P[ \mV_{\pi} \{T_n^{\pi} | \mathcal{X}_n\} ] } \Big\}.
\end{align*}
Then by applying Markov and Chebyshev's inequalities, it is straightforward to see that 
\begin{align} \label{Eq: Concentration Bounds}
\mP_P(\mathcal{A}_1^c) \leq \beta/3 \quad \text{and} \quad \mP_P(\mathcal{A}_2^c) \leq \beta/3.
\end{align}

\paragraph{Step 3.} Here, building on the first two steps, we conclude the result. We begin by upper bounding the type II error of the permutation test as
\begin{align*}
\mP_P (T_n \leq c_{1-\alpha,n}) ~=~ &  \mP_P (T_n \leq c_{1-\alpha,n}, \ \mathcal{A}_1 \cup \mathcal{A}_2)  + \mP_P (T_n \leq c_{1-\alpha,n}, \ \mathcal{A}_1^c \cap \mathcal{A}_2^c) \\[.5em]
\leq ~ & \mP_P (T_n \leq \omega'(P)) +  \mP_P (\mathcal{A}_1^c \cap \mathcal{A}_2^c),
\end{align*}
where, for simplicity, we write
\begin{align*}
\omega'(P) := \mE_P[\mE_{\pi}\{T_n^{\pi} | \mathcal{X}_n\}] + \sqrt{3\beta^{-1}\mV_P [ \mE_{\pi}\{T_n^{\pi} | \mathcal{X}_n \}]}  +  \sqrt{3\alpha^{-1}\beta^{-1} \mE_P[ \mV_{\pi} \{T_n^{\pi} | \mathcal{X}_n\}]}.
\end{align*}
One may check that the type II error of $\ind\{T_n > \omega'(P)\}$ is controlled by $\beta/3$ as long as $\omega'(P) + \sqrt{3\mV_P[T_n]/\beta} \leq \mE_P[T_n]$ from the inequality~(\ref{Eq: Condition of a critical value}) in Step 1. However, this sufficient condition is ensured by condition (\ref{Eq: Sufficient Condition}) of Lemma~\ref{Lemma: Two Moments Method}. Furthermore, the probability of the intersection of the two bad events $\mathcal{A}_1^c \cap \mathcal{A}_2^c$ is also bounded by $2\beta/3$ due to the concentration results in (\ref{Eq: Concentration Bounds}). Hence, by taking the supremum over $P \in \mathcal{P}_1$, we may conclude that 
\begin{align*}
\sup_{P \in \mathcal{P}_1} \mP_P (T_n \leq c_{1-\alpha,n}) \leq \beta.
\end{align*}
This completes the proof of Lemma~\ref{Lemma: Two Moments Method}.

\section{Proof of Theorem~\ref{Theorem: Two-Sample U-statistic}} \label{Section: Proof of Theorem: Two-Sample U-statistic}
We proceed the proof by verifying the sufficient condition in Lemma~\ref{Lemma: Two Moments Method}. We first verify that the expectation of $U_{n_1,n_2}^\pi$ is zero under the permutation law. Let us recall the permuted $U$-statistic $U_{n_1,n_2}^\pi$ in (\ref{Eq: permuted two-sample U-statistic}). In fact, by the linearity of expectation, it suffices to prove 
\begin{align*}
	\mE_{\pi} [h_{\text{ts}}(X_{\pi_{1}},X_{\pi_{2}}; X_{\pi_{n_1+ 1}},X_{\pi_{n_1+2}}) | \mathcal{X}_n] = 0.
\end{align*}
This is clearly the case by recalling the definition of kernel $h_{\text{ts}}$ in (\ref{Eq: Two-Sample kernel}) and noting that the expectation $\mE_{\pi}[g(X_{\pi_i}, X_{\pi_j}) | \mathcal{X}_n]$ is invariant to the choice of $(i,j) \in \mathbf{i}_2^n$, which leads to $\mE_{\pi}[U_{n_1,n_2}^\pi | \mathcal{X}_n] = 0$. Therefore we only need to verify the simplified condition (\ref{Eq: Sufficient Condition II}) under the given assumptions in Theorem~\ref{Theorem: Two-Sample U-statistic}. 

The rest of the proof is divided into two parts. In each part, we prove the following conditions separately,
\begin{align} \label{Eq: Condition I}
	&\mE_P[U_{n_1,n_2}] \geq 2\sqrt{\frac{2\mV_P[U_{n_1,n_2}]}{\beta}} \quad \text{and} \\[.5em] \label{Eq: Condition II}
	& \mE_P[U_{n_1,n_2}] \geq 2\sqrt{\frac{2\mE_P[\mV_{\pi}\{U_{n_1,n_2}^\pi| \mathcal{X}_n \}]}{\alpha\beta}}.
\end{align}
We then complete the proof of Theorem~\ref{Theorem: Two-Sample U-statistic} by noting that (\ref{Eq: Condition I}) and (\ref{Eq: Condition II}) imply the simplified condition (\ref{Eq: Sufficient Condition II}).

\vskip 1em 

\noindent \textbf{Part 1. Verification of condition (\ref{Eq: Condition I}).} In this part, we verify condition~(\ref{Eq: Condition I}). To do so, we state the explicit variance formula of  a two-sample $U$-statistic \citep[e.g. page 38 of][]{lee1990u}. Following the notation of \cite{lee1990u}, we let $\check{\sigma}_{i,j}^2$ denote the variance of a conditional expectation given as
\begin{align*}
	\check{\sigma}_{i,j}^2 = \mV_P [\mE_P \{\overline{h}_{\text{ts}}(y_1,\ldots,y_i,Y_{i+1},\ldots,Y_2;z_1,\ldots,z_j,Z_{j+1},\ldots,Z_{2})\}]
\end{align*} 
for $0 \leq i,j \leq 2$. Then the variance of $U_{n_1,n_2}$ is given by
\begin{align}
	\mV_P[U_{n_1,n_2}] = \sum_{i=0}^2 \sum_{j=0}^2 \binom{2}{i}\binom{2}{j}\binom{n_1-2}{2-i}\binom{n_2-2}{2-j}\binom{n_1}{2}^{-1}\binom{n_2}{2}^{-1} \check{\sigma}_{i,j}^2. \label{Eq: variance expression}
\end{align}
By the law of total variance, one may see that $\check{\sigma}_{i,j}^2 \leq \check{\sigma}_{2,2}^2$ for all $0 \leq i,j \leq 2$. This leads to an upper bound for $\mV_P[U_{n_1,n_2}]$ as
\begin{align*}
	\mV_P[U_{n_1,n_2}]  \leq C_1 \frac{\check{\sigma}_{1,0}^2}{n_1} + C_2 \frac{\check{\sigma}_{0,1}^2}{n_2} + C_3 \left( \frac{1}{n_1} + \frac{1}{n_2} \right)^2 \check{\sigma}_{2,2}^2. 
\end{align*}
Now applying Jensen's inequality, repeatedly, yields
\begin{align*}
	\check{\sigma}_{2,2}^2 \leq \mE_P [ \overline{h}_{\text{ts}}^2(Y_1,Y_2;Z_1,Z_2)] \leq \mE_P [ h_{\text{ts}}^2(Y_1,Y_2;Z_1,Z_2)] \leq C_4 \psi_{YZ,2}(P).
\end{align*}
Then by noting that $\check{\sigma}_{1,0}^2$ and $\check{\sigma}_{0,1}^2$ correspond to the notation $\psi_{Y,1}(P)$ and $\psi_{Z,1}(P)$, respectively,
\begin{align*}
	\mV_P[U_{n_1,n_2}]  \leq C_1 \frac{\psi_{Y,1}(P)}{n_1} + C_2 \frac{\psi_{Z,1}(P)}{n_2} + C_5 \left( \frac{1}{n_1} + \frac{1}{n_2} \right)^2 \psi_{YZ,2}(P). 
\end{align*}
Hence condition~(\ref{Eq: Condition I}) is satisfied by taking the constant $C$ in Theorem~\ref{Theorem: Two-Sample U-statistic} sufficiently large. 

\vskip 1em

\noindent \sloppy \textbf{Part 2. Verification of condition (\ref{Eq: Condition II}).} In this part, we verify condition~(\ref{Eq: Condition II}). Intuitively, the permuted $U$-statistic behaves similarly as the unconditional $U$-statistic under a certain null model. This means that the variance of $U_{n_1,n_2}^\pi$ should have a similar convergence rate as $(n_1^{-1} + n_2^{-1})^2\psi_{YZ,2}(P)$ since $\psi_{Y,1}(P)$ and $\psi_{Z,1}(P)$ are zero under the null hypothesis. We now prove that this intuition is indeed correct. Since $U_{n_1,n_2}^\pi$ is centered under the permutation law, it is enough to study $\mE_P[\mE_{\pi}\{(U_{n_1,n_2}^\pi)^2| \mathcal{X}_n \}]$. Let us write a set of indices $\mathsf{I}_{\text{total}} := \{ (i_1,i_2,j_1,j_2,i_1',i_2',j_1',j_2') \in \mathbb{N}^8_{+}: (i_1,i_2) \in \mathbf{i}_2^{n_1}, (j_1,j_2) \in \mathbf{i}_2^{n_2}, (i_1',i_2') \in \mathbf{i}_2^{n_1}, (j_1',j_2') \in \mathbf{i}_2^{n_2}\}$ and define $\mathsf{I}_{\text{A}} = \{ (i_1,i_2,j_1,j_2,i_1',i_2',j_1',j_2') \in \mathsf{I}_{\text{total}} : \#|\{i_1,i_2\} \cap \{i_1',i_2'\}| + \#|\{j_1,j_2\} \cap \{j_1',j_2'\}|  \leq 1\}$ and its complement $\mathsf{I}_{\text{A}^c} = \{ (i_1,i_2,j_1,j_2,i_1',i_2',j_1',j_2') \in \mathsf{I}_{\text{total}} : \#|\{i_1,i_2\} \cap \{i_1',i_2'\}| + \#|\{j_1,j_2\} \cap \{j_1',j_2'\}| >  1\}$. Here $\#|B|$ denotes the cardinality of a set $B$. Based on this notation and the linearity of expectation, 
\begin{align*}
	& \mE_{\pi}[(U_{n_1,n_2}^\pi)^2| \mathcal{X}_n] \\[.5em]
	~=~ &\frac{1}{(n_1)_{(2)}^2(n_2)_{(2)}^2} \sum_{(i_1,\ldots,j_2') \in \mathsf{I}_{\text{total}}} \mE_{\pi} \Big[ h_{\text{ts}}(X_{\pi_{i_1}},X_{\pi_{i_2}}; X_{\pi_{n_1+ j_1}},X_{\pi_{n_1+j_2}}) \\ 
	& ~~~~~~~~~~~~~~~~~~~~~~~~~~~~~~~~~~~~ \times h_{\text{ts}}(X_{\pi_{i_1'}},X_{\pi_{i_2'}}; X_{\pi_{n_1+ j_1'}},X_{\pi_{n_1+j_2'}}) ~ \Big| \mathcal{X}_n \Big] \\[.5em]
	= ~& (I) + (II), 
\end{align*}
where
\begin{align*}
	(I) ~:=~ &\frac{1}{(n_1)_{(2)}^2(n_2)_{(2)}^2} \sum_{(i_1,\ldots,j_2') \in \mathsf{I}_{\text{A}}} \mE_{\pi} \Big[ h_{\text{ts}}(X_{\pi_{i_1}},X_{\pi_{i_2}}; X_{\pi_{n_1+ j_1}},X_{\pi_{n_1+j_2}}) \\ 
	& ~~~~~~~~~~~~~~~~~~~~~~~~~~~~~~~~~~\times h_{\text{ts}}(X_{\pi_{i_1'}},X_{\pi_{i_2'}}; X_{\pi_{n_1+ j_1'}},X_{\pi_{n_1+j_2'}}) ~ \Big| \mathcal{X}_n \Big],  \\[.5em]
	(II) ~:=~ &\frac{1}{(n_1)_{(2)}^2(n_2)_{(2)}^2} \sum_{(i_1,\ldots,j_2') \in \mathsf{I}_{\text{A}^c}} \mE_{\pi} \Big[ h_{\text{ts}}(X_{\pi_{i_1}},X_{\pi_{i_2}}; X_{\pi_{n_1+ j_1}},X_{\pi_{n_1+j_2}}) \\ 
	& ~~~~~~~~~~~~~~~~~~~~~~~~~~~~~~~~~~\times h_{\text{ts}}(X_{\pi_{i_1'}},X_{\pi_{i_2'}}; X_{\pi_{n_1+ j_1'}},X_{\pi_{n_1+j_2'}}) ~ \Big| \mathcal{X}_n \Big]. 
\end{align*}
We now claim that the first term $(I) = 0$. This is the key observation that makes the upper bound for the variance of the permuted $U$-statistic depend on $(n_1^{-1} + n_2^{-1})^{2}$ rather than a slower rate $(n_1 + n_2)^{-1}$. First consider the case where $\#|\{i_1,i_2\} \cap \{i_1',i_2'\}| + \#|\{j_1,j_2\} \cap \{j_1',j_2'\}|  = 0$, that is, all indices are distinct. Let us focus on the summands of $(I)$. By symmetry, we may assume the set of indices $(i_1,i_2,n_1+j_1,n_1+j_2,i_1',i_2',n_1+j_1',n_1+j_2')$ to be $(1,\ldots,8)$ and observe that
\begin{align*}
	& \mE_{\pi} \Big[ h_{\text{ts}}(X_{\pi_{1}},X_{\pi_{2}}; X_{\pi_{3}},X_{\pi_{4}})h_{\text{ts}}(X_{\pi_{5}},X_{\pi_{6}}; X_{\pi_{7}},X_{\pi_{8}}) \Big| \mathcal{X}_n \Big] \\[.5em]
	\overset{(i)_1}{=} ~ & \mE_{\pi} \Big[ h_{\text{ts}}(X_{\pi_{3}},X_{\pi_{2}}; X_{\pi_{1}},X_{\pi_{4}})h_{\text{ts}}(X_{\pi_{5}},X_{\pi_{6}}; X_{\pi_{7}},X_{\pi_{8}}) \Big| \mathcal{X}_n \Big] \\[.5em]
	\overset{(ii)_1}{=} ~ & - \mE_{\pi} \Big[ h_{\text{ts}}(X_{\pi_{1}},X_{\pi_{2}}; X_{\pi_{3}},X_{\pi_{4}})h_{\text{ts}}(X_{\pi_{5}},X_{\pi_{6}}; X_{\pi_{7}},X_{\pi_{8}}) \Big| \mathcal{X}_n \Big] \\[.5em]
	\overset{(iii)_1}{=} ~ & 0,
\end{align*}
where $(i)_1$ holds since the distribution of the product kernels does not change even after $\pi_1$ and $\pi_3$ are switched and $(ii)_1$ uses the fact that $h_{\text{ts}}(y_1,y_2;z_1,z_2) = - h_{\text{ts}}(z_1,y_2;y_1,z_2)$. $(iii)_1$ follows directly by comparing the first line and the third line of the equations. 

Next consider the case where $\#|\{i_1,i_2\} \cap \{i_1',i_2'\}| + \#|\{j_1,j_2\} \cap \{j_1',j_2'\}| = 1$. In this case, we also argue that the expectation is zero. To verify this, assume that $i_1 = i_1'$ (and the other cases similarly follow). Then, by symmetry again, we have 
\begin{align*}
	& \mE_{\pi} \Big[ h_{\text{ts}}(X_{\pi_{1}},X_{\pi_{2}}; X_{\pi_{3}},X_{\pi_{4}})h_{\text{ts}}(X_{\pi_{1}},X_{\pi_{5}}; X_{\pi_{6}},X_{\pi_{7}}) \Big| \mathcal{X}_n \Big] \\[.5em]
	\overset{(i)_2}{=} ~ &  \mE_{\pi} \Big[ h_{\text{ts}}(X_{\pi_{1}},X_{\pi_{4}}; X_{\pi_{3}},X_{\pi_{2}})h_{\text{ts}}(X_{\pi_{1}},X_{\pi_{5}}; X_{\pi_{6}},X_{\pi_{7}}) \Big| \mathcal{X}_n \Big] \\[.5em]
	\overset{(ii)_2}{=} ~ & -  \mE_{\pi} \Big[ h_{\text{ts}}(X_{\pi_{1}},X_{\pi_{2}}; X_{\pi_{3}},X_{\pi_{4}})h_{\text{ts}}(X_{\pi_{1}},X_{\pi_{5}}; X_{\pi_{6}},X_{\pi_{7}}) \Big| \mathcal{X}_n \Big] \\[.5em]
	\overset{(iii)_2}{=} ~ & 0,
\end{align*}
where $(i)_2$ follows by the same reasoning for $(i)_2$, and $(ii)_2$ holds since $h_{\text{ts}}(y_1,y_2;z_1,z_2) = - h_{\text{ts}}(y_1,z_2;y_1,y_2)$. Then $(iii)_2$ is obvious by comparing the first line and the third line of the equations. Hence, for any choice of indices $(i_1,\ldots,j_2') \in \mathsf{I}_{\text{A}}$, the summands of $(I)$ becomes zero, which leads to $(I) = 0$. 

Now turning to the second term $(II)$, for any $1 \leq i_1\neq i_2,i_3 \neq i_4 \leq n$, we have 
\begin{align*}
	& \big| \mE_P \big[\mE_{\pi}\{ g(X_{\pi_{i_1}},X_{\pi_{i_2}})g(X_{\pi_{i_3}},X_{\pi_{i_4}}) | \mathcal{X}_n \} \big] \big| \\[.5em]
	~\overset{(i)_3}{=}~ & \big|\mE_{\pi} \big[\mE_P\{ g(X_{\pi_{i_1}},X_{\pi_{i_2}})g(X_{\pi_{i_3}},X_{\pi_{i_4}}) | \pi_{i_1},\ldots,\pi_{i_4} \} \big] \big| \\[.5em]
	\overset{(ii)_3}{\leq} ~ & \frac{1}{2} \mE_{\pi} \big[ \mE_P\{ g^2(X_{\pi_{i_1}},X_{\pi_{i_2}}) | \pi_{i_1}, \pi_{i_2} \} \big] + \frac{1}{2} \mE_{\pi} \big[ \mE_P\{ g^2(X_{\pi_{i_3}},X_{\pi_{i_4}}) | \pi_{i_3}, \pi_{i_4} \} \big] \\[.5em]
	\overset{(iii)_3}{\leq} ~ & \psi_{YZ,2}(P),
\end{align*}
where $(i)_3$ uses the the law of total expectation, $(ii)_3$ uses the basic inequality $xy \leq x^2/2 + y^2/2$ and $(iii)_3$ clearly holds by recalling the definition of $\psi_{YZ,2}(P)$. Using this observation, it is not difficult to see that for any $(i_1,\ldots,j_2') \in \mathsf{I}_{\text{total}}$,
\begin{align*}
	& \big| \mE_P\big[ \mE_{\pi} \{ h_{\text{ts}}(X_{\pi_{i_1}},X_{\pi_{i_2}}; X_{\pi_{n_1+ j_1}},X_{\pi_{n_1+j_2}})  \\[.5em]
	& ~~~~~~~~ \times h_{\text{ts}}(X_{\pi_{i_1'}},X_{\pi_{i_2'}}; X_{\pi_{n_1+ j_1'}},X_{\pi_{n_1+j_2'}}) \big| \mathcal{X}_n \} \big] \big| \\[.5em]
	\leq ~ & C_6 \psi_{YZ,2}(P).
\end{align*}
Therefore, by counting the number of elements in $\mathsf{I}_{\text{A}^c}$,
\begin{align*}
	\mE_P[\mV_{\pi}\{U_{n_1,n_2}^\pi| \mathcal{X}_n \}] = \mE_P[(II)] ~\leq~ & C_6 \psi_{YZ,2}(P) \times \frac{1}{(n_1)_{(2)}^2(n_2)_{(2)}^2} \sum_{(i_1,\ldots,j_2') \in \mathsf{I}_{\text{A}^c}}1  \\[.5em]
	\leq ~ & C_7 \psi_{YZ,2}(P) \left( \frac{1}{n_1} + \frac{1}{n_2} \right)^2.
\end{align*}
Hence condition~(\ref{Eq: Condition II}) is satisfied by taking the constant $C$ in Theorem~\ref{Theorem: Two-Sample U-statistic} sufficiently large. This completes the proof of Theorem~\ref{Theorem: Two-Sample U-statistic}.

\section{Proof of Proposition~\ref{Proposition: Multinomial Two-Sample Testing}} \label{Section: Proof of Proposition: Multinomial Two-Sample Testing}
As discussed in the main text, we start proving that the three inequalities in (\ref{Eq: Three conditions for two-sample multinomials}) are fulfilled. Focusing on the first one, we want to show that 
\begin{align*}
\psi_{Y,1}(P) \leq C_1 \sqrt{b_{(1)}} \| p_Y - p_Z \|_2^2 \quad \text{for some $C_1 >0$.}
\end{align*}
By denoting the $k$th component of $p_Y$ and $p_Z$ by $p_Y(k)$ and $p_Z(k)$, respectively, note that 
\begin{align*}
\mE_P[\overline{h}_{\text{ts}}(Y_1,Y_2;Z_1,Z_2) | Y_1]  = \sum_{k=1}^d [\ind(Y_1=k) - p_Z(k)][p_Y(k) - p_Z(k)] 
\end{align*}
and so $\psi_{Y,1}(P)$, which is the variance of the above expression, becomes
\begin{align*}
\psi_{Y,1}(P) = \mE_P \bigg[ \bigg( \sum_{k=1}^d [\ind(Y_1=k) - p_Y(k)][p_Y(k) - p_Z(k)]  \bigg)^2  \bigg].
\end{align*}
Furthermore, observe that 
\begin{align*}
\psi_{Y,1}(P) ~ \overset{(i)}{\leq}~ & 2 \mE_P \bigg[ \bigg( \sum_{k=1}^d \ind(Y_1=k) [p_Y(k) - p_Z(k)]  \bigg)^2 \bigg] + 2  \bigg( \sum_{k=1}^d p_Y(k)[p_Y(k) - p_Z(k)]  \bigg)^2 \\[.5em]
= ~ & 2 \sum_{k=1}^d p_Y(k) [p_Y(k) - p_Z(k)]^2  + 2  \bigg( \sum_{k=1}^d p_Y(k)[p_Y(k) - p_Z(k)]  \bigg)^2 \\[.5em]
\overset{(ii)}{\leq}  ~ &2 \sqrt{\sum_{k=1}^d p_Y^2(k)}\sqrt{\sum_{k=1}^d  [p_Y(k) - p_Z(k)]^4} + 2 \sum_{k=1}^d p_Y^2(k) \sum_{k=1}^d  [p_Y(k) - p_Z(k)]^2 \\[.5em]
\overset{(iii)}{\leq} ~ & 4\sqrt{b_{(1)}} \| p_Y - p_Z \|_2^2,
\end{align*} 
where $(i)$ is based on $(x+y)^2 \leq 2x^2 + 2y^2$, $(ii)$ uses Cauchy-Schwarz inequality and $(iii)$ uses the monotonicity of $\ell_p$ norm (specifically, $\ell_4 \leq \ell_2$) as well as the fact that $\|p_Y\|_2^2 \leq \|p_Y\|_2$. By symmetry, we can also have that 
\begin{align*}
\psi_{Z,1}(P) ~ \leq ~ 4 \sqrt{b_{(1)}} \| p_Y - p_Z \|_2^2. 
\end{align*}
Now focusing on the third line of the claim~(\ref{Eq: Three conditions for two-sample multinomials}), recall that 
\begin{align*}
\psi_{YZ,2}(P) &:= \max \{ \mE_P[g_{\text{Multi}}^2(Y_1,Y_2)], \  \mE_P[g_{\text{Multi}}^2(Y_1,Z_1)], \ \mE_P[g_{\text{Multi}}^2(Z_1,Z_2)] \}
\end{align*}
and by noting that $g_{\text{Multi}}(x,y)$ is either one or zero,
\begin{align*}
& \mE_P[g_{\text{Multi}}^2(Y_1,Y_2)]  = \sum_{k=1}^d p_Y^2(k), \\[.5em]
& \mE_P[g_{\text{Multi}}^2(Z_1,Z_2)]  = \sum_{k=1}^d p_Z^2(k) \quad \text{and}  \\[.5em]
&  \mE_P[g_{\text{Multi}}^2(Y_1,Z_1)]  = \sum_{k=1}^d p_Y(k) p_Z(k) \leq \frac{1}{2} \sum_{k=1}^d p_Y^2(k) +  \frac{1}{2} \sum_{k=1}^d p_Z^2(k),
\end{align*}
where the last inequality uses $xy \leq x^2/2 + y^2/2$. This clearly shows that $\psi_{YZ,2} \leq b_{(1)}$, which confirms the claim (\ref{Eq: Three conditions for two-sample multinomials}). Since the expectation of $U_{n_1,n_2}$ is $\|p_Y - p_Z\|_2^2$, one may see that 
\begin{align*}
\mE_P[U_{n_1,n_2}]  \geq \epsilon_{n_1,n_2}^2 & \geq C_1 \frac{\sqrt{b_{(1)}}}{\alpha^{1/2} \beta} \left( \frac{1}{n_1} + \frac{1}{n_2} \right) \\[.5em]
& \geq C_2 \sqrt{\max \Bigg\{ \frac{\psi_{Y,1}(P)}{\beta n_1}, \ \frac{\psi_{Z,1}(P)}{\beta n_2}, \  \frac{\psi_{YZ,2}(P)}{\alpha\beta}\left( \frac{1}{n_1} + \frac{1}{n_2}\right)^2 \Bigg\}}.
\end{align*} 
Now we apply Theorem~\ref{Theorem: Two-Sample U-statistic} and finish the proof of Proposition~\ref{Proposition: Multinomial Two-Sample Testing}.

\section{Proof of Proposition~\ref{Proposition: Minimum Separation for Two-Sample Multinomial Testing}} \label{Section: Proposition: Minimum Separation for Two-Sample Multinomial Testing}
We first note that Proposition~\ref{Proposition: Multinomial Two-Sample Testing} establishes an upper bound for the minimum separation as $\epsilon_{n_1,n_2}^\dagger \lesssim b_{(1)}^{1/4}n_1^{-1/2}$ where $n_1 \leq n_2$. Hence once we identify a lower bound such that $\epsilon_{n_1,n_2}^\dagger  \gtrsim b_{(1)}^{1/4}n_1^{-1/2}$, the proof is completed. As briefly explained in the main text, our strategy to prove this result is to consider the one-sample problem, which is conceptually easier than the two-sample problem, and establish the matching lower bound. In the one-sample problem, we assume that $p_Z$ is known and observe $n_1$ samples from the other distribution $p_Y$. Based on these $n_1$ samples, we want to test whether $p_Y = p_Z$ or $\|p_Y - p_Z\|_2 \geq \epsilon_{n_1}$. As formalized by \cite{arias2018remember} (see their Lemma 1), the one-sample problem can be viewed as a special case of the two-sample problem where one of the sample sizes is taken to be infinite and thus the minimum separation for the one-sample problem is always smaller than or equal to that for the two-sample problem. This means that if the minimum separation for the one-sample problem, denoted by $\epsilon_{n_1}^\dagger$, satisfies $\epsilon_{n_1}^\dagger \gtrsim b_{(1)}^{1/4} n_1^{-1/2}$, then we also have that $\epsilon_{n_1,n_2}^\dagger  \gtrsim b_{(1)}^{1/4}n_1^{-1/2}$. In the end, it suffices to verify $\epsilon_{n_1}^\dagger \gtrsim b_{(1)}^{1/4} n_1^{-1/2}$ to complete the proof. We show this result based on the standard lower bound technique due to \cite{ingster1987minimax,ingster1993asymptotically}.

\paragraph{$\bullet$ Ingster's method for the lower bound.} Let us recall from Section~\ref{Section: Minimax optimality} that the minimax type II error is given by
\begin{align*}
R_{n,\epsilon_n}^\dagger := \inf_{\phi \in \Phi_{n,\alpha}} \sup_{P \in \mathcal{P}_1(\epsilon_n)} \mP_P^{(n)} (\phi = 0).
\end{align*}
For $P_1,\ldots,P_N \in \mathcal{P}_1(\epsilon_n)$, define a mixture distribution $Q$ given by 
\begin{align*}
Q(A) = \frac{1}{N} \sum_{i=1}^N P_i^n(A).
\end{align*}
Given $n$ i.i.d.~observations $X_1,\ldots,X_n$, we denote the likelihood ratio between $Q$ and the null distribution $P_0$ by
\begin{align*}
L_n = \frac{dQ}{dP_0^n} = \frac{1}{N} \sum_{i=1}^N \prod_{j=1}^n \frac{p_i(X_j)}{p_0(X_j)}. 
\end{align*}
Then one can relate the variance of the likelihood ratio to the minimax type II error as follows. 
\begin{lemma}[Lower bound] \label{Lemma: Ingster's method}
	Let $0 < \beta < 1 - \alpha$. If 
	\begin{align*}
	\mE_{P_0}[L_n^2]  \leq 1 + 4(1 - \alpha - \beta)^2,
	\end{align*}
	then $R_{n,\epsilon_n}^\dagger  \geq \beta$. 
	\begin{proof}
		We present the proof of this result only for completeness. Note that $\mP_{P_0}^n(\phi = 1) \leq \alpha$ for $\phi \in \Phi_{n,\alpha}$. Thus
		\begin{align*}
		R_{n,\epsilon_n}^\dagger \geq ~&\inf_{\phi \in \Phi_{n,\alpha}} \mP_Q(\phi = 0) = \inf_{\phi \in \Phi_{n,\alpha}} \big[ \mP_{P_0}(\phi = 0) + \mathbb{P}_Q(\phi = 0) - \mathbb{P}_{P_0}(\phi = 0)  \big] \\[.5em]
		\overset{(i)}{\geq} ~ &  1 - \alpha + \inf_{\phi \in \Phi_{n,\alpha}} \big[ \mathbb{P}_Q(\phi = 0) - \mathbb{P}_{P_0}(\phi = 0)  \big]  \\[.5em]
		\overset{(ii)}{\geq} ~ & 1 - \alpha - \sup_{A} \big| \mathbb{P}_Q(A) -  \mathbb{P}_{P_0}(A) \big| \\[.5em]
		\overset{(iii)}{=} ~ & 1 - \alpha - \frac{1}{2} \big\| Q - P_0^n \big\|_1. 
		\end{align*}
		where $(i)$ uses the fact that $\mP_{P_0}^n(\phi = 1) \leq \alpha$, $(ii)$ follows by taking the supremum over all measurable sets, $(iii)$ uses the alternative expression for the total variation distance in terms of $L_1$-distance. The result then follows by noting that
		\begin{align*}
		\big\| Q - P_0^n \big\|_1 = \mE_{P_0} [ |L_n(X_1,\ldots,X_n) - 1| ] \leq \sqrt{\mE_{P_0}[L_n^2(X_1,\ldots,X_n)] - 1}.
		\end{align*}
		This proves Lemma~\ref{Lemma: Ingster's method}.
	\end{proof}
\end{lemma}
Next we apply this method to find a lower bound for $\epsilon_{n_1}^\dagger$. To apply Lemma~\ref{Lemma: Ingster's method}, we first construct $Q$ and $P_0$. 

\paragraph{$\bullet$ Construction of $Q$ and $P_0$.}
Suppose that $p_Z$ is the uniform distribution over $\mathbb{S}_d$, that is $p_Z(k) = 1/d$ for $k=1,\ldots,d$. Let $\widetilde{\zeta} = \{\widetilde{\zeta}_1,\ldots,\widetilde{\zeta}_d\}$ be dependent Rademacher random variables uniformly distributed over $\{-1,1\}^d$ such that $\sum_{i=1}^d \widetilde{\zeta}_i = 0$ where $d$ is assumed to be even. More formally we define such a set by
\begin{align} \label{Eq: definition of hypercube}
\mathcal{M}_d:= \big\{ x \in \{-1,1\}^d : \sum_{i=1}^d x_i =0 \big\}.
\end{align}
If $d$ is odd, then we set $\widetilde{\zeta}_d = 0$ and the proof follows similarly. Given $\widetilde{\zeta} \in \mathcal{M}_d$, let us define a distribution $p_{\widetilde{\zeta}}$ as
\begin{align*}
p_{\widetilde{\zeta}}(k)  := p_Z(k) + \delta \sum_{i=1}^d \widetilde{\zeta}_i \ind(k=i),
\end{align*}
where $\delta$ is specified later but $\delta \leq 1/d$. There are $N$ such distributions where $N$ is the cardinality of $\mathcal{M}_d$ and we denote them by $p_{\widetilde{\zeta}(1)}, \ldots, p_{\widetilde{\zeta}(N)}$. By construction we make three observations. First $p_{\widetilde{\zeta}}$ is a proper distribution as each component $p_{\widetilde{\zeta}}(k)$ is non-negative and $\sum_{k=1}^d p_{\widetilde{\zeta}}(k)=1$. Second the $\ell_2$ distance between $p_{\widetilde{\zeta}}$ and $p_Z$ is 
\begin{align} \label{Eq: l_2 norm}
\| p_{\widetilde{\zeta}} - p_Z \|_2 = \delta \sqrt{d}.
\end{align}
Third we see that $b_{(1)} = \max\{\|p_Z\|_2^2, \|p_{\widetilde{\zeta}}\|_2^2 \}$ is lower and upper bounded by
\begin{align} \label{Eq: bounds for b_1}
\frac{1}{d} \leq b_{(1)} \leq \frac{2}{d},
\end{align}
which can be verified based on Cauchy-Schwarz inequality and the fact that $\delta \leq 1/d$. Finally we denote the uniform mixture of $p_{\widetilde{\zeta}(1)},\ldots,p_{\widetilde{\zeta}(N)}$ by
\begin{align*}
Q := \frac{1}{N} \sum_{i=1}^N  p_{\widetilde{\zeta}(i)}
\end{align*}
and let $P_0 = p_Z$. Having $Q$ and $P_0$ at hand, we are now ready to compute the expected value of the squared likelihood ratio.

\paragraph{$\bullet$ Calculation of $\mE_{P_0}[L_n^2] $.}
For each $\widetilde{\zeta}_{(i)} \in \mathcal{M}_d$ and $i=1,\ldots,N$, let us denote the components of $\widetilde{\zeta}_{(i)}$ by $\{\widetilde{\zeta}_{1,(i)},\ldots,\widetilde{\zeta}_{d,(i)}\}$. Based on this notation as well as the definition of $Q$ and $P_0$, the squared the likelihood ratio $L_n^2$ can be written as
\begin{align*}
L_n^2 ~=~ & \frac{1}{N^2} \sum_{i_1=1}^N \sum_{i_2=1}^N \prod_{j=1}^{n_1} \frac{p_{\widetilde{\zeta}(i_1)}(X_j)p_{\widetilde{\zeta}(i_2)}(X_j)}{p_0(X_j)p_0(X_j)} \\[.5em]
=~ &  \frac{1}{N^2} \sum_{i_1=1}^N \sum_{i_2=1}^N \prod_{j=1}^{n_1} \frac{\{  1/d + \delta \sum_{k=1}^d \widetilde{\zeta}_{k(i_1)} \ind(X_j=k) \} \{ 1/d + \delta \sum_{k=1}^d \widetilde{\zeta}_{k(i_2)} \ind(X_j=k) \} }{1/d^2} \\[.5em]
=~ & \frac{1}{N^2} \sum_{i_1=1}^N \sum_{i_2=1}^N \prod_{j=1}^{n_1} \bigg\{ 1 + d\delta \sum_{k=1}^d \widetilde{\zeta}_{k(i_1)} \ind(X_j=k) \bigg\} \bigg\{ 1 + d\delta \sum_{k=1}^d \widetilde{\zeta}_{k(i_2)} \ind(X_j=k) \bigg\}.
\end{align*}
Now by taking the expectation under $P_0$, it can be seen that
\begin{align*}
\mE_{P_0}[L_n^2] ~=~ & \frac{1}{N^2} \sum_{i_1=1}^N \sum_{i_2=1}^N \Bigg( 1 + d \delta^2 \sum_{k=1}^d  \widetilde{\zeta}_{k(i_1)} \widetilde{\zeta}_{k(i_2)}  \Bigg)^{n_1} \\[.5em]
\leq ~ & \frac{1}{N^2} \sum_{i_1=1}^N \sum_{i_2=1}^N  \exp \Bigg( n_1 d \delta^2 \sum_{k=1}^d  \widetilde{\zeta}_{k(i_1)} \widetilde{\zeta}_{k(i_2)}\Bigg),
\end{align*}
where the inequality uses $1 + x \leq e^x$ for all $x \in \mathbb{R}$. By letting $\widetilde{\zeta}^\ast$ be i.i.d.~copy of $\widetilde{\zeta}$, we may see that 
\begin{align*}
\frac{1}{N^2} \sum_{i_1=1}^N \sum_{i_2=1}^N  \exp \Bigg( n_1 d \delta^2 \sum_{k=1}^d  \widetilde{\zeta}_{k(i_1)} \widetilde{\zeta}_{k(i_2)}\Bigg) = \mE_{\widetilde{\zeta},\widetilde{\zeta}^\ast} \Big[ \exp \Big( n_1 d \delta^2 \big\langle \widetilde{\zeta}, \widetilde{\zeta}^\ast \big\rangle \Big) \Big].
\end{align*}
Moreover $\{\widetilde{\zeta}_1,\ldots\widetilde{\zeta}_d\}$ are negatively associated \citep{dubhashi1998balls}. Hence applying Lemma 2 of \cite{dubhashi1998balls} yields
\begin{align*}
\mE_{P_0}[L_n^2] ~ \leq ~ & \mE_{\widetilde{\zeta},\widetilde{\zeta}^\ast} \Big[ \exp \Big( n_1 d \delta^2 \big\langle \widetilde{\zeta}, \widetilde{\zeta}^\ast \big\rangle \Big) \Big] \\[.5em]
\leq ~ & \prod_{i=1}^d  \mE_{\widetilde{\zeta}_i,\widetilde{\zeta}^\ast_i} \Big[ \exp \Big( n_1 d \delta^2 \widetilde{\zeta}_i \widetilde{\zeta}^\ast_i \Big) \Big] = \prod_{i=1}^d \text{cosh} (n_1 d \delta^2) \\[.5em]
\overset{(i)}{\leq} ~ & \prod_{i=1}^d e^{n_1^2 d^2 \delta^4/2} = e^{n_1^2 d^3 \delta^4/2}, 
\end{align*}
where $(i)$ uses the inequality $\text{cosh}(x) \leq e^{x^2/2}$ for all $x \in \mathbb{R}$.

\paragraph{$\bullet$ Completion of the proof.} 
Based on this upper bound, we have from Lemma~\ref{Lemma: Ingster's method} that if 
\begin{align*}
\delta \leq \frac{1}{\sqrt{n_1}d^{3/4}}\big[\log \big\{ 1 + 4 (1-\alpha-\beta)^2 \big\}\big]^{1/4}
\end{align*}
the minimax type II error is lower bounded by $\beta$. Furthermore, based on the expression for the $\ell_2$ norm in (\ref{Eq: l_2 norm}) and the bound for $b_{(1)}$ in (\ref{Eq: bounds for b_1}). The above condition is further implied by
\begin{align*}
\epsilon_{n_1} \leq \frac{b_{(1)}^{1/4}}{\sqrt{n_1}} \big[\log \big\{ 1 + 4 (1-\alpha-\beta)^2 \big\}\big]^{1/4}.
\end{align*}
This completes the proof of Proposition~\ref{Proposition: Minimum Separation for Two-Sample Multinomial Testing}.

\section{Proof of Proposition~\ref{Proposition: Two-Sample Testing for Holder Densities}} \label{Section: Proposition: Two-Sample Testing for Holder Densities}
The proof of Proposition~\ref{Proposition: Two-Sample Testing for Holder Densities} is fairly straightforward based on Proposition~\ref{Proposition: Multinomial Two-Sample Testing} and Lemma 3 of \cite{arias2018remember}. For two vectors $\mathbf{v} = (v_1,\ldots,v_d) \in \mathbb{R}^d$ and $\mathbf{w} =(w_1,\ldots,w_d) \in \mathbb{R}^d$ where $v_i \leq w_i$ for all $i$, we borrow the notation from \cite{arias2018remember} and denote the hyperrectangle by 
\begin{align*}
\boldsymbol{[}\mathbf{v},\mathbf{w}\boldsymbol{]} = \prod_{i=1}^d [v_i,w_i]. 
\end{align*}
Recall that $\kappa_{(1)} = \floor{n_1^{2/(4s+d)}}$ and define $\boldsymbol{H}_{\boldsymbol{\ell}}:= \boldsymbol{[}(\boldsymbol{\ell}-1)/\kappa_{(1)}, \boldsymbol{\ell} /\kappa_{(1)}\boldsymbol{]}$ where $\boldsymbol{\ell} \in \{1,2,\ldots, \kappa_{(1)}\}^d$, 
\begin{align*}
p_{Y}(\boldsymbol{\ell}) := \int_{\boldsymbol{H}_{\boldsymbol{\ell}}} f_Y(t)dt \quad \text{and} \quad p_{Z}(\boldsymbol{\ell}) := \int_{\boldsymbol{H}_{\boldsymbol{\ell}}} f_Z(t)dt. 
\end{align*} 
Since both $f_Y$ and $f_Z$ are in H\"{o}lder's density class $\mathcal{P}_{\text{H\"{o}lder}}^{(d,s)}$ where $|\!|\!|f_Y|\!|\!|_{\infty} \leq L$ and $|\!|\!|f_Z|\!|\!|_{\infty} \leq L$, it is clear to see that 
\begin{align*}
p_Y(\boldsymbol{\ell}) \leq |\!|\!|f_Y|\!|\!|_{\infty} \kappa_{(1)}^{-d} \leq L \kappa_{(1)}^{-d} \quad \text{and} \quad p_Z(\boldsymbol{\ell}) \leq |\!|\!|f_Z|\!|\!|_{\infty} \kappa_{(1)}^{-d} \leq L \kappa_{(1)}^{-d} \quad \text{for all $\boldsymbol{\ell}$.} 
\end{align*}
This gives 
\begin{align} \label{Eq: b_1 bound}
b_{(1)} = \max\{\|p_Y\|_2^2, \|p_Z\|_2^2\} \leq L \kappa_{(1)}^{-d}.
\end{align}
Based on Lemma 3 of \cite{arias2018remember}, one can find a constant $C_1 > 0$ such that 
\begin{align} \label{Eq: L2 norm lower bound}
\|p_Y - p_Z\|_2^2 \geq C_1 \kappa_{(1)}^{-d} \epsilon_{n_1,n_2}^2, 
\end{align}
where $\epsilon_{n_1,n_2}$ is the lower bound for $\|f_Y - f_Z\|_{L_2}$. By combining (\ref{Eq: b_1 bound}) and (\ref{Eq: L2 norm lower bound}), the condition of Proposition~\ref{Proposition: Multinomial Two-Sample Testing} is satisfied when 
\begin{align*}
\kappa_{(1)}^{-d} \epsilon_{n_1,n_2}^2 ~\geq~ C_2 \frac{L^{1/2} \kappa_{(1)}^{-d/2}}{\alpha^{1/2} \beta} \left( \frac{1}{n_1} + \frac{1}{n_2} \right).
\end{align*}
Equivalently, 
\begin{align*}
 \epsilon_{n_1,n_2} ~ \geq ~ C_3 \frac{L^{1/4} \kappa_{(1)}^{d/4}}{\alpha^{1/4}\beta^{1/2}} \left( \frac{1}{n_1} + \frac{1}{n_2} \right)^{1/2}.
\end{align*}
Since $\kappa_{(1)} = \floor{n_1^{2/(4s+d)}}$ and we assume $n_1 \leq n_2$, the above inequality is further implied by
\begin{align*}
\epsilon_{n_1,n_2} ~ \geq ~ \frac{C_4}{\alpha^{1/4}\beta^{1/2}} \left( \frac{1}{n_1} + \frac{1}{n_2}\right)^{\frac{2s}{4s+d}},
\end{align*}
where $C_4$ is a constant that may depend on $s,d,L$. This completes the proof of Proposition~\ref{Proposition: Two-Sample Testing for Holder Densities}.

\section{Proof of Theorem~\ref{Theorem: U-statistic for independence testing}} \label{Section: Theorem: U-statistic for independence testing}
The proof of Theorem~\ref{Theorem: U-statistic for independence testing} is similar to that of Theorem~\ref{Theorem: Two-Sample U-statistic}. First we verify that the permuted $U$-statistic $U_n^{\pi}$, which can be recalled from (\ref{Eq: permuted U-statistic}), has zero expectation. By the linearity of expectation, the problem boils down to showing 
\begin{align*}
	\mE_{\pi} [h_{\text{in}}\{(Y_1,Z_{\pi_1}), (Y_2,Z_{\pi_2}), (Y_3,Z_{\pi_3}), (Y_4,Z_{\pi_4})\} | \mathcal{X}_n] = 0.
\end{align*}
Since $Y_1,\ldots,Y_4$ are constant under permutations, it further boils down to proving
\begin{align*}
	\mE_{\pi} [g_Z(Z_{\pi_1},Z_{\pi_2}) + g_Z(Z_{\pi_3},Z_{\pi_4}) - g_Z(Z_{\pi_1},Z_{\pi_3}) - g_Z(Z_{\pi_2},Z_{\pi_4})  | \mathcal{X}_n] = 0.
\end{align*}
In fact, this equality is clear by noting that $\mE_{\pi}[g_Z(Z_{\pi_i},Z_{\pi_j})]$ is invariant to the choice of $(i,j) \in \mathbf{i}_2^n$, which leads to $\mE_{\pi}[U_n^{\pi}|\mathcal{X}_n] = 0$. Therefore we can focus on the simplified condition~(\ref{Eq: Sufficient Condition II}) to proceed.

The rest of the proof is split into two parts. In each part, we prove the following conditions separately,
\begin{align} \label{Eq: Condition I'}
	&\mE_P[U_{n}] \geq 2\sqrt{\frac{2\mV_P[U_{n}]}{\beta}} \quad \text{and} \\[.5em] \label{Eq: Condition II'}
	& \mE_P[U_{n}] \geq 2\sqrt{\frac{2\mE_P[\mV_{\pi}\{U_{n}^\pi| \mathcal{X}_n \}]}{\alpha\beta}}.
\end{align}
We then complete the proof of Theorem~\ref{Theorem: U-statistic for independence testing} by noting that (\ref{Eq: Condition I'}) and (\ref{Eq: Condition II'}) imply the simplified condition (\ref{Eq: Sufficient Condition II}).

\vskip 1em

\noindent \textbf{Part 1. Verification of condition~(\ref{Eq: Condition I'}).} This part verifies condition~(\ref{Eq: Condition I'}). The main ingredient of this part of the proof is the explicit variance formula of a $U$-statistic \citep[e.g.~page 12 of][]{lee1990u}. Following the notation of \cite{lee1990u}, we define $\check{\sigma}_i^2$ to be the variance of the conditional expectation by 
\begin{align*}
	\check{\sigma}_i^2 := \mV_P [ \mE_P \{ \overline{h}_{\text{in}} (x_1,\ldots,x_i, X_{i+1},\ldots, X_{4}) \}] \quad \text{for $1 \leq i \leq 4$.}
\end{align*}
Then the variance of $U_{n}$ is given by
\begin{align*}
	\mV_P[U_n] = \sum_{i=1}^4  \binom{4}{i} \binom{n-4}{4-i} \binom{n}{4}^{-1} \check{\sigma}_i^2.
\end{align*}
By the law of total variance, it can be seen that $\check{\sigma}_i^2 \leq \check{\sigma}_4^2$ for all $1 \leq i \leq 4$, which leads to an upper bound for $\mV_P[U_n] $ as
\begin{align*}
	\mV_P[U_n] \leq C_1 \frac{\check{\sigma}_1^2}{n} + C_2 \frac{\check{\sigma}_4^2}{n^2}.
\end{align*}  
Now applying Jensen's inequality, repeatedly, yields
\begin{align*}
	\check{\sigma}_{4}^2 \leq \mE_P [ \overline{h}_{\text{in}}^2(X_1,X_2,X_3,X_4)] \leq \mE_P [ h_{\text{in}}^2(X_1,X_2,X_3,X_4)] \leq C_3 \psi_{2}'(P).
\end{align*}
Then by noting that $\check{\sigma}_{1}^2$ corresponds to the notation $\psi_{1}'(P)$, we have that
\begin{align*}
	\mV_P[U_{n_1,n_2}]  \leq C_1 \frac{\psi_{1}'(P)}{n} + C_4 \frac{\psi_{2}'(P)}{n^2}.
\end{align*}
Therefore condition~(\ref{Eq: Condition I'}) is satisfied by taking constant $C$ sufficiently large in Theorem~\ref{Theorem: U-statistic for independence testing}. 

\vskip 1em

\noindent \textbf{Part 2. Verification of condition~(\ref{Eq: Condition II'}).} This part verifies condition~(\ref{Eq: Condition II'}). As mentioned in the main text, the permuted $U$-statistic $U_{n}^{\pi}$ mimics the behavior of $U_n$ under the null hypothesis. Hence one can expect that the variance of $U_{n}^{\pi}$ is similarly bounded by $\psi_{2}'(P)n^{-2}$ up to some constant as $\psi_{1}'(P)$ becomes zero under the null hypothesis. To prove this statement, we first introduce some notation. Let us define a set of indices $\mathsf{J}_{\text{total}}:= \{(i_1,i_2,i_3,i_4,i_1',i_2',i_3',i_4') \in \mathbb{N}^8_+ : (i_1,i_2,i_3,i_4) \in \mathbf{i}_4^n, (i_1',i_2',i_3',i_4') \in \mathbf{i}_4^n \}$ and let $\mathsf{J}_{\text{A}} := \{ (i_1,i_2,i_3,i_4,i_1',i_2',i_3',i_4') \in \mathsf{J}_{\text{total}} : \#|\{i_1,i_2,i_3,i_4\}\cap \{i_1',i_2',i_3',i_4'\} | \leq 1 \}$ and $\mathsf{J}_{\text{A}^c} := \{ (i_1,i_2,i_3,i_4,i_1',i_2',i_3',i_4') \in \mathsf{J}_{\text{total}} : \#|\{i_1,i_2,i_3,i_4\}\cap \{i_1',i_2',i_3',i_4'\} | > 1 \}$. By construction, it is clear that $\mathsf{J}_{\text{total}} = \mathsf{J}_{\text{A}} \cup \mathsf{J}_{\text{A}^c}$. To shorten the notation, we simply write 
\begin{align*}
	h_{\text{in}}(x_1,x_2,x_3,x_4) = h_{\text{in},Y}(y_1,y_2,y_3,y_4) h_{\text{in},Z}(z_1,z_2,z_3,z_4), 
\end{align*}
where $h_{\text{in},Y}(y_1,y_2,y_3,y_4):= g_Y(y_1,y_2) + g_Y(y_3,y_4) - g_Y(y_1,y_3) - g_Y(y_2,y_4)$ and $h_{\text{in},Z}(z_1,z_2,z_3,z_4) := g_Z(z_1,z_2) + g_Z(z_3,z_4) - g_Z(z_1,z_3) - g_Z(z_2,z_4)$. Since $U_{n}^{\pi}$ is centered, our interest is in bounding $\mE_P[\mE_{\pi}\{ (U_n^{\pi})^2 | \mathcal{X}_n \}]$. Focusing on the conditional expectation inside, observe that 
\begin{align*}
	\mE_{\pi} [(U_n^{\pi})^2 | \mathcal{X}_n] &~=~  \frac{1}{n_{(4)}^2}  \sum_{(i_1,\ldots,i_4') \in \mathsf{J}_{\text{total}}} h_{\text{in},Y}(Y_{i_1},Y_{i_2},Y_{i_3},Y_{i_4})h_{\text{in},Y}(Y_{i_1'},Y_{i_2'},Y_{i_3'},Y_{i_4'}) \\
	&\times \mE_{\pi} \big[ h_{\text{in},Z}(Z_{\pi_{i_1}},Z_{\pi_{i_2}},Z_{\pi_{i_3}},Z_{\pi_{i_4}})h_{\text{in},Z}(Z_{\pi_{i_1'}},Z_{\pi_{i_2'}},Z_{\pi_{i_3'}},Z_{\pi_{i_4'}}) | \mathcal{X}_n \big] \\[.5em]
	&~= ~ (I') + (II'),
\end{align*}
where 
\begin{align*}
	(I') ~:=~ & \frac{1}{n_{(4)}^2}  \sum_{(i_1,\ldots,i_4') \in \mathsf{J}_{\text{A}}} h_{\text{in},Y}(Y_{i_1},Y_{i_2},Y_{i_3},Y_{i_4})h_{\text{in},Y}(Y_{i_1'},Y_{i_2'},Y_{i_3'},Y_{i_4'}) \\
	&\times \mE_{\pi} \big[ h_{\text{in},Z}(Z_{\pi_{i_1}},Z_{\pi_{i_2}},Z_{\pi_{i_3}},Z_{\pi_{i_4}})h_{\text{in},Z}(Z_{\pi_{i_1'}},Z_{\pi_{i_2'}},Z_{\pi_{i_3'}},Z_{\pi_{i_4'}}) | \mathcal{X}_n \big], \\[.8em]
	(II') ~:=~ & \frac{1}{n_{(4)}^2}  \sum_{(i_1,\ldots,i_4') \in \mathsf{J}_{\text{A}^c}} h_{\text{in},Y}(Y_{i_1},Y_{i_2},Y_{i_3},Y_{i_4})h_{\text{in},Y}(Y_{i_1'},Y_{i_2'},Y_{i_3'},Y_{i_4'}) \\
	&\times \mE_{\pi} \big[ h_{\text{in},Z}(Z_{\pi_{i_1}},Z_{\pi_{i_2}},Z_{\pi_{i_3}},Z_{\pi_{i_4}})h_{\text{in},Z}(Z_{\pi_{i_1'}},Z_{\pi_{i_2'}},Z_{\pi_{i_3'}},Z_{\pi_{i_4'}})  | \mathcal{X}_n \big].
\end{align*}
We now claim that the first term $(I') = 0$, which is critical to obtain a faster rate $n^{-2}$ rather than $n^{-1}$ in the bound
\begin{align*} 
	\mE_P[\mV_{\pi}\{U_n^{\pi} | \mathcal{X}_n\}] \leq C_2 \frac{\psi_2'(P)}{n^2}.
\end{align*}
However we have already proved in the second part of the proof of Theorem~\ref{Theorem: Two-Sample U-statistic} that 
\begin{align*}
	\mE_{\pi} \big[ h_{\text{in},Z}(Z_{\pi_{i_1}},Z_{\pi_{i_2}},Z_{\pi_{i_3}},Z_{\pi_{i_4}})h_{\text{in},Z}(Z_{\pi_{i_1'}},Z_{\pi_{i_2'}},Z_{\pi_{i_3'}},Z_{\pi_{i_4'}}) | \mathcal{X}_n \big] = 0,
\end{align*}
whenever $(i_1,\ldots,i_4') \in \mathsf{J}_{A}$. This concludes $(I') = 0$ and so $\mE_{\pi} [(U_n^{\pi})^2 | \mathcal{X}_n] = (II')$. To bound $\mE_P[(II')]$, we make an observation that for any $1 \leq i_1 \neq i_2, i_1' \neq i_2' \leq n$, 
\begin{align*}
	& \big| \mE_P \big[ g_Y(Y_{i_1},Y_{i_2}) g_Y(Y_{i_1'},Y_{i_2'}) \mE_{\pi} \big\{  g_Z(Z_{\pi_{i_1}},Z_{\pi_{i_2}}) g_Z(Z_{\pi_{i_1'}},Z_{\pi_{i_2'}})  | \mathcal{X}_n \big\}  \big] \big| \\[.5em]
	\overset{(i)}{=} ~ &  \big| \mE_\pi \big[ \mE_P \big\{ g_Y(Y_{i_1},Y_{i_2}) g_Y(Y_{i_1'},Y_{i_2'}) g_Z(Z_{\pi_{i_1}},Z_{\pi_{i_2}}) g_Z(Z_{\pi_{i_1'}},Z_{\pi_{i_2'}}) | \pi \big\} \big] \big|  \\[.5em]
	\overset{(ii)}{\leq} ~ & \frac{1}{2}\mE_\pi \big[ \mE_P \big\{ g_Y^2(Y_{i_1},Y_{i_2}) g_Z^2(Z_{\pi_{i_1}},Z_{\pi_{i_2}}) | \pi \big\} \big]  + \frac{1}{2}\mE_\pi \big[ \mE_P \big\{ g_Y^2(Y_{i_1'},Y_{i_2'}) g_Z^2(Z_{\pi_{i_1'}},Z_{\pi_{i_2'}}) | \pi \big\} \big]  \\[.5em]
	\overset{(iii)}{\leq} ~ & \psi_{2}'(P),
\end{align*}
where $(i)$ uses the law of total expectation, $(ii)$ uses the basic inequality $xy \leq x^2/2 + y^2/2$ and $(iii)$ follows by the definition of $\psi_2'(P)$. Based on this observation, it is not difficult to see that for any $(i_1,\ldots, i_4') \in \mathsf{J}_{\text{total}}$, 
\begin{align*}
	\big |\mE_P [\mE_{\pi}\{ h_{\text{in}}(X_{\pi_{i_1}}, X_{\pi_{i_2}},X_{\pi_{i_3}},X_{\pi_{i_4}} )h_{\text{in}}(X_{\pi_{i_1'}}, X_{\pi_{i_2'}},X_{\pi_{i_3'}},X_{\pi_{i_4'}} ) | \mathcal{X}_n \} ] \big| \leq C_1 \psi_{2}'(P).
\end{align*}
Therefore, by counting the number of elements in $\mathsf{J}_{\text{A}^c}$, 
\begin{align*}
	\mE_P[\mV_{\pi}\{U_n^\pi | \mathcal{X}_n \}] ~ = ~ & \mE_P[(II')]  \\[.5em]
	\leq ~ & C_1 \psi_{2}'(P) \frac{1}{n_{(4)}^2}  \sum_{(i_1,\ldots,i_4') \in \mathsf{J}_{\text{A}^c}} 1 \\[.5em]
	\leq ~ & C_2 \frac{\psi_{2}'(P)}{n^2}.
\end{align*}

Now by taking constant $C$ in Theorem~\ref{Theorem: U-statistic for independence testing} sufficiently large, one may see that condition~(\ref{Eq: Condition II'}) is satisfied. This completes the proof of Theorem~\ref{Theorem: U-statistic for independence testing}.

\section{Proof of Proposition~\ref{Proposition: Multinomial Independence Testing}} \label{Section: Proposition: Multinomial Independence Testing}
To prove Proposition~\ref{Proposition: Multinomial Independence Testing}, it suffices to verify that the two inequalities~(\ref{Eq: Aim of the proof of multinomial independence testing}) hold. Then the result follows by Theorem~\ref{Theorem: U-statistic for independence testing}. To start with the first inequality in (\ref{Eq: Aim of the proof of multinomial independence testing}), we want to upper bound $\psi_{1}'(P)$ as $\psi_{1}'(P) \leq C_1 \sqrt{b_{(2)}}\|p_{YZ} - p_Yp_Z\|_2^2$. A little algebra shows that 
\begin{align*}
& \mE_P[\overline{h}_{\text{in}}(X_1,X_2,X_3,X_4)|X_1]  - 4 \|p_{YZ} - p_Yp_Z\|_2^2 \\[.5em]
= ~ & 2 \sum_{k=1}^{d_1} \sum_{k'=1}^{d_2}  \big[ \ind(Y_1=k) \ind(Z_1=k') - p_{YZ}(k,k') \big] \big[p_{YZ}(k,k') - p_Y(k)p_Z(k') \big] \\[.5em]
 - & 2  \sum_{k=1}^{d_1} \sum_{k'=1}^{d_2}  \big[ \ind(Y_1=k) - p_Y(k) \big] p_Z(k') \big[p_{YZ}(k,k') - p_Y(k)p_Z(k') \big] \\[.5em]
  - & 2  \sum_{k=1}^{d_1} \sum_{k'=1}^{d_2}  \big[ \ind(Z_1=k') - p_Z(k') \big] p_Y(k) \big[p_{YZ}(k,k') - p_Y(k)p_Z(k') \big] \\[.5em]
:= ~ &  2 (I) - 2 (II) - 2(III) \quad \text{(say).}
\end{align*}
Then by recalling the definition of $\psi_{1}'(P)$ in (\ref{Eq: definition of psi prime functions}) and based on the elementary inequality $(x_1 + x_2 + x_3)^2 \leq 3 x_1^2 + 3 x_2^2 + 3x_3^2$, we have 
\begin{align*}
\psi_{1}'(P) \leq 12 \mE_P[(I)^2] + 12 \mE_P[(II)^2] + 12 \mE_P[(III)^2].
\end{align*}
For convenience, we write $\Delta_{k,k'} := p_{YZ}(k,k') - p_Y(k) p_Z(k')$. Focusing on the first expectation in the above upper bound, the basic inequality $(x+y)^2 \leq x^2/ + y^2/2$ gives
\begin{align*}
\mE_P[(I)^2] ~ \leq ~& \frac{1}{2} \mE_P \bigg[ \bigg\{ \sum_{k=1}^{d_1} \sum_{k'=1}^{d_2} \ind(Y_1=k) \ind(Z_1=k')\Delta_{k,k'}  \bigg\}^2 \bigg] + \frac{1}{2}  \bigg\{ \sum_{k=1}^{d_1} \sum_{k'=1}^{d_2} p_{YZ}(k,k') \Delta_{k,k'}  \bigg\}^2 \\[.5em]
\overset{(i)}{\leq} ~ & \frac{1}{2} \sum_{k=1}^{d_1} \sum_{k'=1}^{d_2} p_{YZ}(k,k') \Delta_{k,k'}^2 + \frac{1}{2} \sum_{k=1}^{d_1} \sum_{k'=1}^{d_2} p_{YZ}^2(k,k') \sum_{k=1}^{d_1} \sum_{k'=1}^{d_2}\Delta_{k,k'}^2 \\[.5em]
\overset{(ii)}{\leq} ~ & \frac{1}{2} \sqrt{\sum_{k=1}^{d_1} \sum_{k'=1}^{d_2} p_{YZ}^2(k,k')}\sqrt{\sum_{k=1}^{d_1} \sum_{k'=1}^{d_2} \Delta_{k,k'}^4}  + \frac{1}{2} \sum_{k=1}^{d_1} \sum_{k'=1}^{d_2} p_{YZ}^2(k,k') \sum_{k=1}^{d_1} \sum_{k'=1}^{d_2}\Delta_{k,k'}^2 \\[.5em]
\overset{(iii)}{\leq} ~ & \sqrt{b_{(2)}} \|p_{YZ} - p_Yp_Z\|_2^2,
\end{align*}
where $(i)$ and $(ii)$ use Cauchy-Schwarz inequality and the monotonicity of $\ell_p$ norm (specifically, $\ell_4 \leq \ell_2$). $(iii)$ follows by the definition of $b_{(2)}$ in (\ref{Eq: definition of b_2}) and the fact that $\|p_{YZ}\|_2^2 \leq \|p_{YZ}\|_2$. Turning to the second term $(II)$, one may see that 
\begin{align*}
\mE_P[(II)^2] ~ \leq ~ & \frac{1}{2} \mE_P \bigg[ \bigg\{ \sum_{k=1}^{d_1} \sum_{k'=1}^{d_2} \ind(Y_1=k) p_Z(k')\Delta_{k,k'}  \bigg\}^2 \bigg] + \frac{1}{2}  \bigg\{ \sum_{k=1}^{d_1} \sum_{k'=1}^{d_2} p_{Y}(k)p_{Z}(k') \Delta_{k,k'}  \bigg\}^2 \\[.5em]
= ~ & \frac{1}{2}(II)_{a} + \frac{1}{2}(II)_{b} \quad \text{(say).}
\end{align*}
Using the fact that $\ind(Y_1 = k_1) \ind(Y_1 = k_2) =\ind(Y_1 = k_1)\ind(k_1 = k_2)$, we may upper bound $(II)_{a}$ by
\begin{align*}
\mE_P[(II)_{a}] ~=~ & \sum_{k=1}^{d_1} p_Y(k) \Bigg[ \sum_{k'=1}^{d_2} p_Z(k') \Delta_{k,k'} \Bigg]^2 \\[.5em]
\overset{(i)}{\leq} ~ & \sqrt{\sum_{k=1}^{d_1} p_Y^2(k)} \sqrt{\sum_{k=1}^{d_1} \Bigg( \sum_{k'=1}^{d_2} p_Z(k') \Delta_{k,k'} \Bigg)^4}  \\[.5em]
\overset{(ii)}{\leq} ~ & \sqrt{\sum_{k=1}^{d_1} p_Y^2(k)} \sqrt{\sum_{k=1}^{d_1} \Bigg( \sum_{k'=1}^{d_2} p_Z^2(k')   \sum_{k{''}=1}^{d_2}\Delta_{k,k{''}}^2 \Bigg)^2} \\[.5em]
\overset{(iii)}{\leq} ~ & \sqrt{\sum_{k=1}^{d_1} p_Y^2(k)\sum_{k'=1}^{d_2} p_Z^2(k')} \sqrt{\sum_{k=1}^{d_1}\Bigg( \sum_{k{'}=1}^{d_2}\Delta_{k,k{'}}^2 \Bigg)^2} \\[.5em]
\overset{(iv)}{\leq} ~ & \sqrt{b_{(2)}} \| p_{YZ} - p_Yp_Z\|_2^2, 
\end{align*}
where both $(i)$ and $(ii)$ use Cauchy-Schwarz inequality, $(iii)$ uses $\|p_Z\|_2^2 \leq \|p_Z\|_2$ and $(iii)$ follows by the monotonicity of $\ell_p$ norm (specifically, $\ell_2 \leq \ell_1$) and the definition of $b_{(2)}$ in (\ref{Eq: definition of b_2}). The second term $(II)_{b}$ is bounded similarly by Cauchy-Schwarz inequality and $\|p_Y\|_2^2 \leq \|p_Y\|_2$ and $\|p_Z\|_2^2 \leq \|p_Z\|_2$. In particular,
\begin{align*}
\mE_P[(II)_b] ~ \leq ~ \sum_{k=1}^{d_1} \sum_{k'=1}^{d_2} p_Y^2(k) p_Z^2(k') \|p_{YZ} - p_Yp_Z\|_2^2 ~\leq~ \sqrt{b_{(2)}} \| p_{YZ} - p_Yp_Z\|_2^2.
\end{align*}
By symmetric, $\mE_P[(III)^2]$ is also upper bounded by $\sqrt{b_{(2)}} \| p_{YZ} - p_Yp_Z\|_2^2$. Hence, putting things together, we have
$\psi_{1}'(P) \leq C_1 \sqrt{b_{(2)}}\| p_{YZ} - p_Yp_Z\|_2^2$.

Next we show that the second inequality of (\ref{Eq: Aim of the proof of multinomial independence testing}), which is $\psi_{(2)}'(P) \leq C_2 b_{(2)}$, holds. By recalling the definition of $\psi_{(2)}'(P)$ in (\ref{Eq: definition of psi prime functions}) and noting that $g_Y^2(Y_1,Y_2) = g_Y(Y_1,Y_2)$ and $g_Z^2(Z_1,Z_2) = g_Z(Z_1,Z_2)$, we shall see that
\begin{align*}
 \mE_P[g_Y(Y_1,Y_2)g_Z(Z_1,Z_2)] &~=~  \sum_{k=1}^{d_1} \sum_{k'=1}^{d_2} p_{YZ}^2(k,k') \leq b_{(2)} , \\[.5em]
 \mE_P[g_Y(Y_1,Y_2)g_Z(Z_1,Z_3)] &~=~  \sum_{k=1}^{d_1} \sum_{k'=1}^{d_2} p_{YZ}(k,k') p_Y(k) p_Z(k') \\[.5em]
 & ~\leq~  \frac{1}{2} \sum_{k=1}^{d_1} \sum_{k'=1}^{d_2} p_{YZ}^2(k,k') + \frac{1}{2}\sum_{k=1}^{d_1} \sum_{k'=1}^{d_2}  p_Y^2(k) p_Z^2(k') \leq b_{(2)}, \\[.5em]
 \mE_P[g_Y(Y_1,Y_2)g_Z(Z_3,Z_4)] &~=~ \sum_{k=1}^{d_1} \sum_{k'=1}^{d_2}  p_Y^2(k) p_Z^2(k')  \leq b_{(2)}.
\end{align*}
Hence both conditions in (\ref{Eq: Aim of the proof of multinomial independence testing}) are satisfied under the assumption in Proposition~\ref{Proposition: Multinomial Independence Testing}. This concludes Proposition~\ref{Proposition: Multinomial Independence Testing}.

\section{Proof of Proposition~\ref{Proposition: Minimum Separation for Multinomial Independence Testing}} \label{Section: Proposition: Minimum Separation for Multinomial Independence Testing}
As in the proof of Proposition~\ref{Proposition: Minimum Separation for Two-Sample Multinomial Testing}, we properly construct a mixture distribution $Q$ and a null distribution $P_0$ and apply Lemma~\ref{Lemma: Ingster's method} to prove the result. To start we consider $P_0$ to be the product of the uniform discrete distributions given by
\begin{align*}
P_0(k_1,k_2) := p_Y(k_1)p_Z(k_2)= \frac{1}{d_1d_2} \quad \text{for all $k_1=1,\ldots,d_1$ and $k_2 = 1,\ldots,d_2$.}
\end{align*}
Let $\widetilde{\zeta} = \{\widetilde{\zeta}_1,\ldots,\widetilde{\zeta}_{d_1}\}$ and $\widetilde{\xi} = \{\widetilde{\xi}_1,\ldots,\widetilde{\xi}_{d_2}\}$ be dependent Rademacher random variables uniformly distributed over $\mathcal{M}_{d_1}$ and $\mathcal{M}_{d_2}$, respectively, where $\mathcal{M}_{d_1}$ and $\mathcal{M}_{d_2}$ are hypercubes defined in (\ref{Eq: definition of hypercube}). Assume that $\widetilde{\zeta}$ and $\widetilde{\xi}$ are independent. Let us denote the cardinality of $\mathcal{M}_{d_1}$ and $\mathcal{M}_{d_2}$ by $N_1$ and $N_2$, respectively. Given $\widetilde{\zeta} \in \mathcal{M}_{d_1}$ and $\widetilde{\xi} \in \mathcal{M}_{d_2}$, we define a distribution $p_{\widetilde{\zeta}, \widetilde{\xi}}$ such that
\begin{align*}
p_{\widetilde{\zeta}, \widetilde{\xi}}(k_1,k_2) := \frac{1}{d_1d_2} + \delta  \sum_{i_1=1}^{d_1} \sum_{i_2=1}^{d_2} \widetilde{\zeta}_{i_1} \widetilde{\xi}_{i_2}  \ind(k_1=i_1) \ind(k_2=i_2),
\end{align*}
where $\delta \leq 1/(d_1d_2)$ and thus $\| p_{\widetilde{\zeta}, \widetilde{\xi}} \|_2^2 \leq 2/(d_1d_2)$. Since $\widetilde{\zeta} \in \mathcal{M}_{d_1}$ and $\widetilde{\xi} \in \mathcal{M}_{d_2}$, it is straightforward to check that 
\begin{align*}
& \sum_{k_1=1}^{d_1} p_{\widetilde{\zeta}, \widetilde{\xi}}(k_1,k_2) = \frac{1}{d_2} +  \delta \Bigg\{\sum_{i_1=1}^{d_1} \widetilde{\zeta}_{i_1} \Bigg\} \times \Bigg\{\sum_{i_2=1}^{d_2}  \widetilde{\xi}_{i_2}  \ind(k_2=i_2) \Bigg\} = \frac{1}{d_2}, \\[.5em]
& \sum_{k_2=1}^{d_2} p_{\widetilde{\zeta}, \widetilde{\xi}}(k_1,k_2) = \frac{1}{d_1} +  \delta \Bigg\{ \sum_{i_1=1}^{d_1}  \widetilde{\zeta}_{i_1}  \ind(k_1=i_1) \Bigg\}  \times \Bigg\{ \sum_{i_2=1}^{d_2} \widetilde{\xi}_{i_2} \Bigg\}  = \frac{1}{d_1} \quad \text{and} \\[.5em]
& \sum_{k_1=1}^{d_1}\sum_{k_2=1}^{d_2}  p_{\widetilde{\zeta}, \widetilde{\xi}}(k_1,k_2) = 1.
\end{align*}
Therefore $p_{\widetilde{\zeta}, \widetilde{\xi}}$ is a joint discrete distribution whose marginals are equivalent to those of the product distribution. Let us denote such distributions by $p_{\widetilde{\zeta}(1), \widetilde{\xi}(1)}, \ldots, p_{\widetilde{\zeta}(N_1), \widetilde{\xi}(N_2)}$. We then consider the uniform mixture $Q$ given by
\begin{align*}
Q:= \frac{1}{N_1N_2} \sum_{i=1}^{N_1} \sum_{j=1}^{N_2} p_{\widetilde{\zeta}(i), \widetilde{\xi}(j)}.
\end{align*}
Note that both $\{\widetilde{\zeta}_1,\ldots,\widetilde{\zeta}_{d_1}\}$ and $\{\widetilde{\xi}_1,\ldots,\widetilde{\xi}_{d_2}\}$ are negatively associated and these two sets are mutually independent by construction. Hence, following Proposition 7 of \cite{dubhashi1998balls}, the pooled random variables $\{\widetilde{\zeta}_1,\ldots,\widetilde{\zeta}_{d_1}, \widetilde{\xi}_1,\ldots,\widetilde{\xi}_{d_2}\}$ are also negatively associated. Having this observation at hand, the remaining steps are exactly the same as those in the proof of Proposition~\ref{Proposition: Minimum Separation for Two-Sample Multinomial Testing}. This together with Proposition~\ref{Proposition: Multinomial Independence Testing} completes the proof of Proposition~\ref{Proposition: Minimum Separation for Multinomial Independence Testing}.

\section{Proof of Proposition~\ref{Proposition: Independence Testing for Holder Densities}} \label{Section: Proposition: Independence Testing for Holder Densities}
The proof of Proposition~\ref{Proposition: Independence Testing for Holder Densities} is based on Proposition~\ref{Proposition: Multinomial Independence Testing} and similar to that of Proposition~\ref{Proposition: Two-Sample Testing for Holder Densities}. By recalling the notation from Appendix~\ref{Section: Proposition: Two-Sample Testing for Holder Densities} and $\kappa_{(2)} = \floor{n^{2/(4s+d_1+d_2)}}$, we define $\boldsymbol{H}_{\boldsymbol{\ell}_Y}:= \boldsymbol{[}(\boldsymbol{\ell}_Y-1)/\kappa_{(2)}, \boldsymbol{\ell}_Y /\kappa_{(2)}\boldsymbol{]}$ and $\boldsymbol{H}_{\boldsymbol{\ell}_Z}:= \boldsymbol{[}(\boldsymbol{\ell}_Z-1)/\kappa_{(2)}, \boldsymbol{\ell}_Z /\kappa_{(2)}\boldsymbol{]}$ where $\boldsymbol{\ell}_Y \in \{1,2,\ldots, \kappa_{(2)}\}^{d_1}$ and $\boldsymbol{\ell}_Z \in \{1,2,\ldots, \kappa_{(2)}\}^{d_2}$. Then we denote the joint and product discretized distributions by 
\begin{align*}
& p_{YZ}(\boldsymbol{\ell_Y,\ell_Z}) := \int_{\boldsymbol{H}_{\boldsymbol{\ell}_Y} \times \boldsymbol{H}_{\boldsymbol{\ell}_Z}} f_{YZ}(t_Y,t_Z)dt_Ydt_Z \quad \text{and} \\[.5em]
& p_{Y}p_{Z}(\boldsymbol{\ell_Y,\ell_Z}) := \int_{\boldsymbol{H}_{\boldsymbol{\ell}_Y} \times \boldsymbol{H}_{\boldsymbol{\ell}_Z}} f_{Y}(t_Y) f_{Z}(t_Z) dt_Ydt_Z. 
\end{align*} 
Since both $f_{YZ}$ and $f_Yf_Z$ are in H\"{o}lder's density class $\mathcal{P}_{\text{H\"{o}lder}}^{(d_1+d_2,s)}$ where $|\!|\!|f_Yf_Z|\!|\!|_{\infty} \leq L$ and $|\!|\!|f_{YZ}|\!|\!|_{\infty} \leq L$, it is clear to see that 
\begin{align*}
& p_{YZ}(\boldsymbol{\ell_Y,\ell_Z}) \leq |\!|\!|f_{YZ}|\!|\!|_{\infty} \kappa_{(2)}^{-(d_1+d_2)} \leq L \kappa_{(2)}^{-(d_1+d_2)} \quad \text{and} \\[.5em]  
& p_{Y}p_{Z}(\boldsymbol{\ell_Y},\boldsymbol{\ell_Z})  \leq |\!|\!|f_Yf_Z|\!|\!|_{\infty} \kappa_{(2)}^{-(d_1+d_2)} \leq L \kappa_{(2)}^{-(d_1+d_2)} \quad \text{for all $\boldsymbol{\ell}_Y,\boldsymbol{\ell}_Z$.} 
\end{align*}
This leads to
\begin{align} \label{Eq: b_2 bound}
b_{(2)} = \max\{ \|p_{YZ} \|_2^2, \|p_Yp_Z\|_2^2 \} \leq L \kappa_{(2)}^{-(d_1+d_2)}.
\end{align}
Furthermore, based on Lemma 3 of \cite{arias2018remember}, one can find a constant $C_1 > 0$ such that 
\begin{align} \label{Eq: L2 norm lower bound II}
\|p_{YZ} - p_Yp_Z\|_2^2 \geq C_1 \kappa_{(2)}^{-(d_1+d_2)} \epsilon_{n}^2, 
\end{align}
where $\epsilon_{n}$ is the lower bound for $\|f_{YZ} - f_Yf_Z\|_{L_2}$. By combining (\ref{Eq: b_2 bound}) and (\ref{Eq: L2 norm lower bound II}), the condition of Proposition~\ref{Proposition: Multinomial Independence Testing} is satisfied when 
\begin{align*}
\kappa_{(2)}^{-(d_1+d_2)} \epsilon_{n}^2 ~\geq~ C_2 \frac{L^{1/2} \kappa_{(2)}^{-(d_1+d_2)/2}}{\alpha^{1/2} \beta n}.
\end{align*}
By putting $\kappa_{(2)} = \floor{n^{2/(4s+d_1+d_2)}}$ and rearranging the terms, the above inequality is equivalent to 
\begin{align*}
\epsilon_n ~ \geq ~ \frac{C_3}{\alpha^{1/4} \beta^{1/2}} \left( \frac{1}{n} \right)^{\frac{2s}{4s+d_1+d_2}}, 
\end{align*}
where $C_3$ is a constant that may depend on $s,d_1,d_2,L$. This completes the proof of Proposition~\ref{Proposition: Independence Testing for Holder Densities}.

\section{Proof of Proposition~\ref{Proposition: Minimum Separation for Independence Testing for Holder Densities}} \label{Section: Proposition: Minimum Separation for Independence Testing for Holder Densities}
The proof of Proposition~\ref{Proposition: Minimum Separation for Independence Testing for Holder Densities} is standard based on Ingster's method in Lemma~\ref{Lemma: Ingster's method}. In particular we closely follow the proof of Theorem 1 in \cite{arias2018remember} which builds on \cite{ingster1987minimax}. Let us start with the construction of a mixture distribution $Q$ and a null distribution $P_0$.

\paragraph{$\bullet$ Construction of $Q$ and $P_0$.}

Let $f_Y$ and $f_Z$ be the uniform density functions on $[0,1]^{d_1}$ and $[0,1]^{d_2}$, respectively. Then the density function of the baseline product distribution $P_0$ is defined by
\begin{align*}
f_{0}(y,z) := f_Y(y) f_Z(z) = 1 \quad \text{for all $(y,z) \in [0,1]^{d_1+d_2}$.}
\end{align*}
We let 
$\varphi_Y: \mathbb{R}^{d_1} \mapsto \mathbb{R}$ and $\varphi_Z: \mathbb{R}^{d_2} \mapsto \mathbb{R}$ be infinitely differentiable functions supported on $[0,1]^{d_1}$ and $[0,1]^{d_2}$ respectively. Furthermore these two functions satisfy 
\begin{align*}
& \int_{[0,1]^{d_1}} \varphi_Y(y) dy = \int_{[0,1]^{d_2}} \varphi_Z(z) dz = 0 \quad \text{and} \\[.5em]
& \int_{[0,1]^{d_1}} \varphi_Y^2(y) dy = \int_{[0,1]^{d_2}} \varphi_Z^2(z) dz = 1.
\end{align*}
For $\boldsymbol{i} \in \mathbb{Z}^{d_1}$, $\boldsymbol{j} \in \mathbb{Z}^{d_2}$ and a positive integer $\kappa$, we write $\varphi_{Y,\boldsymbol{i}}(x) = \kappa^{d_1/2} \varphi_Y(\kappa x - \boldsymbol{i} + 1 )$ and $\varphi_{Z,\boldsymbol{j}}(x) = \kappa^{d_2/2} \varphi_Z(\kappa x - \boldsymbol{j} + 1 )$ where $\varphi_{Y,\boldsymbol{i}}$ and $\varphi_{Z,\boldsymbol{j}}$ are supported on $\boldsymbol{[} (\boldsymbol{i}-1)/\kappa, \boldsymbol{i} / \kappa \boldsymbol{]}$ and $\boldsymbol{[} (\boldsymbol{j}-1)/\kappa, \boldsymbol{j} / \kappa \boldsymbol{]}$. By construction, it can be seen that 
\begin{align*}
 & \int_{[0,1]^{d_1}} \varphi_{Y,\boldsymbol{i}}^2 (y) dy  =  \int_{[0,1]^{d_2}} \varphi_{Z,\boldsymbol{j}}^2 (z) dz = 1, \\[.5em]
 &  \int_{[0,1]^{d_1}} \varphi_{Y,\boldsymbol{i}}(y)dy  = \int_{[0,1]^{d_2}} \varphi_{Z,\boldsymbol{j}}(z)dz= 0 \quad \text{and}  \\[.5em] & \int_{[0,1]^{d_1}} \varphi_{Y,\boldsymbol{i}}(y)\varphi_{Y,\boldsymbol{i}'}(y) dy  = \int_{[0,1]^{d_2}} \varphi_{Z,\boldsymbol{j}}(z)\varphi_{Z,\boldsymbol{j}'}(z) dz =0,
\end{align*}
for $\boldsymbol{i} \neq \boldsymbol{i}'$ and $\boldsymbol{j} \neq \boldsymbol{j}'$. We denote by $\zeta_{\boldsymbol{k}} \in \{0,1\}$ an i.i.d.~sequence of Rademacher variables where $\boldsymbol{k}:=(\boldsymbol{i},\boldsymbol{j}) \in \boldsymbol{[}\kappa \boldsymbol{]}^{d_1+d_2}$. Now for $\rho >0$ specified later, let us define the density function of a mixture distribution $Q$ by
\begin{align*}
f_{\zeta} (y,z) :=  f_{0}(y,z)  + \rho \sum_{\boldsymbol{k} \in \boldsymbol{[}\kappa \boldsymbol{]}^{d_1+d_2}}  \zeta_{\boldsymbol{k}}  \varphi_{Y,\boldsymbol{i}}(y)\varphi_{Z,\boldsymbol{j}}(z).
\end{align*}
By letting $\rho$ such that $\rho \kappa^{(d_1+d_2)/2} |\!|\!|\varphi_{Y,Z}|\!|\!|_{\infty} \leq 1$ where $\varphi_{Y,Z}(y,z) := \varphi_Y(y) \varphi_Z(z)$, it is seen that $f_{\zeta}$ is a proper density function supported on $[0,1]^{d_1+d_2}$ such that 
\begin{align*}
\int_{[0,1]^{d_1}} f_{\zeta} (y,z) dy = \int_{[0,1]^{d_2}} f_{\zeta} (y,z) dz = \int_{[0,1]^{d_1 + d_2}} f_{\zeta} (y,z) dy dz =1.
\end{align*}
Therefore $f_{\zeta}$ has the same marginal distributions as the product distribution $f_0$. Furthermore when $\rho \kappa^{(d_1+d_2)/2+s} M/L \leq 1$ where $M := \max \big\{ 4 |\!|\!|\varphi_{Y,Z}^{(\floor{s})}|\!|\!|_{\infty},2 |\!|\!|\varphi_{Y,Z}^{(\floor{s}+1)} |\!|\!|_{\infty} \big\}$, it directly follows from the proof of Theorem 1 in \cite{arias2018remember} that $f_{\zeta} \in \mathcal{P}_{\text{H\"{o}lder}}^{(d_1+d_2,s)}$. Having these two densities $f_0$ and $f_{\zeta}$ such that
\begin{align*}
|\!|\!|f_{\zeta} - f_0|\!|\!|_{L_2}^2 = \rho^2 \kappa^{d_1 + d_2} = \epsilon_n^2,
\end{align*}
we next compute $\mE_{P_0}[L_n^2]$.

\paragraph{$\bullet$ Calculation of $\mE_{P_0}[L_n^2]$.} By recalling that $f_0(y,z) = 1$ for $(y,z) \in [0,1]^{d_1+d_2}$, let us start by writing $L_n^2$ as
\begin{align*}
L_n^2 ~=~ & \frac{1}{2^{2\kappa^{d_1+d_2}}} \sum_{\zeta,\zeta' \in \{-1,1\}^{\kappa^{d_1+d_2}}} \prod_{i=1}^{n} f_{\zeta}(Y_i,Z_i)f_{\zeta'}(Y_i,Z_i).
\end{align*}
We then use the orthonormal property of $\varphi_{Y,\boldsymbol{i}}$ and $\varphi_{Z,\boldsymbol{j}}$ to see that
\begin{align*}
\mE_{P_0}[L_n^2] ~=~ & \frac{1}{2^{2\kappa^{d_1+d_2}}} \sum_{\zeta,\zeta' \in \{-1,1\}^{\kappa^{d_1+d_2}}} \prod_{i=1}^{n}\mE_0\Bigg[ 1  + \rho^2 \sum_{\boldsymbol{k} \in \boldsymbol{[}\kappa \boldsymbol{]}^{d_1+d_2}}  \zeta_{\boldsymbol{k}} \zeta_{\boldsymbol{k}}'  \varphi_{Y,\boldsymbol{i}}^2(Y_i)\varphi_{Z,\boldsymbol{j}}^2(Z_i) \Bigg] \\[.5em]
=~ & \frac{1}{2^{2\kappa^{d_1+d_2}}} \sum_{\zeta,\zeta' \in \{-1,1\}^{\kappa^{d_1+d_2}}} \Bigg[ 1  + \rho^2 \sum_{\boldsymbol{k} \in \boldsymbol{[}\kappa \boldsymbol{]}^{d_1+d_2}}  \zeta_{\boldsymbol{k}} \zeta_{\boldsymbol{k}}'  \Bigg]^n \\[.5em]
\leq ~ & \mE_{\zeta,\zeta'} \Big[e^{n\rho^2 \langle \zeta, \zeta'  \rangle} \Big],
\end{align*}
where the last inequality uses $(1+x)^n \leq e^{nx}$. Based on the independence among the components of $\zeta$ and $\zeta'$, we further observe that 
\begin{align*}
 \mE_{\zeta,\zeta'} \Big[e^{n\rho^2 \langle \zeta, \zeta'  \rangle} \Big] ~ = ~ \big\{ \text{cosh}(n\rho^2) \big\}^{\kappa^{d_1+d_2}} ~ \leq ~ \exp \left(\kappa^{d_1+d_2} n^2\rho^4/2\right)
\end{align*}
where the last inequality follows by $\text{cosh}(x) \leq e^{x^2/2}$ for all $x \in \mathbb{R}$.

\paragraph{$\bullet$ Completion of the proof.} 
We invoke Lemma~\ref{Lemma: Ingster's method} to finish the proof. From the previous step, we know that 
\begin{align*}
\mE_{P_0}[L_n^2] \leq \exp \left(\kappa^{d_1+d_2} n^2\rho^4/2\right).
\end{align*}
Therefore the condition in Lemma~\ref{Lemma: Ingster's method} is fulfilled when 
\begin{align*}
\kappa^{d_1+d_2} n^2 \rho^4 \leq 2 \log \{1 + 4 (1 - \alpha - \beta)^2 \}.
\end{align*}
Now by setting $\kappa = \floor{n^{2/(4s+d_1+d_2)}}$ and $\rho = c n^{-(2s+d_1+d_2)/(4s+ d_1+d_2)}$, the above condition is further implied by 
\begin{align*}
c \leq 2 \log \{1 + 4 (1 - \alpha - \beta)^2 \}.
\end{align*}
Previously we also use the assumptions that $\rho \kappa^{(d_1+d_2)/2} |\!|\!|\varphi_{Y,Z}|\!|\!|_{\infty} \leq 1$ and $\rho \kappa^{(d_1+d_2)/2+s} M/L \leq 1$. These are satisfied by taking $c$ sufficiently small. This means that when 
\begin{align*}
\epsilon_n \leq c e^{-2s/(4s+d_1+d_2)},
\end{align*} 
for a small $c>0$, the minimax type II error is less than $\beta$. Therefore, combined with Proposition~\ref{Proposition: Independence Testing for Holder Densities}, we complete the proof of Proposition~\ref{Proposition: Minimum Separation for Independence Testing for Holder Densities}.

\section{Proof of Theorem~\ref{Theorem: Two-Sample Concentration}} \label{Section: Theorem: Two-Sample Concentration}
We continue the proof of Theorem~\ref{Theorem: Two-Sample Concentration} from the last line of (\ref{Eq: Chernoff Bound}). First we view $\widetilde{U}_{n_1,n_2}^{\pi,L,\zeta}$ as a quadratic form of $\zeta$ conditional on $\pi$ and $L$. We then borrow the proof of Hanson--Wright inequality \citep[see e.g.][]{rudelson2013hanson,vershynin2018high} to proceed. To do so, let us denote $a_{k_1,k_2}(\pi,L) =  h_{\text{ts}}(X_{\pi_{k_1}},X_{\pi_{k_2}}; X_{\pi_{n_1+\ell_{k_1}}},X_{\pi_{n_1+\ell_{k_2}}})$ for $1 \leq k_1 \neq k_2 \leq n$ and $a_{k_1,k_2}(\pi,L) = 0$ for $1 \leq k_1 = k_2 \leq n$. Let $\mathbf{A}_{\pi,L}$ be the $n \times n$ matrix whose elements are $a_{k_1,k_2}(\pi,L)$. By following the proof of Theorem 1.1 in \cite{rudelson2013hanson}, we can obtain
\begin{align*}
e^{-\lambda t} \mathbb{E}_{\pi,L,\zeta} \big[ \exp \big( \lambda \widetilde{U}_{n_1,n_2}^{\pi,L,\zeta} \big) |  \mathcal{X}_n \big] \leq  \mE_{\pi,L} \big[ \exp \big( -\lambda t + C\lambda^2 \| \mathbf{A}_{\pi,L} \|_F^2 \big)  \big],
\end{align*}
which holds for $0 \leq \lambda \leq  c / \|\mathbf{A}_{\pi,L}\|_{\text{op}}$. Here, $\| \mathbf{A}_{\pi,L} \|_F$ and $\| \mathbf{A}_{\pi,L} \|_{\text{op}}$ denote the Frobenius norm and the operator norm of $\mathbf{A}_{\pi,L}$, respectively. By optimizing over $0 \leq \lambda  \leq c / \|A_{\pi,L}\|_{\text{op}}$, we have that
\begin{align*}
\mP_{\pi} (U_{n_1,n_2}^\pi \geq t ~| \mathcal{X}_{n}) \leq \mE_{\pi,L} \left[  \exp \bigg\{ - C_1 \min \left( \frac{t^2}{\|\mathbf{A}_{\pi,L}\|_F^2}, \ \frac{t}{\|\mathbf{A}_{\pi,L}\|_{\text{op}}} \right)  \bigg\} \right].
\end{align*} 
The proof of Theorem~\ref{Theorem: Two-Sample Concentration} is completed by noting that $\|\mathbf{A}_{\pi,L}\|_{\text{op}} \leq \|\mathbf{A}_{\pi,L}\|_F \leq  C_2 \Sigma_{n_1,n_2}$.

\section{Proof of Corollary~\ref{Corollary: Dependent Rademacher Chaos}} \label{Section: Corollary: Dependent Rademacher Chaos}
Note that the following equality holds:
\begin{align*}
\sum_{(i,j) \in \mathbf{i}_2^n} \widetilde{\zeta}_{i} \widetilde{\zeta}_{j}\left(a_{i, j}-\overline{a}\right) ~ = ~ & \frac{1}{4(n-1)(n-2)} \sum_{\left(i_{1}, i_{2}, i_{3}, i_{4}\right) \in \mathbf{i}_{4}^{n}} \Big\{ \\
&  ( \widetilde{\zeta}_{i_{1}}\widetilde{\zeta}_{i_{3}} + \widetilde{\zeta}_{i_{2}} \widetilde{\zeta}_{i_{4}} - \widetilde{\zeta}_{i_{1}} \widetilde{\zeta}_{i_{4}} - \widetilde{\zeta}_{i_{2}} \widetilde{\zeta}_{i_{3}} ) \left(a_{i_{1}, i_{3}}+a_{i_{2}, i_{4}}-a_{i_{1}, i_{4}}-a_{i_{2}, i_{3}}\right) \Big\},
\end{align*}
which can be verified by expanding the summation on the right-hand side. We also note that $\{\widetilde{\zeta}_{1},\ldots,\widetilde{\zeta}_{n} \} \overset{d}{=} \{b_{\pi_1}, \ldots, b_{\pi_n} \}$ where $b_i = 1$ for $i=1,\ldots, n/2$ and $b_i = -1$ for $i=n/2+1,\ldots,n$. Therefore, we can apply Theorem~\ref{Theorem: Concentration inequality for independence U-statistic} with the bound of $\Sigma_n^2$ in (\ref{Eq: A bound of Sigma}). To be clear, $a_{i,j}$ does not need to be symmetric in its arguments. Theorem~\ref{Theorem: Concentration inequality for independence U-statistic} still holds as long as $g_Z$ is symmetric ($g_Y$ is not necessarily symmetric), which is the case for this application. Alternatively, one can work with the symmetrized version of $a_{i,j}$, i.e.~$\widetilde{a}_{i,j} := (a_{i,j} + a_{j,i})/2$ by observing that $\overline{a} = \overline{\widetilde{a}} := n_{(2)}^{-1} \sum_{(i_1,i_2) \in \mathbf{i}_2^{n}} \widetilde{a}_{i_1,i_2}$ and 
\begin{align*}
\sum_{(i,j) \in \mathbf{i}_2^n} \widetilde{\zeta}_{i} \widetilde{\zeta}_{j}\left(a_{i, j}-\overline{a}\right) = \sum_{(i,j) \in \mathbf{i}_2^n} \widetilde{\zeta}_{i} \widetilde{\zeta}_{j}\left(\widetilde{a}_{i, j}-\overline{\widetilde{a}}\right).
\end{align*}
This completes the proof of Corollary~\ref{Corollary: Dependent Rademacher Chaos}.

\section{Proof of Theorem~\ref{Theorem: Concentration inequality for independence U-statistic II}} \label{Section: Theorem: Concentration inequality for independence U-statistic II}
Continuing our discussion from the main text, we prove Theorem~\ref{Theorem: Concentration inequality for independence U-statistic II} in two steps. In the first step, we replace two independent permutations $\pi,\pi'$ in $\widetilde{U}_{n}^{\pi,\pi',\zeta}$ with their i.i.d.~counterparts $\widetilde{\pi},\widetilde{\pi}'$. Once this decoupling step is done, the resulting statistic can be viewed as a usual degenerate $U$-statistic of i.i.d.~random variables conditional on $\mathcal{X}_n$. This means that we can apply the concentration inequalities for degenerate $U$-statistics in \cite{de1999decoupling} to finish the proof. This shall be done in the second step. For notational convenience, we write 
\begin{equation}
\begin{aligned} \label{Eq: kernel in short}
&h_{\pi,\pi'}(i_1,i_2,i_1+m,i_2+m) \\[.5em]
:= ~ & h_{\text{in}} \{(Y_{\pi_{i_1}'},Z_{\pi_{i_1}}),(Y_{\pi_{i_2}'},Z_{\pi_{i_2}}),(Y_{\pi_{i_2+m}'},Z_{\pi_{\pi_{i_2+m}}}),(Y_{\pi_{i_1+m}'},Z_{\pi_{i_1+m}})\},
\end{aligned}
\end{equation}
throughout this proof. 

\paragraph{1. Decoupling.}
We start with the decoupling part. Let $\widetilde{U}_{n}^{\widetilde{\pi},\widetilde{\pi}',\zeta}$ be defined similarly as $\widetilde{U}_{n}^{\pi,\pi',\zeta}$ but with decoupled permutations $(\widetilde{\pi}, \widetilde{\pi}')$ instead of the original permutations $(\pi,\pi')$. Our goal here is to bound
\begin{align} \label{Eq: Decopuling}
\mE_{\pi,\pi',\zeta} \big[\Psi\big(\lambda \widetilde{U}_{n}^{\pi,\pi',\zeta}\big) | \mathcal{X}_n \big] \leq \mE_{\widetilde{\pi},\widetilde{\pi}',\zeta} \big[\Psi\big(C_n \lambda \widetilde{U}_{n}^{\widetilde{\pi},\widetilde{\pi}',\zeta}\big) | \mathcal{X}_n\big],
\end{align}
where $c < C_n <C$ is some deterministic sequence depending on $n$ with some positive constants $c,C>0$. The way how we associate the original statistic $\widetilde{U}_{n}^{\pi,\pi',\zeta}$ with the decoupled couterpart $\widetilde{U}_{n}^{\widetilde{\pi},\widetilde{\pi}',\zeta}$ is as follows. First, we construct a random subset $K$ of $\{1,\ldots,n\}$ such that $\{\pi\}_{i \in K}$ and $\{\widetilde{\pi}\}_{i \in K}$ have the same distribution so that two test statistics based on $\{\pi\}_{i \in K}$ and $\{\widetilde{\pi}\}_{i \in K}$, respectively, shall have the same distribution. The remainder of the proof is devoted to replacing the subset of permutations $\{\pi_i\}_{i \in K}$ and $\{\widetilde{\pi}\}_{i \in K}$ with the entire set of permutations $\{\pi_i\}_{i=1}^n$ and $\{\widetilde{\pi}_i\}_{i=1}^n$. As far as we know, this idea was first employed by \cite{duembgen1998symmetrization} to decouple the simple linear permuted statistic.  

Let us make this decoupling idea more precise. To do so, we define $K$ to be a random subset of $\{1,\ldots,n\}$ independent of everything else except $\widetilde{\pi}$. Specifically, we assume that the conditional distribution of $K$ given $\widetilde{\pi}$ has the uniform distribution on the set of all $J \in \{1,\ldots,n\}$ such that 
\begin{align*}
& \{ \widetilde{\pi}_i : 1 \leq i \leq n \} = \{ \widetilde{\pi}_i : i \in J \} \quad \text{and} \quad \# |\{ \widetilde{\pi}_i : 1 \leq i \leq n \}| = \# |J|,
\end{align*}
where $\# |A|$ denotes the cardinality of a set $A$. Then as noted in \cite{duembgen1998symmetrization}, $\{\pi_i\}_{i \in K} \overset{d}{=} \{ \widetilde{\pi}_i \}_{i \in K}$ follows. In the same way, define another random subset $K'$ of $\{1,\ldots,n\}$ only depending on $\widetilde{\pi}'$ such that $\{\pi_i'\}_{i \in K'} \overset{d}{=} \{ \widetilde{\pi}_i' \}_{i \in K'}$; note that, by construction, $K$ and $K'$ are independent. Furthermore, we let $\mathcal{B}_{K,n}(i_1,i_2,i_1+m,i_2+m)$ be the event such that all of $\{i_1,i_2,i_1+m,i_2+m\}$ are in the random subset $K$. Then, as $\{\pi_i\}_{i \in K} \overset{d}{=} \{ \widetilde{\pi}_i \}_{i \in K}$ and $\{\pi_i'\}_{i \in K'} \overset{d}{=} \{ \widetilde{\pi}_i' \}_{i \in K'}$, we may observe that
\begin{align} \label{Eq: Equal in distribution K and K'}
\widetilde{U}_n^{\pi,\pi',\zeta}(K,K') ~\overset{d}{=} ~ \widetilde{U}_n^{\widetilde{\pi},\widetilde{\pi}',\zeta}(K,K'),
\end{align}
where
\begin{align*}
\widetilde{U}_n^{\pi,\pi',\zeta}(K,K') :=~& \frac{1}{m_{(2)}} \sum_{(i_1,i_2) \in \mathbf{i}_2^{m}} \zeta_{i_1} \zeta_{i_2} \zeta_{i_1+m} \zeta_{i_2+m}h_{\pi,\pi'}(i_1,i_2,i_1+m,i_2+m) \times \\
&  \ind\{ \mathcal{B}_{K,n}(i_1,i_2,i_1+m,i_2+m)\} \ind\{  \mathcal{B}_{K',n}(i_1,i_2,i_1+m,i_2+m)\}, \\[.5em]
\widetilde{U}_n^{\widetilde{\pi},\widetilde{\pi}',\zeta}(K,K') := ~ & \frac{1}{m_{(2)}} \sum_{(i_1,i_2) \in \mathbf{i}_2^{m}} \zeta_{i_1} \zeta_{i_2} \zeta_{i_1+m} \zeta_{i_2+m}h_{\widetilde{\pi},\widetilde{\pi}'}(i_1,i_2,i_1+m,i_2+m) \times \\
& \ind\{ \mathcal{B}_{K,n}(i_1,i_2,i_1+m,i_2+m)\} \ind\{  \mathcal{B}_{K',n}(i_1,i_2,i_1+m,i_2+m)\}.
\end{align*}
Next we calculate the probability of $\mathcal{B}_{K,n}(i_1,i_2,i_1+m,i_2+m)$. By symmetry, we may assume that $i_1=1,i_2=2,i_1+m=3,i_2+m=4$. In fact, this probability is the same as the probability that all of the first four urns are not empty when one throws $n$ balls independently into $n$ urns (here, each urn is equally likely to be selected). Based on the inclusion--exclusion formula, this probability can be computed as
\begin{align*}
B_n:= \mP\{\mathcal{B}_{K,n}(1,2,3,4)\} = 1 - 4 \left(1 - \frac{1}{n}\right)^n + 6 \left(1 - \frac{2}{n}\right)^n - 4 \left( 1 - \frac{3}{n} \right)^n + \left( 1 - \frac{4}{n}\right)^n.
\end{align*}
Indeed, $B_n$ is monotone increasing for all $n \geq 4$. Hence we have that $\ell \leq B_n \leq u$ for any $n \geq 4$ where $\ell = 1 - 4(3/4)^4 + 6 \left( 1/2 \right)^4 - 4\left( 1/4\right)^4 =  0.09375$ and $u = 1 - 4e^{-1} + 6 e^{-2} - 4 e^{-3} + e^{-4} \approx 0.1597$.
In the next step, we replace the subset of permutations $\{\pi_i\}_{i \in K}$ with the entire set of permutations $\{\pi_i\}_{i=1}^n$ as follows:
\begin{align*}
\mE_{\pi,\pi',\zeta} \big[\Psi\big(\lambda \widetilde{U}_{n}^{\pi,\pi',\zeta}\big) | \mathcal{X}_n \big] ~ \overset{(i)}{\leq} ~ &  \mE_{\pi,\pi',\zeta,K,K'} \big[\Psi\big\{B_n^{-2} \lambda \widetilde{U}_{n}^{\pi,\pi',\zeta}(K,K')\big\} | \mathcal{X}_n \big] \\[.5em]
\overset{(ii)}{=}~ &  \mE_{\widetilde{\pi},\widetilde{\pi}',\zeta,K,K'} \big[\Psi\big\{B_n^{-2} \lambda \widetilde{U}_{n}^{\widetilde{\pi},\widetilde{\pi}',\zeta}(K,K')\big\} | \mathcal{X}_n \big] \\[.5em]
\overset{(iii)}{\leq} ~ &  \mE_{\widetilde{\pi},\widetilde{\pi}',\zeta} \big[\Psi\big(B_n^{-2} \lambda \widetilde{U}_{n}^{\widetilde{\pi},\widetilde{\pi}',\zeta}\big) | \mathcal{X}_n \big],
\end{align*}
where $(i)$ holds by Jensen's inequality with $\mE_{K,K'}[\widetilde{U}_n^{\pi,\pi',\zeta}(K,K')] = B_n^2\widetilde{U}_n^{\pi,\pi',\zeta}$, $(ii)$ is due to the relationship (\ref{Eq: Equal in distribution K and K'}) and $(iii)$ uses Jensen's inequality again with 
\begin{align*}
\widetilde{U}_{n}^{\widetilde{\pi},\widetilde{\pi}',\zeta}(K,K') = \mE_\zeta \big[\widetilde{U}_{n}^{\widetilde{\pi},\widetilde{\pi}',\zeta} ~\big|~ \big\{\zeta_i \big\}_{i \in K}, \big\{\zeta_i \big\}_{i \in K'}, K, K', \mathcal{X}_n, \widetilde{\pi},\widetilde{\pi}' \big].
\end{align*}
This proves the decoupling inequality in (\ref{Eq: Decopuling}). 

\paragraph{2. Concentration.} Having established the decoupled bound in (\ref{Eq: Decopuling}), we are now ready to obtain the main result of Theorem~\ref{Theorem: Concentration inequality for independence U-statistic II}. This part of the proof is largely based on Chapter 4.1.3 of \cite{de1999decoupling}. Recall that 
\begin{align*}
\widetilde{U}_{n}^{\widetilde{\pi},\widetilde{\pi}',\zeta} ~\overset{d}{=}~  \frac{1}{m_{(2)}} \sum_{(i_1,i_2) \in \mathbf{i}_2^{m}} \zeta_{i_1} \zeta_{i_2} h_{\widetilde{\pi},\widetilde{\pi}'}(i_1,i_2,i_1+m,i_2+m)
\end{align*}
and $h_{\widetilde{\pi},\widetilde{\pi}'}(i_1,i_2,i_1+m,i_2+m)$ is given in (\ref{Eq: kernel in short}). Let us write $\mathbf{Q}_{i_1} = ((Y_{\widetilde{\pi}_{i_1}'},Z_{\widetilde{\pi}_{i_1}}), (Y_{\widetilde{\pi}_{i_1+m}'},Z_{\widetilde{\pi}_{i_1+m}}))$ and $\mathbf{Q}_{i_2} = ((Y_{\widetilde{\pi}_{i_2}'},Z_{\widetilde{\pi}_{i_2}}), (Y_{\widetilde{\pi}_{i_2+m}'},Z_{\widetilde{\pi}_{i_2+m}}))$, which are random vectors with four main components. Note that $\mathbf{Q}_1,\ldots,\mathbf{Q}_m$ are independent and identically distributed conditional on $\mathcal{X}_n$. Define
\begin{align*}
h(\mathbf{Q}_{i_1},\mathbf{Q}_{i_2}) := h_{\widetilde{\pi},\widetilde{\pi}'}(i_1,i_2,i_1+m,i_2+m).
\end{align*}
Then $\widetilde{U}_{n}^{\widetilde{\pi},\widetilde{\pi}',\zeta}$ can be viewed as a randomized $U$-statistic with the bivariate kernel $h(\mathbf{Q}_{i_1},\mathbf{Q}_{i_2})$. 
To summarize, we have established that 
\begin{align*}
\mE_{\pi} [ \Psi(\lambda U_n^{\pi}) | \mathcal{X}_n] ~\leq~ \mE_{\zeta,Q} \Bigg[ \Psi \Bigg(B_n^{-2} \lambda \frac{1}{m_{(2)}} \sum_{(i_1,i_2) \in \mathbf{i}_2^{m}} \zeta_{i_1} \zeta_{i_2} h(\mathbf{Q}_{i_1},\mathbf{Q}_{i_2}) \Bigg) ~\Bigg|~ \mathcal{X}_n \Bigg].
\end{align*}
Here, by letting $h^\ast(\mathbf{Q}_{i_1},\mathbf{Q}_{i_2}) = h(\mathbf{Q}_{i_1},\mathbf{Q}_{i_2})/2 + h(\mathbf{Q}_{i_2},\mathbf{Q}_{i_1})/2$, we may express the right-hand side of the above inequality with the symmetrized kernel as 
\begin{align*}
\mE_{\zeta,Q} \Bigg[ \Psi \Bigg(B_n^{-2} \lambda \frac{2}{m_{(2)}} \sum_{1 \leq i_1 < i_2 \leq m} \zeta_{i_{1}} \zeta_{i_{2}} h^\ast(\mathbf{Q}_{i_1},\mathbf{Q}_{i_2})\Bigg) ~\Bigg|~ \mathcal{X}_n \Bigg].
\end{align*}
The rest of the proof follows exactly the same line of that of Theorem 4.1.12 in \cite{de1999decoupling} based on $(i)$ Chernoff bound, $(ii)$ convex modification, $(iii)$ Bernstein's inequality, $(iv)$ hypercontractivity of Rademacher chaos variables and $(v)$ Hoeffding's average \citep{hoeffding1963probability}. In the end, we obtain
\begin{align*}
\mP_{\pi}(nU_n^{\pi} \geq t ~ |\mathcal{X}_n) ~\leq~ C_1 \exp \left( - \lambda t^{2/3} + C_2 \lambda^3 \Lambda_n^2 + \frac{16C_2^2\Lambda_n^2 M_n^2 \lambda^6}{n - (16/3)C_2M_n^2 \lambda^3} \right),
\end{align*}
for $n > (4/3) C_2M_n^2 \lambda^3$, which corresponds to Equation~(4.1.27) of \cite{de1999decoupling}. We complete the proof of Theorem~\ref{Theorem: Concentration inequality for independence U-statistic II} by optimizing the right-hand side over $\lambda$ as detailed in \cite{de1999decoupling}.

\section{Proof of Proposition~\ref{Proposition: Adaptive two-sample test}} \label{Section: Proposition: Adaptive two-sample test}
The proof of this result is motivated by \cite{ingster2000adaptive,arias2018remember} and follows similarly as theirs. First note that type I error control of the adaptive test is trivial by the union bound. Hence we focus on the type II error control. Note that by construction 
\begin{align*}
\left( \frac{n_1}{\log \log n_1} \right)^{\frac{2}{4s + d}} \leq 2^{\gamma_{\text{max}}}.
\end{align*}
Therefore there exists an integer $j \in \{ 1,\ldots, \gamma_{\text{max}} \}$ such that
\begin{align*}
2^{j -1} <   \left( \frac{n_1}{\log \log n_1} \right)^{\frac{2}{4s + d}} \leq 2^{j}.
\end{align*}
We take such $j$ and define $\kappa^\ast := 2^{j} \in \mathsf{K}$. In the rest of the proof, we show that under the given condition, $\phi_{\kappa^\ast,\alpha/\gamma_{\text{max}}}$ has the type II error at most $\beta$. If this is the case, then the proof is completed since $\mP_P(\phi_{\text{adapt}} = 0) \leq \mP_P(\phi_{\kappa^\ast,\alpha/\gamma_{\text{max}}} = 0) \leq \beta$. To this end, let us start by improving Proposition~\ref{Proposition: Multinomial Two-Sample Testing} based on Lemma~\ref{Lemma: Two-Sample U-statistic Improved Version}. Using (\ref{Eq: Three conditions for two-sample multinomials}) and Lemma~\ref{Lemma: Two-Sample U-statistic Improved Version}, one can verify that Proposition~\ref{Proposition: Multinomial Two-Sample Testing} holds if 
\begin{align} \label{Eq: sufficient condition for adaptive test}
\|p_Y - p_Z\|_2^2 \geq \frac{C}{\beta} \log \left(\frac{1}{\alpha} \right) \frac{\sqrt{b_{(1)}}}{n_1},
\end{align}
for some large constant $C > 0$ and $n_1 \asymp n_2$. Hence the multinomial test $\phi_{\kappa^\ast,\alpha/\gamma_{\text{max}}}$ has the type II error at most $\beta$ if condition (\ref{Eq: sufficient condition for adaptive test}) is fulfilled by replacing $\alpha$ with $\alpha / \gamma_{\text{max}}$. Following the proof of Proposition~\ref{Proposition: Two-Sample Testing for Holder Densities} but with $\kappa^\ast$ instead of $\kappa_{(1)}$, we can see that 
\begin{align*}
& b_{(1)} = \max\{\|p_Y\|_2^2, \|p_Z\|_2^2\} \leq L (\kappa^\ast)^{-d} \quad \text{and} \\[.5em]
& \|p_Y - p_Z\|_2^2 \geq C_1(s,d,L) (\kappa^\ast)^{-d} \epsilon_{n_1,n_2}^2.
\end{align*}
Therefore condition (\ref{Eq: sufficient condition for adaptive test}) with $\alpha/ \gamma_{\text{max}}$ is satisfied when 
\begin{align*}
\epsilon_{n_1,n_2}^2 \geq \frac{C_2(s,d,L)}{\beta} \log \left( \frac{\gamma_{\text{max}}}{\alpha} \right) \frac{L^{1/2} \kappa^{\ast d/2}}{n_1}.
\end{align*}
Based on the definition of $\gamma_{\text{max}}$ and $\kappa^\ast$, the above inequality is further implied by
\begin{align*}
\epsilon_{n_1,n_2}^2 \geq C(s,d,L,\alpha,\beta) \left( \frac{\log\log n_1}{n_1} \right)^{\frac{4s}{4s+d}}.
\end{align*}
This completes the proof of Proposition~\ref{Proposition: Adaptive two-sample test}.

\section{Proof of Proposition~\ref{Proposition: Adaptive independence test}}  \label{Section: Proposition: Adaptive independence test}
The proof is almost identical to that of Proposition~\ref{Proposition: Adaptive two-sample test} once we establish the following lemma which is an improvement of Proposition~\ref{Proposition: Multinomial Independence Testing} in size $\alpha$ under additional conditions. 
\begin{lemma}[Multinomial independence testing] \label{Lemma: Multinomial Independence Testing}
Let $Y$ and $Z$ be multinomial random vectors in $\mathbb{S}_{d_1'}$ and $\mathbb{S}_{d_2'}$, respectively. Consider the multinomial problem setting in Proposition~\ref{Proposition: Multinomial Independence Testing} with additional assumptions that $\frac{1}{n} \log (1/\alpha) \leq C_1 /(d_1'd_2')$ for some constant $C_1 > 0$ and $\alpha \leq e^{-1}$, i.e.
\begin{align*}
	e^{-C_1 \frac{n}{d_1'd_2'}} \leq \alpha \leq e^{-1}.
\end{align*}
Suppose that under the alternative hypothesis,
\begin{align*}
	\|p_{YZ} - p_Yp_Z\|_2 \geq \frac{C_2}{\beta^{1/2}} \sqrt{\log\left( \frac{1}{\alpha} \right)} \frac{b_{(2)}^{1/4}}{n^{1/2}},
\end{align*}
for a sufficiently large $C_2>0$ (depending on $C_1$). Then the permutation test in Proposition~\ref{Proposition: Multinomial Independence Testing} has the type II error at most $\beta$. 
\begin{proof}
	Following the proofs of Lemma~\ref{Lemma: Two-Sample U-statistic Improved Version} and Proposition~\ref{Proposition: Multinomial Independence Testing}, we only need to show that the $1-\beta/2$ quantile of the permutation critical value $c_{1-\alpha,n}$ of $U_n$, denoted by $q_{1-\beta/2,n}$, is bounded as
	\begin{align} \label{Eq: quantile bound for the independence testing}
		q_{1-\beta/2,n} \leq \frac{C_3}{\beta} \log\left( \frac{1}{\alpha} \right) \frac{b_{(2)}^{1/2}}{n}.
	\end{align} 
	To establish this result, we first use the concentration bound in Theorem~\ref{Theorem: Concentration inequality for independence U-statistic II} to have
	\begin{align*}
		c_{1-\alpha, n} \leq C_4 \max \Bigg\{ \frac{\Lambda_n}{n} \log \left( \frac{1}{\alpha}\right), \ \frac{1}{n^{3/2}} \log^{3/2} \left( \frac{1}{\alpha} \right)  \Bigg\},
	\end{align*}
	where we use the fact that $M_n \leq 1$. Hence, by Markov's inequality as in Lemma~\ref{Lemma: Two-Sample U-statistic Improved Version}, it can be seen that the quantile $q_{1-\beta/2,n}$ is bounded by 
	\begin{align*}
		q_{1-\beta/2,n} \leq C_4 \max \Bigg\{ \frac{\sqrt{2\mE[\Lambda_n^2]}}{\beta^{1/2} n} \log \left( \frac{1}{\alpha}\right), \ \frac{1}{n^{3/2}} \log^{3/2} \left( \frac{1}{\alpha} \right)  \Bigg\}.
	\end{align*}
	On the other hand, one can easily verify that
	\begin{align*}
		\mE_P[\Lambda_n^2] = \frac{1}{n^4} \sum_{1 \leq i_1,i_2 \leq n} \sum_{1 \leq j_1,j_2 \leq n} \mE \big[ g_Y^2(Y_{i_1},Y_{i_2})g_Z^2(Z_{j_1},Z_{j_2}) \big] \leq b_{(2)} + \frac{C_5}{n}. 
	\end{align*}
	Furthermore, Cauchy--Schwarz inequality shows that
	\begin{align} \label{Eq: CS inequality}
		b_{(2)} = \max\{\|p_{YZ}\|_2^2, \|p_Yp_Z\|_2^2 \} \geq \frac{1}{d_1'd_2'} \geq \frac{1}{C_1 n},
	\end{align}
	where the last inequality uses the assumption $C_1/(d_1'd_2') \geq n^{-1} \log(1/\alpha) \geq n^{-1}$. Therefore we have $\mE_P[\Lambda_n^2] \leq C_6 b_{(2)}$. This further implies that 
	\begin{align*}
		q_{1-\beta/2,n} \leq ~ & C_7 \max \Bigg\{ \frac{\sqrt{2\mE[\Lambda_n^2]}}{\beta^{1/2} n} \log \left( \frac{1}{\alpha}\right), \ \frac{1}{n^{3/2}} \log^{3/2} \left( \frac{1}{\alpha} \right)  \Bigg\} \\[.5em]
		= ~ &  C_7 \max \Bigg\{ \frac{\sqrt{2\mE[\Lambda_n^2]}}{\beta^{1/2} n} \log \left( \frac{1}{\alpha}\right), \ \frac{1}{n} \log \left( \frac{1}{\alpha} \right) \times \sqrt{\frac{1}{n} \log \left( \frac{1}{\alpha} \right)}  \Bigg\} \\[.5em]
		\leq ~ & \frac{C_8}{\beta} \log \left( \frac{1}{\alpha} \right) \frac{b_{(2)}^{1/2}}{n},
	\end{align*}
	where the last inequality uses $\beta \leq \beta^{1/2}$ and
	\begin{align*}
		\sqrt{\frac{1}{n} \log \left(\frac{1}{\alpha} \right)} \leq \sqrt{\frac{C_1}{d_1'd_2'}} \leq \sqrt{C_1 b_{(2)}^{1/2}}.
	\end{align*}
	Hence the quantile is bounded as (\ref{Eq: quantile bound for the independence testing}). This completes the proof of Lemma~\ref{Lemma: Multinomial Independence Testing}.
\end{proof}
\end{lemma}

Let us come back to the proof of Proposition~\ref{Proposition: Adaptive independence test}. Since type I error control is trivial by the union bound, we only need to show the type II error control of the adaptive test. As in the proof of Proposition~\ref{Proposition: Adaptive two-sample test}, we know that there exists an integer $j \in \{1,\ldots, \gamma_{\text{max}}^\ast\}$ such that 
\begin{align} \label{Eq: sandwich bound}
	2^{j -1} <   \left( \frac{n}{\log \log n} \right)^{\frac{2}{4s + d_1 + d_2}} \leq 2^{j}.
\end{align}
We take such $j$ and define $\kappa^\ast := 2^{j} \in \mathsf{K}^\dagger$. Since $\mP_P(\phi_{\text{adapt}}^\dagger = 0) \leq \mP_P(\phi_{\kappa^\ast, \alpha/\gamma_{\text{max}}^\ast}^\dagger = 0)$, it suffices to show that the resulting multinomial test $\phi_{\kappa^\ast, \alpha/\gamma_{\text{max}}^\ast}^\dagger$ controls the type II error by $\beta$ under the given condition. To this end, we invoke Lemma~\ref{Lemma: Multinomial Independence Testing}. Note that there are $(\kappa^{\ast})^{d_1 + d_2}$ number of bins for $\phi_{\kappa^\ast, \alpha/\gamma_{\text{max}}^\ast}^\dagger$, which is bounded by
\begin{align} \label{Eq: kappa bound}
	\left( \kappa^{\ast}\right)^{d_1 + d_2} \overset{(i)}{\leq} 2^{d_1 + d_2} \left( \frac{n}{\log \log n} \right)^{\frac{2(d_1 + d_2)}{4s + d_1 + d_2}} \overset{(ii)}{\leq} 2^{d_1 + d_2} \left( \frac{n}{\log \log n} \right), 
\end{align} 
where $(i)$ follows by the bound (\ref{Eq: sandwich bound}) and $(ii)$ follows since $4s \geq d_1 +d_2$.

Next, in order to apply Lemma~\ref{Lemma: Multinomial Independence Testing}, we need to check the condition in the lemma, namely $\frac{1}{n} \log (1/\alpha) \leq C_1 /(d_1'd_2')$, is satisfied. Here we work with size $\alpha/\gamma_{\mathrm{max}}^\ast$ and the number of bins $d_1'd_2' = \left( \kappa^{\ast}\right)^{d_1 + d_2}$. Therefore the given condition becomes
\begin{align*}
	\frac{1}{n} \log \left( \frac{\gamma_{\mathrm{max}}^\ast}{\alpha} \right) \leq \frac{C_1}{\left( \kappa^{\ast}\right)^{d_1 + d_2}}.
\end{align*}
Recall that 
\begin{align*}
	\gamma_{\text{max}}^\ast := \bigg\lceil\frac{2}{d_1+ d_2} \log_2 \left( \frac{n}{\log \log n}\right)\bigg\rceil \leq \frac{4}{(d_1+ d_2)\log 2} \log \left( \frac{n}{\log \log n}\right)
\end{align*}
Using this upper bound along with the bound~\eqref{Eq: kappa bound}, it can be seen that the condition is fulfilled if
\begin{align*}
	& \frac{1}{n} \biggl[ \log\biggl( \frac{4}{\alpha (d_1 + d_2) \log 2} \biggr) + \log \log \biggl( \frac{n}{\log \log n} \biggr) \biggr] \leq \frac{C_1}{2^{d_1+d_2}} \frac{\log \log n}{n} \\[.8em]
	\Longleftrightarrow ~~ & 2^{d_1+d_2} \frac{\log\bigl( \frac{4}{\alpha (d_1 + d_2) \log 2} \bigr) + \log \log \bigl( \frac{n}{\log \log n} \bigr) }{\log \log n} \leq C_1,
\end{align*}
which is implied by 
\begin{align} \label{Eq: condition of C1}
	2^{d_1+d_2} \biggl[ \log\biggl( \frac{4}{\alpha (d_1 + d_2) \log 2} \biggr)  + 1  \biggr] \leq C_1.
\end{align}
Since $d_1,d_2,\alpha$ are treated as constants in Proposition~\ref{Proposition: Adaptive independence test}, the condition in Lemma~\ref{Lemma: Multinomial Independence Testing} is fulfilled with any constant $C_1$ satisfying the above inequality~(\ref{Eq: condition of C1}). Therefore we can apply Lemma~\ref{Lemma: Multinomial Independence Testing}. 

From the proof of Proposition~\ref{Proposition: Independence Testing for Holder Densities}, we know that
\begin{align*}
	& b_{(2)} = \max\{ \|p_{YZ} \|_2^2, \|p_Yp_Z\|_2^2 \} \leq L (\kappa^\ast)^{-(d_1+d_2)} \quad \text{and} \\[.5em]
	& \|p_{YZ} - p_Yp_Z\|_2^2 \geq C_1(s,L,d_1,d_2) (\kappa^\ast)^{-(d_1+d_2)} \epsilon_{n}^2,
\end{align*}
where $\epsilon_n$ is the lower bound for $|\!|\!|f_{YZ} - f_Yf_Z|\!|\!|_{L_2}$. Combining this observation with Lemma~\ref{Lemma: Multinomial Independence Testing} shows that $\phi_{\kappa^\ast, \alpha/\gamma_{\text{max}}^\ast}^\dagger$ has non-trivial power when
\begin{align*}
	\epsilon_{n}^2 \geq \frac{C_2(s,L,d_1,d_2)}{\beta} \log\left( \frac{\gamma_{\text{max}}^\ast}{\alpha} \right) \frac{L^{1/2} \cdot (\kappa^\ast)^{(d_1+d_2)/2}}{n}.
\end{align*}
By the definition of $\gamma_{\text{max}}^\ast$ and $\kappa^\ast$, this inequality is further implied by 
\begin{align*}
	\epsilon_n^2 \geq C_3(s,L,d_1,d_2,\alpha,\beta) \left( \frac{\log \log n}{n} \right)^{\frac{4s}{4s + d_1 + d_2}}.
\end{align*}
This completes the proof of Proposition~\ref{Proposition: Adaptive independence test}.

\section{Proof of Theorem~\ref{Theorem: Two-Sample testing under Poisson sampling}} \label{Section: Theorem: Two-Sample testing under Poisson sampling}
We use the quantile approach described in Section~\ref{Section: A general strategy with two moments} to prove the result \citep[see also][]{fromont2013two}. More specifically we let $q_{1-\beta/2,n}$ denote the quantile of the permutation critical value $c_{1-\alpha,n}$ of $T_{\chi^2}$. Then as shown in the proof of Lemma~\ref{Lemma: Two Moments Method}, if 
\begin{align} \label{Eq: Goal of Poisson two-sample testing}
\mE_P[T_{\chi^2}] ~ \geq ~ q_{1-\beta/2,n} + \sqrt{\frac{2\mV_P[T_{\chi^2}]}{\beta}},
\end{align}
then the type II error of the permutation test is controlled as 
\begin{align*}
\sup_{P \in \mathcal{P}_1} \mP_P( T_{\chi^2} \leq c_{1-\alpha,n}) ~\leq~ & \sup_{P \in \mathcal{P}_1} \mP_P( T_{\chi^2} \leq q_{1-\beta/2,n})  + \sup_{P \in \mathcal{P}_1} \mP_P( q_{1-\beta/2,n} < c_{1-\alpha,n})  \\[.5em]
\leq ~ &\beta. 
\end{align*}
Therefore we only need to show that the inequality (\ref{Eq: Goal of Poisson two-sample testing}) holds under the condition given in Theorem~\ref{Theorem: Two-Sample testing under Poisson sampling}. Note that \cite{chan2014optimal} present a lower bound for $\mE_P[T_{\chi^2}]$ as 
\begin{equation}
\begin{aligned} \label{Eq: Bounds for M}
\mE_P [T_{\chi^2}] ~ =~ &   \sum_{k=1}^d \frac{\{p_Y(k) - p_Z(k) \}^2}{p_Y(k) + p_Z(k)} n \Bigg( 1 -  \frac{1 - e^{-n\{p_Y(k) + p_Z(k)\}}}{n \{p_Y(k) + p_Z(k)\}} \Bigg) \\[.5em]
\geq ~ & \frac{n^2}{4d + 2n} \| p_Y - p_Z\|_1^2,
\end{aligned}
\end{equation}
 and an upper bound for $\mV_P[T_{\chi^2}]$ by
\begin{align} \label{Eq: Bounds for V}
\mV_P[T_{\chi^2}] ~  \leq ~ 2 \min\{n,d\} + 5n\sum_{k=1}^d \frac{\{p_Y(k) - p_Z(k) \}^2}{p_Y(k) + p_Z(k)}.
\end{align}
In the rest of the proof, we show that for some constant $C_1>0$,
\begin{align} \label{Eq: Sufficient condition for Poisson testing}
q_{1-\beta/2,n}  \leq \frac{C_1}{\beta} \log\left( \frac{1}{\alpha}\right)  \sqrt{\min\{n,d\} }.
\end{align}
Building on these three observations (\ref{Eq: Bounds for M}), (\ref{Eq: Bounds for V}) and (\ref{Eq: Sufficient condition for Poisson testing}), we can verify that the sufficient condition~(\ref{Eq: Goal of Poisson two-sample testing}) is satisfied under the assumption made in Theorem~\ref{Theorem: Two-Sample testing under Poisson sampling}. Although it can be done by following~\cite{chan2014optimal}, their proof may be too concise for some readers (also there is a typo in their algorithm in Section 2 --- the critical value should be $C \sqrt{\min\{n,d\}}$ instead of $C\sqrt{n}$) and so we decide to give detailed explanations in Appendix~\ref{Section: Details on the verification}. Hence all we need to show is condition~(\ref{Eq: Sufficient condition for Poisson testing}).

\subsection{Verification of condition~(\ref{Eq: Sufficient condition for Poisson testing})}
Recall the permuted chi-square statistic $T_{\chi^2}^{\pi}$ given as
\begin{align*}
T_{\chi^2}^{\pi} = \sum_{k=1}^d \frac{(\sum_{i=1}^n X_{\pi_i,k} - \sum_{i=1}^n X_{\pi_{i+n},k} )^2 - V_k - W_k}{V_k + W_k} \ind(V_k+W_k > 0).
\end{align*}
For simplicity, let us write $\omega_k := V_k + W_k$ for $k=1,\ldots,d$. Note that $\omega_1,\ldots, \omega_d$ are permutation invariant and they should be constant under the permutation law. Having this observation in mind, we split the permuted statistic into two parts:
\begin{align*}
T_{\chi^2}^{\pi} ~=~ & \sum_{(i,j) \in \mathbf{i}_2^n} \sum_{k=1}^d \frac{(X_{\pi_i,k} - X_{\pi_{i+n},k})(X_{\pi_j,k} - X_{\pi_{j+n},k})}{\omega_k} \ind(\omega_k >0) \\[.5em]
+ & \sum_{i=1}^n \sum_{k=1}^d \frac{(X_{\pi_i,k} - X_{\pi_{i+n},k})^2}{\omega_k} \ind(\omega_k >0) - \sum_{k=1}^d \ind(\omega_k > 0) \\[.5em]
=~ & T_{\chi^2,a}^{\pi} + T_{\chi^2,b}^{\pi}  \quad \text{(say).}
\end{align*}
Let us first compute an upper bound for the $1-\alpha$ critical value of $T_{\chi^2}^{\pi}$. To do so, recall that $\xi_1,\ldots,\xi_n$ are i.i.d.~Rademacher random variables. From the same reasoning made in Section~\ref{Section: Degenerate two-sample U-statistics}, one can see that $T_{\xi^2,a}^{\pi}$ have the same distribution as
\begin{align*}
\sum_{(i,j) \in \mathbf{i}_2^n} \xi_i \xi_j \Bigg[ \sum_{k=1}^d \frac{(X_{\pi_i,k} - X_{\pi_{i+n},k})(X_{\pi_j,k} - X_{\pi_{j+n},k})}{\omega_k} \ind(\omega_k >0) \Bigg]. 
\end{align*}
Then following the same line of the proof of Theorem~\ref{Theorem: Two-Sample Concentration} with the trivial bound in (\ref{Eq: trivial bound}), we have that for any $t>0$,
\begin{align} \label{Eq: Poisson Result 1}
\mP_{\pi} \big(T_{\chi^2,a}^{\pi} \geq t  ~| \mathcal{X}_n \big) \leq \exp \bigg\{ - C_2 \min \left( \frac{t^2}{\Sigma_{n,\text{pois}}^2 }, \ \frac{t}{\Sigma_{n,\text{pois}}} \right) \bigg\},
\end{align}
where
\begin{align} \label{Eq: definition of Sigma n}
\Sigma_{n,\text{pois}}^2 ~ :=  \sum_{(i,j) \in \mathbf{i}_2^{2n}} \Bigg\{ \sum_{k=1}^d \frac{X_{i,k}X_{j,k}}{\omega_k} \ind(\omega_k > 0) \Bigg\}^2
\end{align}
and $\{X_{1,k},\ldots,X_{2n,k}\} := \{Y_{1,k},\ldots,Y_{n,k},Z_{1,k},\ldots,Z_{n,k} \}$. Also note that
\begin{equation} 
\begin{aligned}  \label{Eq: Poisson Result 2}
T_{\chi^2,b}^{\pi} ~=~ & \sum_{i=1}^n \sum_{k=1}^d \frac{(X_{\pi_i,k} - X_{\pi_{i+n},k})^2}{\omega_k} \ind(\omega_k >0) - \sum_{k=1}^d \ind(\omega_k > 0) \\[.5em]
\leq ~ & \sum_{i=1}^{2n} \sum_{k=1}^d \frac{X_{i,k}^2}{\omega_k} \ind(\omega_k >0) - \sum_{k=1}^d \ind(\omega_k > 0) \\[.5em]
:= ~ & T_{\chi^2,b, \text{up}}
\end{aligned}
\end{equation}
where $T_{\chi^2,b, \text{up}}$ is independent of $\pi$. Furthermore, since each $X_{i,k}$ can have a nonnegative integer and $\omega_k = \sum_{i=1}^{2n} X_{i,k}$, it is clear that $\sum_{i=1}^{2n} X_{i,k}^2 /\omega_k \geq 1$ whenever $\omega_k > 0$. This means that $T_{\chi^2,b, \text{up}}$ is nonnegative. Combining the results (\ref{Eq: Poisson Result 1}) and (\ref{Eq: Poisson Result 2}), for any $t >0$,
\begin{align*}
\mP_{\pi} \big( T_{\chi^2}^{\pi} \geq t + T_{\chi^2,b,\text{up}} ~|\mathcal{X}_n \big) \leq \mP_{\pi}\big( T_{\chi^2,a}^{\pi}  \geq t   ~| \mathcal{X}_n \big) \leq   \exp \bigg\{ - C_3 \min \left( \frac{t^2}{\Sigma_{n,\text{pois}}^2 }, \ \frac{t}{\Sigma_{n,\text{pois}}} \right) \bigg\}.
\end{align*}
By setting the upper bound to be $\alpha$ and assuming $\alpha < e^{-1}$, it can be seen that 
\begin{align*}
c_{1-\alpha,n} \leq C_4 \Sigma_{n,\text{pois}} \log \left( \frac{1}{\alpha} \right) + T_{\chi^2,b,\text{up}}.
\end{align*}
Let $q^\ast_{1-\beta/2,n}$ be the $1-\beta/2$ quantile of the above upper bound, which means that $q_{1-\beta/2,n} \leq q^\ast_{1-\beta/2,n}$. For now, we take the following two bounds for granted:
\begin{align} \label{Eq: Bounds on Poisson Moments}
\mE_P[\Sigma_{n,\text{pois}}^2] \leq C_5 \min\{n,d\} \quad \text{and} \quad \mE_P[T_{\chi^2,b,\text{up}}] \leq C_6,
\end{align}
which are formally proved in Appendix~\ref{Section: Verification of bounds}. Then by using Markov's inequality, for any $t_1,t_2>0$ and $t=t_1+t_2$,
\begin{align*}
\mP_P \big[ C_5 \Sigma_{n,\text{pois}} \log (\alpha^{-1})+ T_{\chi^2,b,\text{up}} \geq t  \big] ~ \leq~ &  \mP_P \big[C_5 \Sigma_{n,\text{pois}} \log (\alpha^{-1}) \geq t_1 \big] + \mP_P \big[ T_{\chi^2,b,\text{up}} \geq t_2 \big] \\[.5em]
\leq ~ & C_6\frac{\mE_P[\Sigma_{n,\text{pois}}^2] \{\log(\alpha^{-1})\}^2}{t_1^2} + C_7 \frac{\mE_P[T_{\chi^2,b,\text{up}}]}{t_2} \\[.5em]
\leq ~ & C_7 \frac{\min(n,d) \{\log(\alpha^{-1})\}^2}{t_1^2} + \frac{C_8}{t_2}.
\end{align*}
Then by setting the upper bound to be $\beta/2$, one may see that for sufficiently large $C_9>0$,
\begin{align*}
q^\ast_{1-\beta/2,n} \leq \frac{C_9}{\beta^{1/2}} \log \left( \frac{1}{\alpha} \right) \sqrt{\min\{n,d\}} + \frac{C_{10}}{\beta},
\end{align*}
which in turn shows that condition~(\ref{Eq: Sufficient condition for Poisson testing}) is satisfied. 

\subsection{Verification of two bounds in (\ref{Eq: Bounds on Poisson Moments})} \label{Section: Verification of bounds}
This section proves the bounds in (\ref{Eq: Bounds on Poisson Moments}), namely, $(a)~\mE_P[\Sigma_{n,\text{pois}}^2] \leq C_1 \min\{n,d\}$ and $(b)~\mE_P[T_{\chi^2,b,\text{up}}] \leq C_2$. 

\paragraph{$\bullet$ Bound (a).} We start by proving $\mE_P[\Sigma_{n,\text{pois}}^2] \leq C_1 \min\{n,d\}$. By recalling the definition of $\Sigma_{n,\text{pois}}$ in (\ref{Eq: definition of Sigma n}), note that
\begin{align*}
\Sigma_{n,\text{pois}}^2 ~ = ~ & \sum_{(i,j) \in \mathbf{i}_2^{n}} \Bigg\{ \sum_{k=1}^d \frac{Y_{i,k}Y_{j,k}}{\omega_k} \ind(\omega_k > 0) \Bigg\}^2 + \sum_{(i,j) \in \mathbf{i}_2^{n}} \Bigg\{ \sum_{k=1}^d \frac{Z_{i,k}Z_{j,k}}{\omega_k} \ind(\omega_k > 0) \Bigg\}^2 \\
  + 2& \sum_{1 \leq i,j \leq n} \Bigg\{ \sum_{k=1}^d \frac{Y_{i,k}Z_{j,k}}{\omega_k} \ind(\omega_k > 0) \Bigg\}^2 \\[.5em]
:= ~ &  \Sigma_{n,Y}^2 + \Sigma_{n,Z}^2 + 2 \Sigma_{n,YZ}^2 \quad \text{(say).}
\end{align*}
Given $1\leq i \neq j \leq n$, expand the first squared term as
\begin{align*}
\Bigg\{ \sum_{k=1}^d \frac{Y_{i,k}Y_{j,k}}{\omega_k} \ind(\omega_k > 0) \Bigg\}^2 ~=~  &  \sum_{k=1}^d \omega_k^{-2} Y_{i,k}^2 Y_{j,k}^2 \ind(\omega_k > 0) \\[.5em] 
+ & \sum_{(k_1,k_2) \in \mathsf{i}_{2}^d} \omega_{k_1}^{-1}\omega_{k_2}^{-1} Y_{i,k_1} Y_{j,k_1} Y_{i,k_2} Y_{j,k_2} \ind(\omega_{k_1} > 0)\ind(\omega_{k_2} > 0) \\[.5em]
=~& (I) + (II) \quad \text{(say).}
\end{align*}
Let us first look at the expectation of $(I)$. Suppose that $Q_1,\ldots,Q_n$ are independent Poisson random variables with parameters $\lambda_1,\ldots,\lambda_n$, respectively. To calculate the above expectation, we use the fact that conditional on the event $\sum_{i=1}^n Q_i = N$, $(Q_1,\ldots,Q_n)$ has a multinomial distribution as
\begin{align*}
(Q_1,\ldots,Q_n) \sim \text{Multinomial}\left(N, ~ \Bigg\{ \frac{\lambda_1}{\sum_{i=1}^n\lambda_i}, \ldots, \frac{\lambda_n}{\sum_{i=1}^n\lambda_i} \Bigg\}  \right).
\end{align*}
Therefore, conditioned on $\omega_k = N$, we observe that 
\begin{align*}
& (Y_{i,k},Y_{j,k}, \omega_k - Y_{i,k} - Y_{j,k}) \\
\sim ~ & \text{Multinomial}\left(N, ~ \Bigg[ \frac{p_Y(k)}{n\{p_Y(k) + p_Z(k)\}}, \ \frac{p_Y(k)}{n\{p_Y(k) + p_Z(k)\}}, \ 1 -  \frac{2p_Y(k)}{n\{p_Y(k) + p_Z(k)\}} \Bigg] \right).
\end{align*}
Using this property and the moment generating function (MGF) of a multinomial distribution (see Appendix~\ref{Section: Multinomial Moments}),
\begin{align*}
\mathbb{E} \big[Y_{i,k}^2 Y_{j,k}^2 ~|~ \omega_k = N\big] ~=~ & N(N-1)(N-2)(N-3) \widetilde{p}_{k,n}^4  \\
+ ~& 2N(N-1)(N-2)\widetilde{p}_{k,n}^3 + N(N-1)\widetilde{p}_{k,n}^2,
\end{align*}
where 
\begin{align*}
\widetilde{p}_{k,n} :=  \frac{p_Y(k)}{n\{p_Y(k)+ p_Z(k)\}}.
\end{align*}
This gives
\begin{align*}
\mE \Bigg[ \frac{Y_{i,k}^2 Y_{j,k}^2 }{\omega_k^{2} }\ind(\omega_k > 0)  \Bigg] ~=~&  \mathbb{E}_N \left[\frac{Y_{i,k}^2 Y_{j,k}^2 }{\omega_k^{2} }\ind(\omega_k > 0) ~\Bigg|~ \omega_k = N \Bigg\} \right] \\[.5em]
\leq ~& \mathbb{E}_N \left[ N^2 \widetilde{p}_{k,n}^4 \ind(N > 0) \right] + 2\mathbb{E}_N \left[ N \widetilde{p}_{k,n}^3 \ind(N > 0) \right] + \mathbb{E}_N \left[ \widetilde{p}_{k,n}^2 \ind(N > 0) \right].
\end{align*}
By noting that $N \sim \text{Poisson}(n\{p_Y(k)+p_Z(k)\})$,
\begin{align*}
\mathbb{E}_N \left[ N^2 \widetilde{p}_{k,n}^4 \ind(N > 0) \right] ~=~& \widetilde{p}_{k,n}^4 \mathbb{E}_N \left[ N^2 \ind(N > 0) \right] \\[.5em] 
= ~& \left(\frac{p_Y(k)}{n\{p_Y(k)+p_Z(k)\}}\right)^4 (n\{p_Y(k)+p_Z(k)\})^2 \\[.5em]
+& \left(\frac{p_Y(k)}{n\{p_Y(k)+p_Z(k)\}}\right)^4  n\{p_Y(k)+p_Z(k)\} \\[.5em]
\leq ~ & \frac{p_Y(k)^4}{(n\{p_Y(k)+p_Z(k)\})^2} + \frac{p_Y(k)^4}{(n\{p_Y(k)+p_Z(k)\})^3},
\end{align*}
and 
\begin{align*}
\mathbb{E}_N \left[ N \widetilde{p}_{k,n}^3 \ind(N > 0) \right] ~=~   & \widetilde{p}_{k,n}^3 \mathbb{E}_N \left[ N  \ind(N > 0) \right] \\[.5em]
~=~&  \frac{p_Y(k)^3}{(n\{p_Y(k)+p_Z(k)\})^2}, \\[.5em]
\mathbb{E}_N \left[ \widetilde{p}_{k,n}^2 \ind(N > 0) \right] ~=~ & \widetilde{p}_{k,n}^2 \mathbb{E}_N \left[ \ind(N > 0) \right] \\[.5em] 
~=~& \left(\frac{p_Y(k)}{n\{p_Y(k)+p_Z(k)\}}\right)^2 \times \left( 1 - e^{-n\{p_Y(k)+p_Z(k)\}} \right).
\end{align*}
Putting these together, 
\begin{align} \nonumber
\mE[(I)]~=~ & \sum_{k=1}^d \mathbb{E} \left[ \frac{Y_{i,k}^2 Y_{j,k}^2}{\omega_k^2} \ind(\omega_k > 0) \right] \\[.5em]  \nonumber
~\leq~ & \sum_{k=1}^d \frac{p_Y(k)^4}{(n\{p_Y(k)+p_Z(k)\})^2} + \sum_{k=1}^d \frac{p_Y(k)^4}{(n\{p_Y(k)+p_Z(k)\})^3}  + \sum_{k=1}^d  \frac{p_Y(k)^3}{(n\{p_Y(k)+p_Z(k)\})^2} \\[.5em] \nonumber
+ & \sum_{k=1}^d  \left(\frac{p_Y(k)}{n\{p_Y(k)+p_Z(k)\}}\right)^2 \times \left( 1 - e^{-n\{p_Y(k)+p_Z(k)\}} \right) \\[.5em]
\leq ~ & \frac{1}{n^2} + \frac{1}{n^3} + \frac{1}{n^2} + \frac{1}{n^2} \min\{d, 2n\},  \label{Eq: Bound on (I)}
\end{align}
where the last inequality uses $1 - e^{-x} \leq \min\{x,1\}$.

Next moving onto the expected value of $(II)$, the independence between Poisson random variables gives
\begin{align*}
& \mathbb{E} \left[\frac{Y_{i,k_1} Y_{j,k_1}}{\omega_{k_1}} \frac{Y_{i,k_2} Y_{j,k_2}}{\omega_{k_2}} \ind(\omega_{k_1} > 0)\ind(\omega_{k_2} > 0)\right] \\[.5em]
= ~ &  \mathbb{E} \left[\frac{Y_{i,k_1} Y_{j,k_1}}{\omega_{k_1}} \ind(\omega_{k_2} > 0)\right]  \mathbb{E} \left[ \frac{Y_{i,k_2} Y_{j,k_2}}{\omega_{k_2}} \ind(\omega_{k_1} > 0)\right]. 
\end{align*}
Again, $(Y_{i,k_1}, Y_{j,k_1},\omega_{k_1} - Y_{i,k_1}- Y_{j,k_1})$ has a multinomial distribution conditional on $\omega_{k_1} = N$. Based on this property, we have
\begin{align*}
\mathbb{E}[Y_{i,k_1}Y_{j,k_1}|\omega_{k_1} = N] = N(N-1)  \widetilde{p}_{k_1,n}^2.
\end{align*}
Thus 
\begin{align*}
\mathbb{E}\left[ \frac{Y_{i,k_1} Y_{j,k_1}}{\omega_{k_1}} \ind(\omega_{k_1} > 0) \right] ~=~ & \mathbb{E}_N \left[ \mathbb{E} \Bigg\{ \frac{Y_{i,k_1} Y_{j,k_1}}{\omega_{k_1}} \ind(\omega_{k_1} > 0) \Bigg| \omega_{k_1} = N \Bigg\} \right] \\[.5em]
=~ & \widetilde{p}_{k_1,n}^2 \mathbb{E}_N \left[ (N-1) \ind(N>0)\right] \\[.5em]
=~ & \widetilde{p}_{k_1,n}^2 \bigg[ n\{p_Y(k_1) + p_Z(k_1) \} - 1 + e^{-n\{p_Y(k_1) + p_Z(k_1) \}} \bigg] \\[.5em]
\leq ~ & \frac{p_Y^2(k_1)}{n\{p_Y(k_1) + p_Z(k_1)\}}.
\end{align*}
This gives
\begin{align} \nonumber 
\mE[(II)] ~=~ & \mE\Bigg[ \sum_{(k_1,k_2) \in \mathsf{i}_{2}^d} \omega_{k_1}^{-1}\omega_{k_2}^{-1} Y_{i,k_1} Y_{j,k_1} Y_{i,k_2} Y_{j,k_2} \ind(\omega_{k_1} > 0)\ind(\omega_{k_2} > 0) \Bigg]  \\[.5em] 
~ \leq ~ & \left( \sum_{i_1 =1}^d \frac{p_{i_1}^2}{n\{p_{i_1} + q_{i_1}\}} \right) \cdot \left( \sum_{i_2 =1}^d \frac{p_{i_2}^2}{n\{p_{i_2} + q_{i_2}\}} \right)  ~\leq~  \frac{1}{n^2}. \label{Eq: Bound on (II)}
\end{align}
Therefore based on (\ref{Eq: Bound on (I)}) and (\ref{Eq: Bound on (II)}), it is clear that $\mE_P[\Sigma_{n,Y}^2] \leq C_2 \min\{n,d\}$. The same analysis further shows that $\mE_P[\Sigma_{n,Z}^2] \leq C_3 \min\{n,d\}$ and $\mE_P[\Sigma_{n,YZ}^2] \leq C_4 \min\{n,d\}$, which leads to $\mE_P[\Sigma_{n,\text{pois}}^2] \leq C_1 \min\{n,d\}$ as desired.

\paragraph{$\bullet$ Bound (b).} Next we prove that $\mE_P[T_{\chi^2,b,\text{up}}] \leq C_2$. Recall that $T_{\chi^2,b,\text{up}}$ is a nonnegative random variable defined in (\ref{Eq: Poisson Result 2}). Since $\omega_k \sim \text{Poisson}(n\{p_Y(k) + p_Z(k)\})$, the second term of $T_{\chi^2,b,\text{up}}$ satisfies
\begin{align*}
\sum_{k=1}^d \mathbb{E}[\ind(\omega_k > 0 )] ~=~ \sum_{k=1}^d \left(1 - e^{-n\{p_Y(k) + p_Z(k)\}}\right).
\end{align*}
Next consider the first term of $T_{\chi^2,b,\text{up}}$:
\begin{align*}
\sum_{i=1}^{2n} \sum_{k=1}^d \frac{X_{i,k}^2}{\omega_k} \ind(\omega_k > 0).
\end{align*}
Note that based on the moments of a multinomial distribution (see Appendix~\ref{Section: Multinomial Moments}), one can compute
\begin{align*}
&\mathbb{E}\big[Y_{i,k}^2 | \omega_k = N\big] = N(N-1) \widetilde{p}_{k,n}^2 + n\widetilde{p}_{k,n}, \\[.5em]
& \mathbb{E}\big[Z_{i,k}^2 | \omega_k = N\big] =  N(N-1) \widetilde{q}_{k,n}^2 + n\widetilde{q}_{k,n},
\end{align*}
where $\widetilde{p}_{k,n} := p_Y(k) /\{n(p_Y(k) + p_Z(k)\}$ and $\widetilde{q}_{k,n} := p_Z(k) /\{n(p_Y(k) + p_Z(k))\}$. Therefore, by the law of total expectation,
\begin{align*}
\mathbb{E}\left[ \frac{Y_{i,k}^2}{\omega_k} \ind(\omega_k > 0) \right] ~ = ~ &\mathbb{E}_N \left[ \mathbb{E} \Bigg\{ \frac{Y_{i,k}^2}{\omega_k} \ind(\omega_k > 0) \Bigg| \omega_k = N \Bigg\} \right] \\[.5em]
= ~ & \widetilde{p}_{k,n}^2 \mathbb{E}_N \left[ (N-1) \ind(N>0) \right] + \widetilde{p}_{k,n} \mathbb{E}_N \left[ \ind(N>0) \right] \\[.5em]
= ~ & \widetilde{p}_{k,n}^2  \left( n\{p_Y(k) + p_Z(k)\} - 1 + e^{-n\{p_Y(k) + p_Z(k)\}} \right) \\[.5em] 
 + & \widetilde{p}_{k,n} \left(1 - e^{-n\{p_Y(k) + p_Z(k)\}}\right).
\end{align*}
Similarly, one can compute
\begin{align*}
\mathbb{E}\left[ \frac{Z_{i,k}^2}{\omega_k} \ind(\omega_k > 0) \right] ~ = ~ & \widetilde{q}_{k,n}^2  \left( n\{p_Y(k) + p_Z(k)\} - 1 + e^{-n\{p_Y(k) + p_Z(k)\}} \right) \\[.5em] 
 +& \widetilde{q}_{k,n} \left(1 - e^{-n\{p_Y(k) + p_Z(k)\}}\right).
\end{align*}
Based on the definition of $\widetilde{p}_{k,n}$ and $\widetilde{q}_{k,n}$, we have the identity
\begin{align*}
&\sum_{i=1}^n \sum_{k=1}^d \widetilde{p}_{k,n} \left(1 - e^{-n\{p_Y(k) + p_Z(k)\}}\right) + \sum_{i=1}^n \sum_{k=1}^d \widetilde{q}_{k,n} \left(1 - e^{-n\{p_Y(k) + p_Z(k)\}}\right) \\[.5em] 
= ~ & \sum_{k=1}^d \left(1 - e^{-n\{p_Y(k) + p_Z(k)\}}\right),
\end{align*}
which is the expected value of $\sum_{k=1}^d \ind(\omega_k > 0)$. Putting everything together, 
\begin{align*}
\mE \big[T_{\chi^2,b,\text{up}}\big] ~=~ & \sum_{i=1}^{2n} \sum_{k=1}^d \mE\Bigg[\frac{X_{i,k}^2}{\omega_k} \ind(\omega_k > 0) \Bigg] -\sum_{k=1}^d \mathbb{E}[\ind(\omega_k > 0 )] \\[.5em]
\leq ~ & \sum_{i=1}^n \sum_{k=1}^d  \frac{p_Y^2(k)}{n\{p_Y(k) + p_Z(k)\}}  + \sum_{i=1}^n \sum_{k=1}^d  \frac{p_Z^2(k)}{n\{p_Y(k) + p_Z(k)\}} \\[.5em]
 ~\leq~ & 2.
\end{align*}
This proves the bound $\mE[T_{\chi^2,b,\text{up}}] \leq C_2$.

\subsection{Details on verifying the sufficient condition (\ref{Eq: Goal of Poisson two-sample testing})} \label{Section: Details on the verification}

First assume that $n < d$. Then the variance (\ref{Eq: Bounds for V}) is dominated by the first term and thus condition~(\ref{Eq: Goal of Poisson two-sample testing}) is fulfilled when
\begin{align*}
\mE_P [T_{\chi^2}] ~\overset{(i)}{\geq}~ & \frac{n^2}{6d}\|p_Y - p_Z\|_1^2 ~\overset{(ii)}{\geq} ~ \frac{n^2}{6d} \epsilon_n^2 ~ \overset{(iii)}{\geq} ~ \frac{C_2}{\beta^{1/2}} \log\left( \frac{1}{\alpha}\right)  \sqrt{n} \\[.5em]
~ \geq ~ &  q_{1-\beta/2,n} + \sqrt{\frac{2\mV_P[T_{\chi^2}]}{\beta}},
\end{align*}
where $(i)$ follows by the bound (\ref{Eq: Bounds for M}), $(ii)$ uses $\|p_Y - p_Z\|_1 \geq \epsilon_n$ and $(iii)$ holds from the bounds (\ref{Eq: Bounds for V}) and (\ref{Eq: Sufficient condition for Poisson testing}) and the condition on $\epsilon_n$, i.e.
\begin{align*}
\epsilon_n ~\geq~ \frac{C_3}{\beta^{1/2}} \sqrt{\log \left( \frac{1}{\alpha} \right)} \frac{d^{1/2}}{n^{3/4}},
\end{align*}
for some large constant $C_3>0$.

Next assume that $n \geq d$. For convenience, let us write
\begin{align*}
\varphi_k := 1 -  \frac{1 - e^{-n p_Y(k) - np_Z(k)}}{np_Y(k) + np_Z(k)} \quad \text{for $k=1,\ldots,d$.}
\end{align*}
We define $\mathsf{I}_d := \{ k \in \{1,\ldots,d\} : 2 \varphi_k  \geq 1 \}$ and denote its complement by $\mathsf{I}_d^c$. Note that $np_Y(k) + np_Z(k) < 2$ for $k \in \mathsf{I}_d^c$ and thus 
\begin{align*}
n\sum_{k=1}^d \frac{\{p_Y(k) - p_Z(k) \}^2}{p_Y(k) + p_Z(k)} ~ = ~ & n\sum_{k \in \mathsf{I}_d} \frac{\{p_Y(k) - p_Z(k) \}^2}{p_Y(k) + p_Z(k)} + n\sum_{k \in \mathsf{I}_d^c} \frac{\{p_Y(k) - p_Z(k) \}^2}{p_Y(k) + p_Z(k)} \\[.5em]
\leq ~ & n\sum_{k \in \mathsf{I}_d} \frac{\{p_Y(k) - p_Z(k) \}^2}{p_Y(k) + p_Z(k)} + 2d.
\end{align*}
Based on this observation with $n \geq d$, the variance of $T_{\chi^2}$ can be further bounded by
\begin{align} \label{Eq: Bound on V}
\mV_P[T_{\chi^2}] ~  \leq ~ 4d + 5n\sum_{k \in \mathsf{I}_d} \frac{\{p_Y(k) - p_Z(k) \}^2}{p_Y(k) + p_Z(k)}.
\end{align}
Let us make one more observation that $n^2\|p_Y - p_Z\|_1^2/(4d+2n) \geq C_4 \beta^{-1}$ for some large constant $C_4>0$, which holds under the assumption on $\epsilon_n$ in Theorem~\ref{Theorem: Two-Sample testing under Poisson sampling} and $n \geq d$. Based on this, the expectation of $T_{\chi^2}$ is bounded by
\begin{equation}
\begin{aligned} \label{Eq: Lower bound of the expectation}
\mE_P [T_{\chi^2}] ~ \geq ~ & \sqrt{\sum_{k=1}^d \frac{\{p_Y(k) - p_Z(k) \}^2}{p_Y(k) + p_Z(k)} n \varphi_k } \sqrt{\frac{n^2}{4d + 2n} \| p_Y - p_Z\|_1^2} \\[.5em]
\geq ~ &  \sqrt{ \frac{C_4n}{2\beta} \sum_{k \in \mathsf{I}_d} \frac{\{p_Y(k) - p_Z(k) \}^2}{p_Y(k) + p_Z(k)}}, 
\end{aligned} 
\end{equation} 
where the last inequality uses the definition of $\mathsf{I}_d$. This gives
\begin{align*}
\mE_P [T_{\chi^2}] ~ \overset{(i)}{\geq} ~ & \frac{1}{2}\mE_P [T_{\chi^2}]  + \frac{1}{2}\sqrt{ \frac{C_4n}{2\beta} \sum_{k \in \mathsf{I}_d} \frac{\{p_Y(k) - p_Z(k) \}^2}{p_Y(k) + p_Z(k)}} \\[.5em]
\overset{(ii)}{\geq} ~ & \frac{n}{12} \|p_Y - p_Z\|_1^2 + \frac{1}{2}\sqrt{ \frac{C_4 n}{2\beta} \sum_{k \in \mathsf{I}_d} \frac{\{p_Y(k) - p_Z(k) \}^2}{p_Y(k) + p_Z(k)}} \\[.5em]
\overset{(iii)}{\geq} ~ & \frac{n}{12} \epsilon_n^2  + \frac{1}{2}\sqrt{ \frac{C_4 n}{2\beta} \sum_{k \in \mathsf{I}_d} \frac{\{p_Y(k) - p_Z(k) \}^2}{p_Y(k) + p_Z(k)}} \\[.5em]
\overset{(iv)}{\geq} ~ & \frac{C_5}{\beta} \log\left( \frac{1}{\alpha}\right)d^{1/2}  + \frac{1}{2}\sqrt{ \frac{C_4 n}{2\beta} \sum_{k \in \mathsf{I}_d} \frac{\{p_Y(k) - p_Z(k) \}^2}{p_Y(k) + p_Z(k)}} \\[.5em]
\overset{(v)}{\geq} ~ & \frac{C_1}{\beta} \log\left( \frac{1}{\alpha}\right)d^{1/2} + \sqrt{\frac{8d}{\beta}} + \sqrt{\frac{10n}{\beta}\sum_{k \in \mathsf{I}_d} \frac{\{p_Y(k) - p_Z(k) \}^2}{p_Y(k) + p_Z(k)}} \\[.5em]
\overset{(vi)}{\geq} ~ & q_{1-\beta/2,n} + \sqrt{\frac{2\mV_P[T_{\chi^2}]}{\beta}},
\end{align*}
where $(i)$ uses the lower bound (\ref{Eq: Lower bound of the expectation}), $(ii)$ and $(iii)$ follow by the bound (\ref{Eq: Bounds for M}) and $\|p_Y - p_Z\|_1 \geq \epsilon_n$, respectively, $(iv)$ follows from the lower bound for $\epsilon_n$ in the theorem statement, $(v)$ holds by choosing $C_4,C_5$ large and lastly $(vi)$ uses (\ref{Eq: Sufficient condition for Poisson testing}) and (\ref{Eq: Bound on V}).

\subsection{Multinomial Moments} \label{Section: Multinomial Moments}
This section collects some moments of a multinomial distribution that are used in the proof of Theorem~\ref{Theorem: Two-Sample testing under Poisson sampling}. Suppose that $\mathbf{X} = (X_1,\ldots,X_d)$ has  a multinomial distribution with the number of trials $n$ and probabilities $(p_1,\ldots,p_d)$. The MGF of $\mathbf{X}$ is given by
\begin{align*}
M_{\mathbf{X}}(\mathbf{t}) = \left(\sum_{i=1}^d p_i e^{t_i} \right)^n.
\end{align*}
We collect some of partial derivatives of the MGF. 
\begin{align*}
\frac{\partial}{\partial t_i}M_{\mathbf{X}}(\mathbf{t}) ~=~ &n \left(\sum_{i=1}^d p_i e^{t_i} \right)^{n-1}p_i e^{t_i}, \\[.5em]
\frac{\partial^2}{\partial t_i\partial t_j}M_{\mathbf{X}}(\mathbf{t}) ~=~& n(n-1) \left(\sum_{i=1}^d p_i e^{t_i} \right)^{n-2}p_i e^{t_i} p_j e^{t_j},\\[.5em]
\frac{\partial^2}{\partial t_i^2}M_{\mathbf{X}}(\mathbf{t}) ~=~& n(n-1) \left(\sum_{i=1}^d p_i e^{t_i} \right)^{n-2}p_i^2 e^{2t_i} + n \left(\sum_{i=1}^d p_i e^{t_i} \right)^{n-1}p_i e^{t_i},  \\[.5em]
\frac{\partial^3}{\partial t_i^2 \partial t_j}M_{\mathbf{X}}(\mathbf{t}) ~=~& n(n-1)(n-2) \left(\sum_{i=1}^d p_i e^{t_i} \right)^{n-3}p_i^2 p_j  e^{2t_i} e^{t_j}  \\
& + n(n-1) \left(\sum_{i=1}^d p_i e^{t_i} \right)^{n-2}p_i p_j e^{t_i} e^{t_j},  \\[.5em]
\frac{\partial^4}{\partial t_i^2 \partial t_j^2}M_{\mathbf{X}}(\mathbf{t}) ~=~& n(n-1)(n-2)(n-3) \left(\sum_{i=1}^d p_i e^{t_i} \right)^{n-4}p_i^2 p_j^2  e^{2t_i} e^{2t_j}  \\
& + n(n-1)(n-2) \left(\sum_{i=1}^d p_i e^{t_i} \right)^{n-3}p_i^2 p_j  e^{2t_i} e^{t_j} \\
& + n(n-1)(n-2) \left(\sum_{i=1}^d p_i e^{t_i} \right)^{n-3}p_i p_j^2  e^{t_i} e^{2t_j} \\
& + n(n-1) \left(\sum_{i=1}^d p_i e^{t_i} \right)^{n-2}p_i p_j e^{t_i} e^{t_j}.
\end{align*}
By setting $\mathbf{t} =0$, for $i \neq j$,
\begin{align*}
\mathbb{E}[X_i] ~=~ &  np_i, \\
\mathbb{E}[X_i^2] ~=~ &  n(n-1)p_i^2 + np_i, \\
\mathbb{E}[X_iX_j] ~=~ &  n(n-1)p_ip_j, \\ 
\mathbb{E}[X_i^2X_j] ~=~ &  n(n-1)(n-2)p_i^2p_j + n(n-1)p_ip_j, \\
\mathbb{E}[X_i^2X_j^2] ~=~ &  n(n-1)(n-2)(n-3)p_i^2p_j^2 + n(n-1)(n-2)p_i^2p_j \\
+ & n(n-1)(n-2)p_ip_j^2 + n(n-1)p_ip_j. 
\end{align*}

\section{Proof of Proposition~\ref{Proposition: Multinomial L1 Testing}} \label{Section: Proof of Proposition: Multinomial L1 Testing}
Recall that the test is carried out via sample-splitting and the critical value of the permutation test is obtained by permuting the labels within $\mathcal{X}_{2n_1}^{\text{split}} = \{Y_1,\ldots,Y_{n_1},Z_1,\ldots,Z_{n_1}\}$. Nevertheless, the distribution of the test statistic is invariant to any partial permutation under the null hypothesis. Based on this property, it can be shown that type I error control of the permutation test via sample-splitting is also guaranteed \citep[see e.g.~Theorem 15.2.1 of][]{lehmann2006testing}. Hence we focus on the type II error control. Note that conditional on $w_1,\ldots,w_d$, the test statistic $U_{n_1,n_2}^{\text{split}}$ can be viewed as a $U$-statistic with kernel $g_{\text{Multi},w}(x,y)$ given in (\ref{Eq: weighted kernel}). Moreover this $U$-statistic is based on the two samples of equal size, which allows us to apply Lemma~\ref{Lemma: Two-Sample U-statistic Improved Version}. Based on this observation, we first study the performance of the test conditioning on $w_1,\ldots,w_d$. We then remove this conditioning part using Markov's inequality and conclude the result. 

\paragraph{$\bullet$ Conditional Analysis.} In this part, we investigate the type II error of the permutation test conditional on $w_1,\ldots,w_d$. As noted earlier, $U_{n_1,n_2}^{\text{split}}$ can be viewed as a $U$-statistic and so we can apply Lemma~\ref{Lemma: Two-Sample U-statistic Improved Version} to proceed. To do so, we need to lower bound the conditional expectation of $U_{n_1,n_2}^{\text{split}}$ and upper bound $\psi_{Y,1}(P)$, $\psi_{Z,1}(P)$ and $\psi_{YZ,2}(P)$. On the one hand, the conditional expectation of $U_{n_1,n_2}^{\text{split}}$ is lower bounded by the squared $\ell_1$ distance as
\begin{align*}
\mE_P \big[U_{n_1,n_2}^{\text{split}}|w_1,\ldots,w_n \big] = \sum_{k=1}^d \frac{[ p_Y(k) - p_Z(k) ]^2}{w_k} \geq \| p_Y - p_Z\|_1^2,
\end{align*}
where the inequality follows by Cauchy-Schwarz inequality and $\sum_{k=1}^d w_k = 1$. On the other hand, $\psi_{Y,1}(P)$, $\psi_{Z,1}(P)$ and $\psi_{YZ,2}(P)$ are upper bounded by
\begin{equation}
\begin{aligned} \label{Eq: Bounding psi functions}
\psi_{Y,1}(P) & ~\leq~  4 \sqrt{\sum_{k=1}^d \frac{p_Y^2(k)}{w_k^2}} \sum_{k=1}^d \frac{[p_Y(k) - p_Z(k)]^2}{w_k}\\[.5em]
\psi_{Z,1}(P) & ~\leq~  4 \sqrt{\sum_{k=1}^d \frac{p_Z^2(k)}{w_k^2}} \sum_{k=1}^d \frac{[p_Y(k) - p_Z(k)]^2}{w_k}\\[.5em]
\psi_{YZ,2}(P) & ~\leq~ \max \Bigg\{\sum_{k=1}^d \frac{p_Y^2(k)}{w_k^2},  \ \sum_{k=1}^d \frac{p_Z^2(k)}{w_k^2}\Bigg\}.
\end{aligned}
\end{equation}
The details of the derivations are presented in Section~\ref{Section: Details on Eq: Bounding psi functions}. Further note that 
\begin{equation}
\begin{aligned} \label{Eq: Bounding p/w}
\sum_{k=1}^d \frac{p_Y^2(k)}{w_k^2} \overset{(i)}{\leq} ~ & 2 \sum_{k=1}^d \frac{[p_Y(k) - p_Z(k)]^2}{w_k^2} + 2 \sum_{k=1}^d \frac{p_Z^2(k)}{w_k^2} \\[.5em]
\overset{(ii)}{\leq} ~ & 4d \sum_{k=1}^d \frac{[p_Y(k) - p_Z(k)]^2}{w_k} + 2 \sum_{k=1}^d \frac{p_Z^2(k)}{w_k^2},
\end{aligned}
\end{equation}
where $(i)$ uses $(x+y)^2 \leq 2 x^2 + 2 y^2$ and $(ii)$ follows since $w_k \geq 1/ (2d)$ for $k=1,\ldots,d$. For notational convenience, let us write 
\begin{align*}
& \|p_Y - p_Z\|_w^2 := \sum_{k=1}^d \frac{[p_Y(k) - p_Z(k)]^2}{w_k}  \quad \text{and} \\[.5em]
& \|p_Z/w\|_2^2 := \sum_{k=1}^d \frac{p_Z^2(k)}{w_k^2}.
\end{align*}
Having this notation in place, we see that the condition (\ref{Eq: Two-Sample U-statistic Sufficient Condition 2}) in Lemma~\ref{Lemma: Two-Sample U-statistic Improved Version} is fulfilled when
\begin{align*}
& \sqrt{\frac{\psi_{Y,1}(P)}{\beta n_1}} \leq C_1 \|p_Y - p_Z\|_w^2, \\[.5em]
& \sqrt{\frac{\psi_{Z,1}(P)}{\beta n_1}} \leq C_2 \|p_Y - p_Z\|_w^2 \quad \text{and} \\[.5em]
& \sqrt{\frac{\psi_{YZ,2}(P)}{\beta}} \log \left( \frac{1}{\alpha} \right) \frac{1}{n_1} \leq C_3 \|p_Y - p_Z\|_w^2.
\end{align*}
Based on the results in (\ref{Eq: Bounding psi functions}) and (\ref{Eq: Bounding p/w}), it can be shown that these three inequalities are implied by
\begin{align*}
& \| p_Y- p_Z \|_w^2 ~\geq~ \frac{C_4}{\beta} \log \left( \frac{1}{\alpha} \right) \frac{\|p_Z/w\|_2}{n_1} \quad \text{and} \\[.5em]
& \| p_Y - p_Z \|_w^2 ~\geq~ \frac{C_5}{\beta^2} \log^2 \left( \frac{1}{\alpha} \right) \frac{d}{n_1^2}. 
\end{align*} 
Moreover, using the lower bound of the conditional expectation $\| p_Y - p_Z \|_w^2 \geq \| p_Y - p_Z\|_1^2 \geq \epsilon_{n_1,n_2}^2$ and the boundedness of $\ell_1$ norm so that $\epsilon_{n_1,n_2}^2 \leq 4$, the above two inequalities are further implied by
\begin{align*}
& \epsilon_{n_1,n_2}^2 ~ \geq ~ \frac{C_4}{\beta} \log \left( \frac{1}{\alpha} \right) \frac{\|p_Z/w\|_2}{n_1} \quad \text{and}  \\[.5em]
& \epsilon_{n_1,n_2}^2 ~ \geq ~ \frac{2\sqrt{C_5}}{\beta} \log \left( \frac{1}{\alpha} \right) \frac{d^{1/2}}{n_1}.
\end{align*}
In other words, for a sufficiently large $C_6>0$, the type II error of the permutation test is at most $\beta$ when
\begin{align} \label{Eq: sufficient condition on eps}
\epsilon_{n_1,n_2}^2 \geq \frac{C_6}{\beta} \log \left( \frac{1}{\alpha} \right) \max\Bigg\{ \frac{\|p_Z/w\|_2}{n_1}, \ \frac{d^{1/2}}{n_1} \Bigg\}. 
\end{align}
Note that the above condition is not deterministic as $w_1,\ldots,w_d$ are random variables. Next we remove this randomness. 

\paragraph{$\bullet$ Unconditioning $w_1,\ldots,w_d$.} 
Recall that $m = \min\{n_2,d\}$ and thus $w_k$ is clearly lower bounded by
\begin{align*}
w_k = \frac{1}{2d} + \frac{1}{2m} \sum_{i=1}^m \ind(Z_{i+n_2} = k) \geq \frac{1}{2d} \bigg[ 1 + \sum_{i=1}^m \ind(Z_{i+n_2} = k) \bigg].
\end{align*}
Based on this bound, one can see that $\|p_Z/w\|_2^2$ has the expected value upper bounded by
\begin{align*}
\mE_P\left[ \sum_{k=1}^d \frac{p_{Z}^2(k)}{w_k^2} \right] ~\leq~ & 4d^2 \sum_{k=1}^d \mE_P \left[ \frac{p_Z^2(k)}{\{1 + \sum_{i=1}^m \ind(Z_{i+n_2}  = k)\}^2} \right] \\[.5em]
\leq~ & 4d^2 \sum_{k=1}^d \mE_P \left[ \frac{p_Z^2(k)}{1 + \sum_{i=1}^m \ind(Z_{i+n_2}  = k)} \right] \\[.5em]
\overset{(i)}{\leq}~& 4d^2 \sum_{k=1}^d \frac{p_Z^2(k)}{(m+1) p_Z(k)} \\[.5em]
\leq~ & \frac{4d^2}{m},
\end{align*}
where $(i)$ uses the fact that when $X \sim \text{Binominal}(n,p)$, we have
\begin{align} \label{Eq: inverse binomial moment}
\mE\left[ \frac{1}{1+ X} \right] = \frac{1 - (1-p)^{n+1}}{(n+1)p} \leq \frac{1}{(n+1)p}.
\end{align}
Using this upper bound of the expected value, Markov's inequality yields
\begin{align*}
\mP_P \left( \sqrt{\sum_{k=1}^d \frac{p_{Z}^2(k)}{w_k^2}} \geq t \right) \leq \frac{1}{t^2}\mE \left[ \sum_{k=1}^d \frac{p_{Z}^2(k)}{w_k^2} \right] \leq \frac{4d^2}{mt^2}.
\end{align*}
By letting the right-hand side be $\beta$ and $\mathcal{A}$ be the event such that $\mathcal{A}:= \{ \|p_Z/w\|_2 < 2d/\sqrt{m\beta} \}$, we know that $\mP_P(\mathcal{A}) \geq 1 -\beta$. Under this good event $\mathcal{A}$, the sufficient condition (\ref{Eq: sufficient condition on eps}) is fulfilled when 
\begin{equation}
\begin{aligned} \label{Eq: unconditional sufficient condition}
\epsilon_{n_1,n_2}^2 \geq~& \frac{C_7}{\beta^{3/2}} \log \left( \frac{1}{\alpha} \right) \max\Bigg\{ \frac{d}{\sqrt{m}n_1}, \ \frac{d^{1/2}}{n_1} \Bigg\} \\[.5em]
= ~& \frac{C_7}{\beta^{3/2}} \log \left( \frac{1}{\alpha} \right) \max\Bigg\{ \frac{d}{n_1\sqrt{n_2}}, \ \frac{d^{1/2}}{n_1} \Bigg\}.
\end{aligned}
\end{equation}

\paragraph{$\bullet$ Completion of the proof.} To complete the proof, let us denote the critical value of the permutation test by $c_{1-\alpha,n_1,n_2}$. Then the type II error of the permutation test is bounded by
\begin{align*}
\mP_P(U_{n_1,n_2}^{\text{split}} \leq c_{1-\alpha,n_1,n_2}) =~ & \mP_P(U_{n_1,n_2}^{\text{split}} \leq c_{1-\alpha,n_1,n_2}, \mathcal{A}) + \mP_P(U_{n_1,n_2}^{\text{split}} \leq c_{1-\alpha,n_1,n_2}, \mathcal{A}^c) \\[.5em] 
\leq ~ & \mP_P(U_{n_1,n_2}^{\text{split}} \leq c_{1-\alpha,n_1,n_2}, \mathcal{A}) + \mP_P(\mathcal{A}^c).
\end{align*}
As shown before, the type II error under the event $\mathcal{A}$ is bounded by $\beta$, which leads to $\mP_P(U_{n_1,n_2}^{\text{split}} \leq c_{1-\alpha,n_1,n_2}, \mathcal{A}) \leq \beta$. Also we have $\mP_P(\mathcal{A}^c) \leq \beta$ proved by Markov's inequality. Thus the unconditional type II error is bounded by $2\beta$. Notice that condition~(\ref{Eq: unconditional sufficient condition}) is equivalent to condition (\ref{Eq: sufficient condition for L1 two-sample testing}) given in Proposition~\ref{Proposition: Multinomial L1 Testing}. Hence the proof is completed by letting $2\beta = \beta'$.

\subsection{Details on Equation~(\ref{Eq: Bounding psi functions})} \label{Section: Details on Eq: Bounding psi functions}
We start with bounding $\psi_{Y,1}(P)$. Following the proof of Proposition~\ref{Proposition: Multinomial Two-Sample Testing}, it can be seen that
\begin{align*}
\psi_{Y,1}(P) ~=~&  \mE_P \bigg[ \bigg( \sum_{k=1}^d w_k^{-1} [\ind(Y_1=k) - p_Y(k)][p_Y(k) - p_Z(k)]  \bigg)^2  \bigg| w_1,\ldots,w_d \bigg] \\[.5em]
\leq ~ & 2 \sum_{k=1}^d w_k^{-2} p_Y(k) [p_Y(k) - p_Z(k)]^2  + 2  \bigg( \sum_{k=1}^d w_k^{-1} p_Y(k)[p_Y(k) - p_Z(k)]  \bigg)^2 \\[.5em]
:= ~ & 2 (I) + 2 (II).
\end{align*}
For the first term $(I)$, we apply Cauchy-Schwarz inequality to have 
\begin{align*}
\sum_{k=1}^d w_k^{-2} p_Y(k) [p_Y(k) - p_Z(k)]^2 \leq ~ & \sqrt{\sum_{k=1}^d \frac{p_Y^2(k)}{w_k^2}} \sqrt{\sum_{k=1}^d \frac{[p_Y(k) - p_Z(k)]^4}{w_k^2}} \\[.5em]
\leq ~ &  \sqrt{\sum_{k=1}^d \frac{p_Y^2(k)}{w_k^2}} \sum_{k=1}^d \frac{[p_Y(k) - p_Z(k)]^2}{w_k},
\end{align*}
where the second inequality follows by the monotonicity of $\ell_p$ norm. For the second term $(II)$, we apply Cauchy-Schwarz inequality repeatedly to have 
\begin{align*}
\bigg( \sum_{k=1}^d w_k^{-1} p_Y(k)[p_Y(k) - p_Z(k)]  \bigg)^2 \leq ~ & \sum_{k=1}^d \frac{p_Y^2(k)}{w_k} \sum_{k=1}^d \frac{[p_Y(k) - p_Z(k)]^2}{w_k} \\[.5em]
\leq  ~ & \sqrt{\sum_{k=1}^d \frac{p_Y^2(k)}{w_k^2}} \sqrt{\sum_{k=1}^d p_Y^2(k)} \sum_{k=1}^d \frac{[p_Y(k) - p_Z(k)]^2}{w_k} \\[.5em]
\overset{(i)}{\leq} ~ & \sqrt{\sum_{k=1}^d \frac{p_Y^2(k)}{w_k^2}} \sum_{k=1}^d \frac{[p_Y(k) - p_Z(k)]^2}{w_k}
\end{align*}
where $(i)$ uses $\sum_{k=1}^d p_Y^2(k) \leq 1$. Combining the results yields
\begin{align*}
\psi_{Y,1}(P) \leq 4 \sqrt{\sum_{k=1}^d \frac{p_Y^2(k)}{w_k^2}} \sum_{k=1}^d \frac{[p_Y(k) - p_Z(k)]^2}{w_k}. 
\end{align*}
By symmetry, it similarly follows that 
\begin{align*}
\psi_{Z,1}(P) \leq 4 \sqrt{\sum_{k=1}^d \frac{p_Z^2(k)}{w_k^2}} \sum_{k=1}^d \frac{[p_Y(k) - p_Z(k)]^2}{w_k}. 
\end{align*}
These establish the first two inequalities in (\ref{Eq: Bounding psi functions}). Next we find an upper bound for $\psi_{YZ,2}(P)$. By recalling the definition of $\psi_{YZ,2}(P)$, we have
\begin{align*}
\psi_{YZ,2}(P) := \max \big\{ & \mE_P[g_{\text{Multi},w}^2(Y_1,Y_2)|w_1,\ldots,w_d], \  \mE_P[g_{\text{Multi},w}^2(Y_1,Z_1) | w_1,\ldots,w_d], \\[.5em] &\mE_P[g_{\text{Multi},w}^2(Z_1,Z_2) | w_1,\ldots,w_d] \big\}.
\end{align*}
Moreover each conditional expected value is computed as
\begin{align*}
\mE_P[g_{\text{Multi},w}^2(Y_1,Y_2)|w_1,\ldots,w_d] =~ &\sum_{k=1}^d w_k^{-2} p_Y^2(k), \\
\mE_P[g_{\text{Multi},w}^2(Z_1,Z_2)|w_1,\ldots,w_d] =~ &\sum_{k=1}^d w_k^{-2} p_Z^2(k), \\[.5em]
\mE_P[g_{\text{Multi},w}^2(Y_1,Z_1)|w_1,\ldots,w_d] =~ & \sum_{k=1}^d w_k^{-2} p_Y(k)p_Z(k) \\[.5em]
\leq ~ & \frac{1}{2} \sum_{k=1}^d w_k^{-2} p_Y^2(k) + \frac{1}{2} \sum_{k=1}^d w_k^{-2} p_Z^2(k) \\[.5em] 
\leq ~  & \max \Bigg\{\sum_{k=1}^d w_k^{-2} p_Y^2(k), \ \sum_{k=1}^d w_k^{-2} p_Z^2(k) \Bigg\}.
\end{align*}
This leads to
\begin{align*}
\psi_{YZ,2}(P) \leq  \max \Bigg\{\sum_{k=1}^d w_k^{-2} p_Y^2(k), \  \sum_{k=1}^d w_k^{-2} p_Z^2(k) \Bigg\}.
\end{align*}

\section{Proof of Proposition~\ref{Proposition: Multinomial independence testing in L1 distance}} \label{Section: Proof of Proposition: Multinomial independence testing in L1 distance}
We note that the test statistic considered in Proposition~\ref{Proposition: Multinomial independence testing in L1 distance} is essentially the same as that considered in Proposition~\ref{Proposition: Multinomial L1 Testing} with different weights. Hence following the same line of the proof of Proposition~\ref{Proposition: Multinomial L1 Testing}, we may arrive at the point (\ref{Eq: sufficient condition on eps}) where
the type II error of the considered permutation test is at most $\beta$ when
\begin{align} \label{Eq: sufficient condition 2}
\epsilon_{n}^2 \geq \frac{C}{\beta} \log \left( \frac{1}{\alpha} \right) \max\Bigg\{ \sqrt{\sum_{k_1=1}^{d_1}\sum_{k_2=1}^{d_2} \frac{p_Y^2(k)p_Z^2(k)}{w_{k_1,k_2}^2}}\frac{1}{n}, \ \frac{d_1^{1/2}d_2^{1/2}}{n} \Bigg\}. 
\end{align}
Similarly as before, let us remove the randomness from $w_{1,1},\ldots,w_{d_1,d_2}$ by applying Markov's inequality. First recall that $m_1 = \min\{n/2,d_1\}$ and $m_2 = \min\{n/2,d_2\}$ and thus
\begin{align*}
w_{k_1,k_2} = ~ & \Bigg[ \frac{1}{2d_1} + \frac{1}{2m_1} \sum_{i=1}^{m_1} \ind(Y_{3n/2+i} = k_1)  \Bigg] \times \Bigg[ \frac{1}{2d_2} + \frac{1}{2m_2} \sum_{j=1}^{m_2} \ind(Z_{5n/2+i} = k_2) \Bigg] \\[.5em]
\leq ~ &   \frac{1}{4d_1d_2} \Bigg[ 1 + \sum_{i=1}^{m_1} \ind(Y_{3n/2+i} = k_1)   \Bigg] \times \Bigg[ 1 + \sum_{j=1}^{m_2} \ind(Z_{5n/2+i} = k_2) \Bigg]. 
\end{align*}
Based on this observation, we have
\begin{align*}
& \mE_P \Bigg[ \sum_{k_1=1}^{d_1}\sum_{k_2=1}^{d_2} \frac{p_Y^2(k)p_Z^2(k)}{w_{k_1,k_2}^2} \Bigg] \\[.5em] 
\leq ~ & 16d_1^2d_2^2 \sum_{k_1=1}^{d_1}\sum_{k_2=1}^{d_2} \mE_P \bigg[ \frac{p_Y^2(k)p_Z^2(k)}{ \{1 + \sum_{i=1}^{m_1} \ind(Y_{3n/2+i} = k_1)\}^2 \{ 1 + \sum_{j=1}^{m_2} \ind(Z_{5n/2+i} = k_2) \}^2 } \bigg] \\[.5em]
\leq ~ & 16d_1^2d_2^2 \sum_{k_1=1}^{d_1}\sum_{k_2=1}^{d_2} \mE_P \bigg[ \frac{p_Y^2(k)p_Z^2(k)}{ \{1 + \sum_{i=1}^{m_1} \ind(Y_{3n/2+i} = k_1)\}\{ 1 + \sum_{j=1}^{m_2} \ind(Z_{5n/2+i} = k_2) \} } \bigg] \\[.5em]
\overset{(i)}{\leq} ~ & 16d_1^2d_2^2 \sum_{k_1=1}^{d_1}\sum_{k_2=1}^{d_2} \frac{p_Y^2(k)p_Z^2(k)}{(m_1+1)(m_2+1)p_Y(k)p_Z(k)} \\[.5em]
\leq ~ & \frac{16d_1^2d_2^2}{(m_1+1)(m_2+1)},
\end{align*}
where $(i)$ uses the independence between $\{Y_{3n/2+1},\ldots,Y_{3n/2+m_1}\}$ and $\{ Z_{5n/2+1}, \ldots, Z_{5n/2+m_2} \}$ and also the inverse binomial moment in (\ref{Eq: inverse binomial moment}). Therefore Markov's inequality yields
\begin{align*}
\mP_P\left( \sqrt{\sum_{k_1=1}^{d_1}\sum_{k_2=1}^{d_2} \frac{p_Y^2(k)p_Z^2(k)}{w_{k_1,k_2}^2}} \geq t \right) \leq ~ & \frac{1}{t^2} \mE_P \Bigg[ \sum_{k_1=1}^{d_1}\sum_{k_2=1}^{d_2} \frac{p_Y^2(k)p_Z^2(k)}{w_{k_1,k_2}^2} \Bigg] \\[.5em]
\leq ~ & \frac{16d_1^2d_2^2}{t^2(m_1+1)(m_2+1)}.
\end{align*}
This implies that with probability at least $1 - \beta$, we have
\begin{align*}
\sqrt{\sum_{k_1=1}^{d_1}\sum_{k_2=1}^{d_2} \frac{p_Y^2(k)p_Z^2(k)}{w_{k_1,k_2}^2}} \leq \frac{4d_1d_2}{\sqrt{\beta(m_1+1)(m_2+1)}}.
\end{align*} 
Under this event, condition~(\ref{Eq: sufficient condition 2}) is implied by 
\begin{align*}
\epsilon_{n}^2 \geq \frac{C_1}{\beta^{3/2}} \log \left( \frac{1}{\alpha} \right) \max\Bigg\{ \frac{d_1d_2}{m_1^{1/2}m_2^{1/2}n}, \ \frac{d_1^{1/2}d_2^{1/2}}{n} \Bigg\}. 
\end{align*}
By putting the definition of $m_1 = \min\{n/2,d_1\}$ and $m_2 = \min\{n/2,d_2\}$ where $d_1 \leq d_2$ and noting that $\epsilon_n^2 \leq 4$, the condition is further implied by 
\begin{align*}
\epsilon_n^2 \geq  \frac{C_2}{\beta^{3/2}} \log \left( \frac{1}{\alpha} \right) \max\Bigg\{ \frac{d_1^{1/4} d_2^{1/2}}{n^{3/4}}, \ \frac{d_1^{1/2}d_2^{1/2}}{n} \Bigg\},
\end{align*}
for a sufficiently large $C_2>0$. The remaining steps are exactly the same as those in the proof of Proposition~\ref{Proposition: Multinomial L1 Testing}. This completes the proof of Proposition~\ref{Proposition: Multinomial independence testing in L1 distance}.

\section{Proof of Proposition~\ref{Proposition: Gaussian MMD}} \label{Section: Proof of Proposition: Gaussian MMD}

The proof of Proposition~\ref{Proposition: Gaussian MMD} is motivated by \cite{meynaoui2019aggregated} who study the uniform separation rate for the HSIC test. In contrast to \cite{meynaoui2019aggregated} who use the critical value based on the (theoretical) null distribution, we study the permutation test base on the MMD statistic. The structure of the proof is as follows. We first upper bound $\psi_{Y,1}(P)$, $\psi_{Z,1}(P)$ and $\psi_{YZ,2}(P)$ to verify the sufficient condition given in Lemma~\ref{Lemma: Two-Sample U-statistic Improved Version}. We then provide a connection between the expected value of the MMD statistic and $L_2$ distance $|\!|\!| f_Y - f_Z |\!|\!|_{L_2}$. Finally, we conclude the proof based on the previous results. Throughout the proof, we write the Gaussian kernel $K_{\lambda_1,\ldots,\lambda_d,d}(x-y)$ in (\ref{Eq: Gaussian kernel}) as $K_{\lambda,d}(x-y)$ so as to simplify the notation.

\paragraph{$\bullet$ Verification of condition (\ref{Eq: Two-Sample U-statistic Sufficient Condition 2}).}
In this part of the proof, we find upper bounds for $\psi_{Y,1}(P)$, $\psi_{Z,1}(P)$ and $\psi_{YZ,2}(P)$. Let us start with $\psi_{Y,1}(P)$. Recall that $\psi_{Y,1}(P)$ is given as
\begin{align*}
\psi_{Y,1}(P) =  \text{Var}_P\{\mE_P[\overline{h}_{\text{ts}}(Y_1,Y_2;Z_1,Z_2)|Y_1] \},
\end{align*}
where $\overline{h}_{\text{ts}}$ is the symmetrized kernel (\ref{Eq: symmetrized kernel}). Using the definition, it is straightforward to see that
\begin{align*}
\psi_{Y,1}(P) = ~ & \mV_P \{ \mE_P[g_{\text{Gau}}(Y_1,Y_2)|Y_1] -  \mE_P[g_{\text{Gau}}(Y_1,Z_1)|Y_1] \} \\[.5em]
\leq ~ & \mE_P [ \{ \mE_P[g_{\text{Gau}}(Y_1,Y_2)|Y_1] -  \mE_P[g_{\text{Gau}}(Y_1,Z_1)|Y_1] \}^2 ].
\end{align*}
Let us denote the convolution $f_Y - f_Z$ and $K_{\lambda,d}$ by
\begin{align*}
(f_Y - f_Z) \ast K_{\lambda,d}(x) = \int_{\mathbb{R}^d}  [f_Y(t) - f_Z(t)]K_{\lambda,d}(x-t)dt,
\end{align*}
where $K_{\lambda,d}$ can be recalled from (\ref{Eq: Gaussian kernel}). Then the upper bound of $\psi_{Y,1}(P)$ is further bounded by
\begin{align*}
\mE_P [ \{ \mE_P[g_{\text{Gau}}(X_1,X_2)|X_1] -  \mE_P[g_{\text{Gau}}(X_1,Y_1)|X_1] \}^2 ] = ~ & \int_{\mathbb{R}^d} f_Y(x) \big[ (f_Y - f_Z) \ast K_{\lambda,d}(x) \big]^2dx \\[.5em]
\leq ~ & |\!|\!|f_Y|\!|\!|_{\infty}  |\!|\!| (f_Y - f_Z) \ast K_{\lambda,d} |\!|\!|_{L_2}^2.
\end{align*}
By symmetry, $\psi_{Z,1}(P)$ can be similarly bounded. Thus
\begin{equation}
\begin{aligned} \label{Eq: bound 1}
\psi_{Y,1}(P) \leq & ~ |\!|\!|f_Y|\!|\!|_{\infty}  |\!|\!| (f_Y - f_Z) \ast K_{\lambda,d} |\!|\!|_{L_2}^2, \\[.5em]
\psi_{Z,1}(P) \leq & ~ |\!|\!|f_Z|\!|\!|_{\infty}  |\!|\!| (f_Y - f_Z) \ast K_{\lambda,d} |\!|\!|_{L_2}^2.
\end{aligned}
\end{equation}
Moving onto $\psi_{YZ,2}(P)$, we need to compute $\mE_P[g_{\text{Gau}}^2(Y_1,Y_2)]$, $\mE_P[g_{\text{Gau}}^2(Z_1,Z_2)]$ and $\mE_P[g_{\text{Gau}}^2(Y_1,Z_1)]$. Note that 
\begin{align*}
K_{\lambda,d}^2(x) = \frac{1}{(4\pi)^{d/2} \lambda_1 \cdots \lambda_d} K_{\lambda/\sqrt{2},d}(x), 
\end{align*}
where $K_{\lambda/\sqrt{2},d}(x)$ is the Gaussian density function (\ref{Eq: Gaussian kernel}) with scale parameters $\lambda_1/\sqrt{2}, \ldots, \lambda_d/\sqrt{2}$. Therefore it can be seen that
\begin{align*}
\mE_P[g_{\text{Gau}}^2(Y_1,Y_2)] ~= ~ & \int_{\mathbb{R}^d}  \int_{\mathbb{R}^d} K_{\lambda,d}^2(y_1-y_2) f_Y(y_1)f_Y(y_2) dy_1 d y_2 \\[.5em]
= ~ & \frac{1}{(4\pi)^{d/2} \lambda_1 \cdots \lambda_d} \int_{\mathbb{R}^d}  \int_{\mathbb{R}^d} K_{\lambda/\sqrt{2},d}(y_1-y_2) f_Y(y_1)f_Y(y_2) dy_1 d y_2 \\[.5em]
\leq  ~ &  \frac{|\!|\!|f_Y|\!|\!|_\infty}{(4\pi)^{d/2} \lambda_1 \cdots \lambda_d} \int_{\mathbb{R}^d}  \bigg[ \int_{\mathbb{R}^d} K_{\lambda/\sqrt{2},d}(y_1-y_2)  dy_1 \bigg]  f_Y(y_2)d y_2  \\[.5em]
\leq ~ &  \frac{M_{f,d}}{(4\pi)^{d/2} \lambda_1 \cdots \lambda_d},
\end{align*}
where $\max\{|\!|\!|f_Y|\!|\!|_\infty, |\!|\!|f_Z|\!|\!|_\infty \} \leq M_{f,d}$. The other two terms $\mE_P[g_{\text{Gau}}^2(Z_1,Z_2)]$ and $\mE_P[g_{\text{Gau}}^2(Y_1,Z_1)]$ are similarly bounded. Thus we have
\begin{align} \label{Eq: bound 2}
\psi_{YZ,2}(P) \leq \frac{M_{f,d}}{(4\pi)^{d/2} \lambda_1 \cdots \lambda_d}. 
\end{align}
Given bounds (\ref{Eq: bound 1}) and (\ref{Eq: bound 2}), Lemma~\ref{Lemma: Two-Sample U-statistic Improved Version} shows that the type II error of the considered permutation test is at most $\beta$ when 
\begin{equation}
\begin{aligned} \label{Eq: condition for the variance part}
\mE_P[U_{n_1,n_2}] ~\geq~ & C_1(M_{f,d},d)  \sqrt{\frac{|\!|\!| (f_Y - f_Z) \ast K_{\lambda,d} |\!|\!|_{L_2}^2}{\beta} \left( \frac{1}{n_1} + \frac{1}{n_2} \right)}  \\[.5em]
+ ~ &  \frac{C_2(M_{f,d},d)}{\sqrt{\lambda_1 \cdots \lambda_d}} \frac{1}{\sqrt{\beta}} \log \left( \frac{1}{\alpha} \right)    \left( \frac{1}{n_1} + \frac{1}{n_2} \right).
\end{aligned}
\end{equation}

\paragraph{$\bullet$ Relating $\mE_P[U_{n_1,n_2}]$ to $L_2$ distance.}
Next we related the expected value of $U_{n_1,n_2}$ to $L_2$ distance between $f_Y$ and $f_Z$. Based on the unbiasedness property of a $U$-statistic, one can easily verify that 
\begin{equation}
\begin{aligned} \label{Eq: condition for the mean part}
\mE_P[U_{n_1,n_2}]  ~ = ~ & \int_{\mathbb{R}^d} \int_{\mathbb{R}^d} K_{\lambda,d} (t_1-t_2)[f_Y(t_1) - f_Z(t_1)][f_Y(t_2) - f_Z(t_2)] dt_1 dt_2 \\[.5em]
= ~ & \int_{\mathbb{R}^d}  [f_Y(t_2) - f_Z(t_2)]  (f_Y - f_Z) \ast K_{\lambda,d} (t_2) dt_2 \\[.5em]
= ~ & \frac{1}{2}|\!|\!| f_Y - f_Z |\!|\!|_{L_2}^2  + \frac{1}{2} |\!|\!| (f_Y -f_Z) \ast K_{\lambda,d} |\!|\!|_{L_2}^2  \\[.5em]
- & \frac{1}{2} |\!|\!| (f_Y - f_Z)  - (f_Y - f_Z) \ast K_{\lambda,d} |\!|\!|_{L_2}^2.
\end{aligned}
\end{equation}
where the last equality uses the fact that $2xy = x^2 + y^2 - (x-y)^2$. 

\paragraph{$\bullet$ Completion of the proof.}
We now combine the previous results (\ref{Eq: condition for the variance part}) and (\ref{Eq: condition for the mean part}) to conclude the result. To be more specific, based on equality (\ref{Eq: condition for the mean part}), it is seen that condition~(\ref{Eq: condition for the variance part}) is equivalent to
\begin{equation}
\begin{aligned} \label{Eq: sufficient condition}
|\!|\!| f_Y - f_Z |\!|\!|_{L_2}^2  ~\geq~ & |\!|\!| (f_Y - f_Z)  - (f_Y - f_Z) \ast K_{\lambda,d} |\!|\!|_{L_2}^2  \\[.5em] 
- ~ &  |\!|\!| (f_Y -f_Z) \ast K_{\lambda,d} |\!|\!|_{L_2}^2 \\[.5em]
+ ~ & C_3(M_{f,d},d)  \sqrt{\frac{|\!|\!| (f_Y - f_Z) \ast K_{\lambda,d} |\!|\!|_{L_2}^2}{\beta} \left( \frac{1}{n_1} + \frac{1}{n_2} \right)}  \\[.5em]
+ ~ &  \frac{C_4(M_{f,d},d)}{\sqrt{\lambda_1 \cdots \lambda_d}} \frac{1}{\sqrt{\beta}} \log \left( \frac{1}{\alpha} \right)  \left( \frac{1}{n_1} + \frac{1}{n_2} \right).
\end{aligned}
\end{equation}
Based on the basic inequality $\sqrt{xy} \leq x + y$ for $x,y \geq 0$, we can upper bound the third line of the above equation as
\begin{align*}
C_3(M_{f,d},d)  \sqrt{\frac{|\!|\!| (f_Y - f_Z) \ast K_{\lambda,d} |\!|\!|_{L_2}^2}{\beta} \left( \frac{1}{n_1} + \frac{1}{n_2} \right)} ~ \leq ~ &  \frac{C_5(M_{f,d},d)}{\beta} \left( \frac{1}{n_1} + \frac{1}{n_2} \right) \\[.5em]
+~ &   |\!|\!| (f_Y - f_Z) \ast K_{\lambda,d} |\!|\!|_{L_2}^2.
\end{align*}
Therefore the previous inequality (\ref{Eq: sufficient condition}) is implied by 
\begin{align*}
\epsilon_{n_1,n_2}^2  ~ \geq ~ & |\!|\!| (f_Y - f_Z)  - (f_Y - f_Z) \ast K_{\lambda,d} |\!|\!|_{L_2}^2 \\[.5em]
& + \frac{C(M_{f,d},d)}{\beta \sqrt{\lambda_1 \cdots \lambda_d}}  \log \left( \frac{1}{\alpha} \right) \cdot \left( \frac{1}{n_1} + \frac{1}{n_2} \right),
\end{align*}
where we used the condition $\prod_{i=1}^d\lambda_i \leq 1$. This completes the proof of Proposition~\ref{Proposition: Gaussian MMD}.

\vskip 1em

\section{Proof of Proposition~\ref{Proposition: Two-moment method for MC-based tests}} \label{Section: Proof of Proposition: Two-moment method for MC-based tests}
As stated in Appendix~\ref{Section: Monte Carlo-based permutation tests}, type I error control for the MC-based test is well-known and its proof can be found in, for example, \cite{hemerik2018exact}. Hence we focus on the type II error rate. Recall that the original permutation distribution is defined as
\begin{align*}
	F_{T_n^{\pi}}(t) = \frac{1}{M_n} \sum_{\pi \in \mathbf{\Pi}_n} \ind\{T_n^{\pi} \leq t \},
\end{align*}
and we denote its Monte Carlo approximation by
\begin{align*}
	\widehat{F}_{T_n^{\pi}}(t) :=  \frac{1}{B} \sum_{i=1}^{B}  \mathds{1} \big(T_n^{\pi_i} \leq t \big), 
\end{align*}
where $\pi_1,\ldots,\pi_B$ be random permutations uniformly sampled with replacement from $\mathbf{\Pi}_n$. Let us define an event 
\begin{align*}
	\mathcal{A} := \Bigg\{ \sup_{t \in \mathbb{R}} \big| \widehat{F}_{T_n^{\pi}}(t) - F_{T_n^{\pi}}(t)  \big| \leq \sqrt{\frac{1}{2B} \log \left( \frac{4}{\beta} \right)}  \Bigg\},
\end{align*}
which holds with probability at least $1 - \beta/2$ by Dvoretzky--Kiefer--Wolfowitz inequality \citep{dvoretzky1956asymptotic,massart1990tight}. Then, under the event $\mathcal{A}$, the definition of $\widehat{c}_{1-\alpha, n}$ in (\ref{Eq: MC critical value}) yields
\begin{align*}
	\widehat{c}_{1-\alpha, n} ~ \leq ~ & \inf\bigg\{t: \frac{1}{B+1} \sum_{i=1}^{B}  \mathds{1} \big(T_n^{\pi_i} \leq t \big) \geq 1-\alpha \bigg\} \\[.5em]
	= ~ & \inf\bigg\{t: \widehat{F}_{T_n^{\pi}}(t) \geq \frac{B+1}{B}(1-\alpha) \bigg\} \\[.5em]
	\leq ~ &  \inf\bigg\{t: F_{T_n^{\pi}}(t) \geq \underbrace{\frac{B+1}{B}(1-\alpha) + \sqrt{\frac{1}{2B} \log \left( \frac{4}{\beta} \right)}}_{:= 1 - \alpha^\ast}  \bigg\} \\[.5em]
	= ~ & c_{1-\alpha^\ast,n}.
\end{align*}
Therefore the type II error of the MC-based test is bounded by
\begin{align} \nonumber
	\sup_{P \in \mathcal{P}_1} \mP_P(\widehat{p} > \alpha) ~=~& \sup_{P \in \mathcal{P}_1} \mP_P(T_n \leq \widehat{c}_{1-\alpha,n}, \mathcal{A} \cup \mathcal{A}^c) \\[.5em]
	\leq~ & \sup_{P \in \mathcal{P}_1} \mP_P(T_n \leq c_{1-\alpha^\ast,n}) + \frac{\beta}{2}.  \label{Eq: intermediate step}
\end{align}
After some algebra, it can be verified that when $B \geq 8 \alpha^{-2} \log (4/\beta)$, it holds $1-\alpha^\ast \leq 1 - \alpha/2$. Consequently, we can upper bound the type II error of $\mathds{1}(T_n \leq c_{1-\alpha^\ast,n})$ by  
\begin{align*}
	\sup_{P \in \mathcal{P}_1} \mP_P(T_n \leq c_{1-\alpha^\ast,n}) ~\leq~ & \sup_{P \in \mathcal{P}_1} \mP_P(T_n \leq c_{1-\alpha/2,n}) \\[.5em] 
	\overset{(i)}{\leq} ~ & \beta/2, 
\end{align*}
where the inequality~$(i)$ follows by Lemma~\ref{Lemma: Two Moments Method} under the condition~(\ref{Eq: Sufficient Condition for MC}). Combining this bound with (\ref{Eq: intermediate step}), we obtain the desired type II error rate. This completes the proof of Proposition~\ref{Proposition: Two-moment method for MC-based tests}. 

\end{document}